\documentclass[a4paper,12pt]{article}

\usepackage[all]{xy}
\usepackage[T1]{fontenc}
\usepackage[dvips]{graphicx}
\usepackage{amsmath,amssymb,amsfonts,amsthm,latexsym,mathrsfs,textcomp,verbatim, enumerate, soul, color}
\xyoption{web}

\title{A topos-theoretic approach to\\ Stone-type dualities}
\author{Olivia Caramello \vspace{3 mm}\\ {\small DPMMS, University of Cambridge,}\\{\small Wilberforce Road, Cambridge CB3 0WB, U.K.}\\{\small O.Caramello@dpmms.cam.ac.uk}\thanks{The author gratefully acknowledges the support of a Research Fellowship from Jesus College, Cambridge (U.K.) and of a visiting position from the Centro di Ricerca Matematica Ennio De Giorgi, Pisa (Italy).}}

\date{March 17, 2011}
\begin{document}

\mathcode`\<="4268  
\mathcode`\>="5269  
\mathcode`\.="313A  
\mathchardef\semicolon="603B 
\mathchardef\gt="313E
\mathchardef\lt="313C

\newcommand{\app}
 {{\sf app}}

\newcommand{\Ass}
 {{\bf Ass}}

\newcommand{\ASS}
 {{\mathbb A}{\sf ss}}

\newcommand{\Bb}
{\mathbb}

\newcommand{\biimp}
 {\!\Leftrightarrow\!}

\newcommand{\bim}
 {\rightarrowtail\kern-1em\twoheadrightarrow}

\newcommand{\bjg}
 {\mathrel{{\dashv}\,{\vdash}}}

\newcommand{\bstp}[3]
 {\mbox{$#1\! : #2 \bim #3$}}

\newcommand{\cat}
 {\!\mbox{\t{\ }}}

\newcommand{\cinf}
 {C^{\infty}}

\newcommand{\cinfrg}
 {\cinf\hy{\bf Rng}}

\newcommand{\cocomma}[2]
 {\mbox{$(#1\!\uparrow\!#2)$}}

\newcommand{\cod}
 {{\rm cod}}

\newcommand{\comma}[2]
 {\mbox{$(#1\!\downarrow\!#2)$}}

\newcommand{\comp}
 {\circ}

\newcommand{\cons}
 {{\sf cons}}

\newcommand{\Cont}
 {{\bf Cont}}

\newcommand{\ContE}
 {{\bf Cont}_{\cal E}}

\newcommand{\ContS}
 {{\bf Cont}_{\cal S}}

\newcommand{\cover}
 {-\!\!\triangleright\,}

\newcommand{\cstp}[3]
 {\mbox{$#1\! : #2 \cover #3$}}

\newcommand{\Dec}
 {{\rm Dec}}

\newcommand{\DEC}
 {{\mathbb D}{\sf ec}}

\newcommand{\den}[1]
 {[\![#1]\!]}

\newcommand{\Desc}
 {{\bf Desc}}

\newcommand{\dom}
 {{\rm dom}}

\newcommand{\Eff}
 {{\bf Eff}}

\newcommand{\EFF}
 {{\mathbb E}{\sf ff}}

\newcommand{\empstg}
 {[\,]}

\newcommand{\epi}
 {\twoheadrightarrow}

\newcommand{\estp}[3]
 {\mbox{$#1 \! : #2 \epi #3$}}

\newcommand{\ev}
 {{\rm ev}}

\newcommand{\Ext}
 {{\rm Ext}}

\newcommand{\fr}
 {\sf}

\newcommand{\fst}
 {{\sf fst}}

\newcommand{\fun}[2]
 {\mbox{$[#1\!\to\!#2]$}}

\newcommand{\funs}[2]
 {[#1\!\to\!#2]}

\newcommand{\Gl}
 {{\bf Gl}}

\newcommand{\hash}
 {\,\#\,}

\newcommand{\hy}
 {\mbox{-}}

\newcommand{\im}
 {{\rm im}}

\newcommand{\imp}
 {\!\Rightarrow\!}

\newcommand{\Ind}[1]
 {{\rm Ind}\hy #1}

\newcommand{\iten}[1]
{\item[{\rm (#1)}]}

\newcommand{\iter}
 {{\sf iter}}

\newcommand{\Kalg}
 {K\hy{\bf Alg}}

\newcommand{\llim}
 {{\mbox{$\lower.95ex\hbox{{\rm lim}}$}\atop{\scriptstyle
{\leftarrow}}}{}}

\newcommand{\llimr}
 {{\mbox{$\lower.95ex\hbox{{\rm lim}}$}\atop{\scriptstyle
{\rightarrow}}}{}}

\newcommand{\llimd}
 {\lower0.37ex\hbox{$\pile{\lim \\ {\scriptstyle
\leftarrow}}$}{}}

\newcommand{\Mf}
 {{\bf Mf}}

\newcommand{\Mod}
 {{\bf Mod}}

\newcommand{\MOD}
{{\mathbb M}{\sf od}}

\newcommand{\mono}
 {\rightarrowtail}

\newcommand{\mor}
 {{\rm mor}}

\newcommand{\mstp}[3]
 {\mbox{$#1\! : #2 \mono #3$}}

\newcommand{\Mu}
 {{\rm M}}

\newcommand{\name}[1]
 {\mbox{$\ulcorner #1 \urcorner$}}

\newcommand{\names}[1]
 {\mbox{$\ulcorner$} #1 \mbox{$\urcorner$}}

\newcommand{\nml}
 {\triangleleft}

\newcommand{\ob}
 {{\rm ob}}

\newcommand{\op}
 {^{\rm op}}
 
\newcommand{\palrr}[4]{ 
  \def\labelstyle{\scriptstyle} 
  \xymatrix{ {#1} \ar@<0.5ex>[r]^{#2} \ar@<-0.5ex>[r]_{#3} & {#4} } } 
  
\newcommand{\palrl}[4]{ 
  \def\labelstyle{\scriptstyle} 
  \xymatrix{ {#1} \ar@<0.5ex>[r]^{#2}  &  \ar@<0.5ex>[l]^{#3} {#4} } }  

\newcommand{\pepi}
 {\rightharpoondown\kern-0.9em\rightharpoondown}

\newcommand{\pmap}
 {\rightharpoondown}

\newcommand{\Pos}
 {{\bf Pos}}

\newcommand{\prarr}
 {\rightrightarrows}

\newcommand{\princfil}[1]
 {\mbox{$\uparrow\!(#1)$}}

\newcommand{\princid}[1]
 {\mbox{$\downarrow\!(#1)$}}

\newcommand{\prstp}[3]
 {\mbox{$#1\! : #2 \prarr #3$}}

\newcommand{\pstp}[3]
 {\mbox{$#1\! : #2 \pmap #3$}}

\newcommand{\relarr}
 {\looparrowright}

\newcommand{\rlim}
 {{\mbox{$\lower.95ex\hbox{{\rm lim}}$}\atop{\scriptstyle
{\rightarrow}}}{}}

\newcommand{\rlimd}
 {\lower0.37ex\hbox{$\pile{\lim \\ {\scriptstyle
\rightarrow}}$}{}}

\newcommand{\rstp}[3]
 {\mbox{$#1\! : #2 \relarr #3$}}

\newcommand{\scn}
 {{\bf scn}}

\newcommand{\scnS}
 {{\bf scn}_{\cal S}}

\newcommand{\semid}
 {\rtimes}

\newcommand{\Sep}
 {{\bf Sep}}

\newcommand{\sep}
 {{\bf sep}}

\newcommand{\Set}
 {{\bf Set}}

\newcommand{\Sh}
 {{\bf Sh}}

\newcommand{\ShE}
 {{\bf Sh}_{\cal E}}

\newcommand{\ShS}
 {{\bf Sh}_{\cal S}}

\newcommand{\sh}
 {{\bf sh}}

\newcommand{\Simp}
 {{\bf \Delta}}

\newcommand{\snd}
 {{\sf snd}}

\newcommand{\stg}[1]
 {\vec{#1}}

\newcommand{\stp}[3]
 {\mbox{$#1\! : #2 \to #3$}}

\newcommand{\Sub}
 {{\rm Sub}}

\newcommand{\SUB}
 {{\mathbb S}{\sf ub}}

\newcommand{\tbel}
 {\prec\!\prec}

\newcommand{\tic}[2]
 {\mbox{$#1\!.\!#2$}}

\newcommand{\tp}
 {\!:}

\newcommand{\tps}
 {:}

\newcommand{\tsub}
 {\pile{\lower0.5ex\hbox{.} \\ -}}

\newcommand{\wavy}
 {\leadsto}

\newcommand{\wavydown}
 {\,{\mbox{\raise.2ex\hbox{\hbox{$\wr$}
\kern-.73em{\lower.5ex\hbox{$\scriptstyle{\vee}$}}}}}\,}

\newcommand{\wbel}
 {\lt\!\lt}

\newcommand{\wstp}[3]
 {\mbox{$#1\!: #2 \wavy #3$}}
 
\newcommand{\fu}[2]
{[#1,#2]}


%
%
%
\def\pushright#1{{
   \parfillskip=0pt            
   \widowpenalty=10000         
   \displaywidowpenalty=10000  
   \finalhyphendemerits=0      
  %
   \leavevmode                 
   \unskip                     
   \nobreak                    
   \hfil                       
   \penalty50                  
   \hskip.2em                  
   \null                       
   \hfill                      
   {#1}                        
  %
   \par}}                      

\def\qed{\pushright{$\square$}\penalty-700 \smallskip}

\newtheorem{theorem}{Theorem}[section]

\newtheorem{proposition}[theorem]{Proposition}

\newtheorem{scholium}[theorem]{Scholium}

\newtheorem{lemma}[theorem]{Lemma}

\newtheorem{corollary}[theorem]{Corollary}

\newtheorem{conjecture}[theorem]{Conjecture}

\newenvironment{proofs}%
 {\begin{trivlist}\item[]{\bf Proof }}%
 {\qed\end{trivlist}}

  \newtheorem{rmk}[theorem]{Remark}
\newenvironment{remark}{\begin{rmk}\em}{\end{rmk}}

  \newtheorem{rmks}[theorem]{Remarks}
\newenvironment{remarks}{\begin{rmks}\em}{\end{rmks}}

  \newtheorem{defn}[theorem]{Definition}
\newenvironment{definition}{\begin{defn}\em}{\end{defn}}

  \newtheorem{eg}[theorem]{Example}
\newenvironment{example}{\begin{eg}\em}{\end{eg}}

  \newtheorem{egs}[theorem]{Examples}
\newenvironment{examples}{\begin{egs}\em}{\end{egs}}


\bgroup           
\let\footnoterule\relax  
\maketitle

\begin{abstract}
We present an abstract unifying framework for interpreting Stone-type dualities; several known dualities are seen to be instances of just one topos-theoretic phenomenon, and new dualities are introduced. In fact, infinitely many new dualities between preordered structures and locales or topological spaces can be generated through our topos-theoretic machinery in a uniform way. We then apply our topos-theoretic interpretation to obtain results connecting properties of preorders and properties of the corresponding locales or topological spaces, and we establish adjunctions between various kinds of categories as natural applications of our general methodology. In the last part of the paper, we exploit the theory developed in the previous parts to obtain a topos-theoretic interpetation of the problem of finding explicit descriptions of models of `ordered algebraic theories' presented by generators and relations, and give several examples which illustrate the effectiveness of our methodology. In passing, we provide a number of other applications of our theory to Algebra, Topology and Logic.
\end{abstract} 
\egroup

\newpage

\tableofcontents

\newpage

\section{Introduction}

In this paper we present a general topos-theoretic interpretation of `Stone-type dualities'; by this term we refer, following the standard terminology, to a class of dualities or equivalences between categories of preordered structures and categories of posets, locales or topological spaces, a class which notably includes the classical Stone duality for Boolean algebras (or, more generally, for distributive lattices), the duality between spatial frames and sober spaces, the equivalence between preorders and Alexandrov spaces, the Lindenbaum-Tarski duality between sets and complete atomic Boolean algebras, and the Birkhoff's duality between finite distributive lattices and finite posets. 

We introduce an abstract framework in which all of these dualities are interpreted as instances of just one topos-theoretic phenomenon, and in which several new dualities are introduced. In fact, the known dualities, as well as the new ones, all arise from the application of one `general machinery for generating dualities' to specific `sets of inputs' which vary from case to case. 

In section \ref{locales} we show that, under relatively mild hypotheses, one can naturally identify, through an isomorphism of categories, the opposite of a given category of ordered structures with a subcategory of the category of locales, in general in more than one way; in fact, this result (Theorem \ref{teoabstr}) provides us with an infinite number of dualities between categories of posets and subcategories of the category of locales. Anyway, these subcategories are not in general closed under arbitrary isomorphisms of locales, and in fact, in order to obtain `intrinsic' dualities between categories of preorders and categories of locales whose objects (resp. arrows) can be characterized as the locales (resp. locale homomorphisms) which satisfy some locale-theoretic invariant, we have to enlarge the target subcategory of locales to include all the isomorphic copies of the objects and arrows in it, and to look for a functor defined on this enlarged subcategory which is inverse (up to isomorphism) to the functor from posets to locales forming one half of the original isomorphism of categories. This can be done, under some natural assumptions, by functorially transferring topos-theoretic invariants across two different sites of definition of the same topos, according to the method `toposes as bridges' introduced in \cite{OC10}. For instance, the dualities between a given category $\cal K$ of preorders and a category of locales arise from the process of assigning to each structure $\cal C$ of $\cal K$, equipped with a subcanonical Grothendieck topology $J_{\cal C}$ in such a way that the morphisms in the category $\cal K$ induce morphisms of the associated sites, the locale $Id_{J}({\cal C})$ of $J$-ideals on $\cal C$, and from the inverse process of functorially recovering $\cal C$ from the locale $Id_{J}({\cal C})$ (equivalently, from the topos $\Sh({\cal C}, J) \simeq \Sh(Id_{J}({\cal C}))$) through a topos-theoretic invariant. The covariant equivalences with categories of locales are established in a similar way; the structures need not be equipped with any Grothendieck topology, and one relies on the well-known possibility of assigning a geometric morphism $[{\cal C}, \Set]\to [{\cal D}, \Set]$ to a given functor ${\cal C} \to {\cal D}$ in a canonical way. 

In section \ref{dualtop}, we give a general methodology for `enriching' a given duality (resp. equivalence) between a category $\cal K$ of preorders and a category of locales to a duality (resp. equivalence) between $\cal K$ and a category of topological spaces; this methodology relies on an appropriate choice of points of the toposes corresponding to the structures.          
       
In the following sections of the paper, we investigate further consequences of the topos-theoretic perspective introduced in the previous sections, again in light of the method `toposes as bridges' of \cite{OC10}. For example, we apply this method to establish various results connecting properties of preordered structures with properties of the corresponding locales or topological spaces, and we obtain a number of adjunctions between categories of these kinds. 

The theory developed in the present paper provides a unified perspective on the subject of Stone-type dualities, in that the well-known dualities are easily recovered as applications of it. Anyway, what we consider to be the main interest of our topos-theoretic machinery is, apart from the conceptual enlightenment that it brings into the world of classical dualities, its inherent technical flexibility. In fact, one can generate infinitely many new dualities by applying it; examples are provided in the paper to illustrate how to do this in practice, and the reader will be able to use the method to generate his or her favorite applications.   

The different `ingredients' that our `machinery' for generating dualities with categories of locales or topological spaces takes as `inputs' are: the initial category $\cal K$ of preordered structures, the subcanonical Grothendieck topologies $J_{\cal C}$ on the structures $\cal C$ in $\cal K$, the topos-theoretic invariant enabling one to recover a structure $\cal C$ from the topos $\Sh({\cal C}, J_{\cal C})$ and, if a duality with topological spaces is to be generated, appropriate sets of points of the toposes $\Sh({\cal C}, J_{\cal C})$ (and functions between them). In fact, the more general approach of section \ref{generalization} provides us with an additional degree of freedom in the choice of ingredients. Given such ingredients, dualities are generated in an automatic and `uniform' way by the `machine', as different concrete instances of a unique abstract pattern; in this way, the problem of building dualities gets reduced in many important cases to the much easier one of choosing appropriate sets of ingredients for this `machine'.

In connection with the perspectives outlined in \cite{OC10}, we remark that in this paper we have just `brought to the surface' a limited number of results which can be established by means of the methods of \cite{OC10}, carefully selected by virtue of their `representativeness' in illustrating the nature and variety of the insights obtainable by applying our techniques. In fact, the paper contains many general ideas which can be applied in the context of arbitrary, rather than just preordered, categories; for example, the method of section \ref{Mordual} of building dualities or equivalences starting from Morita-equivalences (by equipping structures, regarded as categories, with subcanonical Grothendieck topologies in such a way that the `structure-preserving' maps between the structures yield morphisms of the associated sites, and recovering each of the structures functorially from the corresponding topos through a topos-theoretic invariant), is potentially applicable beyond the preordered context that we have addressed in the present paper (cf. section \ref{conclusions} for a further discussion of this point). Even in the preordered context that we have addressed in the present paper, much remains to be discovered, in the form of new dualities, representation theorems, adjunctions and characterization theorems arising from translations of properties from one side to another of a given duality or equivalence. Anyway, the careful reader will realize that much of all of this can be easily uncovered in a semi-automatic way, by using similar means to those that we have adopted in the paper to generate our examples. 

\subsection{An overview of the paper}

The contents of the present work can be described more in detail as follows. 

In section \ref{subterminaltop}, we introduce a general method for building topological spaces from toposes equipped with a set of points. Specifically, we show that, given any subframe $\Gamma$ of the frame of subterminals of a locally small cocomplete topos $\cal E$ and any set of points of $\cal E$ indexed by a set $I$, we can naturally define a topology on the set $I$, which we call the ($\Gamma$-)subterminal topology, and that this construction can be naturally made functorial. The interest of this notion lies in its level of generality, which encompasses that of classical topology (every topological space arises from this construction in a canonical way), as well as in its formulation as a topos-theoretic invariant admitting a `natural behaviour' with respect to sites. Indeed, as shown in section \ref{exsub}, this notion allows us to recover, with natural choices of sites of definition and of sets of points of toposes, many interesting topological spaces considered in the literature, leaving at the same time enough freedom to construct new ones with particular properties. In section \ref{Sobriety} we give a general criterion for deciding the sobriety of topological spaces built in this way.

In section \ref{general}, we present our general topos-theoretic interpretation of Stone-type dualities; first we discuss dualities with categories of locales, then we introduce a general methodology, based on the notion of subterminal topology, for `enriching' them so to obtain dualities with categories of topological spaces. 

In section \ref{ex}, we discuss various examples of dualities generated by using our method; we recover the classical Stone duality for distributive lattices (and Boolean algebras), the Alexandrov duality between preorders and Alexandrov spaces, the Lindenbaum-Tarski duality, the duality between spatial frames and sober spaces, and we establish new ones, including localic and topological dualities for meet-semilattices, an equivalence between the category of posets and a category of spatial locales (equivalently, a category of sober topological spaces), a localic duality for $k$-frames (for a regular cardinal $k$), and new dualities between specific categories of preordered structures. 

In section \ref{generalization} we further generalize the method of section \ref{locales} for building dualities or equivalences with categories of locales starting from Morita-equivalences of the form $\Sh({\cal C}, J) \simeq \Sh(Id_{J}({\cal C}))$ to general equivalences $\Sh({\cal C}, J)\simeq \Sh({\cal D}, K)$, where $\cal C$ and $\cal D$ are preordered categories. This allows an abstract symmetric definition of the functors yielding the dualities, and provides us with an additional degree of freedom in building dualities or equivalences between categories of preordered structures. Grothendieck Comparison Lemma turns out be an extremely fruitful source of Morita-equivalences to which we can apply our methods; we illustrate this point in section \ref{addex} by generating several new dualities or equivalences. In particular, we establish a duality which naturally generalizes Birkhoff's duality for finite distributive lattices, and a duality which generalizes the well-known duality between algebraic lattices and sup-semilattices.  

In section \ref{adj}, we apply the method `toposes as bridges' of \cite{OC10} to the Morita-equivalences $\Sh({\cal C}, J) \simeq \Sh(Id_{J}({\cal C}))$ and to the other equivalences established in section \ref{generalization} to obtain adjunctions which extend the dualities obtained in the previous sections; in particular, we establish reflections from various categories of preordered structures to the category of frames, as well as reflections between categories of posets satisfying some generalized `distributive law' and full subcategories of them consisting of posets satisfying certain `topological conditions'. In this context, as another application of our general method, we establish adjunctions between categories of toposes paired with points (as defined in section \ref{subterminal}) and categories of topological spaces.

In section \ref{insights}, we prove a number of results connecting properties of preordered structures with properties of the locales or topological spaces corresponding to them via the dualities or equivalences considered in the previous sections. Again, the technique that we employ for performing these `translations' is that of using toposes as `bridges' for transferring properties between their distinct sites of definition; specifically, we consider a number of logically-motivated topos-theoretic invariants, admitting bijective site characterizations, and rephrase them in terms of the two different representations $\Sh({\cal C}, J)$ and $\Sh(Id_{J}({\cal C}))$. We also establish various other results for preordered structures by using the same method. Of course, the topos-theoretic notions completely disappear in the final formulation of the results, they are just instrumental for performing the `automatic translation' (in the sense of \cite{OC10}) of properties from one site of definition into another. 

Section \ref{spacesprop} deals with the problem of giving concrete characterizations of the topological spaces arising by putting the subterminal topology on a given set of points of a localic topos (equivalently, on a given set of models of a propositional geometric theory). To this end, we introduce a general method for building analogues of the Zariski spectrum for structures which can be described as models of propositional theories; in this context, we also provide characterizations of the syntactic categories of these theories as ordered algebraic structures presented by generators and relations. We treat this latter topic in full generality by introducing a generalized notion of first-order mathematical theory and a corresponding notion of syntactic category; we then restrict our attention to the syntactic categories of these general propositional theories and show that they can be characterized as models presented by generators and relations of certain Horn theories, which we call `ordered algebraic theories'. These notions pave the way for a topos-theoretic interpretation of the problem of giving explicit descriptions of models of theories of this kind presented by generators and relations. After discussing the abstract features of this interpretation, we illustrate its effectiveness by discussing several examples, notably including the construction (in section \ref{freecomjoin}) of the free frame on a complete join-semilattice, which solves in particular an open problem stated by P. Resende and S. Vickers in 2003 (cf. \cite{RV}). In section \ref{zariski}, we analyze the classical Zariski spectrum from our general topos-theoretic perspective and extend the definition of the Zariski topology on the collection of prime ideals of a ring to arbitrary collections of subsets of the ring.

In the appendix of the paper, we provide `elementary' proofs (that is, proofs which do not rely on results in Topos Theory) of some of the central results in the paper, including those which constitute the `machinery for generating dualities' of sections \ref{general} and \ref{generalization}; in fact, although we have established most of our results through a natural combination of abstract topos-theoretic techniques in light of the philosophy `toposes as bridges' of \cite{OC10}, most of our results, as well as their proofs, can be directly reformulated in the language of Locale Theory. The comparison between the `abstract' and `concrete' proofs of our results is an interesting one; on one hand, the direct arguments might be judged preferable to the abstract ones because they can be easier to understand to the reader who is not familiar with Topos Theory (in fact, this is the main reason why we have decided to provide them in the appendix of the paper), while on the other hand it is precisely the topos-theoretic arguments that are mostly illuminating from a conceptual perspective, and which can be most naturally generalized beyond the contexts that we have specifically addressed in the paper (cf. section \ref{conclusions} below for a further discussion of this point).

\subsection{Terminology and notation}\label{terminology}

The terminology and notation used in this paper are standard and borrowed from \cite{El}, if not otherwise specified. Conventions that we will frequently employ include the following.

Given a topological space $X$, we denote by $X_{0}$ its underlying set and by ${\cal O}(X)$ its frame of open sets. For a locale $L$, we denote by ${\cal O}(L)$ its underlying frame and, for a morphism $g:L\to L'$ of locales, we denote by $g^{\ast}:{\cal O}(L')\to {\cal O}(L)$ the frame homomorphism corresponding to it. Occasionally, when there is no risk of confusion, we denote the frame underlying a locale $L$ by the same letter.

By a \emph{subframe} of a frame $A$ we mean a subset $B\subseteq A$ such that $B$ is a frame with respect to the order induced by that on $A$ (equivalently, $0_{A}, 1_{A}\in B$ and $B$ is closed under arbitrary joins and finitary meets in $A$). 

By an \emph{indexing function} of a set of points $\textsc{P}$ of a topos $\cal E$ by a set $X$ we mean a surjective function $i:X \to \textsc{P}$. If $\cal E$ has arbitrary set-indexed products, any such indexing $i$ naturally identifies with a geometric morphism $[X, \Set]\to {\cal E}$, denoted by $\tilde{i}$ and defined as follows: $\tilde{i}^{\ast}(A)(x)=i(x)^{\ast}(A)$ for any object $A$ of $\cal E$ and $\tilde{i}^{\ast}(A\to B)(x)=i(x)^{\ast}(A\to B)$ for any arrow $A\to B$ in $\cal E$. For any topos $\cal E$, we call a geometric morphism $\xi:[X, \Set]\to {\cal E}$ an \emph{indexing} of points of $\cal E$ by the set $X$; any such morphism induces an indexing function $i_{\xi}$ with domain $X$ of a set of points of $\cal E$, such that $\tilde{i_{\xi}}\cong \xi$; for any $x\in X$, we denote also by $\xi_{x}$ the point $i_{\xi}(x)$ of $\cal E$.  

We consider points of toposes up to isomorphisms, rather than strictly.

We will denote by \textbf{Loc} the category of locales and by \textbf{Top} the category of topological spaces.

Given a category $\cal C$, we denote by $ob({\cal C})$ the collection of the objects of $\cal C$; for any arrow $f$ in $\cal C$, we denote by $dom(f)$ its domain and by $codom(f)$ its codomain. We sometimes write $c\in {\cal C}$ to mean that $c$ is an object of $\cal C$. The operation of composition of arrows in a category will always be denoted by the symbol $\circ$.

Given a preorder $\cal P$, and an element $a\in {\cal P}$, we denote by $(a)\downarrow$ the set of elements $b\in {\cal P}$ such that $b\leq a$.

The union of a family of subobjects $\{a_{i} \mono a \textrm{ | } i\in I\}$ in a topos will be denoted by $\mathbin{\mathop{\textrm{\huge $\vee$}}\limits_{i\in I}}a_{i}\mono a$.

The powerset of a given set $A$ will be denoted by ${\mathscr{P}}(A)$. The collection of all the finite subsets of a set $A$ will be denoted by ${\mathscr{P}}_{fin}(A)$. 

By an atom in a category $\cal C$ with an initial object $0_{\cal C}$ we mean an object $c\ncong 0_{{\cal C}}$ of $\cal C$ such that for any monomorphism $b\mono c$, either $b\cong 0_{\cal C}$ or $b\mono c$ is an isomorphism; by an atomic subobject in $\cal C$ we mean a subobject whose domain is an atom in $\cal C$.

When we define `concrete' categories in this paper by only specifying their objects and arrows, we tacitly assume that the definition of the identities and composition rule is straightforward, i.e. that the identity arrow on any object is given by the identity function on its underlying set and the composition of arrows in the category is given by the set-theoretic composition of the underlying functions; that is, we consider all these `concrete' categories as categories structured over the category of sets.

\section{Topos-theoretic topologies}\label{subterminaltop}

In this section we introduce a method for building topological spaces starting from toposes equipped with a set of points. The resulting notions will be central for our purposes, in that, as we shall see in section \ref{dualtop}, they will allow us to naturally `lift` a given duality (resp. equivalence) between a category $\cal K$ of preorders and a category of locales to a duality (resp. equivalence) between $\cal K$ and a category of topological spaces.  

\subsection{Spaces of points of a topos}\label{subterminal}

First, let us recall the following standard definition of the \emph{space of points} of a locale (cf. for example p. 491 \cite{El}).

A point of a locale $X$ is defined to be a locale morphism $p: 1 \to X$, where $1$ denotes the locale corresponding to the one-point space; equivalently, it is a frame homomorphism $p^{\ast}: {\cal O}(X) \to {\cal O}(1)\cong \{0,1\}$. 

\begin{definition}
Let $X$ be a locale. The \emph{space of points} of $X$ is the set $X_{p}$ of all the points of $X$ equipped with the topology given by the image of the frame homomorphism $\phi_{X}: {\cal O}(X) \to {\mathscr{P}}(X_{p})$ defined by:
\[
\phi_{X}(U)=\{p\in X_{p} \textrm{ | } p^{\ast}(U)=1\},
\] 
for any $U\in {\cal O}(X)$.
\end{definition}

From now on we will consider $X_{p}$ as a space equipped with this topology. 

We note that we can interpret this definition topos-theoretically, as follows. The set $X_{p}$ of points of $X$ corresponds bijectively with the set $\textsc{P}$ of (isomorphism classes of) points of the topos $\Sh(X)$ (cf. Proposition C1.4.5 \cite{El}), and, for any point $p$ of $X$, the map $p^{\ast}:{\cal O}(X)\to \{0,1\}$ corresponds precisely to the action on subterminals of the inverse image $f_{p}^{\ast}:\Sh(X) \to \Set$ of the geometric morphism $f_{p}:\Set \to \Sh(X)$ corresponding to the point $p$. In these terms, the frame homomorphism $\phi_{X}$ acquires the following expression: 
\[
\phi_{X}(U)\cong\{f_{p} \in \textsc{P} \textrm{ | } f_{p}^{\ast}(U)\cong 1_{\Set}\}.
\]  

This remark naturally leads to the following more general definition.
  
\begin{definition}\label{funddef}
Let $\cal E$ be a locally small cocomplete topos, $\Gamma$ be a subframe of $\Sub_{\cal E}(1)$ and $i:X \to \textsc{P}$ be an indexing function of a set $\textsc{P}$ of points of $\cal E$ by a set $X$. The \emph{$\Gamma$-subterminal} topology ${\tau}^{\cal E}_{\Gamma, i}$ on the set $X$ is the image of the function $\phi_{\Gamma, {\cal E}}: \Gamma \to {\mathscr{P}}(X)$ given by 
\[
\phi_{\Gamma, {\cal E}}(u)=\{x\in X \textrm{ | } \xi(x)^{\ast}(u) \cong 1_{\Set}\}.
\]
In other words, the subsets in ${\tau}^{\cal E}_{\Gamma, i}$ are precisely those of the form $\phi_{\Gamma, {\cal E}}(u)$ where $u$ ranges among the subterminals in $\Gamma$.
\end{definition}

If $\Gamma=\Sub_{{\cal E}}(1)$, we simply call ${\tau}^{\cal E}_{\Gamma, i}$ the \emph{subterminal topology}, and denote it by ${\tau}^{\cal E}_{i}$. If $\xi:[X, \Set]\to {\cal E}$ is an indexing of points of $\cal E$ by the set $X$ and $i_{\xi}$ is the corresponding indexing function of a set of points of $\cal E$, we denote ${\tau}^{\cal E}_{\Gamma, i}$ (resp. ${\tau}^{\cal E}_{i}$) also by ${\tau}^{\cal E}_{\Gamma, \xi}$ (resp. ${\tau}^{\cal E}_{\xi}$).

Instances of this notion have occasionally been considered in the literature; for example, the construction for $\Gamma=\Sub_{\cal E}(1)$ is used in the proof of Theorem 7.25 \cite{topos} (in fact, this result is subsumed by our Theorem \ref{topol} below). 

Under the hypotheses of Defininition \ref{funddef}, we say that a collection $\textsc{P}$ of points of $\cal E$ \emph{separates the subterminals in $\Gamma$} if for any non-isomorphic subterminals $u, v$ in $\Gamma$ there exists a point $p \in \textsc{P}$ such that $p^{\ast}(u) \not \cong p^{\ast}(v)$. Note that if $\textsc{P}$ is a separating set of points for $\cal E$ (i.e. the inverse images of the points in $\textsc{P}$ jointly reflect isomorphisms) then in particular $\textsc{P}$ separates the subterminals in $\Gamma$, for any $\Gamma$.
 
\begin{theorem}\label{topol}
Let $\cal E$ be a locally small cocomplete topos, $\Gamma \subseteq \Sub_{{\cal E}}(1)$ be a subframe of $\Sub_{{\cal E}}(1)$ and $i:X\to \textsc{P}$ be an indexing function of a set of points of $\cal E$. Then
\begin{enumerate}[(i)]
\item The collection ${\tau}^{\cal E}_{\Gamma, i}$ of subsets of $X$ defines a topology on the set $X$;

\item The collection $\textsc{P}$ of points of $\cal E$ separates the subterminals in $\Gamma$ if and only if the frame of open sets of the topology ${\tau}^{\cal E}_{\Gamma, i}$ on $X$ is isomorphic to the frame $\Gamma$ via the map $\phi_{\Gamma, {\cal E}}$. In particular, if $\textsc{P}$ separates the subterminals of $\cal E$ then $\phi_{\Sub_{\cal E}(1), {\cal E}}$ yields an isomorphism between the frame of open sets of the topology ${\tau}^{\cal E}_{i}$ and $\Sub_{{\cal E}}(1)$. 
\end{enumerate}
\end{theorem}

\begin{proofs}
$(i)$ Clearly, it suffices to prove that the function $\phi_{\Gamma, {\cal E}}: \Gamma \to {\mathscr{P}}(X)$ of Definition \ref{funddef} is a frame homomorphism. 

The fact that for any $u,v\in \Gamma$, $\phi_{\Gamma, {\cal E}}(u\wedge v)=\phi_{\Gamma, {\cal E}}(u)\cap \phi_{\Gamma, {\cal E}}(v)$ (where $u\wedge v$ denotes the intersection of $u\mono 1$ and $v\mono 1$ in $\Sub_{\cal E}(1)$ or, equivalently, in $\Gamma$) follows immediately from the fact that each $\xi(x)^{\ast}$ preserves intersections of subobjects (being the inverse image functor of a geometric morphism).   

Similarly, the equality $\mathbin{\mathop{\textrm{\huge $\cup$}}\limits_{i\in I}}\phi_{\Gamma, {\cal E}}(u_{i})=\phi_{\Gamma, {\cal E}}({\mathbin{\mathop{\textrm{\huge $\vee$}}\limits_{i\in I}}u_{i}})$ follows from the fact that each $\xi(x)^{\ast}$ preserves unions of subobjects (being the inverse image functor of a geometric morphism). 

$(ii)$ The function $\phi_{\Gamma, {\cal E}}$ is, by $(i)$, a frame homomorphism and, by definition of ${\tau}^{\cal E}_{\Gamma, i}$, always surjective, so it is an isomorphism precisely when it is injective i.e. when $\textsc{P}$ separates the subterminals in $\Gamma$.
\end{proofs}

Given $\cal E$, $\Gamma$ and $\xi:X\to \textsc{P}$ as in the hypotheses of the theorem, we denote the set $X$ equipped with the topology ${\tau}^{\cal E}_{\Gamma, i}$ (resp. ${\tau}^{\cal E}_{i}$) by $X_{{\tau}^{\cal E}_{\Gamma, i}}$ (resp. $X_{{\tau}^{\cal E}_{i}}$). If $\xi$ is the identity function on a set of points $\textsc{P}$ of $\cal E$, we denote ${\tau}^{\cal E}_{i}$ (resp. $X_{{\tau}^{\cal E}_{i}}$) also by ${\tau}^{\cal E}_{\textsc{P}}$ (resp. $X_{{\tau}^{\cal E}_{\textsc{P}}}$); if moreover $\cal E$ has only a \emph{set} of points (up to isomorphism) and $\textsc{P}$ is the collection of all the points of $\cal E$ we denote ${\tau}^{\cal E}_{\textsc{P}}$ (resp. $X_{{\tau}^{\cal E}_{\textsc{P}}}$) simply by ${\tau}^{\cal E}$ (resp. $X_{{\tau}^{\cal E}}$), and we call $X_{{\tau}^{\cal E}}$ the \emph{space of points} of the topos $\cal E$.  

We note that part $(ii)$ of the theorem generalizes the following well-known result from Locale Theory: if a locale $L$ has enough points then the frame of open sets of the space of points of $L$ is isomorphic to ${\cal O}(L)$.

\begin{remarks}\label{topolinvariant}

\begin{enumerate}[(a)]

\item By Theorem \ref{topol}(ii), if $\textsc{P}$ separates the subterminals of $\cal E$ then $\Sh(X_{{\tau}^{\cal E}_{\textsl{P}}})$ is (equivalent to ) the localic part of the hyperconnected-localic factorization of the unique (up to isomorphism) geometric morphism ${\cal E}\to \Set$. 

\item If $\cal E$ is a localic topos with enough points and $\textsl{P}$ is the collection of all the points of $\cal E$ then the topological space $X_{{\tau}^{\cal E}_{\textsl{P}}}$ is sober (cf. Theorem \ref{sobriety} below in view of Remark \ref{topolinvariant}(a)).

\item  Every topological space is of the form $X_{{\tau}^{\cal E}_{\xi}}$ for some topos $\cal E$ and indexing function $i$ of a set of points of $\cal E$. Indeed, given a topological space $X$, $X$ is homeomorphic to the space $X_{{\tau}^{\cal E}_{i_{X}}}$, where $\cal E$ is the topos $\Sh(X)$ and $i_{X}:X_{0} \to \textsl{P}$ is the indexing function sending a point $x$ of $X_{0}$ to the geometric morphism $i_{X}(x):\Set \to \Sh(X)$ whose inverse image is the stalk functor at the point $x$.

\item The ($\Gamma$-)subterminal topology is a topos-theoretic invariant; that is, if $f:{\cal E}\to {\cal F}$ is an equivalence of toposes then for any indexing function $i:X \to \textsc{P}$ of a set $\textsc{P}$ of points of $\cal E$, denoted by $i':X\to \textsc{Q}$ the indexing function of points of $\cal F$ defined by setting $i'(x)$ equal to the point of $\cal F$ given by the composite of $i(x)$ with the equivalence $f$, the spaces $X_{{\tau}^{\cal E}_{i}}$ and $X_{{\tau}^{\cal F}_{i'}}$ are homeomorphic. In particular, if $\cal E$ and $\cal F$ are two equivalent toposes with the property of having only a \emph{set} of points (for example, if $\cal E$ and $\cal F$ are localic toposes) then the spaces $X_{{\tau}^{\cal E}}$ and $X_{{\tau}^{\cal F}}$ are homeomorphic.

\end{enumerate}
\end{remarks}

The construction of the $\Gamma$-subterminal topology can be naturally made functorial, as follows.

Let us define a category $\mathfrak{Top}_{t}$ whose objects are the triples $({\cal E}, \Gamma, \xi)$ where $\cal E$ is a locally small cocomplete topos, $\Gamma$ is a subframe of $\Sub_{\cal E}(1)$ and $\xi:[X, \Set] \to {\cal E}$ is an indexing of set of points of $\cal E$, and whose arrows $({\cal E}, \Gamma, \xi)\to ({\cal F}, \Gamma', \xi')$, where $\xi:[X, \Set] \to {\cal E}$ and $\xi':[Y, \Set] \to {\cal F}$ are indexings of points respectively of $\cal E$ and of $\cal F$, are the pairs $(f, l)$ where $f:{\cal E}\to {\cal F}$ is a geometric morphism such that $f^{\ast}$ sends the subterminals in $\Gamma'$ to subterminals in $\Gamma$ and $l:X\to Y$ is a function such that, denoted by $E_{l}:[X, \Set]\to [Y, \Set]$ the geometric morphism induced by $l$ as in Example A4.1.4 \cite{El}, the diagram

\[  
\xymatrix {
[X, \Set] \ar[r]^{E_{l}} \ar[d]^{\xi} & [Y, \Set] \ar[d]^{\xi'} \\
{\cal E} \ar[r]^{f} & {\cal F}}
\]\\
commutes (up to isomorphism). Identities and composition in $\mathfrak{Top}_{t}$ are defined componentwise in the obvious way.

Notice that, given an geometric morphism $f:{\cal E}\to {\cal F}$ such that $f^{\ast}$ sends the subterminals in $\Gamma'$ to subterminals in $\Gamma$, $f$ lifts to an arrow $({\cal E}, \Gamma, \xi)\to ({\cal F}, \Gamma' \xi')$ (i.e. there exists $l:X\to Y$ such that $(f, l)$ yields a morphism $({\cal E}, \Gamma, \xi)\to ({\cal F}, \Gamma' \xi')$ in $\mathfrak{Top}_{s}$), where $\xi:[X, \Set] \to {\cal E}$ and $\xi':[Y, \Set] \to {\cal F}$ are indexings of points respectively of $\cal E$ and of $\cal F$, if and only if there is a way of assigning to each $x\in X$ an element $y\in Y$ with the property that the composite $f\circ \xi_{x}:\Set \to {\cal F}$ is equal to $\xi_{y}$. Note that if $i_{\xi'}$ is bijective then there is at most one such `lifting' $(f, l)$; in particular, if $\cal F$ has only a set of points and $i_{\xi'}$ is the identity on this set of points then there is exactly one `lifting' $(f, l_{f})$ for any geometric morphism $f:{\cal E}\to {\cal F}$.  

Every arrow $(f, l):({\cal E}, \Gamma, \xi)\to ({\cal F}, \Gamma', \xi')$ in $\mathfrak{Top}_{t}$ gives rise to a continuous map of topological spaces $X_{{\tau}^{\cal E}_{\Gamma_{\cal E}, \xi}}\to X_{{\tau}^{\cal E}_{\Gamma_{\cal F}, \xi'}}$ with underlying function $l$. Indeed, by the commutativity of the diagram above, for any subterminal $v$ in $\Gamma_{\cal F}$, $l^{-1}(\phi_{\Gamma_{\cal F}, {\cal F}}(v))=\{x\in X \textrm{ | } l(x) \in \phi_{\Gamma_{\cal F}, {\cal F}}(v)\}=\{x\in X \textrm{ | } \xi(x)^{\ast}(f^{\ast}(v)) \cong 1_{\Set}\}=\phi_{\Gamma_{\cal E}, {\cal E}}({f^{\ast}(v)})$. We thus have a functor 
\[
\Theta_{t}: \mathfrak{Top}_{t} \to \textbf{Top}.
\]
which sends an object $({\cal E}, \Gamma, \xi)$ of $\mathfrak{Top}_{t}$ to the topological space $X_{{\tau}^{\cal E}_{\Gamma_{\cal E}, \xi}}$ and an arrow $(f, l):({\cal E}, \Gamma, \xi)\to ({\cal F}, \Gamma', \xi')$ in $\mathfrak{Top}_{t}$ to the continuous map $l:X_{{\tau}^{\cal E}_{\Gamma_{\cal E}, \xi}}\to X_{{\tau}^{\cal E}_{\Gamma_{\cal F}, \xi'}}$. We will discuss various properties of this functor in section \ref{adj}. 

We denote by $\mathfrak{Top}_{p}$ the full subcategory of $\mathfrak{Top}_{t}$ on the objects of the form $({\cal E}, \xi, \Sub_{\cal E}(1))$ , and call it the category of \emph{toposes paired with points}. Note that the objects of $\mathfrak{Top}_{p}$ can be simply identified with the pairs $({\cal E}, \xi)$, where $\cal E$ is a locally small cocomplete topos and $\xi:[X, \Set] \to {\cal E}$ is an indexing of a set of points of $\cal E$.     
We denote the restriction of the functor $\Theta_{t}$ to the category $\mathfrak{Top}_{p}$ simply by $\Theta_{p}:\mathfrak{Top}_{p} \to \textbf{Top}$. 

For any two subframes $\Gamma$ and $\Delta$ of $\Sub_{\cal E}(1)$, if $\Gamma$ is a subframe of $\Delta$ then we have a continuous surjection of topological spaces $X_{{\tau}^{\cal E}_{\Delta}, \xi}\to X_{{\tau}^{\cal E}_{\Gamma}, \xi}$; this surjection induces a geometric surjection of toposes $\Sh(X_{{\tau}^{\cal E}_{\Delta}, \xi}) \to \Sh(X_{{\tau}^{\cal E}_{\Gamma}, \xi})$ (cf. Example A4.2.7(c) \cite{El}).

The usefulness of Theorem \ref{topol} lies in the fact that it allows us to build topological spaces from toposes through a topos-theoretic invariant which has a natural behaviour with respect to sites; indeed, by Diaconescu's equivalence, the points of a topos $\Sh({\cal C}, J)$ correspond precisely to the flat $J$-continuous functors $\cal C \to \Set$ (in particular, if $\cal E$ is the classifying topos of a geometric theory $\mathbb T$ then the points of $\cal E$ correspond precisely to the models of $\mathbb T$ in $\Set$). Moreover, the formulation of the notion of ($\Gamma$-)subterminal topology as a topos-theoretic invariant paves the way, in light of the methodologies introduced in \cite{OC10}, for an effective transfer of properties between topological spaces constructed from two different sites of definition of a given topos. We will see concrete applications of this remark in sections \ref{insights} and \ref{zariski} below.

In passing, we note that the specialization order $\leq$ on $X_{{\tau}^{\cal E}_{\Gamma, \xi}}$ can be naturally characterized as a topos-theoretic invariant: $x\leq x'$ in $X_{{\tau}^{\cal E}_{\Gamma, \xi}}$ if and only if for every subterminal $u$ in $\Gamma$, $x\in \phi_{\Gamma, {\cal E}}(u)$ implies $x'\in \phi_{\Gamma, {\cal E}}(u)$.    

\subsection{Examples}\label{exsub}

Let us begin our list of examples of subterminal topologies by giving an explicit description of the spaces of points of toposes of the form $\Sh({\cal C}, J)$, where $\cal C$ is a preorder category and $J$ is a Grothendieck topology on it. To this end, we introduce the following notions.

\begin{definition}
Let $({\cal C}, J)$ be a site.

\begin{enumerate}[(a)]

\item A \emph{$J$-ideal} on $\cal C$ is a subset $I\subseteq ob({\cal C})$ such that for any arrow $f:b\to a$ in $\cal C$ if $a\in I$ then $b\in I$, and for any $J$-covering sieve $R$ on an object $c$ of $\cal C$, if $dom(f)\in I$ for every $f\in R$ then $c\in I$; we denote by $Id_{J}({\cal C})$ the set of $J$-ideals on $\cal C$, endowed with the subset-inclusion order relation. If $J$ is the trivial topology on $\cal C$, we call the $J$-ideals on $\cal C$ simply \emph{ideals}, and we denote $Id_{J}({\cal C})$ by $Id({\cal C})$. 

\item Given an object $c$ of $\cal C$, we the \emph{principal $J$-ideal} $(c)\downarrow_{J}$ generated by $c$ is the smallest $J$-ideal on $\cal C$ containing the object $c$, that is the collection of all the objects $d\in {\cal C}$ such that there exists a $J$-covering sieve $R$ on $d$ with the property that for every $f\in R$ there exists an arrow $dom(f)\to c$ in $\cal C$.

\item Let $({\cal C}, \leq)$ be a preorder category. A \emph{$J$-prime filter} on $\cal C$ is a subset $F\subseteq ob({\cal C})$ such that $F$ is non-empty, $a\in F$ implies $b\in F$ whenever $a\leq b$ in $\cal C$, for any $a, b\in F$ there exists $c\in F$ such that $c\leq a$ and $c\leq b$, and for any $J$-covering sieve $\{a_{i} \to a \textrm{ | } i\in I\}$ in $\cal C$ if $a\in F$ then there exists $i\in I$ such that $a_{i}\in F$. 

\end{enumerate}

\end{definition}

Notice that if $\cal C$ is cartesian (i.e. a meet-semilattice) then the second condition in the definition of $J$-prime filter can be equivalently replaced by the requirement that for any $a,b\in F$, $a\wedge b\in F$ (where $\wedge$ denotes the meet operation in $\cal C$), while the condition that $F$ should be non-empty can be replaced by the requirement that $1\in F$ (where $1$ is the top element of $\cal C$).

The notion of $J$-ideal defined above makes also sense for a (not necessarily Grothendieck) coverage $J$ on ${\cal C}$, in the sense of Definition C2.1.1 \cite{El}; in fact, we will use this more general notion in section \ref{concreteness} below.

\begin{remarks}\label{rem2}
\begin{enumerate}[(a)]

\item The $J$-ideals on $\cal C$ can be identified with the subterminals of the topos $\Sh({\cal C}, J)$ (cf. Example C1.1.16 \cite{El}); under this identification, the usual order in $\Sub_{\Sh({\cal C}, J)}(1_{\Sh({\cal C}, J)})$ corresponds to the subset-inclusion order on the set of $J$-ideals on $\cal C$, and the principal $J$-ideal $(c)\downarrow_{J}$ on an object $c$ of $\cal C$ corresponds to the $J$-closure of the subobject $S_{c}\mono 1$ in $[{\cal C}^{\textrm{op}}, \Set]$ given by the monic part of the cover-mono factorization in $[{\cal C}^{\textrm{op}}, \Set]$ of the unique arrow $Hom_{\cal C}(-,c)\to 1_{[{\cal C}^{\textrm{op}}, \Set]}$.

\item If $J$ is subcanonical then $(c)\downarrow_{J}$ is equal to the collection of the objects $d$ of $\cal C$ such that there exists an arrow $d\to c$ in $\cal C$. Indeed, since $1_{[{\cal C}^{\textrm{op}}, \Set]}$ is always a sheaf, the action $a_{J}(S_{c})$ of the associated sheaf functor $a_{J}:[{\cal C}^{\textrm{op}}, \Set]\to \Sh({\cal C}, J)$ on the subobject $S_{c}\mono 1$ coincides with its $J$-closure. In particular, if $({\cal C}, \leq)$ is a preorder then $(c)\downarrow_{J}$ is equal to the set $(c)\downarrow$ of elements $d\in {\cal C}$ such that $d\leq c$.        

\end{enumerate}
\end{remarks}

We can now state the following proposition.

\begin{proposition}\label{mslattice}
Let $\cal C$ be a preorder and $J$ be a Grothendieck topology on it. Then the topological space $X_{{\tau}^{\Sh({\cal C}, J)}}$ is homeomorphic to the space which has as set of points the collection ${\cal F}^{J}_{\cal C}$ of the $J$-prime filters on $\cal C$ and as open sets the sets the form
\[
{\cal F}_{I}=\{F\in {\cal F}^{J}_{\cal C} \textrm{ | } F\cap I\neq \emptyset\},
\] 
where $I$ ranges among the $J$-ideals on $\cal C$. In particular, a sub-basis for this topology is given by the sets
\[
{\cal F}_{c}=\{F\in {\cal F}^{J}_{\cal C} \textrm{ | } c\in F\},
\] 
where $c$ varies among the elements of $\cal C$. 
\end{proposition}

\begin{proofs}
As we observed in Remark \ref{rem}(a), the subterminals in $\Sh({\cal C}, J)$ can be identified with the $J$-ideals on $\cal C$.

The points of $X_{{\tau}^{\Sh({\cal C}, J)}}$ are, by definition of $X_{{\tau}^{\Sh({\cal C}, J)}}$, the isomorphism classes of geometric morphisms $\Set \to \Sh({\cal C}, J)$. These correspond, by Diaconescu's equivalence, to the $J$-continuous flat functors $\cal C \to \Set$, and these in turn correspond exactly to the $J$-prime filters on $\cal C$ (cf. section \ref{logical} below).

Now, if $F$ is the $J$-prime filter corresponding to a point $p:\Set \to \Sh({\cal C}, J)$ then, since every $J$-ideal is the union of the principal $J$-ideals generated by the elements belonging to it, $p^{\ast}(I)=\mathbin{\mathop{\textrm{\huge $\cup$}}\limits_{c\in I}}p^{\ast}((c)\downarrow_{J})$ is isomorphic to $1_{\Set}$ if and only if there exists $c\in I$ such that $p^{\ast}((c)\downarrow_{J})\cong 1_{\Set}$ (equivalently, $c\in F$), from which our thesis follows. 
\end{proofs}

Now that we have an explicit description, provided by Proposition \ref{mslattice}, of the spaces of points of toposes of the form $\Sh({\cal C}, J)$, where $\cal C$ is a preorder category, it is natural to wonder if we can also explicitly describe in these terms the continuous map $X_{{\tau}^{f}}:X_{{\tau}^{\Sh({\cal D}, K)}}\to X_{{\tau}^{\Sh({\cal C}, J)}}$ resulting from applying the functor $\Theta_{p}: \mathfrak{Top}_{p} \to \textbf{Top}$ to the morphism $(\dot{f}, l_{f}):(\Sh({\cal D}, K), P_{\Sh({\cal D}, K)}, \to (\Sh({\cal C}, J)), P_{\Sh({\cal D}, K)})$, where $\dot{f}:\Sh({\cal D}, K) \to \Sh({\cal C}, J)$ is the geometric morphism induced by a morphism of sites $f:({\cal C}, J)\to ({\cal D}, K)$, $P_{\Sh({\cal D}, K)}$ (resp. $P_{\Sh({\cal C}, J)}$) is the indexing induced by the identity indexing function on the set $p_{\Sh({\cal D}, K)}$ (resp. $p_{\Sh({\cal C}, J)}$) of all the points of $\Sh({\cal D}, K)$ (resp. $\Sh({\cal C}, J)$) and $l_{f}:p_{\Sh({\cal D}, K)}\to p_{\Sh({\cal C}, J)}$ is the function induced by composition with $f$. The following proposition gives a positive answer to this question.

\begin{proposition}\label{mslatticefun}
With the notation above, if we identify $X_{{\tau}^{\Sh({\cal D}, K)}}$ (resp. $X_{{\tau}^{\Sh({\cal C}, J)}}$) with the set ${\cal F}^{K}_{\cal D}$ (resp. ${\cal F}^{J}_{\cal C}$) of $K$-prime filters on $\cal D$ (resp. of $J$-prime filters on $\cal C$) endowed with the subterminal topology, as in Proposition \ref{mslattice}, the continuous map $X_{{\tau}^{f}}:X_{{\tau}^{\Sh({\cal D}, K)}}\to X_{{\tau}^{\Sh({\cal C}, J)}}$ admits the following description: for any filter $F\in X_{{\tau}^{\Sh({\cal D}, K)}}$, 
\[
X_{{\tau}^{f}}(F)=f^{-1}(F).
\] 
\end{proposition}

\begin{proofs}
This immediately follows from the fact that if a point $p:\Set \to \Sh({\cal D}, K)$ of the topos $\Sh({\cal D}, K)$ corresponds to a $K$-prime filter $F$ on $\cal D$ (as in Proposition \ref{mslattice}) then the point $\dot{f} \circ p$ of the topos $\Sh({\cal C}, J)$ corresponds to the $J$-prime filter $f^{-1}(F)$ on $\cal C$. This fact is in turn easily verified by using the explicit description of Diaconescu's equivalence and observing that the inverse image functor $\dot{f}^{\ast}$ satisfies the property that for any $c\in {\cal C}$, $\dot{f}^{\ast}(y(c))\cong y'(f(c))$, where $y:{\cal C} \to \Sh({\cal C}, J)$ and $y':{\cal D}\to \Sh({\cal D}, K)$ are the composites of the associated sheaf functors $[{\cal C}^{\textrm{op}}, \Set]\to \Sh({\cal C}, J)$ and $[{\cal D}^{\textrm{op}}, \Set]\to \Sh({\cal D}, K)$ respectively with the Yoneda embeddings ${\cal C}\to [{\cal C}^{\textrm{op}}, \Set]$ and ${\cal D}\to [{\cal D}^{\textrm{op}}, \Set]$.   
\end{proofs}

\begin{examples}\label{exa}
\begin{enumerate}[(a)]

Let us now point out several interesting topological spaces which naturally arise from putting the subterminal topology on the domain of some indexing function of set of points of a topos.

\item \emph{The trivial topology}

Given a topos $\cal E$ and an indexing $\xi$ of a set of points of $\cal E$ by a set $X$, if $\Gamma$ is equal to the subframe $\{0_{\cal E},1_{\cal E}\} \subseteq \Sub_{\cal E}(1)$ then the topological space $X_{{\tau}^{\cal E}_{\Gamma, \xi}}$ identifies with the trivial topological space with underlying set $X$.

\item \emph{The Alexandrov topology}

Let $({\cal P}, \leq)$ be a preorder, $\cal E$ be the topos $[{\cal P}, \Set]$, and let $i:{\cal P}\to \textsl{S}$ be the indexing function sending an element $p\in {\cal P}$ to the geometric morphism $e_{p}:\Set\to [{\cal P}, \Set]$ whose inverse image $e_{p}^{\ast}:[{\cal P}, \Set]\to \Set$ is the evaluation functor at the object $p$. The subterminals in $[{\cal P}, \Set]$ can be identified with the subsets $T$ of $\cal P$ such that for any $a\leq b$ in $\cal P$, $a\in T$ implies $b\in T$; clearly, the open set of the topological space ${\cal P}_{{{\tau}^{\cal E}}, i}$ corresponding to a subterminal $T$ coincides with $T$ itself (regarded as a subset of $\cal P$). So, the topological space ${\cal P}_{{{\tau}^{\cal E}}, i}$ is precisely the Alexandrov space associated to the preorder $\cal P$ (i.e. the topological space whose underlying set is $\cal P$ and whose open sets are the upper sets in $\cal P$).    

\item \emph{The Stone topology for distributive lattices}

Let $\cal D$ be a distributive lattice (regarded as a preorder coherent category) and $\cal E$ be the topos $\Sh({\cal D}, J_{coh})$, where $J_{coh}$ is the coherent topology on $\cal D$ (recall that, for any element $d\in {\cal D}$, the $J_{coh}$-covering sieves on $d$ are precisely the sieves on $d$ which contain finite families $\{d_{i}\leq d \textrm{ | } i\in I\}$ with the property that $d=\mathbin{\mathop{\textrm{\huge $\vee$}}\limits_{i\in I}}d_{i}$).

The subterminals in $\cal E$ can be identified with the $J_{coh}$-ideals in $D$ i.e. with the usual ideals of the distributive lattice $\cal D$. The points of $\Sh({\cal D}, J_{coh})$ can be identified with the prime filters on $\cal D$. Using Proposition \ref{mslattice}, we obtain that the subterminal topology ${\tau}^{\cal E}_{X}$ on the collection $X$ of prime filters on $\cal D$ has as basis of open sets the collection of sets of the form $\{F\in X \textrm{ | } d\in F\}$, where $d$ varies among the elements of $\cal D$. We have thus recovered the classical Stone topology on the set of prime filters on $\cal D$; in particular, if $\cal D$ is a Boolean algebra then the prime filters on $\cal D$ coincide with the ultrafilters on $\cal D$ and hence we recover the Stone topology on the set of ultrafilters of the Boolean algebra $\cal D$.   

\item \emph{A topology for meet-semilattices}

Let $\cal M$ be a meet-semilattice (regarded as a preorder cartesian category) and let $\cal E$ be the topos $[{\cal M}^{op}, \Set]$. The points of the topos $\cal E$ can be identified with the filters on $\cal M$. By Proposition \ref{mslattice}, the topology ${\tau}^{\cal E}_{X}$ on the set $X$ of filters on $\cal M$ has as basis of open sets the collection of sets of the form $\{F\in X \textrm{ | } d\in F\}$ where $d$ varies among the elements of $\cal M$.      

\item \emph{The space of points of a locale}

Let $L$ be a locale and $\cal E$ be the localic topos $\Sh(L)$; the points of the topos $\Sh(L)$ correspond bijectively with the points of the locale $L$ and hence there is only a \emph{set} $\textsl{P}$ of such points. Now, the points of a locale $L$ can be identified with the completely prime filters on the frame ${\cal O}(L)$ corresponding to $L$, while the subterminals in $\Sh(L)$ can be identified with the elements $u$ of $L$. Proposition \ref{mslattice} thus yields that the open sets of the subterminal topology ${\tau}^{\cal E}_{X}$ on the set $X$ of completely prime filters on $L$ are precisely the sets of the form $F_{u}=\{F\in X \textrm{ | } u\in F\}$ where $u$ ranges among the elements of $L$. We have thus recovered the usual topology on the space of points of the locale (as described for example in section C1.2 of \cite{El}). 

\item \emph{A logical topology} 

Let $\cal E$ be the classifying topos of a geometric theory $\mathbb T$ over a signature $\Sigma$. The points of $\cal E$ can be identified with the models of $\mathbb T$ in $\Set$. Let $X$ be a collection $X$ of such models. The subterminals in $\cal E$ can be identified with the $\mathbb T$-provable equivalence classes of geometric sentences over $\Sigma$, and the open set $F_{\phi}$ corresponding in the subterminal topology ${\tau}^{\cal E}_{X}$ to such a formula $\phi$ (regarded as a subterminal in $\cal E$) is the collection of all the models in $X$ which satisfy $\phi$. Thus $X_{{\tau}^{\cal E}_{\textsl{P}}}$ yields in this case a `logical topology' on the collection $X$ of models of $\mathbb T$; if every model of $\mathbb T$ in $\Set$ occurs as an element of $X$, we call the resulting topological space the \emph{logical space} associated to the theory $\mathbb T$. Notice that, since every geometric sentence over $\Sigma$ is $\mathbb T$-provably equivalent to a disjunction of coherent sentences, the collection of sets of the form $F_{\phi}$ for $\phi$ \emph{coherent} over $\Sigma$, forms a basis for the topology $X_{{\tau}^{\cal E}}$.    

\item \emph{The Zariski topology}

Let $A$ be a commutative ring with unit, and let $L(A)$ be the distributive lattice generated by symbols $D(a)$, $a \in A$, subject to the relations $D(1_{A})=1_{L(A)}$, $D(a \cdot b) = D(a) \wedge D(b)$, $D(0_{A})=0_{L(A)}$, and $D(a+ b) \leq D(a) \vee D(b)$. If we equip $L(A)$ with the coherent topology then the space $X_{{\tau}^{\Sh(L(A), J)}}$ is homeomorphic to the topological space obtained by equipping the prime spectrum $Spec(A)$ of $A$ with the Zariski topology (cf. section \ref{zariski} below).

\end{enumerate}
\end{examples}

\subsection{Sobriety}\label{Sobriety}

In this section, we present a natural application of the invariant concept of point of a topos to the investigation of the property of sobriety of a topological space; specifically, Theorem \ref{sobriety} below gives a criterion for a topological space built through the subterminal topology to be sober.

\begin{theorem}\label{sobriety}
Let $\cal E$ be a locally small cocomplete topos and let $\textsl{P}$ be a set of points of $\cal E$ which separates the subterminals of $\cal E$. Let $h:{\cal E}\to {\cal F}$ be the hyperconnected part of the hyperconnected-localic factorization of the unique geometric morphism ${\cal E}\to \Set$. Then the topological space $X_{{\tau}^{\cal E}_{\textsl{P}}}$ is sober if and only if the function from $\textsl{P}$ to the collection of points of $\cal F$ which sends a point $s$ in $\textsl{P}$ to the point of $\cal F$ given by the composite geometric morphism $h\circ s$ is a bijection. In particular, if $\cal E$ is localic then the topological space $X_{{\tau}^{\cal E}_{\textsl{P}}}$ is sober if and only if $\textsl{P}$ is the collection of all the points of $\cal E$.
\end{theorem}

\begin{proofs}
Let $X$ be a topological space and let ${\cal F}_{X}$ be the collection of all the completely prime filters on the frame ${\cal O}(X)$ of open sets of $X$. It is well-known (cf. p. 491 \cite{El}) that a topological space $X$ is sober if and only if the map $\eta_{X}:X_{0}\to {\cal F}_{X}$ which sends a point $x\in X_{0}$ to the filter $\eta_{X}(x)$ consisting of all the open sets which contain $x$ is a bijection. Now, the filters in ${\cal F}_{X}$ can be identified with the geometric morphisms $\Set \to \Sh(X)$ and hence we can reformulate the property of sobriety of $X$ by saying that the map $\xi_{X}$ from $X_{0}$ to the collection of all the points of the topos $\Sh(X)$ sending a point $x\in X_{0}$ to the geometric morphism $f_{x}:\Set \to \Sh(X)$ whose inverse image is the stalk functor at $x$ is a bijection. 

Now, since $\textsl{P}$ separates the subterminals in $\cal E$, $\cal F$ is equivalent to the topos $\Sh(X_{{\tau}^{\cal E}_{\textsl{P}}})$ (by Remark \ref{topolinvariant}(a)). It is easy to verify that, under the equivalence ${\cal F}\simeq \Sh(X_{{\tau}^{\cal E}_{\textsl{P}}})$ the map $\xi_{X_{{\tau}^{\cal E}_{\textsl{P}}}}$ defined above corresponds precisely to the map sending a point $s$ in $\textsl{P}$ (regarded as an element of the underlying set of $X_{{\tau}^{\cal E}_{\textsl{P}}}$) to the point of $\cal F$ given by the composite geometric morphism $h\circ s$. From this our thesis follows immediately. 
\end{proofs}
    
It follows at once from the theorem that the Stone spaces associated to distributive lattices, the topological spaces associated to meet-semilattices as in Example \ref{exa}(d) and the logical spaces of Example \ref{exa}(f), are all sober. The following criterion for an Alexandrov space to be sober also follows immediately from our theorem (the equivalence $(i)\biimp (iii)$ in Corollary \ref{alexsober} below is well-known). 
 
\begin{corollary}\label{alexsober}
Let $X$ be an Alexandrov space, and let $X_{\leq}$ be the preorder obtained by equipping the underlying set of $X$ with the specialization order. Then the following conditions are equivalent: 
\begin{enumerate}[(i)]
\item $X$ is sober;

\item Every flat functor $X_{\leq}^{\textrm{op}}\to \Set$ is representable;

\item Every non-empty directed ideal of $X_{\leq}$ is principal.

\end{enumerate}

\end{corollary}

\begin{proofs}
$(i)\biimp (ii)$ If $X$ is an Alexandrov space then $X$ is homeomorphic to the space ${(X_{0})}_{{{\tau}^{\cal E}}, i}$, where $i:X_{0}\to \textsl{S}$ is the indexing function of points of the topos $[X_{\leq}, \Set]$ considered in Example \ref{exa}(b). From Theorem \ref{sobriety} we thus know that $X$ is sober if and only if every point of the topos $[X_{\leq}, \Set]$ is one of the points in $\textsl{S}$. Now, recalling the equivalence between geometric morphisms $\Set \to [X_{\leq}, \Set]$ and flat functors $X_{\leq}^{\textrm{op}}\to \Set$, under which the point of $[X_{\leq}, \Set]$ whose inverse image functor is the evaluation functor at the object $x$ of $X_{\leq}$ corresponds to the flat functor $X_{\leq}^{\textrm{op}}\to \Set$ represented by $x$, we can alternatively reformulate this condition as the requirement that every flat functor $X_{\leq}^{\textrm{op}}\to \Set$ should be representable.  

$(ii)\biimp (iii)$ Every flat functor $F:{\cal P} \to \Set$ from a preorder $\cal P$ to $\Set$ takes values in $\{0, 1\}\hookrightarrow \Set$ and, under the assignment $F \to F^{-1}(1)$, the flat functors $X_{\leq}^{\textrm{op}}\to \Set$ correspond bijectively with the non-empty directed ideals on $X_{\leq}$ (cf. section \ref{logical} below); in these terms, the condition for a flat functor $X_{\leq}^{\textrm{op}}\to \Set$ to be representable amounts precisely to the requirement that the corresponding ideal should be principal, which proves our thesis.      
   
\end{proofs}

\section{The general approach to Stone-type dualities}\label{general}

In this section we present our general topos-theoretic framework for interpreting `Stone-type dualities'. For a classical treatment of these dualities, done from a locale-theoretic and categorical perspective, we refer the reader to the excellent book \cite{stone} by Johnstone, which in fact provided significant inspiration for the present work.

It will be clear from our analysis that the known Stone-type dualities are just a few of a large class of dualities that can be established through our topos-theoretic machinery. 

The main ingredient of our interpretation is the following well-known result from Topos Theory (cf. Example A4.6.2(e) \cite{El}): any topos $\Sh({\cal C}, J)$ of sheaves on a site $({\cal C}, J)$ whose underlying category $\cal C$ is a preorder is localic. 

Now, for any locally small cocomplete topos $\cal E$, we can consider the hyperconnected-localic factorization ${\cal E} \to {\cal F} \to \Set$ of the unique geometric morphism ${\cal E} \to \Set$. Notice that, since $\cal F$ is localic, the Comparison Lemma yields an equivalence ${\cal F}\simeq \Sh(\Sub_{{\cal F}}(1_{\cal F}))$ and hence, by Proposition A4.6.6 \cite{El}, we have an equivalence ${\cal F}\simeq \Sub_{{\cal E}}(1_{\cal E})$. In particular, if ${\cal E}$ is localic then $\cal E$ is equivalent to the localic topos $\Sh(L)$ where $L$ is the locale $\Sub_{\cal E}(1)$ of subterminals in $\cal E$. If we apply this to the topos $\Sh({\cal C}, J)$ of sheaves on a site $({\cal C}, J)$ whose underlying category $\cal C$ is a preorder we thus obtain, in view of Remark \ref{rem}(a), an equivalence $\Sh({\cal C}, J) \simeq \Sh(Id_{J}({\cal C}))$ where $Id_{J}({\cal C})$ is the \emph{locale of $J$-ideals on $\cal C$} i.e. the locale whose corresponding frame consists of the set of $J$-ideals on $\cal C$ equipped with the subset-inclusion order relation. That is, we have the following result.

\begin{theorem}\label{fund}
Let $\cal C$ be a preorder and $J$ be a Grothendieck topology on $\cal C$. Then the toposes $\Sh({\cal C}, J)$ and $\Sh(Id_{J}({\cal C}))$ are equivalent. 
\end{theorem}\qed  

As we shall argue below, this representation theorem, which can be read logically as a Morita-equivalence between two distinct geometric theories (cf. section \ref{logical} below), provides a general framework for analyzing the known Stone-type dualities and extracting new information about them, as well as for generating new dualities. In fact, we shall identify a set of general principles for building dualities starting from Theorem \ref{fund}, and illustrate them in action in several examples, including the classical ones. 

As the reader will have the opportunity to notice, this approach represents a clear implementation of the philosophy `toposes as bridges' introduced in \cite{OC10}. Indeed, the dualities or equivalences arise precisely from the process of `functorializing' a bunch of Morita-equivalences given by Theorem \ref{fund}; for each structure we have a Morita-equivalence, and the relationship between the preordered structure and the corresponding locale or topological space is determined by the expression of a topos-theoretic invariant in terms of the two different sites of definition of the topos. 

We start by making a distinction between

\begin{enumerate}
\item dualities between categories of preorders and categories of \emph{locales} and

\item dualities between categories of preorders and categories of \emph{topological spaces}. 

\end{enumerate}

The first kind of dualities have an essentially constructive nature, while the second class of dualities, which notably includes the classical Stone dualities for distributive lattices and Boolean algebras, may require some form of the axiom of choice. Anyway, the two classes of dualities are strongly interconnected, in that, as we shall see below, it is often possible to extract from a duality of the second kind a duality of the first kind and viceversa.

Broadly speaking, our method for building dualities of the first kind consists in equipping each of the structures $\cal C$ in a given category of preorders with an appropriate Grothendieck topology $J_{\cal C}$ in a such a way that this assignment is `natural' in $\cal C$. Such a choice induces a functor from the category $\cal K$ or its opposite (according to whether the duality is covariant or contravariant) to the category of locales sending each structure $\cal C$ to the locale $Id_{J_{\cal C}}({\cal C})$ of $J_{\cal C}$-ideals in $\cal C$. The operation of `recovering' a structure $\cal C$ from the corresponding topos $\Sh({\cal C}, J_{\cal C})\simeq \Sh(Id_{J_{\cal C}}({\cal C}))$ through a topos-theoretic invariant gives rise to a functor going in the converse direction which yields, together with it, the desired duality or equivalence.  

Once a duality of the first kind is established, we use the notion of subterminal topology introduced in section \ref{subterminaltop} to `enrich' the given duality to a duality with a category of topological spaces. This can be done in various ways, and this process of `enrichment' might require, depending on the case, some form of the axiom of choice.

We thus begin by focusing on the first kind of dualities. 

\subsection{Dualities with categories of locales}\label{locales}

Let us fix a category $\cal K$ of preordered structures. Our aim is to equip each structure $\cal C$ in $\cal K$ (regarded here as a preorder category) with a Grothendieck topology $J_{\cal C}$ on $\cal C$ in such a way that the assignment ${\cal C}\to Id_{J_{\cal C}}({\cal C})$ can be made into a functor (either covariant or contravariant) from $\cal K$ to the category $\textbf{Loc}$ of locales. 

The following notions will be central for our purposes.

Recall from section C2.3 of \cite{El} that a \emph{morphism of sites} $({\cal C}, J_{{\cal C}})\to ({\cal D}, J_{{\cal D}})$ is a flat functor ${\cal C}\to {\cal D}$ (i.e., a functor $F:{\cal C}\to {\cal D}$ such that for any object $d$ of $\cal D$ the category $(d\downarrow F)$ is cofiltered) which is cover-preserving, i.e. which sends every $J_{{\cal C}}$-covering sieve to a family which generates a $J_{{\cal D}}$-covering sieve. In particular, if $\cal C$ and $\cal D$ are meet-semilattices (regarded as cartesian categories) a morphism of sites $({\cal C}, J_{{\cal C}})\to ({\cal D}, J_{{\cal D}})$ is a cover-preserving meet-semilattice homomorphism ${\cal C}\to {\cal D}$. 

For preorder categories $\cal C$ and $\cal D$, the following `concrete' characterization of flat functors ${\cal C}\to {\cal D}$ holds.

\begin{proposition}
Let $({\cal C}, \leq)$ and $({\cal D}, \leq')$ be preorder categories and let $F:{\cal C}\to {\cal D}$ be a functor. Then $F$ is flat if and only if both of the following conditions hold:
\begin{enumerate}[(i)]
\item For any $d\in {\cal D}$ there exists $c\in {\cal C}$ such that $d\leq' F(c)$;

\item For any object $d\in {\cal D}$ and any objects $c,c'\in {\cal C}$ such that $d\leq' F(c)$ and $d\leq' F(c')$ there exists $c''\in {\cal C}$ such that $c''\leq c$, $c''\leq c'$ and $d\leq' F(c'')$.
\end{enumerate}

\end{proposition}

\begin{proofs}
The proposition follows immediately from the definition of cofiltered category.
\end{proofs}

Coming back to our original problem of making the assignment ${\cal C}\to Id_{J_{\cal C}}({\cal C})$ into a functor (either covariant or contravariant) from $\cal K$ to the category $\textbf{Loc}$ of locales, we have to distinguish between the covariant and contravariant case.   

\begin{enumerate}[(i)]
\item If we want to obtain a contravariant functor from $\cal K$ to $\textbf{Loc}$, we equip each structure $\cal C$ in $\cal K$ with a Grothendieck topology $J_{{\cal C}}$ on $\cal C$ in such a way that every arrow $f:{\cal C}\to {\cal D}$ in $\cal K$ gives rise to a morphism of sites $\hat{f}:({\cal C}, J_{{\cal C}})\to ({\cal D}, J_{{\cal D}})$;

\item If we want to obtain a covariant functor ${\cal K} \to \textbf{Loc}$, we do not equip each of the structures in $\cal K$ with any Grothendieck topology.
\end{enumerate}

In case $(i)$, we obtain a functor $A:{\cal K}^{\textrm{op}}\to \textbf{Loc}$ from the opposite of the category $\cal K$ to the category $\textbf{Loc}$ of locales, while in case $(ii)$ we obtain a functor $B:{\cal K}\to \textbf{Loc}$. 

In case $(i)$, the functor $A:{\cal K}^{\textrm{op}}\to \textbf{Loc}$ is defined as follows. Given a structure $\cal C$ in $\cal K$, we put $A({\cal C})=Id_{J_{\cal C}}({\cal C})$. Given an arrow $f:{\cal C}\to {\cal D}$ in ${\cal K}$, the morphism of sites $\hat{f}:({\cal C}, J_{{\cal C}})\to ({\cal D}, J_{{\cal D}})$ corresponding to $f$ gives rise, functorially, to a geometric morphism $\dot{f}:\Sh({\cal D}, J_{{\cal D}}) \to \Sh({\cal C}, J_{{\cal C}})$ (cf. Corollary C2.3.4 \cite{El}), which corresponds, via the equivalences of Theorem \ref{fund}, to a geometric morphism $\Sh(Id_{J_{\cal D}}({\cal D})) \to \Sh(Id_{J_{\cal C}}({\cal C}))$, which in turn corresponds, by Proposition C1.4.5 \cite{El}, to a unique morphism of locales $Id_{J_{\cal D}}({\cal D}) \to Id_{J_{\cal C}}({\cal C})$; we set $A(f):Id_{J_{\cal D}}({\cal D}) \to Id_{J_{\cal C}}({\cal C})$ equal to this morphism. It is easy to see that, concretely, $A(f)$ acts, at the level of frames, as the homomorphism sending a $J_{{\cal C}}$-ideal $I$ on $\cal C$ to the smallest $J_{\cal D}$-ideal on $\cal D$ containing the image of $I$ under $f$. 

In case $(ii)$, the functor $B:{\cal K}\to \textbf{Loc}$ is defined as follows. Given a structure $\cal C$ in $\cal K$, we put $B({\cal C})=Id({{\cal C}^{\textrm{op}}})$. Any arrow $f:{\cal C}\to {\cal D}$ in ${\cal K}$ gives rise, functorially, to a geometric morphism $[{\cal C}, \Set] \to [{\cal D}, \Set]$ (as in Example A4.1.4 \cite{El}), which corresponds via the equivalences of Theorem \ref{fund} to a geometric morphism $\Sh(Id({{\cal C}^{\textrm{op}}})) \to \Sh(Id({{\cal D}^{\textrm{op}}}))$, which in turn corresponds (by Proposition C1.4.5 \cite{El}) to a unique morphism of locales $Id({{\cal C}^{\textrm{op}}}) \to Id({{\cal D}^{\textrm{op}}})$; we set $B(f):Id({\cal C}^{\textrm{op}}) \to Id({{\cal D}^{\textrm{op}}})$ equal to this morphism. Concretely, $B(f)$ acts, at the level of frames, as the homomorphism sending an ideal $I$ on $\cal D$ to the inverse image $f^{-1}(I)$ of $I$ under $f$.

We note that, if we regard both ${\cal K}^{\textrm{op}}$ (resp. $\cal K$) and $\textbf{Loc}$ as preordered $2$-categories in the natural way (i.e., given two arrows $f,g:{\cal C}\to {\cal D}$ in $\cal K$ we set $f\leq g$ in $\cal K$ (or in ${\cal K}^{\textrm{op}}$), if for every $c\in {\cal C}$ $f(c)\leq g(c)$, and similarly for arrows in $\textbf{Loc}$, where the inequalities are considered in the dual category $\textbf{Frm}$), then the functor $A$ (resp. the functor $B$) becomes a $2$-functor which is covariant on $2$-cells (cf. Remark C2.3.5 and Proposition C1.4.5 \cite{El}). Concretely, this amounts precisely to saying that if $f\leq g$ in $\cal K$ then $A(f)\leq A(g)$ (resp. $B(f)\leq B(g)$) in $\textbf{Loc}$. 

So far, we have described a general method for constructing a (either covariant or contravariant) functor from a given category $\cal K$ of preorders to the category $\textbf{Loc}$ of locales; to build, starting from such a functor, an equivalence or duality between $\cal K$ and a subcategory of $\textbf{Loc}$, we have to care about how to go in the other direction. The general strategy is the following (we describe it for the case $(i)$ but our arguments can be trivially adapted to work in the case $(ii)$, cf. section \ref{prealex} below): if we are able to recover a structure ${\cal C}$ in $\cal K$ (uniquely up to canonical isomorphism) from the topos $\Sh({\cal C}, J_{{\cal C}})$ (equivalently, from the locale $A({\cal C})=Id_{J_{\cal C}}({\cal C})$) by means of a topos-theoretic invariant (in the sense of \cite{OC10}) functorially in ${\cal C} \in {\cal K}^{\textrm{op}}$ then we can expect to be able to use the invariant to define a functor $I_{A}:{\cal U}\to {\cal K}^{\textrm{op}}$ on a subcategory $\cal U$ of $\textbf{Loc}$ (namely, the extended image of $A$, cf. Definition \ref{extendedimage} below) which, together with $A$, yields an equivalence of categories ${\cal K}^{\textrm{op}} \simeq {\cal U}$. 

We note that in order to be able to recover a structure $\cal C$ from the topos $\Sh({\cal C}, J_{{\cal C}})$, the topology $J_{\cal C}$ must be `small enough' so that the associated sheaf functor $a_{J_{\cal C}}:[{\cal C}^{\textrm{op}}, \Set]\to \Sh({\cal C}, J_{{\cal C}})$ does not send two distinct representable functors to isomorphic objects. Clearly, if $J_{\cal C}$ is subcanonical then this problem does not subsist (in fact, as we shall see in sections \ref{ex} and \ref{Mordual} below, the classical examples of dualities, as well as our new examples, all arise when the Grothendieck topologies are subcanonical) If $\cal C$ is moreover a poset then the objects of $\cal C$ correspond bijectively with the principal $J_{{\cal C}}$-ideals on $\cal C$ (since, $J_{{\cal C}}$ being subcanonical, the latter are all of the form $(c)\downarrow$ for $c\in {\cal C}$), so the problem reduces to that of characterizing the principal ($J_{{\cal C}}$-)ideals on $\cal C$ among the $J$-ideals on $\cal C$ (i.e. the subterminals of the topos $\Sh({\cal C}, J_{{\cal C}})$) by means of a topos-theoretic invariant. We will discuss this problem in full generality in section \ref{charinv}, and will show that it has a positive solution in many cases of interest.

Recall that a Grothendieck topology $J$ on a small category $\cal C$ is said to be subcanonical if all the representable functors ${\cal C}^{\textrm{op}}\to \Set$ are $J$-sheaves. If $({\cal C}, \leq)$ is a preorder category, the condition for a Grothendieck topology $J$ to be subcanonical admits the following more concrete description: for every $J$-covering sieve $S:=\{c_{i} \leq c \textrm{ | } i\in I\}$ on an object $c\in {\cal C}$, $c$ is the supremum in $\cal C$ of the $c_{i}$ for $i\in I$ (i.e., for any object $c'$ in $\cal C$ such that for every $i\in I$ $c_{i}\leq c'$, $c\leq c'$). We can show this as follows. If $\cal C$ is a preorder then any representable functor $F:{\cal C}^{\textrm{op}}\to \Set$ is a subterminal object of the topos $[{\cal C}^{\textrm{op}}, \Set]$ and hence, the terminal object $1$ of $[{\cal C}^{\textrm{op}}, \Set]$ being always a sheaf, $F$ is a $J$-sheaf if and only if it is $J$-closed as a subobject of $1$. But this condition, for a representable ${\cal C}(-, a)$, amounts precisely to requiring that for any $J$-covering sieve $S:=\{c_{i} \leq c \textrm{ | } i\in I\}$ on an object $c\in {\cal C}$, if $c_{i}\leq a$ for all $i\in I$ then $c\leq a$, from which our thesis follows immediately.  

Recall that, if a functor $F:{\cal C}\to {\cal D}$ between two categories $\cal C$ and $\cal D$ is injective on objects (i.e., for any $c,c'\in {\cal C}$, $F(c)=F(c')$ implies $c=c'$) then we have a subcategory $Im(F)$ of $\cal D$, called the \emph{image} of $F$, whose objects are those of the form $F(c)$ for an object $c$ of $\cal C$ and whose arrows are those of the form $F(f)$ for an arrow $f$ of $\cal C$. Similarly, under the hypothesis that $F$ creates isomorphisms (i.e., for any isomorphism $l:F(c)\cong F(d)$ there exists an isomorphism $u:c\cong d$ such that $F(u)=l$), we can give the following definition.

\begin{definition}\label{extendedimage}
Let $F:{\cal C}\to {\cal D}$ be a functor which creates isomorphisms. The \emph{extended image} $ExtIm(F)$ of $F$ is the subcategory of $\cal D$ having as objects the objects of $\cal D$ which are isomorphic to an object of the form $F(c)$, and as arrows the arrows $f:x\to y$ in $\cal D$ such that there exist objects $c,c'\in {\cal C}$, an arrow $u:c\to c'$ in $\cal C$ and isomorphisms $x\cong F(c)$ and $y\cong F(c')$ such that $F(u)$ is the factorization of $f$ through these isomorphisms.
\end{definition}

The functors $F:{\cal C}\to {\cal D}$ which create isomorphisms enjoy a nice property, namely the fact that whenever they have a categorical left inverse, they yield an equivalence of categories between $\cal C$ and the extended image $ExtIm(F)$ of $F$. We shall appeal to this fact, recorded in the following proposition, in this section and in the next one.

\begin{proposition}\label{extim}
Let $F:{\cal C}\to {\cal D}$ be a functor which creates isomorphisms and let $G:{\cal D}\to {\cal C}$ be a functor such that $G\circ F \cong 1_{\cal C}$. Then the functors $\dot{F}:{\cal C}\to ExtIm(F)$ and $\dot{G}:ExtIm(F) \to {\cal C}$ obtained from $F$ and $G$ by restricting the codomain (resp. the domain) of $F$ (resp. of $G$) to the category $ExtIm(F)$ are categorical inverses to each other, and hence yield an equivalence of categories ${\cal C}\simeq ExtIm(F)$.
\end{proposition} 

\begin{proofs}
Let $\alpha:G\circ F \cong 1_{\cal C}$ be a natural isomorphism; then clearly $\alpha$ yields a natural isomorphism $\dot{G}\circ \dot{F} \cong 1_{\cal C}$. We can construct a natural isomorphism $\beta:\dot{F}\circ \dot{G}\to 1_{ExtIm(F)}$ as follows. Given $x\in ExtIm(F)$ there exists $c\in {\cal C}$ and an isomorphism $r:x\to F(c)$ in $\cal D$; we put $\beta(x):F(G(x))\to x$ equal to the composite $r^{-1}\circ F(\alpha(c)) \circ F(G(r))$. It is easy to see that this assignment does not depend on the choice of the object $c$ and of the isomorphism $r$, and that it is natural in $c\in ExtIm(F)$; therefore, all the arrows $\beta(x)$ being isomorphisms, this assignment defines a natural isomorphism $\beta:\dot{F}\circ \dot{G}\to 1_{ExtIm(F)}$.
\end{proofs}

Note that the proof of Proposition \ref{extim} does not require any form of the axiom of choice.

It is clear from their definitions that our functors $A$ and $B$ are injective on objects. In fact, we can also prove that they are faithful, provided that all the Grothendieck topologies $J_{\cal C}$ are subcanonical and the categories in $\cal K$ are posets. Let us show this for the functors $A$; we will discuss the faithfulness of the functor $B$ in section \ref{prealex}. 

From Lemma C2.3.8 \cite{El} we know that, given a morphism of sites $\hat{f}:({\cal C}, J_{{\cal C}})\to ({\cal D}, J_{{\cal D}})$ as above, if $J_{\cal C}$ and $J_{\cal D}$ are subcanonical then the inverse image of the geometric morphism $\dot{f}:\Sh({\cal D}, J_{{\cal D}}) \to \Sh({\cal C}, J_{{\cal C}})$ factors through the Yoneda embeddings ${\cal C} \hookrightarrow \Sh({\cal C}, J_{{\cal C}})$ and ${\cal D} \hookrightarrow \Sh({\cal D}, J_{{\cal D}})$ yielding a functor ${\cal C}\to {\cal D}$ which is isomorphic to $f$. If the categories $\cal C$ and $\cal D$ are posets, we can conclude that this functor actually coincides with $f$. This argument implies that, under the hypothesis that the categories $\cal C$ are posets and the topologies $J_{\cal C}$ are subcanonical, the functor $A$ is faithful; indeed, for any arrow $f$ in ${\cal K}^{\textrm{op}}$, $f$ is equal to the restriction ${\cal C} \to {\cal D}$ of the factorization $\Sh({\cal C}, J_{{\cal C}})\to \Sh({\cal D}, J_{{\cal D}})$ of the inverse image functor of the geometric morphism $\Sh(A(f)):\Sh(A({\cal D}))\to \Sh(A({\cal C}))$ through the equivalences of Theorem \ref{fund}. Therefore, since $A$ is faithful and injective on objects, it yields an isomorphism of categories between ${\cal K}^{\textrm{op}}$ and the image of $A$.    

We can summarize what we have established as follows.

\begin{theorem}\label{teoabstr}
With the above notation, if all the Grothendieck topologies $J_{\cal C}$ are subcanonical and all the categories in $\cal K$ are posets then the functor $A:{\cal K}^{\textrm{op}}\to \textbf{Loc}$ yields an isomorphism of categories between ${\cal K}^{\textrm{op}}$ and the subcategory of $\textbf{Loc}$ given by the image of $A$.
\end{theorem}

If our categories are just preorders (rather than posets), $f$ remains determined by $A(f)$ only up to isomorphism, and hence $A$ becomes faithful when we factorize the set of arrows in $\cal K$ by the equivalence relation which identifies two arrows $f,g:{\cal C}\to {\cal D}$ if and only if they are isomorphic (i.e., for any $c\in {\cal C}$, $f(c)\cong g(c)$); if we denote the resulting quotient category of $\cal K$ by ${\cal K}_{\cong}$ we thus obtain that the factorization $A_{\cong}:{\cal K}_{\cong}^{\textrm{op}} \to \textbf{Loc}$ of the functor $A$ through the quotient functor ${\cal K}^{\textrm{op}} \to {\cal K}_{\cong}^{\textrm{op}}$ (note that this factorization exists since $A$ is a $2$-functor when ${\cal K}^{\textrm{op}}$ and $\textbf{Loc}$ are regarded as $2$-categories in the natural way) is faithful and hence yields, under the hypothesis that all the Grothendieck topologies $J_{\cal C}$ be subcanonical, an isomorphism of categories between ${\cal K}_{\cong}^{\textrm{op}}$ and the subcategory of $\textbf{Loc}$ given by the image of $A_{\cong}$.

Theorem \ref{teoabstr} provides a categorical equivalence between the opposite of our original category of preorders $\cal K$ and a subcategory of $\textbf{Loc}$. Still, it would be desirable to have a duality of $\cal K$ with a subcategory of $\textbf{Loc}$ whose objects (and arrows) admit an intrinsic description in localic terms; now, since such characterizations can clearly determine the objects (and arrows) in $\textbf{Loc}$ only up to isomorphism, in order to obtain an `intrinsic' duality, we necessarily have to enlarge the target of the functor $A$ to the extended image of $A$. In fact, if we are able to define a functor $I_{A}:ExtIm(A)\to {\cal K}^{\textrm{op}}$ which allows us to recover each structure $\cal C$ in ${\cal K}^{\textrm{op}}$ from $A({\cal C})$ up to natural isomorphism (i.e. such that that the composite functor $I_{A}\circ A$ is naturally isomorphic to the identity functor on $\cal C$) then the composite functor $A\circ I_{A}$ will automatically be naturally isomorphic to the identity functor on $ExtIm(A)$ (by Proposition \ref{extim}). 

From now on we will suppose for simplicity that all the structures $\cal C$ in $\cal K$ are posets and that all the Grothendieck topologies $J_{\cal C}$ are subcanonical. 

Concerning the definition of the functor $I_{A}$, we note that it suffices to specify the action of $I_{A}$ on the objects, its action on the arrows being uniquely determined by it; indeed, for any locale $L$ in the extended image of $A$, $I_{A}(L)$ can be identified with a subset $B_{L}$ of $L$, and the action of $I_{A}$ on an arrow $f:L\to L'$ in $ExtIm(A)$ must be equal to the factorization $B_{L'}\to B_{L}$ of the frame homomorphism corresponding to $f$ through the inclusions $B_{L}\hookrightarrow L$ and $B_{L'}\hookrightarrow L'$. 

The definition of the functor $I_{A}$ on the objects relies on the existence of an appropriate topos-theoretic invariant $U$ which enables one to identify the principal $J_{\cal C}$-ideals on $\cal C$ among the $J_{\cal C}$-ideals on $\cal C$, for any structure $\cal C$ in $\cal K$; indeed, if such an invariant exists then each structure $\cal C$ in $\cal K$ can be recovered, up to isomorphism, as the poset of subterminals of the topos $\Sh({\cal C}, J_{\cal C})$ which satisfy the invariant $U$, and $I_{A}$ is defined as the functor which sends an locale $L$ in $ExtIm(A)$ to the set of subterminals of $\Sh(L)$ (i.e., of elements of $L$) which satisfy the invariant $U$, endowed with the induced order. As we shall see in section \ref{charinv}, the problem of the existence of an appropriate invariant $U$ is strictly related to the possibility of describing the Grothendieck topologies $J_{\cal C}$ `uniformly' by means of a topos-theoretic invariant; indeed, if the topologies $J_{\cal C}$ are `uniformly induced' by an invariant $C$ then the notion of $C$-compact subterminal yields an invariant $U$ which has the required property (cf. section \ref{charinv} below and in particular Theorem \ref{construction}). 

The existence of an appropriate topos-theoretic invariant $U$ is, in a sense, the only element of `non-canonicity' in our machinery building dualities; indeed, it is not \emph{a priori} guaranteed that such an invariant should exist for our category of structures $\cal K$ and for our choice of the Grothendieck topologies $J_{\cal C}$. On the other hand, our notion of Grothendieck topology induced by a topos-theoretic invariant, introduced in section \ref{charinv} below, provides us with a systematic means for building appropriate invariants for many interesting and naturally arising classes collections of Grothendieck topologies. In particular, all the classical examples of dualities fall under this framework, and in fact the level of generality of this notion is such that it allows one to easily build infinitely many new dualities through our machinery.              

Let us now discuss the problem of characterizing the subcategories of $\textbf{Loc}$ which arise as the extended image of a functor $A$ (the case of the functor $B$ is analogous and our arguments below can be easily adapted to work in that case, cf. section \ref{prealex} below) in `topological' terms, as for example in the case of the classical Stone duality. Clearly, one cannot reasonably expect this kind of problems to admit `canonical' solutions. On the other hand, topos-theoretic invariants can be profitably used to establish characterizations of the objects and arrows in the extended image of the functor $A$ starting from properties of the structures in $\cal K$. The notion of basis of a locale (resp. of basis of a topological space) is a particularly convenient one for achieving these characterizations. 

By a \emph{basis} of a frame $L$ (equivalently, of the locale corresponding to it) we mean a subset $B\subseteq L$ such that every element of $L$ is a join of elements in $B$; note that a basis of a topological space $X$ is precisely a basis of the frame ${\cal O}(X)$ of open sets of $X$. Below, when we say that a basis (of a frame or a topological space) is closed under some kind of operation we mean that the result of applying the operation to any set of elements of the basis belonging to its domain is again an element of the basis.

As we shall see in section \ref{charinv}, if the Grothendieck topologies $J_{\cal C}$ can be defined `uniformly' by means of a topos-theoretic invariant $C$ (technically, are all $C$-induced in the sense of Definition \ref{inducedtop} below) which satisfies a natural property (technically, the hypothesis of Theorem \ref{construction}) then the principal $J_{\cal C}$-ideals on $\cal C$ can be characterized among the $J_{\cal C}$-ideals on $\cal C$ (i.e. the subterminals of the topos $\Sh({\cal C}, J_{\cal C})$) as the `$C$-compact' ones (in the sense of Definition \ref{compact}(a) below), and the locales $L$ in the extended image of our functor $A$ can be characterized `intrinsically' as the locales which possess a basis of $C$-compact elements satisfying some specific invariant properties (cf. Theorem \ref{propext}). So, for example, if $\cal K$ consists of all the Boolean algebras, each of which equipped with the coherent topology, one requires the elements of the basis to be closed in the locale under finite meets and joins and to form, with respect to the induced order, a Boolean algebra (see section \ref{ex} below for more examples). 

Concerning the problem of characterizing the arrows in $ExtIm(A)$ `intrinsically', we observe that if the arrows ${\cal C}\to {\cal D}$ in $\cal K$ coincide exactly with the morphisms of sites $({\cal C}, J_{\cal C})\to ({\cal D}, J_{\cal D})$, as in fact it happens in many cases of interest, then, if the principal ideals on any $\cal C$ in $\cal K$ can be characterized as above as the $C$-compact elements of the locale $Id_{J_{\cal C}}({\cal C})$, the arrows in $ExtIm(A)$ can be characterized precisely as the locale morphisms whose associated frame homomorphisms send $C$-compact elements to $C$-compact elements (cf. Theorem \ref{propext}). If the arrows in $\cal K$ induce morphisms of sites $({\cal C}, J_{\cal C})\to ({\cal D}, J_{\cal D})$ but not all such morphisms yield arrows in $\cal K$ then the arrows in $ExtIm(A)$ can still be characterized, although less `intrinsically', as the locale morphisms whose associated frame homomorphisms, when restricted to the relevant bases of the locales, yield a function which, when the latter bases are regarded as objects of $\cal K$, can be identified with an arrow in $\cal K$ between them (cf. Theorem \ref{propext} below).

\subsection{Dualities with categories of topological spaces}\label{dualtop}

So far we have discussed how to build dualities or equivalences between categories of preorders and categories of locales. Our purpose in this section is to show that dualities (resp. equivalences) with categories of topological spaces can be naturally obtained starting from dualities (resp. equivalences) with categories of locales, this latter class being the one which, as we have seen, can be most naturally investigated by means of our topos-theoretic approach.

Our strategy consists in `enriching' a given duality (resp. equivalence) between a category of preorders $\cal K$ and a category of locales established by the method of section \ref{locales} to a duality (resp. equivalence) between $\cal K$ and a category of topological spaces, by means of an appropriate choice of points of the toposes corresponding to the structures in $\cal K$ as in section \ref{locales}.

We denote by $U:\textbf{Top}\to \textbf{Loc}$ the usual functor sending a topological space $X$ to the locale ${\cal O}(X)$ of open sets of $X$. 

Let as assume to start with a duality $A:{\cal K}^{\textrm{op}}\to \textbf{Loc}$ obtained by the method of section \ref{locales}. Suppose that we have assigned to every structure $\cal C$ in $\cal K$ an indexing $\xi_{\cal C}:[X_{\cal C}, \Set] \to \Sh({\cal C}, J_{{\cal C}})$ of a set of points of the topos $\Sh({\cal C}, J_{{\cal C}})$ which separates the subterminals in $\Sh({\cal C}, J_{{\cal C}})$, and to each arrow $f:{\cal C}\to {\cal D}$ in $\cal K$ a function $l_{f}:X_{\cal D} \to X_{\cal C}$ such that, denoted by $\dot{f}$ the geometric morphism $\Sh({\cal D}, J_{{\cal D}}) \to \Sh({\cal C}, J_{{\cal C}})$ induced by the morphism of sites $\hat{f}:({\cal C}, J_{{\cal C}}) \to ({\cal D}, J_{{\cal D}})$, the pair $(\dot{f}, l_{f})$ defines an arrow $(\Sh({\cal D}, J_{{\cal D}}), \xi_{\cal D}) \to (\Sh({\cal C}, J_{{\cal C}}), \xi_{\cal C})$ in the category $\mathfrak{Top}_{p}$ of toposes paired with points (cf. section \ref{subterminal} above). Then we can define a functor $\tilde{A}:{\cal K}^{\textrm{op}}\to \textbf{Top}$ as follows: 
\[
\tilde{A}({\cal C})=\Theta_{p}((\Sh({\cal C}, J_{{\cal C}}), \xi_{\cal C}))={X_{\cal C}}_{{{\tau}^{\Sh({\cal C}, J_{{\cal C}})}_{\xi_{\cal C}}}}
\]
for ${\cal C}\in {\cal K}$, and
\[
\tilde{A}(f)=\Theta_{p}(\dot{f},l_{f})=l_{f}:X_{\cal D} \to X_{\cal C}
\]
for any arrow $f:{\cal C}\to {\cal D}$ in $\cal K$ (cf. section \ref{subterminal} for the definition of the functor $\Theta_{p}$). 

By Theorem \ref{topol}(ii), we have an isomorphism of functors $U\circ \tilde{A}\cong A$. Now, if $\tilde{A}$ creates isomorphisms (note that, since $A$ creates isomorphisms, this condition is automatically satisfied when $U$ is faithful on the image of $\tilde{A}$, for example when all the topological spaces in the image of $\tilde{A}$ are sober) then the functor $\tilde{I_{A}}:ExtIm(\tilde{A}) \to {\cal K}^{\textrm{op}}$ defined as the composite of $I_{A}$ with the restriction of the functor $U$ to $ExtIm(\tilde{A})$ clearly yields an inverse (up to isomorphism) to the functor $\tilde{A}$ (cf. Proposition \ref{extim}). We thus obtain an equivalence between ${\cal K}^{\textrm{op}}$ and the subcategory of $\textbf{Top}$ given by the extended image of $\tilde{A}$, which `lifts' the equivalence between ${\cal K}^{\textrm{op}}$ and the extended image of $A$ from which we started.  

Let us now consider the covariant case. Let as assume to start with a functor $B:{\cal K} \to \textbf{Loc}$ obtained by the method of section \ref{locales}. Similarly as above, suppose that we have assigned to each structure $\cal C$ in $\cal K$ an indexing $\xi_{\cal C}:[X_{\cal C}, \Set] \to [{\cal C}, \Set]$ of a set of points of the topos $[{\cal C}, \Set]$ which separates the subterminals in $[{\cal C}, \Set]$, and to each arrow $f:{\cal C}\to {\cal D}$ in $\cal K$ a function $l_{f}:X_{\cal C} \to X_{\cal D}$ such that, denoted by $\dot{f}$ the geometric morphism $[{\cal C}, \Set] \to [{\cal D}, \Set]$ induced by $f$ as in section \ref{locales} above, the pair $(\dot{f}, l_{f})$ defines an arrow $([{\cal C}, \Set], \xi_{\cal C}) \to ([{\cal D}, \Set], \xi_{\cal D})$ in the category $\mathfrak{Top}_{p}$. Then we can define a functor $\tilde{B}:{\cal K}\to \textbf{Top}$ by putting
\[
\tilde{B}({\cal C})=\Theta_{p}(([{\cal C}, \Set], \xi_{\cal C}))={X_{\cal C}}_{{{\tau}^{[{\cal C}, \Set]}_{\xi_{\cal C}}}}
\]
for ${\cal C}\in {\cal K}$, and 
\[
\tilde{B}(f)=\Theta_{p}(\dot{f}, l_{f})=l_{f}:X_{\cal C} \to X_{\cal D}
\]
for any arrow $f:{\cal C}\to {\cal D}$ in $\cal K$. 

We verify as above that $U\circ \tilde{B}\cong B$; and, provided that $B$ is faithful and creates isomorphisms, if $\tilde{B}$ creates isomorphisms then $\tilde{B}$ yields an equivalence of categories onto its extended image. 

We note that, in the contravariant (resp. covariant) case, if for every structure $\cal C$ in $\cal K$ the indexing $\xi_{\cal C}:[X_{\cal C}, \Set] \to \Sh({\cal C}, J_{{\cal C}})$ (resp. $\xi_{\cal C}:[X_{\cal C}, \Set] \to [{\cal C}, \Set]$) of points of the topos $\Sh({\cal C}, J_{{\cal C}})$ (resp. of the topos $[{\cal C}, \Set]$) is induced by the identity indexing function on the set $\textsc{P}_{\cal C}$ of \emph{all} the points of the topos $\Sh({\cal C}, J_{{\cal C}})$ (resp. of the topos $[{\cal C}, \Set]$) then for every arrow $f$ in $\cal K^{\textrm{op}}$ (resp. in $\cal K$), the definition of $l_{f}:X_{\cal D} \to X_{\cal C}$ (resp. of $l_{f}:X_{\cal C} \to X_{\cal D}$) is automatic; indeed, $l_{f}$ must be equal to the function sending a point $p:\Set \to \Sh({\cal D}, J_{{\cal D}})$ (resp. a point $p:\Set \to [{\cal C}, \Set]$) to the composite $\dot{f}\circ p$.

We note that, while the method of section \ref{locales} provides a `canonical' way of building dualities or equivalences with subcategories of \textbf{Loc}, once the category $\cal K$ of structures and the Grothendieck topologies $J_{\cal C}$ have been fixed, the method for building dualities with categories of topological spaces that we have just explained requires the choice of appropriate sets of points of the toposes involved, and this can be done in general in several different ways, possibly leading to different dualities or equivalences between the category $\cal K$ and a subcategory of $\textbf{Top}$ (cf. section \ref{ex} below for a concrete example of this phenomenon).   

In connection with this, we also remark that some form of the axiom of choice may be necessary in establishing that certain sets of points of a given topos separate the subterminals of the topos (we refer the reader to section \ref{ex} for a detailed analysis of various examples also in light of the axiom of choice).

Concerning the problem of characterizing the subcategories of $\textbf{Top}$ arising as the extended image of a functor $\tilde{A}$ (resp. of a functor $\tilde{B}$) in `topological' terms, we have a particularly satisfying answer when the topological spaces in the extended image of the functor are all sober (cf. Theorem \ref{sobriety} for a topos-theoretic criterion of sobriety); indeed, under this hypothesis the spaces in $\cal V$ can be characterized precisely as the sober topological spaces whose frames of open sets (regarded as locales) belong to the extended image of the functor $A$ (resp. of the functor $B$). Indeed, one direction is obvious while the other one follows from the fact that the assignment $X\to {\cal O}(X)$ defines an embedding of the category of sober spaces into the category of locales (cf. Corollary C1.2.3 \cite{El}). In particular, if the locales in the extended image of $A$ (resp. of $B$) can be characterized by means of a topos-theoretic invariant $C$ as in section \ref{locales} above then the spaces in $ExtIm(\tilde{A})$ (resp. in $ExtIm(\tilde{B})$) can be characterized precisely as the sober topological spaces such that the collection of their $C$-compact open sets forms a basis which, endowed with the subset-inclusion ordering, has the structure of a preorder in $\cal K$ and satisfies some specific invariant properties (as in the characterization of the locales in $ExtIm(A)$ (resp. in $ExtIm(B)$)), while the arrows in $ExtIm(\tilde{A})$ (resp. in $ExtIm(\tilde{B})$) can be characterized precisely as the continuous maps between the spaces in $ExtIm(\tilde{A})$ (resp. in $ExtIm(\tilde{B})$) such that their inverse images send $C$-compact open sets to $C$-compact open sets. 

\subsection{Characterizing principal ideals through invariants}\label{charinv}

In this section, we investigate the problem of recovering a (preorder) category $\cal C$, equipped with a (subcanonical) Grothendieck topology $J$ from the topos $\Sh({\cal C}, J)$ (up to categorical equivalence) by means of a topos-theoretic invariant. Notice that any categorical equivalence between two preorder categories is an isomorphism, so our considerations will address in particular the problem of recovering a preordered structure $\cal C$ in $\cal K$ from the topos $\Sh({\cal C}, J_{{\cal C}})$ (respectively, from the topos $[{\cal C}, \Set]$) uniquely up to isomorphism.

Let us denote by $l_{{\cal C}}^{J}:{\cal C}\to \Sh({\cal C}, J)$ the composite of the Yoneda embedding $Y:{\cal C}\to [{\cal C}^{\textrm{op}}, \Set]$ with the associated sheaf functor $a_{J}:[{\cal C}^{\textrm{op}}, \Set]\to \Sh({\cal C}, J)$. The functor $l_{{\cal C}}^{J}$ maps $\cal C$ into $\Sh({\cal C}, J)$ not faithfully in general; anyway, if $J$ is subcanonical then $l_{{\cal C}}^{J}$ can be identified with the Yoneda embedding and hence it is full and faithful. It is thus natural to ask, under this hypothesis, whether one can characterize the objects of $\Sh({\cal C}, J)$ of the form $l_{{\cal C}}^{J}(c)$ for an object $c$ of $\cal C$ in terms of a topos-theoretic invariant applied to the topos $\Sh({\cal C}, J)$. If ${\cal C}$ is a preorder $({\cal P}, \leq)$ then for an object $c$ of $\cal C$ the unique arrow $l_{{\cal C}}^{J}(c)\to 1_{\Sh({\cal C}, J)}$ is a monomorphism and hence $l_{{\cal C}}^{J}(c)$ can be identified with a $J$-ideal on $\cal C$; in fact, if $J$ is subcanonical then $l_{{\cal C}}^{J}(c)$ is (isomorphic to) the the principal $J$-ideal $(c)\downarrow$ generated by $c$ (i.e. the collection of all the objects $d$ of $\cal P$ such that $d\leq c$). This leads us to investigating the problem of characterizing the subterminals in $\Sh({\cal C}, J)$ which correspond to the principal $J$-ideals on $\cal C$ through a topos-theoretic invariant applied to the topos $\Sh({\cal C}, J)$.

The following definition will be central for our purposes. Recall from section \ref{exsub} that a $J$-ideal $I$ on $\cal C$ is said to be principal if $I=(c)\downarrow_{J}$ for an object $c$ of $\cal C$.

\begin{definition}
Let $({\cal C}, J)$ be a site and $I$ be a $J$-ideal on $\cal C$. We say that $I$ is \emph{$J$-compact} if for every covering $\{I_{k} \textrm{ | } k\in K\}$ of $I$ by subterminals in $\Sh({\cal C}, J)$ there exists a $J$-covering sieve $\{f_{h}:c_{h}\to c \textrm{ | } h\in H\}$ in $\cal C$ such that for every $h\in H$, $(c_{h})\downarrow_{J}$ is contained in some $I_{k}$, and $\{(c_{h})\downarrow_{J} \textrm{ | } h\in H\}$ covers $I$ in $\Sh({\cal C}, J)$.   
\end{definition}  

\begin{theorem}\label{Jcomp}
Let $({\cal C}, J)$ be a site and $I$ be a $J$-ideal on $\cal C$. Then $I$ is a principal $J$-ideal if and only if it is $J$-compact.
\end{theorem}

\begin{proofs}
Let us suppose that $I=(c)\downarrow_{J}$ is a $J$-principal ideal. The union $\mathbin{\mathop{\textrm{\huge $\vee$}}\limits_{k\in K}}I_{k}$ in $\Sh({\cal C}, J)$ of a family $\{I_{k} \textrm{ | } k\in K\}$ of $J$-ideals is the $J$-closure of the union $\mathbin{\mathop{\textrm{\huge $\cup$}}\limits_{k\in K}}I_{k}$ in $[{\cal C}^{\textrm{op}}, \Set]$ of the ideals $I_{k}$ (considered as subterminals in $[{\cal C}^{\textrm{op}}, \Set]$); specifically, $\mathbin{\mathop{\textrm{\huge $\vee$}}\limits_{k\in K}}I_{k}=\{d\in {\cal C} \textrm{ | } \{f:e\to d \textrm{ | } e\in \mathbin{\mathop{\textrm{\huge $\cup$}}\limits_{k\in K}}I_{k}\} \in J(d)\}$. This implies that $I=\mathbin{\mathop{\textrm{\huge $\vee$}}\limits_{k\in K}}I_{k}$ if and only if $\{f:e\to c \textrm{ | } e\in \mathbin{\mathop{\textrm{\huge $\cup$}}\limits_{k\in K}}I_{k}\} \in J(c)$. So, if we put $S=\{f:e\to c \textrm{ | } e\in \mathbin{\mathop{\textrm{\huge $\cup$}}\limits_{k\in K}}I_{k}\}$ then we have that for any $f:e\to c$ in $S$, $(e)\downarrow_{J}$ is contained in some $I_{k}$. To show that the sieve $S$ witnesses the fact that $I$ is $J$-compact, it remains to prove that $I=\mathbin{\mathop{\textrm{\huge $\vee$}}\limits_{f\in S}}(dom(f))\downarrow_{J}$. Since each $(dom(f))\downarrow_{J}$ is contained in $I$ (being included in some $I_{k}$) then the inclusion ``$\supseteq$'' holds. The other inclusion follows from the fact that every $J$-ideal that contains each of the $(dom(f))\downarrow_{J}$ must contain $c$ (since $S$ is $J$-covering) and hence $I$.    

Conversely, suppose that $I$ is $J$-compact. We want to prove that $I$ is principal. If we take as family $\{I_{k} \textrm{ | } k\in K\}$ the singleton family $\{I\}$, the fact that $I$ is $J$-compact implies that there exists a $J$-covering sieve $\{f_{h}:c_{h}\to c \textrm{ | } h\in H\}$ in $\cal C$ such that $I=\mathbin{\mathop{\textrm{\huge $\vee$}}\limits_{h\in H}}(c_{h})\downarrow_{J}$. But $\mathbin{\mathop{\textrm{\huge $\vee$}}\limits_{h\in H}}(c_{h})\downarrow_{J}=(c)\downarrow_{J}$, from which it follows that $I=(c)\downarrow_{J}$ and hence that $I$ is a principal $J$-ideal.
\end{proofs}

The characterization of Theorem \ref{Jcomp} represents a first step for achieving, under appropriate hypotheses, characterizations of the property of being a principal $J$-ideal on a preorder $\cal C$ as a topos-theoretic invariant on the topos $\Sh({\cal C}, J)$. To this end, we introduce the following definitions.

\begin{definition}\label{ade}
\begin{enumerate}[(a)]
\item We say that a family of subobjects $\{a_{i} \mono a \textrm{ |} i\in I\}$ of a given object $a$ in a (locally small cocomplete) topos $\cal E$ is \emph{refined} by a family $\{b_{j} \mono a \textrm{ | } j\in J\}$ of subobjects of $a$ if for every $j\in J$ there exists $i\in I$ such that $b_{j}\mono a$ factors through $a_{i}\mono a$ and the union $\mathbin{\mathop{\textrm{\huge $\vee$}}\limits_{i\in I}}a_{i}$ in $\Sub_{\cal E}(a)$ is equal to the union $\mathbin{\mathop{\textrm{\huge $\vee$}}\limits_{j\in J}}b_{j}$ in $\Sub_{\cal E}(a)$.

\item Let $({\cal C}, J)$ be a site. We say that a topos-theoretic invariant property $C$ of families ${\cal F}$ of subterminals in a (locally small cocomplete) topos $\cal E$ is \emph{$({\cal C}, J)$-adequate} if, when $\cal E$ is the topos $\Sh({\cal C}, J)$ and $\cal F$ is a family of principal $J$-ideals on $\cal C$ (regarded as subterminals in $\Sh({\cal C}, J)$), $\cal F$ has a refinement which satisfies $C$ if and only if there exists a $J$-covering sieve $\{f_{h}:d_{h}\to d \textrm{ | } h\in H\}$ in $\cal C$ such that for any $h\in H$, $(d_{h})\downarrow_{J}$ is contained in some ideal in $\cal F$ and the union in $\cal E$ of the family $\{(d_{h})\downarrow_{J} \textrm{ | } h\in H\}$ is equal to the union in $\cal E$ of the family $\cal F$ (equivalently, the objects of $\cal C$ whose corresponding principal $J$-ideals are contained in some ideal in $\cal F$ are the domains of a family of arrows in $\cal C$ which generates a $J$-covering sieve).
\end{enumerate}
\end{definition}

\begin{remark}\label{grot}
If $({\cal C}, \leq)$ is a poset category and $J$ is subcanonical then the elements in $\cal F$ can be thought of as elements of $\cal C$, and the condition that there exist a $J$-covering sieve $\{f_{h}:d_{h}\to d \textrm{ | } h\in H\}$ in $\cal C$ such that for any $h\in H$, $(d_{h})\downarrow_{J}$ is contained in some ideal in $\cal F$ and the union in $\cal E$ of the family $\{(d_{h})\downarrow_{J} \textrm{ | } h\in H\}$ is equal to the union in $\cal E$ of the family $\cal F$ can be equivalently reformulated as the requirement that the supremum $s$ in $\cal C$ of the set of the objects of $\cal C$ whose corresponding principal $J$-ideals are contained in some ideal in $\cal F$ should exist and the collection of all the arrows in $\cal C$ from such objects to $s$ should generate a $J$-covering sieve.   
\end{remark}

\begin{theorem}\label{equivJcomp}
Let $C$ be an invariant property of families of subterminals in a (locally small cocomplete) topos which is $({\cal C}, J)$-adequate. Then a $J$-ideal on $\cal C$ is $J$-compact (equivalently, by Theorem \ref{Jcomp}, principal) if and only if every covering family for $I$ has a refinement which satisfies $C$.       
\end{theorem}

\begin{proofs}
If a $J$-ideal $I$ is $J$-compact then every covering family for $I$ has a refinement by principal $J$-ideals which satisfies $C$.

Conversely, suppose that every covering family for $I$ has a refinement which satisfies $C$. If we take as covering family the collection of all the principal $J$-ideals of the form $(c)\downarrow_{J}$ for $c\in I$ then we conclude, from the fact that $C$ is $({\cal C}, J)$-adequate, that $I$ is the join of a family of principal $J$-ideals indexed over a $J$-covering sieve, and hence it is principal (equivalently, $J$-compact), as required.
\end{proofs}

Let us now introduce some natural topos-theoretic invariants of families of subobjects in a topos, which we shall use in section \ref{ex} in the context of our examples.

Below, by atomic (resp. supercompact) subobject of a given object $a$ in a locally small cocomplete topos $\cal E$ we mean a subobject $m$ of $a$ such that $m$ is non-zero and does not contain any proper subobject of $a$ in $\cal E$ (resp. such that every covering of $m$ in the frame $\Sub_{{\cal E}}(a)$ contains the identity subobject of $a$). An element $a$ of a frame $L$ (equivalently, of the corresponding locale) is an atom (resp. supercompact) if, denoted by $y:L\to \Sh(L)$ the Yoneda embedding, $y(a)\mono 1_{\Sh(L)}$ is an atomic (resp. supercompact) subobject in the topos $\Sh(L)$. 

\begin{definition}\label{refinsubcover} 
Let $\cal E$ be a locally small cocomplete topos and $\cal F$ be a family of subobjects of a given object in $\cal E$.

\begin{enumerate}[(a)]
\item The family $\cal F$ is said to \emph{have a finite subcover} if there exists a finite subfamily ${\cal F}'$ of $\cal F$ such that the union of the subobjects in ${\cal F}'$ is equal to the union of the subobjects in $\cal F$;

\item The family $\cal F$ is said to \emph{have a singleton subcover} if there exists a single subobject in $\cal F$ which is the union of the subobjects in the family ${\cal F}$;

\item Given a regular cardinal $k$, the family $\cal F$ is said to \emph{have a $k$-subcover} if there exists a subfamily ${\cal F}'$ of $\cal F$ of cardinality $\lt k$ such that the union of the subobjects in ${\cal F}'$ is equal to the union of the subobjects in $\cal F$;

\item The family $\cal F$ is said to \emph{have a disjunctive refinement} (resp. to \emph{have a finite disjunctive refinement}) if there exists a family (resp. a finite family) of pairwise disjoint subobjects which refines $\cal F$;

\item The family $\cal F$ is said to \emph{have an atomic refinement} (resp. to \emph{have a finite atomic refinement} if either $\cal F$ is a singleton or there exists a family (resp. a finite family) of atomic subobjects which refines $\cal F$;

\item The family $\cal F$ is said to \emph{have a supercompact refinement} (resp. to \emph{have a finite supercompact refinement} if either $\cal F$ is a singleton or there exists a family (resp. a finite family) of supercompact subobjects which refines $\cal F$;

\item The family $\cal F$ is said to \emph{have a directed refinement} if there exists a directed family $\cal G$ (i.e. a non-empty family such that for any two subobjects in $\cal G$ there exists a subobject in $\cal G$ which contains both) which refines $\cal F$.
 
\end{enumerate}
\end{definition}

Notice that for $k=2$ (resp. $k=\omega$) the invariant `to have a $k$-subcover' specialize to the invariant `to have a singleton subcover' (resp. `to have a finite subcover').   

\begin{remark}\label{refin}
All the invariants defined above are of the form `to have a refinement satisfying $C$' for a topos-theoretic invariant $C$ of families of subterminals in a topos.
\end{remark}

The definition above naturally leads us to introducing the following notions (the concepts of compactness and supercompactness in the definition below already appear in \cite{El}).

\begin{definition}
Let $\cal E$ be a locally small cocomplete topos and let $L$ be a locale, with corresponding Yoneda embedding $y:L\to \Sh(L)$.
\begin{enumerate}[(a)]
\item An object $A$ of $\cal E$ is said to be \emph{compact} if every covering family of $A$ in $\cal E$ has a finite subcover. 

An element $a$ of a locale $L$ (equivalently, of a frame ${\cal O}(L)$) is said to be compact if the object $y(a)$ is compact in $\Sh(L)$, equivalently if whenever $a=\mathbin{\mathop{\textrm{\huge $\vee$}}\limits_{i\in I}}a_{i}$ in $L$ there exists a finite subset $I_{0}\subseteq I$ such that $a=\mathbin{\mathop{\textrm{\huge $\vee$}}\limits_{i\in I_{0}}}a_{i}$. 

\item An object $A$ of $\cal E$ is said to be \emph{supercompact} if every covering family of $A$ has a singleton subcover. 

An element $a$ of a locale $L$ (equivalently, of a frame ${\cal O}(L)$) is said to be supercompact if the object $y(a)$ is supercompact in $\Sh(L)$, equivalently if whenever $a=\mathbin{\mathop{\textrm{\huge $\vee$}}\limits_{i\in I}}a_{i}$ in $L$ there exists an index $i\in I$ such that $a=a_{i}$.

\item For a regular cardinal $k$, an object $A$ of $\cal E$ is said to be \emph{$k$-compact} if every covering family of $A$ in $\cal E$ has a $k$-subcover. 

An element $a$ of a locale $L$ (equivalently, of a frame ${\cal O}(L)$) is said to be $k$-compact if the object $y(a)$ is $k$-compact in $\Sh(L)$, equivalently if whenever $a=\mathbin{\mathop{\textrm{\huge $\vee$}}\limits_{i\in I}}a_{i}$ in $L$ there exists a subset $I_{0}\subseteq I$ of cardinality $\lt k$ such that $a=\mathbin{\mathop{\textrm{\huge $\vee$}}\limits_{i\in I_{0}}}a_{i}$. 

\item An object $A$ of $\cal E$ is said to be \emph{disjunctively compact} (resp. \emph{infinitarily disjunctively compact}) if every covering family of $A$ in $\cal E$ has a finite disjunctive refinement (resp. a disjunctive refinement). 

An element $a$ of a locale $L$ (equivalently, of a frame ${\cal O}(L)$) is said to be disjunctively compact (resp. infinitarily disjunctively compact) if the object $y(a)$ is disjunctively compact (resp. infinitarily disjunctively compact) in $\Sh(L)$, equivalently if whenever $a=\mathbin{\mathop{\textrm{\huge $\vee$}}\limits_{i\in I}}a_{i}$ in $L$ there exists a finite family (resp. a family) $\{b_{j} \leq a \textrm{ | } j\in J\}$ of elements of $L$ such that for every $j\in J$ there exists $i\in I$ such that $b_{j}\leq a_{i}$, the join $\mathbin{\mathop{\textrm{\huge $\vee$}}\limits_{i\in I}}a_{i}$ in $L$ is equal to the join $\mathbin{\mathop{\textrm{\huge $\vee$}}\limits_{j\in J}}b_{j}$ and for any distinct $j,j'\in J$, $b_{j}\wedge b_{j'}=0$.

\item An object $A$ of $\cal E$ is said to be \emph{atomically compact} (resp. \emph{infinitarily atomically compact}) if every covering family of $A$ in $\cal E$ has a finite atomic refinement (resp. an atomic refinement). 

An element $a$ of a locale $L$ (equivalently, of a frame ${\cal O}(L)$) is said to be \emph{atomically compact} (resp. \emph{infinitarily atomically compact}) if the object $y(a)$ is atomically compact (resp. infinitarily atomically compact) in $\Sh(L)$, equivalently if whenever $a=\mathbin{\mathop{\textrm{\huge $\vee$}}\limits_{i\in I}}a_{i}$ in $L$ either there is $i\in I$ such that $a_{i}=a$ or there exists a finite family (resp. a family) $\{b_{j} \leq a \textrm{ | } j\in J\}$ of elements of $L$ such that the $b_{j}$ (for $j\in J$) are atoms in $L$, for every $j\in J$ there exists $i\in I$ such that $b_{j}\leq a_{i}$, and the join $\mathbin{\mathop{\textrm{\huge $\vee$}}\limits_{i\in I}}a_{i}$ in $L$ is equal to the join $\mathbin{\mathop{\textrm{\huge $\vee$}}\limits_{j\in J}}b_{j}$.

\item An object $A$ of $\cal E$ is said to be \emph{supercompactly compact} (resp. \emph{infinitarily supercompactly compact}) if every covering family of $A$ in $\cal E$ has a finite supercompact refinement (resp. a supercompact refinement). 

An element $a$ of a locale $L$ (equivalently, of a frame ${\cal O}(L)$) is said to be \emph{supercompactly compact} (resp. \emph{infinitarily supercompactly compact}) if the object $y(a)$ is supercompactly compact (resp. infinitarily supercompactly compact) in $\Sh(L)$, equivalently if whenever $a=\mathbin{\mathop{\textrm{\huge $\vee$}}\limits_{i\in I}}a_{i}$ in $L$ either there is $i\in I$ such that $a_{i}=a$ or there exists a finite family (resp. a family) $\{b_{j} \leq a \textrm{ | } j\in J\}$ of elements of $L$ such that the $b_{j}$ (for $j\in J$) are supercompact elements of $L$, for every $j\in J$ there exists $i\in I$ such that $b_{j}\leq a_{i}$, and the join $\mathbin{\mathop{\textrm{\huge $\vee$}}\limits_{i\in I}}a_{i}$ in $L$ is equal to the join $\mathbin{\mathop{\textrm{\huge $\vee$}}\limits_{j\in J}}b_{j}$.

\item An object $A$ of $\cal E$ is said to be \emph{directedly compact} if every covering family of $A$ in $\cal E$ has a directed refinement. 

An element $a$ of a locale $L$ (equivalently, of a frame ${\cal O}(L)$) is said to be \emph{directedly compact} if the object $y(a)$ is directedly compact in $\Sh(L)$, equivalently if whenever $a=\mathbin{\mathop{\textrm{\huge $\vee$}}\limits_{i\in I}}a_{i}$ in $L$ there exists a family $\{b_{j} \leq a \textrm{ | } j\in J\}$ of elements of $L$ such that for any $j, j'\in J$ there exists $j''\in J$ such that $b_{j}\leq b_{j''}$ and $b_{j'}\leq b_{j''}$, for every $j\in J$ there exists $i\in I$ such that $b_{j}\leq a_{i}$, and the join $\mathbin{\mathop{\textrm{\huge $\vee$}}\limits_{i\in I}}a_{i}$ in $L$ is equal to the join $\mathbin{\mathop{\textrm{\huge $\vee$}}\limits_{j\in J}}b_{j}$.

\end{enumerate}
\end{definition}

\begin{remark}\label{loctop}
Topos-theoretic invariants of families of subterminals in a locally small cocomplete topos can be identified with locale-theoretic invariants. Indeed, a topos-theoretic invariant $U$ of families of subterminals in a locally small cocomplete topos yields a locale-theoretic invariant property $L_{U}$ of families of elements of a locale, defined by saying that for any locale $L$, $a\in L$ satisfies $L_{U}$ if and only if $y(a)$ satisfies $U$ in the topos $\Sh(L)$, where $y:L\hookrightarrow \Sh(L)$ is the Yoneda embedding; conversely, given a locale-theoretic invariant property $Q$, we have a topos-theoretic invariant $T_{Q}$ of families of subterminals of a topos, defined by saying that a family of subterminals in a locally small cocomplete topos $\cal E$ satisfies $T_{Q}$ if and only if, regarded as a family of elements of $\Sub_{\cal E}(1)$, it satisfies $Q$. 
\end{remark}

The following abstract definition unifies the notions introduced above and enlightens the link between them and the concepts of Definition \ref{refinsubcover}.

\begin{definition}\label{compact}
\begin{enumerate}[(a)]
\item Given a topos-theoretic invariant property $C$ of families of subobjects of a given object in a (locally small cocomplete) topos, an object $A$ of a topos is said to be \emph{$C$-compact} if and only if any covering family of subobjets of $A$ in $\cal E$ has a refinement by a family of subobjects satisfying the property $C$;

\item Given a topos-theoretic invariant $C$ of families of subterminals in a topos (equivalently, by Remark \ref{loctop}, a locale-theoretic invariant of families of elements of a locale), an element $l$ of a locale $L$ (equivalently, of the corresponding locale) is said to be \emph{$C$-compact} if every covering family of $l$ in $L$ has a refinement satisfying $C$; that is, whenever $a=\mathbin{\mathop{\textrm{\huge $\vee$}}\limits_{i\in I}}a_{i}$ in $L$ there exists a family $\{b_{j} \leq a \textrm{ | } j\in J\}$ of elements of $L$ satisfying $C$ such that for every $j\in J$ there exists $i\in I$ such that $b_{j}\leq a_{i}$, and the join $\mathbin{\mathop{\textrm{\huge $\vee$}}\limits_{i\in I}}a_{i}$ in $L$ is equal to the join $\mathbin{\mathop{\textrm{\huge $\vee$}}\limits_{j\in J}}b_{j}$ in $L$.
\end{enumerate}
\end{definition} 

\begin{examples}\label{excomp}
\begin{enumerate}[(a)]
\item An object is compact if and only if it is $C$-compact where $C$ is the invariant `to be finite';  

\item An object is supercompact if and only if it is $C$-compact where $C$ is the invariant `to be a singleton';
  
\item An object is $k$-compact (for a regular cardinal $k$) if and only if it is $C$-compact wheren $C$ is the invariant `to be of cardinality $\lt k$';  

\item An object is disjunctively compact (resp. infinitarily disjunctively compact) if and only if it is $C$-compact where $C$ is the invariant `to be finite and disjoint' (resp. to be disjoint');  

\item An object is atomically compact (resp. infinitarily atomically compact) if and only if it is $C$-compact where $C$ is the invariant `to be a singleton or to be finite and consisting of atomic subobjects' (resp. to be a singleton or to consist of atomic subobjects');

\item An object is supercompactly compact (resp. infinitarily supercompactly compact) if and only if it is $C$-compact where $C$ is the invariant `to be a singleton or to be finite and consisting of supercompact subobjects' (resp. to be a singleton or to consist of supercompact subobjects')

\item An object is directedly compact if and only if it is $C$-compact where $C$ is the invariant `to be directed'.
  
\end{enumerate}
\end{examples}

We note that on any \emph{disjunctively coherent category} $\cal C$ (i.e., regular category in which the initial object exists and the unions of finite families of pairwise disjoint subobjects exist and are stable under pullback) one can put the \emph{disjunctive topology} i.e. the Grothendieck topology $D_{\cal C}$ on $\cal C$ whose $D_{\cal C}$-covering sieves on any object $c$ of $\cal C$ are exactly the covering sieves on $c$ in $\cal C$ which contain finite families of subobjects of $c$ which are pairwise disjoint from each other. Note that on a Boolean algebra (regarded as a coherent category) the disjunctive topology coincides with the coherent topology, while on a total order (regarded as a coherent category) the disjunctive topology specializes exactly to the trivial topology. 

Similarly, on any \emph{disjunctively geometric category} $\cal C$ (i.e., regular category in which the initial object exists and the unions of set-indexed families of pairwise disjoint subobjects exist and are stable under pullback) one can put the \emph{infinitary disjunctive topology}, i.e. the Grothendieck topology $G_{\cal C}$ on $\cal C$ whose $G_{\cal C}$-covering sieves on any object $c$ of $\cal C$ are exactly the covering sieves on $c$ in $\cal C$ which contain set-indexed families of subobjects of $c$ which are pairwise disjoint from each other. 

Recall that, given a regular cardinal $k$, a \emph{$k$-geometric category} is a regular category in which unions of arbitrary families of $\lt k$ subobjects exist and are stable under pullback; on a $k$-geometric category $\cal C$ we can define the \emph{$k$-covering topology} as the Grothendieck topology on $\cal C$ whose covering sieves on any object $c$ of $\cal C$ are exactly the sieves on $c$ which contain covering families of subobjects of $c$ of cardinality $\lt k$. Note that the $2$-geometric categories are precisely the regular categories, while the $\omega$-geometric categories are precisely the coherent categories.

We shall call a poset a \emph{disjunctively distributive lattice} (resp. an \emph{disjunctively distributive frame}) if is is a disjunctively coherent category (resp. a disjunctively geometric category) when regarded as a preorder category. We define the category $\textbf{DJLat}$ of disjunctively distributive lattices as the category whose objects are the disjunctively distributive lattices and whose arrows are the meet-semilattices homomorphisms between them which preserve (the $0$ and) finite pairwise disjoint joins. Similarly, we define the category $\textbf{DJFrm}$ of disjunctively distributive frames as the category whose objects are the disjunctively distributive frames and whose arrows are the meet-semilattices homomorphisms between them which preserve (the $0$ and) arbitrary pairwise disjoint joins.

We shall call a poset a \emph{$k$-frame} (cf. \cite{kframes}) if it is a $k$-geometric category when regarded as a preorder category, equivalently if it is a meet-semilattice in which there exists the join of any family of $\lt k$-elements and the infinite distributive law of these joins with respect to binary meets holds. We define the category $k\textrm{-}\textbf{Frm}$ of $k$-frames as the category whose objects are the $k$-frames and whose arrows are the meet-semilattices homomorphisms between them which send joins of families of $\lt k$-elements to joins. 

Recall that a \emph{preframe} is a poset with finite meets and directed joins, in which the finite meets distribute over the directed joins. We denote by $\textbf{PFrm}$ the category having as objects the preframes and as arrows the meet-semilattice homomorphisms between them which send directed joins to directed joins. Given a preframe $\cal P$, we define the \emph{directed topology} $J^{dir}_{{\cal P}}$ on $\cal P$ as the Grothendieck topology on $\cal P$ whose covering sieves on any object $p\in {\cal P}$ are the sieves on $p$ which contain a directed sieve on $p$, i.e. a sieve $S$ on $p$ such that for any arrows $a\leq p$ and $b\leq p$ in $S$ there exists $c\leq p$ in $S$ such that $a\leq c$, $b\leq c$.       

Let us define a \emph{weakly atomic meet-semilattice} to be a meet-semilattice in which the bottom element $0$ exists and in which finite joins of atoms always exist and distribute over finite meets. Similarly, we define an \emph{infinitarily weakly atomic meet-semilattice} to be a meet-semilattice in which the bottom element $0$ exists and in which arbitrary joins of atoms always exist and distribute over finite meets. We define the category $\textbf{WAtMSLat}$ of weakly atomic meet-semilattices as the category whose objects are the weakly atomic meet-semilattices and whose arrows are the meet-semilattices homomorphisms between them which preserve the $0$, send atoms to atoms and preserve finite joins of atoms. Similarly, we define the category $\textbf{IWAtMSLat}$ of infinitarily weakly atomic meet-semilattices as the category whose objects are the infinitarily weakly atomic meet-semilattices and whose arrows are the meet-semilattices homomorphisms between them which preserve the $0$, send atoms to atoms and preserve arbitrary joins of atoms. 

On a weakly atomic meet-semilattice $\cal M$ (regarded as a cartesian preorder category) one can define the \emph{atomically generated topology} $J^{at}_{\cal M}$, as follows. For any sieve $S$ in $\cal M$ on an object $m\in {\cal M}$, we set $S\in J^{at}_{\cal M}(m)$ if and only if either $S$ is the maximal sieve or $S$ is a covering sieve in $\cal C$ (i.e. $m$ is the supremum in $\cal M$ of all the domains of the arrows in $S$) generated by a finite family of atoms in $\cal M$. Clearly, on an infinitarily weakly atomic meet-semilattice we can define the \emph{infinitary atomically generated topology}, as the (obvious) infinitary analogue of the atomically generated topology.   

If in the definition of weakly atomic meet-semilattice we replace atoms by supercompact objects we obtain a similar but more general notion: we define a \emph{weakly supercompact meet-semilattice} to be a meet-semilattice in which the bottom element $0$ exists and in which finite joins of supercompact elements always exist and distribute over finite meets. Similarly, we define an \emph{infinitarily weakly supercompact meet-semilattice} to be a meet-semilattice in which the bottom element $0$ exists and in which arbitrary joins of supercompact elements exist and distribute over finite meets. We define the category $\textbf{WScMSLat}$ of weakly supercompact meet-semilattices as the category whose objects are the weakly supercompact meet-semilattices and whose arrows are the meet-semilattices homomorphisms between them which preserve the $0$, send supercompact elements to supercompact elements and preserve finite joins of supercompact elements. Similarly, we define the category $\textbf{IWScMSLat}$ of infinitarily weakly supercompact meet-semilattices as the category whose objects are the infinitarily weakly supercompact meet-semilattice and whose arrows are the meet-semilattices homomorphisms between them which preserve the $0$, send supercompact elements to supercompact elements and preserve arbitrary joins of supercompact elements. 

On a weakly supercompact meet-semilattice $\cal M$ (regarded as a cartesian preorder category) one can define the \emph{supercompactly generated topology} $J^{sc}_{\cal M}$, as follows: for any sieve $S$ in $\cal M$ on an object $m\in {\cal M}$, we set $S\in J^{sc}_{\cal M}(m)$ if and only if either $S$ is the maximal sieve or $S$ is a covering sieve in $\cal C$ (i.e. $m$ is the supremum in $\cal M$ of all the domains of the arrows in $S$) generated by a finite family $\{a_{i}\leq m \textrm{ | } i\in I\}$ where the $a_{i}$ are all supercompact elements of $\cal M$. Clearly, on an infinitarily weakly supercompact meet-semilattice we can define the \emph{infinitary supercompactly generated topology}, as the (obvious) infinitary analogue of the supercompactly generated topology.   

Note that all the Grothendieck topologies defined above are subcanonical.

The following theorem shows that, for several naturally arising subcanonical sites $({\cal C}, J)$, including the ones introduced above, there are topos-theoretic invariants which are $({\cal C}, J)$-adequate.   

\begin{theorem}\label{thmadequate}
Let $({\cal C}, J)$ be a site. Then

\begin{enumerate}[(i)]
\item If $\cal C$ is a coherent category and $J$ is the coherent topology on $\cal C$ then the invariant `to be finite' is $({\cal C}, J)$-adequate;

\item If $\cal C$ is a regular category and $J$ is the regular topology on $\cal C$ then the invariant `to be a singleton' is $({\cal C}, J)$-adequate;

\item If $\cal C$ is any category and $J$ is the trivial topology on $\cal C$ then the invariant `to be a singleton' is $({\cal C}, J)$-adequate;

\item If $\cal C$ is a $k$-geometric category (for a regular cardinal $k$) and $J$ is the $k$-covering topology on $\cal C$ then the invariant `to have cardinality $\lt k$' is $({\cal C}, J)$-adequate;   

\item If $\cal C$ is a disjunctively distributive lattice (resp. a disjunctively distributive frame) and $J$ is the disjunctive topology (resp. the infinitary disjunctive topology) on $\cal C$ then the invariant `to be finite and disjoint' (resp. `to be disjoint') is $({\cal C}, J)$-adequate;

\item If $\cal C$ is a weakly atomic meet-semilattice (resp. an infinitarily weakly atomic meet-semilattice) and $J$ is the atomically generated topology (resp. the infinitary atomically generated topology) on $\cal C$  then the invariant `to be a singleton or to be finite and consisting of atomic subobjects' (resp. `to be a singleton or to consist of atomic subobjects') is $({\cal C}, J)$-adequate;

\item If $\cal C$ is a weakly supercompact meet-semilattice (resp. an infinitarily weakly supercompact meet-semilattice) and $J$ is the supercompactly generated topology (resp. the infinitary atomically generated topology) on $\cal C$  then the invariant `to be a singleton or to be finite and consisting of supercompact subobjects' (resp. `to be a singleton or to consist of supercompact subobjects') is $({\cal C}, J)$-adequate.

\end{enumerate}

\end{theorem}

\begin{proofs}
$(i)$ Given a family ${\cal F}=\{(c_{i})\downarrow_{J} \textrm{ | } i\in I\}$ of principal $J$-ideals on $\cal C$, if this family has a finite subcover $\{(c_{i})\downarrow_{J} \textrm{ | } i\in I_{0}\}$, where $I_{0}$ is a finite subset of $I$, we can construct a $J$-covering sieve $S$ on $\cal C$ such that the $c_{i}$ (for $i\in I_{0}$) are exactly the domains $dom(f)$ of a family of arrows $f$ generating the sieve $S$. Indeed, for any $i\in I_{0}$, consider the union $e\mono 1_{\cal C}$ of the subobjects arising as the monic parts of the cover-mono factorizations of the unique arrows $c_{i}\to 1_{\cal C}$ to the terminal object $1_{\cal C}$ of $\cal C$; we have, for each $i\in I_{0}$, a canonical arrow $c_{i}\to e$, and if we define $S$ to be the sieve on $e$ generated by these arrows $c_{i}\to e$ then $S$ is $J$-covering (by definition of coherent topology) and satisfies the condition in the definition of $({\cal C}, J)$-adequate invariant. 

Conversely, if $\cal F$ satisfies the condition in the definition of $({\cal C}, J)$-adequate invariant then clearly $\cal F$ has a finite subcover.    

$(ii)$, $(iii)$ and $(iv)$ follow from an entirely analogous argument to that for $(i)$.

$(v)$ We only prove the finitary version of the result, the infinitary one being entirely analogous to it. Given a family ${\cal F}$ of principal $J$-ideals on $\cal C$, if this family has a finite disjunctive refinement $\{I_{k} \textrm{ | } k\in K\}$ then for any distinct $k,k'\in K$, $I_{k}\cap I_{k'}=\{0_{\cal C}\}$, and for any $k\in K$ there exists $c_{k}$ in $\cal C$ such that $(c_{k})\downarrow_{J}$ belongs to $\cal F$ and $I_{k}\subseteq (c_{k})\downarrow_{J}$. First, we note that for every $k\in K$, $I_{k}$ is a principal $J$-ideal. Indeed, since $c_{k}\in \mathbin{\mathop{\textrm{\huge $\vee$}}\limits_{}{\cal F}}=\mathbin{\mathop{\textrm{\huge $\vee$}}\limits_{k\in K}I_{k}}$ there exists a $J$-covering sieve $S:=\{a_{w}^{k} \mono c_{k} \textrm{ | } w\in W_{k}\}$ on $c_{k}$ made of pairwise disjoint subobjects of $c_{k}$ such that every $a_{w}^{k}$ belongs to $\mathbin{\mathop{\textrm{\huge $\cup$}}\limits_{k\in K}I_{k}}$. By definition of $J$, $c_{k}$ is equal to the disjoint join $\mathbin{\mathop{\textrm{\huge $\vee$}}\limits_{w\in W_{k}}a^{k}_{w}}$ in $\cal C$. Now, if $b\in I_{k}$ then $b\leq c_{k}$ (since $I_{k}\subseteq (c_{k})\downarrow_{J}=(c_{k})\downarrow$) and hence $b=b\wedge c_{k}=\mathbin{\mathop{\textrm{\huge $\vee$}}\limits_{w\in W_{k}}(b\wedge a^{k}_{w})}=\mathbin{\mathop{\textrm{\huge $\vee$}}\limits_{w\in W_{k} \textrm{ | } a_{w}^{k}\in I_{k}}(b\wedge a^{k}_{w})}=b\wedge ( \mathbin{\mathop{\textrm{\huge $\vee$}}\limits_{w\in W_{k} \textrm{ | } a_{w}^{k}\in I_{k}} a^{k}_{w}})$ (since the $I_{k}$ are pairwise disjoint). But, $I_{k}$ being a $J$-ideal, the element $a_{k}:=\mathbin{\mathop{\textrm{\huge $\vee$}}\limits_{w\in W_{k} \textrm{ | } a_{w}^{k}\in I_{k}} a^{k}_{w}}$ belongs to $I_{k}$, from which it follows that $I_{k}=(a_{k})\downarrow=(a_{k})\downarrow_{J}$ is a principal $J$-ideal. Now, the fact that $l_{{\cal C}}^{J}:{\cal C}\to \Sh({\cal C}, J)$ is full and faithful, preserves the initial object and finite meets, combined with the fact that $(c)\downarrow_{J}=l_{{\cal C}}^{J}(c)$ for any $c\in {\cal C}$, ensures that the sieve $S$ satisfies the condition in the definition of $({\cal C}, J)$-adequate invariant with respect to the family $\cal F$. 

The fact that if a family $\cal F$ of principal $J$-ideals on $\cal C$ satisfies the condition in the definition of $({\cal C}, J)$-adequate invariant then $\cal F$ has a disjunctive refinement follows similarly from the same arguments. 

$(vi)$ We only prove the finitary version of the result, the infinitary one being entirely analogous. Given a family ${\cal F}$ of principal $J$-ideals on $\cal C$, if this family has a finite atomic refinement then either $\cal F$ is a singleton $(c)\downarrow_{J}$, in which case the maximal sieve on $c$ satisfies the condition of $({\cal C}, J)$-adequate invariant with respect to the family $\cal F$, or there exists a family  $\{I_{k} \textrm{ | } k\in K\}$ of $J$-ideals such that for any $k\in K$ $I_{k}$ is an atom in $\Sh({\cal C}, J)$ and $I_{k}\subseteq (c_{k})\downarrow_{J}$ for some $c_{k}\in {\cal C}$ such that $(c_{k})\downarrow_{J}$ belongs to $\cal F$. First, we observe that for every $k\in K$, $I_{k}$ is a principal $J$-ideal; indeed, this follows immediately from the fact that $I_{k}$ is an atom, since every $J$-ideal can be expressed as a union of principal ideals. Now, suppose that, for any $k\in K$, $I_{k}=(a_{k})\downarrow_{J}$. It is immediate to see (by using the fact that $l_{{\cal C}}^{J}:{\cal C}\to \Sh({\cal C}, J)$ is full and faithful and preserves the initial object) that for any object $c$ in $\cal C$, $l_{{\cal C}}^{J}(c)=(c)\downarrow_{J}$ is an atom in $\Sh({\cal C}, J)$ if and only if $c$ is an atom in $\cal C$. Thus the elements $a_{k}$ are atoms in $\cal C$; so the supremum $a$ of the $a_{k}$ in $\cal C$ exists and the sieve $\{a_{k}\mono a \textrm{ | } k\in K \}$ is $J$-covering and satisfies the condition in the definition of $({\cal C}, J)$-adequate invariant with respect to the family $\cal F$. 

The converse implication is clear for the same reasons. 

$(vii)$ The proof is entirely analogous to that of $(vi)$, the key point being that, under the hypotheses of the theorem, $l_{{\cal C}}^{J}(c)=(c)\downarrow_{J}$ is a supercompact object of $\Sh({\cal C}, J)$ if and only if $c$ is an supercompact object of $\cal C$.  

\end{proofs}

Note that in fact part $(iv)$ of Theorem \ref{thmadequate} subsumes both parts $(i)$ and $(ii)$ of the theorem (obtained respectively as the particular cases of $(iv)$ for $k=\omega$ and $k=2$).

\begin{corollary}\label{cor}
Let $({\cal C}, J)$ be a site.

\begin{enumerate}[(i)]
\item If $\cal C$ is a coherent category and $J$ is the coherent topology on $\cal C$ then a $J$-ideal $I$ in $\cal C$ is principal if and only if it is compact in $\Sh({\cal C}, J)$;

\item If $\cal C$ is a regular category and $J$ is the regular topology on $\cal C$ then a $J$-ideal $I$ in $\cal C$ is principal if and only if it is supercompact in $\Sh({\cal C}, J)$; 

\item If $\cal C$ is any category and $J$ is the trivial topology on $\cal C$ then a $J$-ideal $I$ in $\cal C$ is principal if and only if it is supercompact in $\Sh({\cal C}, J)$; 

\item If $\cal C$ is a $k$-geometric category, for a regular cardinal $k$, and $J$ is the $k$-covering topology on $\cal C$ then a $J$-ideal $I$ in $\cal C$ is principal if and only if it is $k$-compact in $\Sh({\cal C}, J)$;

\item If $\cal C$ is a disjunctively distributive lattice (resp. a disjunctively distributive frame) and $J$ is the disjunctive topology (resp. the infinitary disjunctive topology) on $\cal C$ then a $J$-ideal $I$ in $\cal C$ is principal if and only if it is disjunctively compact (resp. infinitary disjunctively compact) in $\Sh({\cal C}, J)$; 

\item If $\cal C$ is a a weakly atomic meet-semilattice (resp. an infinitarily weakly atomic meet-semilattice) and $J$ is the atomically generated topology (resp. the infinitary atomically generated topology) on $\cal C$ then a $J$-ideal $I$ in $\cal C$ is principal if and only if it is atomically compact (resp. infinitarily atomically compact) in $\Sh({\cal C}, J)$. 

\item If $\cal C$ is a a weakly supercompact meet-semilattice (resp. an infinitarily weakly supercompact meet-semilattice) and $J$ is the atomically generated topology (resp. the infinitary supercompactly generated topology) on $\cal C$ then a $J$-ideal $I$ in $\cal C$ is principal if and only if it is supercompactly compact (resp. infinitarily supercompactly compact) in $\Sh({\cal C}, J)$. 

\end{enumerate}

\end{corollary}

\begin{proofs}
This immediately follows from Theorem \ref{thmadequate} and Theorem \ref{equivJcomp} in view of Theorem \ref{Jcomp}, Remark \ref{refin} and Examples \ref{excomp}.     
\end{proofs}

We remark that the particular case of Corollary \ref{cor}(i) for distributive lattices (i.e. preorder coherent categories) already appears in the proof of Proposition II 3.2 \cite{stone}. 

The invariants considered in the Corollary have all the property of being $({\cal C}, J_{\cal C})$-adequate for a given class of sites $({\cal C}, J_{\cal C})$ associated to the structures in a certain category $\cal K$ as in section \ref{locales}. This motivates the following definition.

\begin{definition}
Let $\cal K$ be a category of structures as in section \ref{locales}, where each of the structures $\cal C$ in $\cal K$ is equipped with a Grothendieck topology $J_{\cal C}$, and let $C$ be an invariant property of families of subterminals in a topos. 

The invariant $C$ is said to be \emph{$\cal K$-adequate} is it is $({\cal C}, J_{\cal C})$-adequate for every structure $\cal C$ in $\cal K$.
\end{definition}

\begin{examples}\label{excomp2}
\begin{enumerate}[(a)]

\item The invariant `to be finite' is $\textbf{DLat}$-adequate, where $\textbf{DLat}$ is the category of distributive lattices;
\item The invariant `to be a singleton' is $\textbf{MSLat}$-adequate, where $\textbf{MSLat}$ is the category of meet-semilattices; 
\item The invariant `to be of cardinality $\lt k$' (for a regular cardinal $k$) is $k\textrm{-}\textbf{Frm}$-adequate, where $k\textrm{-}\textbf{Frm}$ is the category of $k$-frames;
\item The invariant `to be finite and disjoint' is $\textbf{DJLat}$-adequate, where\\ $\textbf{DJLat}$ is the category of disjunctively distributive lattices;
\item The invariant `to be disjoint' is $\textbf{DJFrm}$-adequate, where $\textbf{DJFrm}$ is the category of disjunctively distributive frames;
\item The invariant `to be a singleton or finite and consisting of atomic subobjects' is $\textbf{WAtMSLat}$-adequate, where $\textbf{WAtMSLat}$ is the category of weakly atomic meet-semilattices;
\item The invariant `to be a singleton or to consist of atomic subobjects' is $\textbf{IWAtMSLat}$-adequate, where $\textbf{WAtMSLat}$ is the category of infinitarily weakly atomic meet-semilattices;
\item The invariant `to be a singleton or finite and consisting of supercompact subobjects' is $\textbf{WScMSLat}$-adequate, where $\textbf{WScMSLat}$ is the category of weakly supercompact meet-semilattices;
\item The invariant `to be a singleton or to consist of supercompact subobjects' is $\textbf{IWScMSLat}$-adequate, where $\textbf{IWScMSLat}$ is the category of infinitarily weakly supercompact meet-semilattices.
\end{enumerate}
\end{examples}

Let $\cal K$ be a category of poset structures, as in section \ref{locales}.
  
If $C$ is a $\cal K$-adequate invariant then, by Theorem \ref{equivJcomp}, every poset structure $\cal C$ in $\cal K$ can be recovered up to isomorphism from the corresponding topos $\Sh({\cal C}, J_{\cal C})$ as the poset of subterminals which are $C$-compact. In fact, if the invariant $C$ satisfies the property that for any structure $\cal C$ in $\cal K$ and for any family $\cal F$ of principal $J_{\cal C}$-ideals on $\cal C$, $\cal F$ has a refinement satisfying $C$ (if and) only if it has a refinement satisfying $C$ made of principal $J_{\cal C}$-ideals on $\cal C$ (cf. Theorem \ref{construction} below) then we can identify a stronger invariant property which is satisfied by all the embeddings ${\cal C} \hookrightarrow Id_{J_{\cal C}}({\cal C})$ (for $\cal C$ in $\cal K$), namely the fact that every covering of a principal $J_{\cal C}$-ideal on $\cal C$ in the topos $\Sh({\cal C}, J_{\cal C})$ is refined by a covering of principal $J_{\cal C}$-ideals on $\cal C$ which satisfy $C$; this follows by an inspection of the proof of Theorem \ref{Jcomp}, recalling the definition of $({\cal C}, J_{\cal C})$-adequate invariant. Therefore, if an embedding $B_{L}\hookrightarrow L$ is isomorphic to any of the embeddings ${\cal C} \hookrightarrow Id_{J_{\cal C}}({\cal C})$, so that $B_{L}$ is the subset of $C$-compact elements of $L$, then we can conclude that it satisfies the property that every covering in $L$ of an element of $B_{L}$ is refined by a covering made of elements of $B_{L}$ which satisfies the invariant $C$. Notice also that all the embeddings $i_{\cal C}:{\cal C} \hookrightarrow Id_{J_{\cal C}}({\cal C})$ satisfy the property that the $J_{\cal C}$-covering sieves are sent by $i_{\cal C}$ to covering families in $Id_{J_{\cal C}}({\cal C})$. 
    
Given a category $\cal K$ and a functor $A:{\cal K}^{\textrm{op}}\to \textbf{Loc}$ as in section \ref{locales}, in order to use a $\cal K$-adequate invariant $C$ for achieving an intrinsic characterization of the locales in the extended image of $A$, and to define a functor on this subcategory which is an inverse (up to isomorphism) to the functor $A$, we need $C$ to satisfy an additional property, namely the fact that for any locale $L$ with a basis $B_{L}$ of $C$-compact elements which, regarded as a preorder with the induced order, belongs to $\cal K$, and such that the embedding $B_{L}\hookrightarrow L$ possibly satisfies some invariant properties that are known to hold for the embeddings ${\cal C} \hookrightarrow Id_{J_{\cal C}}({\cal C})$, the topology $J^{L}_{can}|_{B_{L}}$ coincides with the topology $J_{B_{L}}$ with which $B$ comes equipped as a structure in $\cal K$, where $J^{L}_{can}$ is the canonical topology on $L$. Indeed, if this is the case then, regarding $B_{L}$ as full preorder subcategory of $L$, the Comparison Lemma yields an equivalence $\Sh(L)\simeq \Sh(B_{L}, J^{L}_{can}|_{B_{L}})$, and hence, $J^{L}_{can}|_{B_{L}}$ being equal to $J_{B_{L}}$, $L$ is isomorphic to the image $A(B_{L})$ of $B_{L}$ under the functor $A$. 
 
This motivates the following definition. Below, given a locale $L$, we denote by $J_{can}^{L}$ the canonical Grothendieck topology on $L$.

\begin{definition}
Let $\cal K$ be a category of poset structures as in section \ref{locales}, where each of the structures $\cal C$ in $\cal K$ is equipped with a Grothendieck topology $J_{\cal C}$, and let $P$ be an invariant property of embeddings $B_{L}\hookrightarrow L$ of a basis $B_{L}$ of a frame $L$ into $L$ which holds for all the canonical embeddings ${\cal C}\hookrightarrow Id_{J_{\cal C}}({\cal C})$, where $\cal C$ is a structure in $\cal K$. 

A topos-theoretic invariant property $C$ of families of subterminals in a locally small cocomplete topos is said to be \emph{$\cal K$-compatible relative to $P$} if for any frame $L$ with a basis $B_{L}$ of $C$-compact elements which, regarded as a poset with the induced order, belongs to $\cal K$, if the embedding $B_{L}\hookrightarrow L$ satisfies 
\begin{enumerate}[(i)]
\item property $P$;
\item the property that every covering in $L$ of an element of $B_{L}$ is refined by a covering satisfying $C$ made of elements of $B_{L}$, and
\item the property that the $J_{B_{L}}$-covering sieves are sent by the embedding $B_{L}\hookrightarrow L$ into covering families in $L$, where $J_{B_{L}}$ is the Grothendieck topology with which $B_{L}$ comes equipped as a structure in $\cal K$ 
\end{enumerate}
then $J_{B_{L}}$ is equal to the induced Grothendieck topology $J_{can}^{L}|_{B_{L}}$ on $B_{L}$.  

An invariant $C$ is said to be \emph{$\cal K$-compatible} if it is $\cal K$-compatible relative to $P$ for some invariant property $P$ as above. 
\end{definition}

It is easy to verify that all the invariants considered above are compatible with respect to the corresponding categories of structures. 

Anyway, there is a systematic way for building, given a category $\cal K$ of poset structures as in section \ref{locales}, invariants which are both $\cal K$-adequate and $\cal K$-compatible. This method is based on a formalization of the vague idea of a Grothendieck topology defined `through a topos-theoretic invariant', as given by the following definition.

\begin{definition}\label{inducedtop}
Let $\cal K$ be a category of poset structures as in section \ref{locales}, where each of the structures $\cal C$ in $\cal K$ is equipped with a Grothendieck topology $J_{\cal C}$, let $P$ be an invariant property of embeddings $B_{L}\hookrightarrow L$ of a basis $B_{L}$ of a frame $L$ into $L$ which holds for all the canonical embeddings ${\cal C}\hookrightarrow Id_{J_{\cal C}}({\cal C})$ (where $\cal C$ is a structure in $\cal K$), and $C$ be a topos-theoretic invariant property of families of subterminals in a locally small cocomplete topos. 

Given a structure $\cal C$ in $\cal K$, the Grothendieck topology $J_{\cal C}$ is said to be \emph{$C$-induced relative to $P$} if for any $J_{can}^{L}$-dense monotone embedding $i:{\cal C} \hookrightarrow L$ into a frame $L$ which satisfies property $P$ and such that the $J_{\cal C}$-covering sieves on $\cal C$ are sent by $i$ to covering families in $L$, for any family $\cal A$ of objects in $\cal C$ there exists a $J_{\cal C}$-covering sieve on an object $c\in {\cal C}$ such that the arrows $a\mono c$ for $a\in {\cal A}$ generate $S$ if and only if the image $i(\cal A)$ of the family $\cal A$ in $L$ has a refinement satisfying $C$ made of objects of the form $i(c')$ (for $c'\in {\cal C}$).

The Grothendieck topology $J_{\cal C}$ is said to be \emph{$C$-induced} if it is $C$-induced relative to $P$ for some invariant property $P$ as above.  
\end{definition}

\begin{remark}\label{chargrot}
Under the hypotheses of Definition \ref{inducedtop}, if $J_{\cal C}$ is $C$-induced then the $J_{\cal C}$-covering sieves on an object $c\in {\cal C}$ can be characterized precisely as the sieves on $c$ in $\cal C$ which contain a family of arrows $\cal A$ with codomain $c$ such that the family $\{i(dom(f)) \textrm{ | } f\in {\cal A}\}$ satisfies $C$.  
\end{remark}

Let us give some examples of topologies induced (in our sense) by topos-theoretic invariants.

\begin{theorem}\label{induced}

Let $({\cal C}, J)$ be a site. Then
\begin{enumerate}[(i)]
\item If $\cal C$ is a distributive lattice and $J$ is the coherent topology on $\cal C$ then $J$ is $C$-induced where $C$ is the property `to be a finite family';

\item If $\cal C$ is any poset and $J$ is the trivial topology on $\cal C$ then $J$ is $C$-induced where $C$ is the property `to be a singleton family';

\item If $\cal C$ is a $k$-frame (for a regular cardinal $k$) and $J$ is the $k$-covering topology on $\cal C$ then $J$ is $C$-induced where $C$ is the property `to be a $k$-family (i.e. a family of cardinality $\lt k$)'; 

\item If $\cal C$ is a disjunctively distributive lattice (resp. a disjunctively distributive frame) and $J$ is the disjunctive topology (resp. the infinitary disjunctive topology) on $\cal C$ then $J$ is $C$-induced where $C$ is the property `to be a finite and disjont (i.e. whose elements are pairwise disjoint subobjects) family' (resp. `to be a disjoint family'); 

\item If $\cal C$ is a weakly atomic meet-semilattice (resp. an infinitarily weakly atomic meet-semilattice) and $J$ is the atomically generated topology (resp. the infinitary atomically generated topology) on $\cal C$ then $J$ is $C$-induced where $C$ is the property `to be a singleton family or a  finite and atomic (i.e. whose elements are atoms) family' (resp. `to be singleton family or an atomic family'); 

\item If $\cal C$ is a weakly supercompact meet-semilattice (resp. an infinitarily weakly supercompact meet-semilattice) and $J$ is the supercompactly generated topology (resp. the infinitary atomically generated topology) on $\cal C$ then $J$ is $C$-induced where $C$ is the property `to be a singleton family or a finite and supercompact (i.e. whose elements are supercompact subobjects) family' (resp. `to be a singleton family or a supercompact family'). 
\end{enumerate}
\end{theorem}

\begin{proofs}
$(i)$, $(ii)$ and $(iii)$. We prove $(iii)$, since $(i)$ and $(ii)$ are particular cases of it. The topology $J$ is clearly $C$-induced relative to $P$ where $P$ is the vacuous property.

$(iv)$ The topology $J$ is easily seen to be $C$-induced relative to $P$ where $P$ is the following property of embeddings $i:B_{L}\hookrightarrow L$: `$i$ preserves finite meets' or equivalently (since the canonical embedding ${\cal C}\hookrightarrow Id_{J}({\cal C})$ satisfies the property that ${\cal C}$ is closed under finite meets in $Id_{J}({\cal C})$) the property `$i(B_{L})$ is closed in $L$ under finite meets'. 

$(v)$ The topology $J$ is easily seen to be $C$-induced relative to $P$ where $P$ is the vacuous property.

$(vi)$ This is entirely analogous to $(v)$; again, $P$ can be taken to be the vacuous property. 
\end{proofs}   

\begin{theorem}\label{construction}
Let $\cal K$ be a category of poset structures as in section \ref{locales}. If all the Grothendieck topologies $J_{\cal C}$ associated to the structures $\cal C$ in $\cal K$ are $C$-induced (for a topos-theoretic invariant $C$ of families of subterminals in a locally small cocomplete topos) and the invariant $C$ satisfies the property that for any structure $\cal C$ in $\cal K$ and for any family $\cal F$ of principal $J_{\cal C}$-ideals on $\cal C$, $\cal F$ has a refinement satisfying $C$ (if and) only if it has a refinement satisfying $C$ made of principal $J_{\cal C}$-ideals on $\cal C$  then the invariant $C$ is both $\cal K$-adequate and $\cal K$-compatible.
\end{theorem}

\begin{proofs}
Let us first show that the invariant $C$ is $\cal K$-adequate. Given a structure $\cal C$ in $\cal K$, consider the canonical embedding ${\cal C}\hookrightarrow Id_{J_{\cal C}}({\cal C})$. Given a family $\cal F$ of principal $J_{\cal C}$-ideals on $\cal C$, we have to prove that $\cal F$ has a refinement satisfying $C$ if and only if there exists a $J$-covering sieve $\{f_{h}:d_{h}\to d \textrm{ | } h\in H\}$ in $\cal C$ such that for any $h\in H$, $(d_{h})\downarrow_{J}$ is contained in some ideal in $\cal F$ and the union in $\cal E$ of the family $\{(d_{h})\downarrow_{J} \textrm{ | } h\in H\}$ is equal to the union in $\cal E$ of the family $\cal F$. Now, by our hypothesis on $C$, $\cal F$ has a refinement satisfying $C$ if and only if $\cal F$ has a refinement satisfying $C$ made of principal $J_{\cal C}$-ideals on $\cal C$, and, $J_{{\cal C}}$ being $C$-induced, this latter condition is equivalent to the existence of a $J$-covering sieve $\{f_{h}:d_{h}\to d \textrm{ | } h\in H\}$ with the required property.   

To prove that $C$ is $\cal K$-adequate, we show that $C$ is $\cal K$-adequate relative to $P$, where $P$ is the intersection of all the properties witnessing the fact that the topologies $J_{\cal C}$ are $C$-induced, so that all the Grothendieck topologies $J_{\cal C}$ are $C$-induced relative to $P$. Let $L$ be a frame with a basis $B_{L}$ of $C$-compact elements which, regarded as a preorder with the induced order, is equal to a structure $\cal C$ in $\cal K$; suppose moreover that the embedding $i:{\cal C}=B_{L}\hookrightarrow L$ satisfies property $P$, the property that every covering in $L$ of an element of $B_{L}$ is refined by a covering made of elements of $B_{L}$ which satisfies the invariant $C$, and the property that the $J_{B_{L}}$-covering sieves are sent by the embedding $B_{L}\hookrightarrow L$ into covering families in $L$ (where $J_{B_{L}}$ is the Grothendieck topology with which $B_{L}$ comes equipped as a structure in $\cal K$); we have to prove that $J_{\cal C}=J^{L}_{can}|_{\cal C}$, where $J^{L}_{can}$ is the canonical topology on $L$. Now, a sieve $S$ on an object $c\in {\cal C}$ is $J^{L}_{can}|_{\cal C}$-covering if and only if the sieve in $L$ generated by the image $i(S)$ of the arrows in $S$ under the embedding $i$ generates a $J^{L}_{can}$-covering sieve on $i(c)$ in $L$; but, by our hypotheses on the embedding $i$, this condition holds if and only if there exists a family of arrows of the form $\{i(c')\leq i(c)\}$ (for $c'\in {\cal C}$) which refines the family $\{i(dom(f))\leq i(c) \textrm{ | } f\in S\}$ and satisfies the invariant $C$. On the other hand, by Remark \ref{chargrot}, for any object $c\in {\cal C}$, the $J_{\cal C}$-covering sieves on $c$ are precisely those which contain a family of arrows $\cal A$ in $\cal C$ with codomain $c$ such that the family $\{i(dom(f)) \textrm{ | } f\in {\cal A}\}$ satisfies $C$. From this, the equality $J_{\cal C}=J^{L}_{can}|_{\cal C}$ is clear.  
\end{proofs}

\begin{remark}\label{rmkprincipal}
All the invariants $C$ appearing in the statement of Theorem \ref{induced} satisfy the property in the statement of Theorem \ref{construction}. This is clear for the invariants in parts $(i)$, $(ii)$ and $(iii)$ of the theorem, and it follows from the proofs of the relevant parts of Theorem \ref{thmadequate} for the other ones.
\end{remark}

Theorem \ref{construction} is important because it assures that, whenever the subcanonical Grothendieck topologies $J_{\cal C}$ corresponding to the structures in a given category $\cal K$ are `uniformly defined through a topos-theoretic invariant' (in the sense of being $C$-induced for an invariant $C$ satisfying the hypotheses of Theorem \ref{construction}), there is a topos-theoretic invariant, namely the property of $C$-compactness, which enables us to recover each structure $\cal C$ from the corresponding topos $\Sh({\cal C}, J_{\cal C})$ up to isomorphism, and to define a functor $I_{A}$ on the extended image of the functor $A:{\cal K}^{\textrm{op}}\to \textbf{Loc}$ of section \ref{locales} which is a categorical inverse to $A$ and hence yields a duality between $\cal K$ and the subcategory of $\textbf{Loc}$ given by the extended image of $A$. Recall that $I_{A}$ is defined as the functor which acts on the objects by sending a locale $L$ in $ExtIm(A)$ to the poset of $C$-compact elements of $L$ and acts on the arrows by sending a locale morphism $f:L\to L'$ in $ExtIm(A)$ to the restriction of its associated frame homomorphism to the subsets of $C$-compact elements of $L$ and $L'$ (cf. section \ref{locales}). Summarizing, we have established the following general `duality theorem'.

\begin{theorem}\label{thmcentral}
Let $\cal K$ be a category of poset structures as in section \ref{locales}, and let $C$ be a topos-theoretic invariant of families of subterminals in a locally small cocomplete topos which satisfies the property that for any structure $\cal C$ in $\cal K$ and for any family $\cal F$ of principal $J_{\cal C}$-ideals on $\cal C$, $\cal F$ has a refinement satisfying $C$ (if and) only if it has a refinement satisfying $C$ made of principal $J_{\cal C}$-ideals on $\cal C$. If all the Grothendieck topologies $J_{\cal C}$ associated to the structures $\cal C$ in $\cal K$ are $C$-induced then the functor $A:{\cal K}^{\textrm{op}}\to \textbf{Loc}$ of section \ref{locales} admits a categorical inverse $I_{A}:ExtIm(A)\to {\cal K}^{\textrm{op}}$ sending a locale $L$ in $ExtIm(A)$ to the poset of $C$-compact elements of $L$ and a locale morphism $f:L\to L'$ in $ExtIm(A)$ to the restriction of its associated frame homomorphism to the subsets of $C$-compact elements of $L$ and $L'$. 
\end{theorem}

The following result provides an `intrinsic' characterization of the category $ExtIm(A)$.

\begin{theorem}\label{propext}
Let us assume that the hypotheses of Theorem \ref{thmcentral} are satisfied, and that $C$ is $\cal K$-compatible relative to a property $P$. Then
\begin{enumerate}[(i)]
\item The locales in the extended image $ExtIm(A)$ of the functor $A:{\cal K}^{\textrm{op}}\to \textbf{Loc}$ of section \ref{locales} are precisely the locales $L$ with a basis $B_{L}$ of $C$-compact elements which, regarded as a poset with the induced order, belongs to $\cal K$, and such that the embedding $B_{L}\hookrightarrow L$ satisfies property $P$, the property that every covering in $L$ of an element of $B_{L}$ has a refinement by a covering made of elements of $B_{L}$ which satisfies the invariant $C$, and the property that the $J_{B_{L}}$-covering sieves are sent by the embedding $B_{L}\hookrightarrow L$ into covering families in $L$, where $J_{B_{L}}$ is the Grothendieck topology with which $B_{L}$ comes equipped as a structure in $\cal K$;
\item The arrows $L\to L'$ in $ExtIm(A)$ are the locale morphisms $L\to L'$ between locales in $ExtIm(A)$ whose associated frame homomorphisms send $C$-compact elements of $L'$ to $C$-compact elements of $L$ in such a way that their restriction to the subsets of $C$-compact elements of $L'$ and $L$ can be identified with an arrow in $\cal K$.
\end{enumerate}  
\end{theorem}

\begin{proofs}
The characterization of the objects of $ExtIm(A)$ follows immediately from the definition of $\cal K$-compatible invariant, in light of Theorem \ref{construction}.

To characterize the arrows in $ExtIm(A)$ we observe that, by Lemma C2.3.8 \cite{El}, since the Grothendieck topologies $J_{\cal C}$ are subcanonical, the geometric morphisms $l:\Sh({\cal D}, J_{{\cal D}}) \to \Sh({\cal C}, J_{{\cal C}})$ which are induced (uniquely up to isomorphism) by a morphism of sites $({\cal C}, J_{{\cal C}}) \to ({\cal D}, J_{{\cal D}})$ can be characterized precisely as those such that the inverse image $l^{\ast}$ sends representables on $\cal C$ to representables on $\cal D$. Now, under the equivalences $\Sh({\cal C}, J_{\cal C}) \simeq \Sh(Id_{J_{\cal C}}({\cal C}))$ and $\Sh({\cal D}, J_{\cal D}) \simeq \Sh(Id_{J_{\cal D}}({\cal D}))$ of Theorem \ref{fund}, these geometric morphisms correspond exactly to the frame homomorphisms ${\cal O}(Id_{J_{\cal C}}({\cal C})) \to {\cal O}(Id_{J_{\cal D}}({\cal D}))$ which send principal ($J_{\cal C}$-)ideals on $\cal C$ to principal ($J_{\cal D}$-)ideals on $\cal D$. Therefore, if $C$ is $\cal K$-compatible then the arrows in $\cal U$ can be characterized precisely as the locale morphisms whose associated frame homomorphisms send $C$-compact elements to $C$-compact elements in such a way that their restriction to the subsets of $C$-compact elements can be identified with an arrow in $\cal K$.         
\end{proofs}

Note that if, for any structures $\cal C$ and $\cal D$ in $\cal K$, the arrows ${\cal C}\to {\cal D}$ in $\cal K$ coincide exactly with the morphisms of sites $({\cal C}, J_{\cal C})\to ({\cal D}, J_{\cal D})$ then the condition `in such a way that their restriction to the subsets of $C$-compact elements of $L'$ and $L$ can be identified with an arrow in $\cal K$' in the statement of Theorem \ref{propext} can be omitted, since it is automatically satisfied (cf. section \ref{locales} above). 

A similar characterization holds for the extended images of the covariant functors $B:{\cal K}\to \textbf{Loc}$ defined in section \ref{locales} (cf. section \ref{prealex} below).

\begin{remark}
A basis $B_{L}\hookrightarrow L$ satisfying the hypotheses of Theorem \ref{propext} is unique if it exists. Indeed, $C$ being $\cal K$-adequate, for any $\cal C$ in $\cal K$ the elements of $Id_{J_{\cal C}}({\cal C})$ which are principal $J_{\cal C}$-ideals on $\cal C$ are precisely those which are $C$-compact, and hence $B_{L}$ must be equal to the subset of $L$ consisting of its $C$-compact elements. 
\end{remark}

The following proposition is often useful for achieving explicit descriptions of the condition in Theorem \ref{propext} that a basis $B_{L}$ of a locale $L$ in $ExtIm(A)$ should have the structure of a poset in $\cal K$ (cf. section \ref{ex} below for some concrete applications of this result).

\begin{proposition}\label{multicomposition}
Let $C$ be a topos-theoretic invariant property of families of subterminals in a locally small cocomplete topos which is multicomposition-stable in the sense that, for any family ${\cal F}:=\{a_{i} \textrm{ | } i\in I\}$ of subterminals in a locally small cocomplete topos $\cal E$ and any collection of families ${\cal F}_{i}:=\{b_{ij} \leq a_{i} \textrm{ | } j\in I_{j}\}$ of subterminals in $\cal E$ indexed over a set $I$, if $\cal F$ and all the ${\cal F}_{i}$ (for $i\in I$) satisfy $C$ then the `multicomposite' family ${\cal F} \ast \{{\cal F}_{i} \textrm{ | } i\in I\}:=\{b_{ij}, i\in I, j\in I_{j}\}$ also satisfies $C$. Then for any family $\cal G$ of subterminals in a locally small cocomplete topos $\cal E$ which satisfies $C$, if all the objects in $\cal G$ are $C$-compact then the union in $\cal E$ of all the subterminals in $\cal G$ is again $C$-compact.  
\end{proposition}  

\begin{proofs}
Let ${\cal G}=\{c_{i} \textrm{ | }\in I\}$ be a family of subterminals in a locally small cocomplete topos $\cal E$ such that all the $c_{i}$ (for $i\in I$) are $C$-compact, and let $c=\mathbin{\mathop{\textrm{\huge $\vee$}}\limits_{i\in I}c_{i}}$ be the union in $\cal E$ of all the subterminals in $\cal G$. Given a covering $\{u_{k} \textrm{ | }k\in K\}$ of $c$ in $\cal E$, we have, for any $i\in I$, $c_{i}=\mathbin{\mathop{\textrm{\huge $\vee$}}\limits_{k\in K}(c_{i}\wedge u_{k})}$. This implies, since the $c_{i}$ are $C$-compact, that there exists, for each $i\in I$, a family ${\cal G}_{i}$ satisfying $C$ which refines the family $\{c_{i}\wedge u_{k} \textrm{ | } k\in K\}$. Then the fact that $C$ is multicomposition-stable implies that the family ${\cal G} \ast \{{\cal G}_{i} \textrm{ | } i\in I\}$ satisfies $C$; but this family clearly refines the family $\{u_{k} \textrm{ | }k\in K\}$, from which our thesis follows.   
\end{proofs}

Note that all the invariants that we have considered in this section are multicomposition-stable.
  
\section{Old and new dualities}\label{ex}

In this section, we give some applications of the `machinery' for generating dualities or equivalences described in the previous sections. Specifically, we discuss how well-known dualities can be naturally recovered from our machinery and we establish new dualities and equivalences as applications of our method.

\subsection{Preorders and Alexandrov topologies}\label{prealex}

The well-known equivalence between preorders and Alexandrov spaces falls under the category of covariant equivalences obtainable by the technique of sections \ref{locales} and \ref{dualtop}. 

Let $\textbf{Pro}$ be the category of preorders and monotone maps between them. The construction of covariant equivalences given in section \ref{locales} yields a functor $B:\textbf{Pro} \to \textbf{Loc}$, which assigns to a preorder $\cal P$ the frame $Id({{\cal P}^\textrm{op}})$ of upper sets in $\cal P$ and to an arrow $f:{\cal P}\to {\cal Q}$ in $\textbf{Pro}$ the frame homomorphism $B(f):Id({{\cal Q}^\textrm{op}}) \to Id({{\cal P}^\textrm{op}})$ sending an upper set $S$ in $\cal Q$ to the upper set in $\cal P$ given by the inverse image $f^{-1}(S)$ of $S$ under $f$. 

We shall call the image $B({\cal P})=Id({{\cal P}^\textrm{op}})$ of a preorder $\cal P$ under the functor $B$ the \emph{Alexandrov locale} associated to $\cal P$. 

By Corollary \ref{cor}(iii), any poset $\cal P$ is isomorphic in $\textbf{Pro}$ to the opposite of the set of supercompact subterminals of the topos $[{\cal P}, \Set]$ endowed with the order induced by that in $\Sub_{[{\cal P}, \Set]}(1)$. In fact, this can be made functorial, as follows. Given posets $\cal P$ and $\cal Q$, and an arrow $f:{\cal P}\to {\cal Q}$ in $\textbf{Pro}$, denoted by $E(f):[{\cal P}, \Set] \to [{\cal Q}, \Set]$ the geometric morphism induced by $f$ as in Example A4.1.4 \cite{El}, we have (cf. the proof of Lemma A4.1.5 \cite{El}) that the left adjoint to the inverse image functor $E(f)^{\ast}:[{\cal Q}, \Set] \to [{\cal P}, \Set]$ of $E(f)$ restricts to the subsets of supercompact subterminals of $[{\cal Q}, \Set]$ and of $[{\cal P}, \Set]$ yielding a map which is isomorphic to the opposite $f^{\textrm{op}}:{\cal P}^{\textrm{op}}\to {\cal Q}^{\textrm{op}}$ of the map $f$. 

Let $\textbf{Pos}$ denote the full subcategory of $\textbf{Pro}$ on the posets. Then the argument above ensures that the restriction $B_{\textbf{Pos}}:\textbf{Pos} \to \textbf{Loc}$ of $B:\textbf{Pro} \to \textbf{Loc}$ is faithful (and injective on objects) and hence yields an isomorphism of categories between $\textbf{Pos}$ and a subcategory of the category $\textbf{Loc}$ of locales, namely the \emph{image} of $B_{\textbf{Pos}}$. This represents the analogue for covariant functors of Theorem \ref{teoabstr}. 

Anyway, $B_{\textbf{Pos}}:\textbf{Pos} \to \textbf{Loc}$ also yields an equivalence of categories onto its extended image $ExtIm(B_{\textbf{Pos}})=ExtIm(B)$, as follows. 

The inverse functor $I_{B}:ExtIm(B) \to \textbf{Pos}$ sends a locale $L$ in $ExtIm(B)$ to the opposite of the poset of supercompact elements of ${\cal O}(L)$ and a morphism $g:L\to L'$ of locales in $ExtIm(B)$ to the opposite of the restriction of the left adjoint $g_{!}:{\cal O}(L) \to {\cal O}(L')$ to $g^{\ast}:{\cal O}(L') \to {\cal O}(L)$ to the posets of supercompact elements of the frames ${\cal O}(L)$ and ${\cal O}(L')$. Note that the left adjoint $g_{!}$ has the following expression in terms of $g^{\ast}$: $g_{!}(l)=\mathbin{\mathop{\textrm{\huge $\wedge$}}\limits_{l\leq g^{\ast}(l')}l'}$. 

We can read the equivalence $\textbf{Pos} \simeq ExtIm(B)$ just established more naturally by composing it with the isomorphism of categories $\textbf{Pos}\cong \textbf{Pos}$ sending every poset $\cal P$ in $\textbf{Pos}$ to the opposite poset ${\cal P}^{\textrm{op}}$, and every arrow $f:{\cal P}\to {\cal Q}$ to the opposite arrow $f^{\textrm{op}}:{\cal P}^{\textrm{op}}\to {\cal Q}^{\textrm{op}}$. The resulting equivalence $\textbf{Pos} \simeq ExtIm(B)$ sends a poset $\cal P$ in $\textbf{Pos}$ to its ideal completion $Id({\cal P})$ (i.e., the set of all the ideals of $\cal P$ with the subset-inclusion order) and a locale $L$ in $ExtIm(B)$ to the subset of supercompact elements of $L$ endowed with the induced natural order.    

Let us now turn to the problem of giving an intrinsic characterization of the category $ExtIm(B)$. We shall call the locales in $ExtIm(B)$ the \emph{Alexandrov locales}.

\begin{proposition}\label{alexlocale}
Let $L$ be a locale. Then the following conditions are equivalent:

\begin{enumerate}[(i)]
\item $L$ is an Alexandrov locale;

\item $L$ has a basis of supercompact elements;

\item $\Sh(L)$ is equivalent to a presheaf topos of the form $[{\cal P}, \Set]$, where $\cal P$ is a preorder (equivalently, a poset) category.
\end{enumerate}
   
\end{proposition}

\begin{proofs}
$(i)\imp (ii)$ It is clear that if $L$ is an Alexandrov locale then it has a basis of supercompact elements; indeed, this property is a locale-theoretic invariant, and if $L$ is of the form $Id({{\cal P}^\textrm{op}})$ for a preorder $\cal P$ then the collection of all the principal upper sets generated by an element of $\cal P$ forms a basis of $L$.

$(ii)\imp (iii)$ If $L$ has a basis $B$ of supercompact elements then $B$, regarded as a full subcategory of $L$, is $J^{L}_{can}$-dense (where $J^{L}_{can}$ is the canonical Grothendieck topology on $L$), and hence $\Sh(L)\simeq \Sh(B, J^{L}_{can}|_{B})$ by the Comparison Lemma. But since every element of $B$ is supercompact then $J^{L}_{can}|_{B}$ is the trivial Grothendieck topology on $\cal B$ and hence $\Sh(L)\simeq \Sh(B, J^{L}_{can}|_{B})\simeq [B^{\textrm{op}}, \Set]$.

$(iii) \imp (i)$ If $\Sh(L)\simeq [{\cal P}, \Set]$ then, since $[{\cal P}, \Set]\simeq \Sh(Id({{\cal P}^\textrm{op}})$ (cf. Theorem \ref{fund}), $L$ is isomorphic to $Id({{\cal P}^\textrm{op}})$ and therefore it is an Alexandrov locale.   
\end{proofs} 

Concerning the arrows in $ExtIm(B)$, it is clear from the discussion preceding Proposition \ref{alexlocale} that they are precisely the morphisms $g:L\to L'$ of locales in $ExtIm(B)$ such that $g^{\ast}:{\cal O}(L') \to {\cal O}(L)$ has a left adjoint $g_{!}:{\cal O}(L) \to {\cal O}(L')$ (equivalently, $g^{\ast}$ preserves arbitrary infima) which sends supercompact elements to supercompact elements. Let us call these morphisms the complete maps. In fact, the condition that $g_{!}$ should send supercompact elements to supercompact elements is superfluous, since it is automatically satisfied (cf. section \ref{comparison} below).

Let $\textbf{AlexLoc}$ denote the category whose objects are the Alexandrov locales and whose arrows are the complete maps between them. We have thus obtained the following result.

\begin{theorem}\label{proalex}
Via the functors $B_{\textbf{Pos}}:\textbf{Pos}\to \textbf{AlexLoc}$ and\\ $I_{B}:\textbf{AlexLoc} \to \textbf{Pos}$ defined above, the category $\textbf{Pos}$ is equivalent to the category $\textbf{AlexLoc}$.
\end{theorem}
 
\begin{remarks}
\begin{enumerate}[(a)]

\item The functor $B_{\textbf{Pos}}$ in the statement of Theorem \ref{proalex} corresponds, under the involution $\textbf{Pos}\to \textbf{Pos}$ sending a poset to its dual, to the functor giving one half of the duality between posets and completely distributive algebraic lattices described in \cite{duality}; in fact, it follows from Theorem \ref{proalex} and the duality theorem in \cite{duality} that the category $\textbf{AlexLoc}$ defined above coincides, when considered as a category of frames, with the category of completely distributive algebraic lattices defined in \cite{duality} (since the two theorems ensure that both these categories are equal to the extended image of the same functor, namely $B_{\textbf{Pos}}:\textbf{Pos}^{\textrm{op}} \to \textbf{Frm}$); anyway, it is worth to note that our definition of the inverse functor $I_{B}:\textbf{AlexLoc} \to \textbf{Pos}$ is completely different from the construction of the inverse functor given in \cite{duality};

\item Clearly, every finite distributive lattice is a frame with a basis of supercompact elements; clearly, the posets corresponding to these frames via the equivalence of Theorem \ref{proalex} are precisely the finite ones. The equivalence of Theorem \ref{proalex} thus restricts to Birkhoff duality between finite distributive lattices and finite posets. 

\end{enumerate}
\end{remarks} 
 
The functor $B:\textbf{Pro} \to \textbf{AlexLoc}$ from preorders to Alexandrov locales can be lifted to an equivalence between $\textbf{Pro}$ and a subcategory of the category $\textbf{Top}$ of topological spaces in various ways, according to the method of section \ref{dualtop}. 

Let us choose, for each preorder ${\cal P}$, as separating set of points of the topos $[{\cal P}, \Set]$, the collection of geometric morphisms $\{e_{p}:\Set\to [{\cal P}, \Set] \textrm{ | } p\in P\}$ indexed by the underlying set $P$ of $\cal P$ as in Example \ref{exa}(b). Clearly, every arrow $f:{\cal P} \to {\cal Q}$ in $\textbf{Pro}$ yields a function $l_{f}:P \to Q$ (given by the action of $f$ on the underlying sets of the preorders) which satisfies the necessary compatibility conditions of section \ref{dualtop}. As we observed in Example \ref{exa}(b), the subterminal topology induced on the underlying set $P$ of a preorder $\cal P$ is precisely the Alexandrov topology on $\cal P$; the method of section \ref{dualtop} thus yields in this case the well-known equivalence between the category $\textbf{Pro}$ of preorders and the subcategory of $\textbf{Top}$ given by the Alexandrov spaces (note that the extended functor $\tilde{B}:\textbf{Pro}\to \textbf{Top}$ is already faithful and hence it is no longer necessary, unlike in the case of Alexandrov locales, to restrict to the category $\textbf{Pos}$ of posets to obtain an equivalence of categories).

Note that the characterization of Alexandrov locales given by Proposition \ref{alexlocale} does not lift to a topological characterization of Alexandov spaces as the sober topological spaces whose locales of open sets are Alexandrov locales, since Alexandrov spaces are not in general sober (cf. Proposition \ref{alexsober}). On the other hand, we can build an equivalence of $\textbf{Pos}$ with a category of sober topological spaces, as follows.

For any preorder $\cal P$, the topos $[{\cal P}, \Set]$, being localic, has only a \emph{set} of points (up to isomorphism). By choosing, for each preorder ${\cal P}$, as set of points of the topos $[{\cal P}, \Set]$ the collection $\textsl{S}_{\cal P}$ of \emph{all} its points and, for each arrow $f:{\cal P}\to {\cal Q}$, the function $l_{f}:\textsl{S}_{\cal P}\to \textsl{S}_{\cal Q}$ sending a point $\Set \to [{\cal P}, \Set]$ to the composite of it with the geometric morphism $[{\cal P}, \Set] \to [{\cal Q}, \Set]$ induced by $f$ as in section \ref{locales} above, we obtain, by the technique of section \ref{dualtop}, a lifting of the equivalence of $\textbf{Pos}$ with the category $ExtIm(B_{\textbf{Pos}})$ of Alexandrov locales to an equivalence $\tilde{B_{\textbf{Pos}}}:\textbf{Pos} \simeq ExtIm(\tilde{B_{\textbf{Pos}}})$ between $\textbf{Pos}$ and a subcategory of $\textbf{Top}$. By Theorem \ref{sobriety}, all the spaces $X$ in $ExtIm(\tilde{B_{\textbf{Pos}}})$ are sober and have the property that their locales of open sets ${\cal O}(X)$ are Alexandrov locales; hence, by the remarks in section \ref{dualtop}, the spaces in $ExtIm(\tilde{B_{\textbf{Pos}}})$ can be characterized precisely as the sober topological spaces whose frames of open sets are, regarded as locales, Alexandrov locales, while the arrows $X\to Y$ in $ExtIm(\tilde{B_{\textbf{Pos}}})$ are exactly the continuous maps $g:X\to Y$ between such spaces.  

The topological space $\tilde{B_{\textbf{Pos}}}({\cal P})$ corresponding to a poset $\cal P$ can be described concretely as the space having as set of points the collection ${\cal F}^{dir}_{\cal P}$ of all the non-empty directed ideals on $\cal P$ and as open sets the subsets of the form
\[
{\cal F}_{U}=\{F\in {\cal F}^{dir}_{\cal P} \textrm{ | } F\cap U\neq \emptyset\},
\] 
where $U$ ranges among the upper sets in $\cal P$ (cf. Proposition \ref{mslattice}). Clearly, this space is precisely the sobrification of the Alexandrov space corresponding to $\cal P$. In these terms, the arrow $\tilde{B_{\textbf{Pos}}}(f)$ corresponding to an arrow $f:{\cal P}\to {\cal Q}$ in $\textbf{Pos}$ via the functor $\tilde{B_{\textbf{Pos}}}$ can be described as the map ${\cal F}^{dir}_{\cal P}\to {\cal F}^{dir}_{\cal Q}$ sending a ideal set $I$ in ${\cal F}^{dir}_{\cal P}$ to the ideal in ${\cal F}^{dir}_{\cal Q}$ generated by the image of $I$ under $f$.  

Let us denote by $\textbf{SSC}$ the subcategory of $\textbf{Top}$ having as objects the sober topological spaces with a basis of supercompact open sets and as arrows the continuous maps $f:X\to Y$ between such spaces. The inverse functor $\tilde{I_{B}}:\textbf{SSC}\to \textbf{Pos}$ sends a topological space $X$ in $\textbf{SSC}$ to the opposite of the poset of its supercompact open sets and a continuous map $f:X\to Y$ in $\textbf{SSC}$ to the opposite of the restriction of the left adjoint $f_{!}:{\cal O}(X) \to {\cal O}(Y)$ to $f^{-1}:{\cal O}(Y) \to {\cal O}(X)$ to the posets of supercompact open sets of $X$ and of $Y$. The technique of section \ref{dualtop} thus yields a functor $\tilde{I_{B_{\textbf{Pos}}}}:\textbf{SSC}\to \textbf{Pos}$ which is a categorical inverse to the functor $\tilde{B_{\textbf{Pos}}}:\textbf{Pos} \to \textbf{SSC}$. Note that $\tilde{I_{B_{\textbf{Pos}}}}$ sends a space $X$ in $\textbf{SSC}$ to the opposite of the poset of supercompact open sets of $X$ and a continuous map $f:X\to Y$ in $\textbf{SSC}$ to the opposite of the restriction of the left adjoint $f_{!}:{\cal O}(X) \to {\cal O}(Y)$ to $f^{-1}:{\cal O}(Y) \to {\cal O}(X)$ to the posets of supercompact open sets of $X$ and of $Y$.

Summarizing, have the following result.

\begin{theorem}
Via the functors $\tilde{B_{\textbf{Pos}}}:\textbf{Pos} \to \textbf{SSC}$ and $\tilde{I_{B_{\textbf{Pos}}}}:\textbf{SSC}\to \textbf{Pos}$ defined above, the category $\textbf{Pos}$ is equivalent to the category $\textbf{SSC}$.
\end{theorem}

\subsection{Stone duality for distributive lattices}\label{stonedist}

This duality falls under the category of contravariant equivalences obtainable by means of the technique of sections \ref{locales} and \ref{dualtop}.
 
Let $\textbf{DLat}$ be the category of distributive lattices and homomorphisms (i.e. maps preserving finite meets and finite joins) between them. We can equip each lattice $\cal D$ in $\textbf{DLat}$, regarded as a preorder coherent category, with the coherent Grothendieck topology $J_{{\cal D}}$. The arrows $f:{\cal D}\to {\cal D}'$ in $\textbf{DLat}$ can be identified with the morphisms of sites $({\cal D}, J_{\cal D}) \to ({\cal D}', J_{\cal D}')$; so, by the results of section \ref{locales}, we have a functor $A:\textbf{DLat}^{\textrm{op}}\to \textbf{Loc}$ sending a lattice $\cal D$ in $\textbf{DLat}$ to the locale $Id_{J_{\cal D}}({\cal D})$. Now, by Corollary \ref{cor}(i), any $\cal D$ in $\textbf{DLat}$ can be recovered from the topos $\Sh({\cal D}, J_{\cal D})$ as the set of its compact subterminals, endowed with the induced natural order. By Theorem \ref{induced}(i) and Remark \ref{rmkprincipal}, the hypotheses of Theorem \ref{construction} are satisfied, and hence, since by Theorem \ref{induced}(i) and Remark \ref{rmkprincipal}, the hypotheses of Theorem \ref{construction} are satisfied, we have an inverse functor $I_{A}:ExtIm(A)\to \textbf{DLat}^{\textrm{op}}$, defined on the extended image of the functor $A:\textbf{DLat}^{\textrm{op}}\to ExtIm(A)$. The functor $I_{A}$ sends a locale in $ExtIm(A)$ to the lattice of compact elements in it, and an arrow $f:L\to L'$ in $ExtIm(A)$ to the restriction of its associated frame homomorphism to the subsets of compact elements of $L$ and $L'$.

By Theorem \ref{propext}, the objects in the extended image $ExtIm(A)$ of the functor $A$ can be characterized as the locales which have a basis of compact elements which form, with respect to the induced natural order, a (distributive) lattice (equivalently, as the locales such that the collection of their compact elements is closed under finite meets and forms a basis of them), while the arrows in $ExtIm(A)$ are precisely the locale morphisms whose associated frame homomorphisms send compact elements to compact elements. The category of such locales is called by Johnstone in \cite{stone} the category of \emph{coherent locales}. The method of section \ref{dualtop} thus yields in this case the well-known duality between the category $\textbf{DLat}$ of distributive lattices and the category of coherent locales (cf. Corollary II3.3 \cite{stone}).

We shall call the locale corresponding to a distributive lattice $\cal D$ via this duality the \emph{Stone locale} associated to $\cal D$.

This duality between $\textbf{DLat}$ and the category of coherent locales can be lifted (under a form of the axiom of choice, necessary to ensure that the toposes $\Sh({\cal D}, J_{\cal D})$ have enough points) to a duality between $\textbf{DLat}$ and a subcategory of the category $\textbf{Top}$ of topological spaces, by means of the technique of section \ref{dualtop} above. Indeed, we can choose, for each lattice ${\cal D}$ in $\textbf{DLat}$, as set of points $P_{\cal D}$ of the topos $\Sh({\cal D}, J_{\cal D})$ the collection of all its points, and define the action on arrows $f:{\cal D}\to {\cal D}'$ in $\textbf{DLat}$ accordingly, as specified in section \ref{dualtop}.  

In this way we obtain an equivalence between $\textbf{DLat}^{\textrm{op}}$ and a subcategory $ExtIm(\tilde{A})$ of $\textbf{Top}$, whose objects can be characterized as the sober topological spaces such that the collection of their compact open sets is closed under finite intersection and forms a basis for the topology, and whose arrows are the continuous maps between these spaces such that the inverse image of any compact open set is a compact open set. The technique of section \ref{dualtop} thus yields in this case the classical Stone duality between the category $\textbf{DLat}$ of distributive lattices and the category of spectral topological spaces (cf. Example \ref{exa}(c)). We shall call the topological space corresponding to a distributive lattice $\cal D$ via this duality the \emph{Stone space} associated to $\cal D$.

If we restrict the duality between $\textbf{DLat}$ and $ExtIm(A)$ to the full subcategory $\textbf{Boole}$ of $\textbf{DLat}$ on the Boolean algebras, we obtain an equivalence between $\textbf{Boole}^{\textrm{op}}$ and a full subcategory of $ExtIm(A)$, whose objects are the locales which have a a basis of compact elements which, with respect to the induced order, forms a Boolean algebra (equivalently, on the locales such that the collection of all their complemented elements forms a basis for the locale), and whose arrows are the locale morphisms whose associated frame homomorphisms send complemented elements to complemented elements. By lifting this duality to a topological duality as above, we obtain a duality between $\textbf{Bool}$ and the subcategory of $\textbf{Top}$ whose objects are the sober topological spaces which have a basis of clopen subsets and whose arrows are the continuous maps between such spaces; that is, we recover precisely the classical Stone duality for Boolean algebras.

\subsection{The duality between spatial frames and sober spaces} 

Let $\textbf{Frm}_{sp}$ be the category of spatial frames and frame homomorphisms between them. We recall that a frame $\cal F$ is spatial if there exists a topological space $X$ such that $\cal F$ is isomorphic to the frame ${\cal O}(X)$ of open sets of $X$. Obviously, the opposite of the category $\textbf{Frm}_{sp}$ trivially identifies with a subcategory $\cal U$ of $\textbf{Loc}$, and this identification has the form of a functor $A:{\textbf{Frm}_{sp}}^{\textrm{op}}\to {\cal U}=ExtIm(A)$ with inverse $I_{A}:{\cal U}\to {\textbf{Frm}_{sp}}^{\textrm{op}}$, as in section \ref{locales} above. Indeed, one can naturally equip each frame $\cal F$, regarded as a preorder geometric category, with the canonical geometric Grothendieck topology $J_{\cal F}$; morphisms of frames ${\cal F} \to {\cal F}'$ yield morphisms of sites $({\cal F}, J_{\cal F}) \to ({\cal F}', J_{{\cal F}'})$, and for any frame $\cal F$, the locale $Id_{J_{\cal F}}({\cal F})$ is precisely the locale with underlying frame $\cal F$.

Starting with this trivial duality, if we take for any frame $\cal F$ in $\textbf{Frm}_{sp}$ as set of points of the topos $\Sh({\cal F}, J_{\cal F})$ the collection of \emph{all} the points of $\Sh({\cal F}, J_{\cal F})$ (equivalently, the completely prime filters on $\cal F$), then this set of points is separating for the topos (since $\cal F$ is spatial) and hence the technique of section \ref{dualtop} yields precisely the well-known duality between the category $\textbf{Frm}_{sp}$ and the category of sober spaces and continuous maps between them (cf. Example \ref{exa}(e)). 

\subsection{Lindenbaum-Tarski duality}\label{tarski}

Recall that an element $a$ of a frame $F$ is said to be an atom if $a\neq 0$ and for any $b\leq a$, either $b=0$ or $b=a$. A frame $F$ is said to be \emph{atomic} it it has a basis of atoms, i.e. if every element of $F$ can be expressed as a join of atoms; accordingly, we say that a topological space has a \emph{basis of atomic open subsets} if the frame ${\cal O}(X)$ has a basis of atoms, equivalently if every open set can be expressed as a union of non-empty open sets which do not contain any proper non-empty open set. 

Let $\textbf{AtFrm}$ be the category having as objects the atomic frames and as arrows the frame homomorphisms between them. We can equip each frame $\cal C$ in $\textbf{AtFrm}$ with the canonical topology $J_{{\cal C}}$. For any $\cal C$ in $\textbf{AtFrm}$, the collection $At({\cal C})$ of atoms of $\cal C$, regarded as a full subcategory of $\cal C$, is $J_{\cal C}$-dense and the induced Grothendieck topology $J_{\cal C}|_{At({\cal C})}$ is the trivial one; the Comparison Lemma thus yields an equivalence of toposes 
\[
\Sh({\cal C}, J_{\cal C})\simeq [At({\cal C}), \Set].
\]
On the other hand, the category $At({\cal C})$ is a discrete preorder and hence the upper sets in $At({\cal C})$ coincides precisely with the subsets of $At({\cal C})$; therefore, Theorem \ref{fund} yields an equivalence of toposes
\[
[At({\cal C}), \Set] \simeq \Sh({\mathscr{P}(At({\cal C}))}),  
\]
where ${\mathscr{P}(At({\cal C}))}$ is the locale given by the full powerset of $At({\cal C})$ endowed with the subset-inclusion order.

Composing the two equivalences we obtain an equivalence of toposes
\[
\Sh({\cal C}, J_{\cal C})\simeq \Sh({\mathscr{P}(At({\cal C}))}).
\]
This equivalence in turn entails an isomorphism of frames 
\[
{\cal C}\cong {\mathscr{P}(At({\cal C}))},
\]
which is precisely the isomorphism given by the Lindenbaum-Tarski representation theorem (generalized from atomic complete Boolean algebras to atomic frames).  

We note that the opposite $\textbf{AtFrm}^{\textrm{op}}$ of the category $\textbf{AtFrm}$ is naturally identified with a subcategory of $\textbf{Loc}$; so, for any subcategory $\cal A$ of $\textbf{AtFrm}$, the restriction $i|_{\cal A}$ to $\cal A$ of the inclusion functor $i:\textbf{AtFrm}^{\textrm{op}} \hookrightarrow \textbf{Loc}$ yields a duality between $\cal A$ and the subcategory of $\textbf{Loc}$ given by the image of $i|_{\cal A}$. In fact, such identification can be seen as one induced by the method of section \ref{locales}, since every arrow $f:{\cal C}\to {\cal D}$ in $\textbf{AtFrm}$ yields a morphism of sites $({\cal C}, J_{\cal C})\to ({\cal D}, J_{\cal D})$, and for any $\cal C$ in $\textbf{AtFrm}$, the frame of $J_{\cal C}$-ideals on $\cal C$ is isomorphic to $\cal C$.  

To lift this trivial duality of atomic frames with locales to a topological duality, we select as set of points of the topos $\Sh({\cal C}, J_{\cal C})\simeq [At({\cal C}), \Set]$ the collection $\xi^{A}_{\cal C}$ of the points of $[At({\cal C}), \Set]$ corresponding to the elements of $At({\cal C})$, as in section \ref{prealex}. To make this assignment functorial, we have, according to the indications given in section \ref{dualtop}, to assign to an arrow $f:{\cal C}\to {\cal D}$ in $\textbf{AtFrm}$ a function $l_{f}:At({\cal D})\to At({\cal C})$ in such a way that the pair $(\dot{f}, l_{f})$ defines an arrow $(\Sh({\cal D}, J_{{\cal D}}), \xi^{A}_{\cal D}) \to (\Sh({\cal C}, J_{{\cal C}}), \xi^{A}_{\cal C})$ in the category $\mathfrak{Top}_{p}$ of toposes paired with points (cf. section \ref{subterminal}). This condition corresponds precisely, under the equivalences $\Sh({\cal C}, J_{\cal C})\simeq [At({\cal C}), \Set]$ and $\Sh({\cal D}, J_{\cal D})\simeq [At({\cal D}), \Set]$, to the requirement that $\dot{f}:[At({\cal D}), \Set] \to [At({\cal C}), \Set]$ be induced by the arrow $l_{f}:At({\cal D})\to At({\cal C})$ as in Example A4.1.4 \cite{El}. Now, by Lemma A4.1.5 \cite{El} and the remarks preceding it (combined with the fact that every discrete category is Cauchy-complete), a geometric morphism $[At({\cal D}), \Set] \to [At({\cal C}), \Set]$ is induced by a function $At({\cal D})\to At({\cal C})$ as in Example A4.1.4 \cite{El} if and only if it is \emph{essential}, i.e. if its inverse image functor has a left adjoint (in fact, the function inducing the morphism can be recovered, up to isomorphism, from this left adjoint as its restriction to the full subcategories on the representable functors). This property is a topos-theoretic invariant, and so, following the principles of \cite{OC10}, we attempt to reformulate it in terms of the different representations of our toposes established above. 

Let us start from the first representation of our morphism $\dot{f}$, as the geometric morphism $\Sh({\cal D}, J_{{\cal D}}) \to \Sh({\cal C}, J_{{\cal C}})$ induced by a frame homomorphism $f:{\cal C}\to {\cal D}$. Recall that the assignment $L\to \Sh(L)$ defines a full and faithful $2$-functor from the $2$-category of locales to the $2$-category of Grothendieck toposes (cf. Proposition C1.4.5 \cite{El}); therefore, bearing in mind the possibility of characterizing adjoint functors `equationally' in terms of their unit and counit, we conclude that $\dot{f}$ is essential if and only if $f$ has a left adjoint $f_{!}:{\cal D}\to {\cal C}$ (where $\cal C$ and $\cal D$ are regarded as poset categories); note that, by the Special Adjoint Functor Theorem, this condition holds if and only if $f$ is complete, i.e. preserves arbitrary meets as well as arbitrary joins.

This argument shows that the pair $(\dot{f}, l_{f})$ considered above defines an arrow $(\Sh({\cal D}, J_{{\cal D}}), \xi^{A}_{\cal D}) \to (\Sh({\cal C}, J_{{\cal C}}), \xi^{A}_{\cal C})$ in the category of toposes paired with points if and only if the arrow $f$ is complete.

Now, let us turn to the second representation of our morphism $\dot{f}$ as a geometric morphism $s:[At({\cal D}), \Set] \to [At({\cal C}), \Set]$. By the remarks above, $s$ is essential if and only if there exists a function $l:At({\cal D})\to At({\cal C})$ such that, denoted by $\tilde{l}:[At({\cal D}), \Set]\to [At({\cal C}), \Set]$ the geometric morphism induced by $l$ as in Example A4.1.4 \cite{El}, $s$ is isomorphic to $\tilde{l}$. Now, since the inverse image functor of $\tilde{l}$ is given by composition with $l$, this condition says precisely that $s$ corresponds, under the equivalences $\Sh({\cal C}, J_{\cal C})\simeq [At({\cal C}), \Set]$ and $\Sh({\cal D}, J_{\cal D})\simeq [At({\cal D}), \Set]$, to the geometric morphism $\Sh({\cal D}, J_{\cal D}) \to \Sh({\cal C}, J_{\cal C})$ induced by the frame homomorphism $\mathscr{P}(l):{\mathscr{P}(At({\cal C}))} \to {\mathscr{P}(At({\cal D}))}$ taking inverse images of subsets along $l$. 

By putting together the two pieces of information just collected, we can conclude that a frame homomorphism $r:\mathscr{P}(At({\cal C})) \to \mathscr{P}(At({\cal D}))$ is of the form $\mathscr{P}(l)$ for some function $l:At({\cal D})\to At({\cal C})$ if and only if $r$ is complete; notice that this is the other key ingredient of Lindenbaum-Tarski duality. In fact, our discussion above provides us with a concrete description of the function $l$ in terms of the frame homomorphism $r:{\cal C}\to {\cal D}$; indeed, we know that the function $l$ inducing the morphism can be recovered as the action of the left adjoint $r_{!}:{\cal D}\to {\cal C}$ to $r$ on the atoms (the fact that this left adjoint sends atoms to atoms is a consequence of the general topos-theoretic fact that the left adjoints to the inverse image functors of essential geometric morphisms sends representables to representables). Explicitly, for an element $d\in {\cal D}$, $r_{!}(d)$ is equal to the infimum in $\cal C$ of the set $\{c\in {\cal C} \textrm{ | } r(c)\geq d\}$.          

Coming back to our original aim of building a topological duality for atomic frames, we can infer from our considerations above that restricting to the subcategory $\textbf{CAtFrm}$ of $\textbf{AtFrm}$ whose objects are the atomic frames and whose arrows are the complete frame homomorphisms between them is precisely the condition which allows us to make the assignment ${\cal C}\to \xi^{A}_{\cal C}$ (of a set of points of the topos $\Sh({\cal C}, J)$ to a frame $\cal C$ in $\textbf{CAtFrm}$) functorial in the sense of section \ref{dualtop}. We thus have a `lifting' functor $\tilde{A}:\textbf{CAtFrm} \to \textbf{Top}$ which sends a frame $\cal C$ in $\textbf{CAtFrm}$ to the space obtained by equipping the set of points $At({\cal C})$ with the subterminal topology (on the topos $[At({\cal C}), \Set]$), in other words the discrete topological space with underlying set $At({\cal C})$ (cf. Example \ref{exa}(b)). In fact, since every such space is sober, the general results of section \ref{dualtop} yield the following characterization of the extended image of the functor $\tilde{A}:\textbf{CAtFrm} \to \textbf{Top}$: the objects of $ExtIm(\tilde{A})$ are the sober topological spaces which have a basis of atomic open subsets, while the arrows of $ExtIm(\tilde{A})$ are precisely the continuous maps between them. 

On the other hand, we have just characterized concretely the objects in $ExtIm(\tilde{A})$ as the spaces which are homeomorphic to the discrete topological spaces on the set of atoms of an atomic frame, and it is clear, from the equivalence ${\cal C}\cong {\mathscr{P}(At({\cal C}))}$ established above, that every discrete topological space is of that form. So $ExtIm(\tilde{A})$ can alternatively be described as the subcategory of $\textbf{Top}$ whose objects are the discrete topological spaces and whose arrows are the (continuous) functions between them. Since this category is (trivially) isomorphic to the category $\Set$ of sets, we conclude that the functor $\tilde{A}$ yields a duality between the category $\textbf{CAtFrm}$ and the opposite $\Set^{\textrm{op}}$ of the category of sets. Concretely, the functor $\tilde{A}$ sends an atomic frame $\cal C$ to the set $At({\cal C})$ of its atoms and an arrow $r:{\cal C}\to {\cal D}$ in $\textbf{CAtFrm}$ to the restriction $At({\cal D})\to At({\cal C})$ of the left adjoint to $r$ to the sets of atoms of $\cal C$ and $\cal D$. The inverse $\tilde{I_{A}}:ExtIm(\tilde{A})\simeq {\Set}^{\textrm{op}} \to \textbf{CAtFrm}$ of this functor is obtained, as indicated in section \ref{dualtop}, by recovering each frame $\cal C$ in $\textbf{CAtFrm}$ from the topos $\Sh({\mathscr{P}(At({\cal C}))}) \simeq [At({\cal C}), \Set]$ (equivalently, from the set $At({\cal C})$) through a topos-theoretic invariant, functorially in $\textbf{CAtFrm}$. In this case, the topos-theoretic invariant (of families of subterminals in a topos) is the vacuous one, and hence $\tilde{I_{A}}$ sends a set $S$ to the powerset $\mathscr{P}(S)$ and a function $f:S\to T$ to the function $\mathscr{P}(f):\mathscr{P}(S) \to \mathscr{P}(T)$ sending a subset to the inverse image of it under $f$.

We have thus recovered the Lindenbaum-Tarski duality (cf. for example \cite{stone} for a classical treatment of it).   
 
Note that, in passing, we have established the following intrinsic characterization of discrete spaces. 

\begin{proposition}\label{discrete}
Let $X$ be a topological space. Then $X$ is discrete if and only if it is sober and has a basis of atomic open subsets.
\end{proposition}\qed

Finally, we remark that, in the same way we established Lindenbaum-Tarski duality for atomic frames, we can establish a more general topological duality between, on one hand, the category having as objects the frames with a basis of supercompact elements and as arrows the complete frame homomorphisms between them and, on the other hand, the category of Alexandrov spaces whose underlying preorder is a poset; this latter duality specializes to Lindenbaum-Tarski duality when sets are regarded as discrete preorders (in fact, it is essentially the same as the equivalence of Theorem \ref{proalex}).

\subsection{A duality for meet-semilattices}\label{meet}

Let $\textbf{MSLat}$ be the category of meet-semilattices and homomorphisms (i.e. maps preserving finite meets) between them. 

Any meet-semilattice can be considered as a cartesian preorder category; as such, when equipped with the trivial Grothendieck topology, it gives rise to a cartesian site, and the arrows in $\textbf{MSLat}$ can be identified with the morphisms of these associated sites. The method of section \ref{locales} thus yields a faithful functor $A:\textbf{MSLat}^{\textrm{op}}\to \textbf{Loc}$ sending an object $\cal M$ of $\textbf{MSLat}$ to the frame $Id({{\cal M}})$ of lower sets in $\cal M$ and an arrow $f:{\cal M}\to {\cal N}$ to the frame homomorphism $A(f):Id({{\cal M}})\to Id({{\cal N}})$ which assigns to a lower set $S$ in $\cal M$ the lower set $A(f)(S)$ in $\cal N$ generated by the image $f(S)$ of $S$ in $\cal N$ under $f$. 

By Corollary \ref{cor}(iii), every $\cal M$ in $\textbf{MSLat}$ is isomorphic to the subset of supercompact subterminals in the topos $[{\cal M}^{\textrm{op}}, \Set]$, endowed with the natural order between subterminals in $[{\cal M}^{\textrm{op}}, \Set]$. By Theorem \ref{induced}(ii) and Remark \ref{rmkprincipal}, the hypotheses of Theorem \ref{construction} are satisfied, and hence we have an inverse functor $I_{A}:\textbf{SCLoc}\to \textbf{MSLat}^{\textrm{op}}$ to $A$, where $\textbf{SCLoc}$ is the subcategory of $\textbf{Loc}$ given by the extended image of the functor $A$. By Theorem \ref{propext}, the locales in $\textbf{SCLoc}$ can be characterized precisely as those which have a basis of supercompact elements which is closed under finite meets (equivalently, as those such that the collection of their supercompact elements is closed under finite meets and forms a basis for them), while the arrows in $\textbf{SCLoc}$ are the locale morphisms whose associated frame homomorphisms send supercompact elements to supercompact elements. The functor $I_{A}:\textbf{SCLoc}\to \textbf{MSLat}^{\textrm{op}}$ sends a locale $L$ in $\textbf{SCLoc}$ to the collection of its supercompact elements (endowed with the natural order) and an arrow $f:L\to L'$ in $\textbf{SCLoc}$ to the restriction of its associated frame homomorphism to the subsets of supercompact elements of $L$ and $L'$.  

We have thus obtained a duality between $\textbf{MSLat}$ and a category of locales:

\begin{theorem}\label{meetsm}
The functors 
\[
A:\textbf{MSLat}^{\textrm{op}} \to \textbf{SCLoc}
\]
and 
\[
I_{A}:\textbf{SCLoc}\to \textbf{MSLat}^{\textrm{op}}
\] 
defined above yield a duality between the category $\textbf{MSLat}$ and the category $\textbf{SCLoc}$.  
\end{theorem}

We shall call the locale corresponding to a meet-semilattice $\cal M$ via this duality the \emph{ideal locale} associated to $\cal M$.

Note that, at the level of objects, the duality of Theorem \ref{meetsm} can be read as the assertion that every meet-semilattice $\cal M$ is isomorphic to the subset of supercompact elements of the locale of ideals $Id({{\cal M}})$ of $\cal M$ (endowed with the natural induced order), and every locale with a basis of supercompact elements which is closed under finite meets is isomorphic to the locale of ideals $Id({{\cal M}})$ of the meet-semilattice $\cal M$ consisting of its supercompact elements.

The duality between $\textbf{MSLat}$ and $\textbf{SCLoc}$ can be lifted to a duality between $\textbf{MSLat}$ and a subcategory of the category $\textbf{Top}$ of topological spaces by using the method of section \ref{dualtop}. Let us choose, for each $\cal M$ in $\textbf{MSLat}$, as set of points of the topos $[{\cal M}^{\textrm{op}}, \Set]$, the set of \emph{all} the points of the topos $[{\cal M}^{\textrm{op}}, \Set]$, and define the function corresponding to the arrows in $\textbf{MSLat}$ accordingly, as specified in section \ref{dualtop}. This defines a `lifting functor' $\tilde{A}:\textbf{MSLat}^{\textrm{op}} \to \textbf{Top}$ having an inverse defined on its extended image, which we denote by $\textbf{SCTop}$. The spaces in $\textbf{SCTop}$ can be characterized precisely as the sober spaces which have a basis of supercompact open sets which is closed under finite intersections, while the arrows in $\textbf{SCTop}$ can be characterized as the continuous maps of topological spaces in $\textbf{SCTop}$ such that the inverse image of any supercompact open set is supercompact. 

Concretely, by Proposition \ref{mslattice}, the extended functor $\tilde{A}:\textbf{MSLat}^{\textrm{op}} \to \textbf{Top}$ sends to a meet-semilattice $\cal M$ the topological space having as set of points the collection ${\cal F}_{\cal M}$ of filters on $\cal M$ and as open sets the sets of the form
\[
{\cal F}_{I}=\{F\in {\cal F}_{\cal M} \textrm{ | } F\cap I\neq \emptyset\},
\] 
where $I$ ranges among the ideals of $\cal M$. In particular, a basis for this topological space is given by the sets
\[
{\cal F}_{m}=\{F\in {\cal F}_{\cal M} \textrm{ | } m\in {\cal M}\},
\] 
where $m$ varies among the elements of $\cal M$. We shall call this topological space the \emph{ideal topological space} associated to $\cal M$.

By Proposition \ref{mslatticefun}, the continuous map ${\cal F}_{\cal N}\to {\cal F}_{\cal M}$ corresponding to a meet-semilattice homomorphism $f:{\cal M}\to {\cal N}$ via the functor $\tilde{A}$ is the map sending a filter $F$ in ${\cal F}_{\cal N}$ to the inverse image $f^{-1}(F)$.   

The inverse functor $\tilde{I_{A}}:\textbf{SCTop}\to \textbf{MSLat}^{\textrm{op}}$ sends a topological space $X$ in $\textbf{SCTop}$ to the poset $\textbf{SC}({{\cal O}(X)})$ of its supercompact open sets, and a continuous map $f:X\to Y$ in $\textbf{SCTop}$ to the restriction $f^{-1}:\textbf{SC}({{\cal O}(Y)}) \to \textbf{SC}({{\cal O}(X)})$ of the inverse image $f^{-1}$ to the posets of supercompact open sets of $X$ and $Y$. 

Summmarizing, we have the following result.

\begin{theorem}\label{meetdual}
Via the functors 
\[
\tilde{A}:\textbf{MSLat}^{\textrm{op}} \to \textbf{SCTop}
\]
and 
\[
\tilde{I_{A}}:\textbf{SCTop}\to \textbf{MSLat}^{\textrm{op}}
\]
defined above, the categories $\textbf{MSLat}$ and $\textbf{SCTop}$ are dual to each other.
\end{theorem}

\begin{remark}
\begin{enumerate}[(a)]
\item The functor $\tilde{A}:\textbf{MSLat}^{\textrm{op}} \to \textbf{SCTop}$ coincides with the functor from meet-semilattices to topological spaces giving one half of the topological duality for meet-semilattices established by M. A. Moshier and P. Jipsen in \cite{jipsen}. In fact, it follows from Theorem \ref{meetdual} and the duality theorem in \cite{jipsen} that the category $\textbf{SCTop}$ coincides with the subcategory of $\textbf{Top}$ corresponding to the category of meet-semilattices via the duality in \cite{jipsen} (since the two theorems ensure that both these categories are equal to the extended image of the same functor, namely $\tilde{A}:\textbf{MSLat}^{\textrm{op}} \to \textbf{SCTop}$); anyway, it is worth to note that our definition of the inverse functor $\tilde{I_{A}}:\textbf{SCTop}\to \textbf{MSLat}^{\textrm{op}}$ is completely different from the construction of the inverse functor given in \cite{jipsen};

\item The functor $\tilde{A}:\textbf{MSLat}^{\textrm{op}} \to \textbf{SCTop}$ can also be identified with the functor from meet-semilattices to topological meet-semilattices (regarded here as topological spaces) giving one half of the duality between meet-semilattices and compact zero-dimensional meet-semilattices established in \cite{semilattices}. 
\end{enumerate}
\end{remark}

\subsection{Other dualities}\label{comparison}

We have seen in section \ref{stonedist} that the well-known Stone duality for distributive lattices is recovered through our method by equipping each  distributive lattice with the coherent topology. Anyway, there are several other interesting subcanonical Grothendieck topologies which one can put on distributive lattices, or on any other kind of preordered structures; for each of these choices, we can try to use our machinery to build dualities which, although being similar in spirit to Stone's one, are different enough to capture other aspects of the structures involved. We have already presented some examples above, of both old and new dualities obtained through the application of our methodology; in this section we discuss some further ones, of more `exotic' nature, which we hope should illuminate the flexibility of our general machinery. 

First, let us build a localic duality for $k$-frames (for a regular cardinal $k$), which specializes to the localic Stone duality for distributive lattices and to the duality of Theorem \ref{meetsm}. If we equip each $k$-frame ${\cal F}$ with the $k$-covering topology $J^{k}_{{\cal F}}$ we obtain that the morphisms ${\cal F}\to {\cal F}'$ of $k$-frames are precisely the morphisms of sites $({\cal F}, J^{k}_{{\cal F}}) \to ({\cal F}', J^{k}_{{\cal F}'})$. Therefore, by the technique of section \ref{locales}, we have a functor $A:{k\textrm{-}\textbf{Frm}}^{\textrm{op}}\to \textbf{Loc}$ which sends a poset $\cal F$ in ${k\textrm{-}\textbf{Frm}}$ to the locale of $J^{k}_{{\cal F}}$-ideals on $\cal F$ and an arrow $f:{\cal F}\to {\cal F}'$ in ${k\textrm{-}\textbf{Frm}}$ to the frame homomorphism $Id_{J^{k}_{{\cal F}}}({\cal F}) \to Id_{J^{k}_{{\cal F}'}}({\cal F}')$ sending an ideal $I$ in $Id_{J^{k}_{{\cal F}}}({\cal F})$ to the $J^{k}_{{\cal F}'}$-ideal on ${\cal F}'$ generated by the image $f(I)$ of $I$ under $f$. 

By Theorem \ref{induced}(iii) and Remark \ref{rmkprincipal}, the hypotheses of Theorem \ref{construction} are satisfied, and hence we obtain, by the results of section \ref{locales}, that the functor $A$ admits an inverse $I_{A}$ defined on its extended image of, which we denote by $\textbf{Loc}_{k}$. This latter subcategory has as objects the locales which have a basis of $k$-compact elements which is closed under finite meets and as arrows the locale morphisms whose associated frame homomorphisms send $k$-compact elements to $k$-compact elements. The inverse functor $I_{A}: \textbf{Loc}_{k} \to {k\textrm{-}\textbf{Frm}}^{\textrm{op}}$ sends a locale in $\textbf{Loc}_{k}$ to the poset of its $k$-compact elements and a locale morphism $L\to L'$ to the restriction of its associated frame homomorphism to the subsets of $k$-compact elements of $L$ and $L'$. Summarizing, we have the following result.

\begin{theorem}
Via the functors $A:k\textrm{-}\textbf{Frm}^{\textrm{op}} \to \textbf{Loc}_{k}$ and $I_{A}: \textbf{Loc}_{k} \to {k\textrm{-}\textbf{Frm}}^{\textrm{op}}$ defined above, the categories ${k\textrm{-}\textbf{Frm}}$ and $\textbf{Loc}_{k}$ are dual to each other. 
\end{theorem} 

Now, let us build dualities for disjunctively distributive lattices and for weakly atomic meet-semilattices. To this end, we recall from section \ref{charinv} the definition of disjunctive topology and atomically generated topology. 

Given a disjunctively distributive lattice $\cal D$, the \emph{disjunctive topology} $J_{{\cal D}}^{dj}$ on $\cal D$ is the Grothendieck topology on $\cal D$ (regarded as a preorder category) whose $J_{{\cal D}}^{dj}$-covering sieves on any element $d$ are exactly the sieves in $\cal D$ on $d$ which contain a finite family $\{d_{i} \textrm{ | } i\in I\}$ of elements $d_{i}\leq d$ such that for each pair of distinct $i, j\in I$, $d_{i}\wedge d_{j}=0$, and $\mathbin{\mathop{\textrm{\huge $\vee$}}\limits_{i\in I}} d_{i}=d$. Similarly, on an disjunctively distributive frame one can consider the infinitary disjunctive topology. 

Given a weakly atomic meet-semilattice $\cal M$, the \emph{atomically generated top-}\\ \emph{ology} $J_{{\cal M}}^{at}$ on $\cal M$ is defined as follows: the $J_{{\cal M}}^{at}$-covering sieves on any element $d$ are the maximal sieves and the sieves in $\cal M$ on $d$ which contain a finite family $\{d_{i} \textrm{ | } i\in I\}$ of elements $d_{i}\leq d$ such that $\mathbin{\mathop{\textrm{\huge $\vee$}}\limits_{i\in I}} d_{i}=d$ and each $d_{i}$ is an atom in $\cal M$ (i.e. for any $e\in {\cal M}$, if $e\leq d$ then either $e=d$ or $e=0$). Similarly, on any infinitarily weakly atomic meet-semilattice one can consider the infinitary atomically generated topology.

The definition of supercompactly generated topology on a weakly supercompact meet-semilattice is of course entirely analogous.  

First, let us focus on disjunctively distributive lattices.

Any arrow ${\cal D}\to {\cal D}'$ in $\textbf{DJLat}$ (cf. section \ref{charinv} for the definition of this category) yields a morphism of sites $({\cal D}, J_{{\cal D}}^{dj})\to ({\cal D}', J_{{\cal D}'}^{dj})$. So, by the technique of section \ref{locales}, we have a faithful functor $A:\textbf{DJLat}^{\textrm{op}} \to \textbf{Loc}$, which sends a poset $\cal D$ in $\textbf{DJLat}$ to the locale of $J_{{\cal D}}^{dj}$-ideals on $\cal D$ and an arrow $f:{\cal D}\to {\cal D}'$ in $\textbf{DJLat}$ to the frame homomorphism $Id_{J_{{\cal D}}^{dj}}({\cal D}) \to Id_{J_{{\cal D}'}^{at}}({\cal D}')$ which sends an ideal $I$ in $Id_{J_{{\cal D}}^{dj}}({\cal D})$ to the $J_{{\cal D}'}^{dj}$-ideal on ${\cal D}'$ generated by the image $f(I)$ of $I$ in ${\cal D}'$ under $f$. 

By Corollary \ref{cor}(v), every $\cal D$ in $\textbf{DJLat}^{\textrm{op}}$ can be recovered, up to isomorphism, as the poset of disjunctively compact subterminals in the topos $\Sh({\cal D}, J_{{\cal D}}^{dj})$. By Theorem \ref{induced}(iv) and Remark \ref{rmkprincipal}, the hypotheses of Theorem \ref{construction} are satisfied, and hence we obtain, by the results of section \ref{locales}, that the functor $A$ admits an inverse $I_{A}$ defined on its extended image of, which we denote by $\textbf{Loc}_{dj}$. By Theorem \ref{propext} and Proposition \ref{multicomposition}, we can characterize the objects in $\textbf{Loc}_{dj}$ as the locales which have a basis of disjunctively compact elements which is closed under finite meets and satisfies the property that any covering of an element of the basis has a disjunctively compact refinement by elements of the basis, while the arrows in $\textbf{Loc}_{dj}$ are precisely the locale morphisms whose associated frame homomorphisms send disjunctively compact elements to disjunctively compact elements. The functor $I_{A}:\textbf{Loc}_{dj}\to \textbf{DJLat}^{\textrm{op}}$ sends a locale in $\textbf{Loc}_{dj}$ to the poset of its disjunctively compact elements and a locale morphism $L\to L'$ in $\textbf{Loc}_{dj}$ to the restriction of its associated frame homomorphism to the sets of disjunctively compact elements of $L$ and $L'$. Summarizing, we have the following result:

\begin{theorem}
Via the functors $A:\textbf{DJLat}^{\textrm{op}}\to \textbf{Loc}_{dj}$ and $I_{A}:\textbf{Loc}_{dj} \to \textbf{DJLat}^{\textrm{op}}$ defined above, the categories $\textbf{DJLat}$ and $\textbf{Loc}_{dj}$ are dual to each other.
\end{theorem}\qed

Note that the category $\textbf{DLat}$ of distributive lattices can be identified as a full subcategory of the category $\textbf{DJLat}$, so this duality restricts in particular to a duality between $\textbf{DLat}$ and a subcategory of $\textbf{Loc}_{dj}$.  

Since all the toposes involved are coherent (and hence, under a form of the axiom of choice, have enough points by Deligne's theorem), this duality can be lifted to a duality with a category of topological spaces, by assigning to each $\cal D$ in $\textbf{DJLat}$ the set of points of the topos $\Sh({\cal D}, J_{{\cal D}}^{dj})$ (equivalently, the collection of \emph{disjunctive filters} $F$ on $\cal D$, i.e. the filters on $\cal D$ with the property that for any $J_{{\cal D}}^{dj}$-covering sieve $S$ on an object $d$ of $\cal D$ if $d\in F$ then there exists an element $f\in S$ such that $dom(f)$ belongs to $F$) topologized with the subterminal topology, and assigning function to arrows accordingly, as specified in section \ref{dualtop}. The functor $\tilde{A}:\textbf{DJLat}^{\textrm{op}} \to \textbf{Top}$ sends a poset $\cal D$ in $\textbf{DJLat}$ to the topological space $\tilde{A}({\cal D})$ having as underlying set the set ${\cal D}_{dj}$ of disjunctive filters on $\cal D$ and as basic open sets those of the form ${\cal F}_{c}=\{F\in {\cal D}_{dj} \textrm{ | } d\in F\}$ for $d\in {\cal D}$ (cf. Proposition \ref{mslattice}) and an arrow $f:{\cal D}\to {\cal D}'$ in $\textbf{DJLat}$ to the continuous map ${{\cal D}'}_{dj}\to {\cal D}_{dj}$ sending a filter $F$ in ${{\cal D}'}_{dj}$ to the inverse image $f^{-1}(F)$ (cf. Proposition \ref{mslatticefun}). The extended image of the functor $\tilde{A}$ is the subcategory $\textbf{Top}_{dj}$ of $\textbf{Top}$ whose objects are the sober topological spaces with a basis of disjunctively compact open sets which is closed under finite intersections and satisfies the property that any covering of a basic open set has a disjunctively compact refinement by basic open sets, and whose arrows are the continuous maps between such spaces such that the inverse image of any disjunctively compact open set is a disjunctively compact open set. The inverse functor $\tilde{I_{A}}:\textbf{Top}_{dj}\to \textbf{DJLat}^{\textrm{op}}$ sends a topological space $X$ in $\textbf{Top}_{dj}$ to the poset $\textbf{DJ}{{\cal O}(X)}$ of its disjunctively compact open sets and a continuous map $f:X\to Y$ in $\textbf{Top}_{dj}$ to the restriction $f^{-1}:\textbf{DJ}{\cal O}(Y) \to \textbf{DJ}{\cal O}(X)$ of the inverse image $f^{-1}$ to the posets of disjunctively compact open sets of $X$ and $Y$. 

Summarizing, we have the following result.

\begin{theorem}
Via the functors 
\[
\tilde{A}:\textbf{DJLat}^{\textrm{op}} \to \textbf{Top}_{dj}
\]
and 
\[
\tilde{I_{A}}:\textbf{Top}_{dj}\to \textbf{DJLat}^{\textrm{op}}
\]
defined above, the categories $\textbf{DJLat}$ and $\textbf{Top}_{dj}$ are dual to each other. 
\end{theorem}\qed

Similarly, by using Theorem \ref{induced}(v) and Remark \ref{rmkprincipal}, one obtains a duality between the category $\textbf{WAtMSLat}$ of weakly atomic meet-semilattices (cf. section \ref{charinv}) and the subcategory $\textbf{Loc}_{WAt}$ of \textbf{Loc} whose objects are the locales which have a basis of atomically compact elements which is closed under finite meets and whose arrows are the locale maps whose associated frame homomorphisms send atoms to atoms and atomically compact elements to atomically compact elements. The functor $A:\textbf{WAtMSLat}^{\textrm{op}}\to \textbf{Loc}_{WAt}$ sends a poset $\cal M$ in $\textbf{WAtMSLat}$ to the locale of $J_{{\cal M}}^{at}$-ideals on $\cal M$ and an arrow $f:{\cal M}\to {\cal N}$ in $\textbf{WAtMSLat}$ to the frame homomorphism $Id_{J_{{\cal M}}^{at}}({\cal M}) \to Id_{J_{{\cal N}}^{at}}({\cal N})$ which sends an ideal $I$ in $Id_{J_{{\cal M}}^{at}}({\cal M})$ to the $J_{{\cal N}}^{at}$-ideal on $\cal N$ generated by the image $f(I)$ of $I$ in ${\cal N}$ under $f$. The inverse functor $I_{A}:\textbf{Loc}_{WAt} \to \textbf{WAtMSLat}^{\textrm{op}}$ sends a locale $L$ in $\textbf{Loc}_{WAt}$ to the poset of atomically compact elements of $L$ and an arrow $f:L\to L'$ in $\textbf{Loc}_{WAt}$ to the restriction of its associated frame homomorphism to the subsets of atomically compact elements of $L$ and of $L'$.    

Summarizing, we have the following result.

\begin{theorem}
Via the functors 
\[
A:\textbf{WAtMSLat}^{\textrm{op}}\to \textbf{Loc}_{WAt}
\]
and
\[ 
I_{A}:\textbf{Loc}_{WAt} \to \textbf{WAtMSLat}^{\textrm{op}}
\]
defined above, the categories $\textbf{WAtMSLat}$ and $\textbf{Loc}_{WAt}$ are dual to each other.
\end{theorem}

Of course, such a duality admits an infinitary version for infinitarily weakly atomic meet-semilattices.

Since all the toposes involved are coherent (and hence, under a form of the axiom of choice, have enough points by Deligne's theorem), this duality can also be lifted to a topological duality, as in the example above. For instance, one can assign to each $\cal M$ in $\textbf{WAtMSLat}$ the set of points of the topos $\Sh({\cal M}, J_{{\cal M}}^{at})$ (i.e. the collection of \emph{atomic filters} $F$ on $\cal M$ i.e. the filters $F$ on $\cal M$ with the property that if a join of atoms belongs to $F$ then at least one of these atoms belongs to $F$) topologized with the subterminal topology, and defining the action on arrows accordingly. The resulting functor $\tilde{A}:\textbf{WAtMSLat}^{\textrm{op}} \to \textbf{Top}$ sends a poset $\cal M$ in $\textbf{WAtMSLat}$ to the topological space $\tilde{A}({\cal M})$ having as underlying set the set ${\cal M}_{at}$ of atomic filters on $\cal M$ and as basic open sets those of the form ${\cal F}_{c}=\{F\in {\cal M}_{at} \textrm{ | } d\in F\}$ for $d\in {\cal M}$ (cf. Proposition \ref{mslattice}), and an arrow $f:{\cal M}\to {\cal N}$ in $\textbf{WAtMSLat}$ to the continuous map ${{\cal N}}_{at}\to {\cal M}_{at}$ sending a filter $F$ in ${{\cal N}}_{at}$ to the inverse image $f^{-1}(F)$ (cf. Proposition \ref{mslatticefun}). The extended image of the functor $\tilde{A}$ is the subcategory $\textbf{Top}_{at}$ of $\textbf{Top}$ whose objects are the sober topological spaces with a basis of atomically compact open sets which is closed under finite intersections, and whose arrows are the continuous maps between such spaces such that the inverse image of any atomically compact open set is an atomically compact open set, and the inverse image of any atomic open set is an atomic open set. The inverse functor $\tilde{I(A)}:\textbf{Top}_{at}\to \textbf{WAtMSLat}^{\textrm{op}}$ sends a topological space $X$ in $\textbf{Top}_{at}$ to the poset $\textbf{AT}({{\cal O}(X)})$ of its atomically compact open sets and a continuous map $f:X\to Y$ in $\textbf{Top}_{at}$ to the restriction $f^{-1}:\textbf{AT}({\cal O}(Y)) \to \textbf{AT}({\cal O}(X))$ of the inverse image $f^{-1}$ to the posets of atomically compact open sets of $X$ and $Y$. Summarizing, we have the following result.

\begin{theorem}
Via the functors 
\[
\tilde{A}:\textbf{WAtMSLat}^{\textrm{op}}\to \textbf{Top}_{at}
\]
and 
\[
\tilde{I_{A}}:\textbf{Top}_{at} \to \textbf{WAtMSLat}^{\textrm{op}}
\]
defined above, the categories $\textbf{WAtMSLat}$ and $\textbf{Top}_{at}$ are dual to each other.
\end{theorem}

Similarly, one can obtain localic and topological dualities for weakly supercompact meet-semilattices, and a localic duality for infinitarily weakly supercompact meet-semilattices.

We have limited ourselves to presenting just a few examples of application of our machinery. Anyway, the reader should be convinced at this point that the generality and flexibility of the method described in sections \ref{locales} and \ref{dualtop} enables one to easily establish many other new dualities or equivalences. In particular, `natural' dualities, such as the ones that we have obtained above, can be generated by choosing the subcanonical topologies $J_{\cal C}$ in such a way to capture the operations defining the structure of the corresponding poset $\cal C$; for example, in the classical case of Stone duality, the joins in a distributive lattice can be directly defined in terms of the coherent topology on it, while the finite joins of pairwise disjoint elements in a disjunctively distributive lattice can be characterized in terms of the disjunctive topology on it. In connection with this, it is worth to remark that if $\cal C$ is a meet-semilattice then the Yoneda embedding $y:{\cal C} \hookrightarrow \Sh({\cal C}, J_{\cal C})\simeq \Sh(Id_{J_{\cal C}}({\cal C}))$ preserves finite meets and sends $J_{\cal C}$-covering sieves to covering families; so, meets in $\cal C$ correspond to meets in $Id_{{J_{\cal C}}}({\cal C})$, while operations on $\cal C$ defined in terms of $J$ correspond via the embedding ${\cal C}\hookrightarrow Id_{{J_{\cal C}}}({\cal C})$ to operations on $Id_{{J_{\cal C}}}({\cal C})$ involving joins (in particular, if $\cal C$ possesses a bottom element $0$ which is $J_{\cal C}$-covered by the empty sieve then $y$ sends the initial object $0$ of $\cal C$ to the initial object of $\Sh({\cal C}, J_{\cal C})$, that is to the bottom element of the frame $Id_{{J_{\cal C}}}({\cal C})$). 

\section{A further generalization}\label{generalization}

We note that the fundamental ingredient that we have taken as a starting point in section \ref{locales} to generate dualities consists of bunches of Morita-equivalences which are instances of Theorem \ref{fund}.  

More generally, we can expect to be able to extract representation theorems or dualities between partially ordered structures starting from general Morita-equivalences of the form $\Sh({\cal C}, J)\simeq \Sh({\cal D}, K)$, where $\cal C$ and $\cal D$ are poset categories. Indeed, if $J$ is $U$-induced for a topos-theoretic invariant $U$ of families of subterminals in a topos satisfying the hypothesis of Theorem \ref{construction}, $\cal C$ can be represented as the poset consisting of the elements of the locale $Id_{K}({\cal D})$ which are $U$-compact. And if we have a bunch of categorical equivalences $\Sh({\cal C}, J)\simeq \Sh({\cal D}, K)$ for a collection of pairs of poset structures $({\cal C}, {\cal D})$, we can expect to be able to build a duality or equivalence between categories having as objects respectively the structures of the form $\cal C$ and the structures of the form $\cal D$.

\subsection{From Morita-equivalences to dualities}\label{Mordual}

To formalize the idea described above, suppose that we have two collections ${\cal K}$ and $\cal H$ of poset structures, each of which equipped with a Grothendieck topology (we denote by $J_{\cal C}$ (resp. by $K_{\cal D}$) the Grothendieck topology associated to a structure $\cal C$ (resp. $\cal D$) in $\cal K$ (resp. in $\cal H$)), two functions $f:{\cal K}\to {\cal H}$ and $g:{\cal H}\to {\cal K}$ which are inverse to each other up to isomorphisms, and a categorical equivalence $\Sh({\cal C}, J_{\cal C})\simeq \Sh(f({\cal C}), K_{f({\cal C})})$ for each $\cal C$ in $\cal K$ (equivalently, a categorical equivalence $\Sh({\cal D}, K_{\cal D})\simeq \Sh(g({\cal D}), J_{g({\cal D})})$ for each $\cal D$ in $\cal H$). 

We shall adopt the following conventions:

\begin{enumerate} 
\item Given a morphism of posites $f:({\cal A}, J)\to ({\cal B}, K)$, we denote by $A(f)$ the frame homomorphism $Id_{J}({\cal A})\to Id_{K}({\cal B})$ which sends an ideal $I\in Id_{J}({\cal A})$ to the $K$-ideal on $\cal B$ generated by the image $f(I)$ of $I$ under $f$. Recall from section \ref{locales} that $A(f)$ can be identified with the restriction to the subterminals of the inverse image functor of the geometric morphism $\Sh({\cal B}, K) \to \Sh({\cal A}, J)$ induced by the morphism of sites $f$;

\item Given a functor $f:{\cal A}\to {\cal B}$ between posets, we denote by $B(f):Id({\cal B}^\textrm{op}) \to Id({\cal A}^\textrm{op})$ the frame homomorphism sending a subset $I\subseteq B$ in $Id({\cal B}^\textrm{op})$ to the inverse image $f^{-1}(I)$ of $I$ under $f$. Recall from section \ref{prealex} that $B(f)$ can be identified with the restriction to the subterminals of the inverse image functor of the geometric morphism $[{\cal A}, \Set] \to [{\cal B}, \Set]$ induced by the functor $f$ as in A4.1.4 \cite{El}. 
\end{enumerate}

Let us suppose that $U$ (resp. $V$) is a topos-theoretic invariant satisfying the hypothesis of Theorem \ref{thmcentral} with respect to the category $\cal K$ (resp. the category $\cal H$), so that all the Grothendieck topologies $J_{\cal C}$ (resp. $K_{\cal D}$) are $U$-induced (resp. $V$-induced). Given a frame $L$, we denote by $J_{can}^{L}$ the canonical topology on $L$ (i.e. the Grothendieck topology on $L$ whose covering sieves are exactly those which contain a covering family in $L$) and by $U\textrm{-comp}(L)$ (resp. $V\textrm{-comp}(L)$) the collection of the elements of $L$ which are $U$-compact (resp. $V$-compact); the Grothendieck topology on $U\textrm{-comp}(L)$ (resp. on $V\textrm{-comp}(L)$) induced by the canonical topology $J_{can}^{L}$ on $L$ by regarding $U\textrm{-comp}(L)$ (resp. $V\textrm{-comp}(L)$) as a full subcategory of $L$ via the inclusion $U\textrm{-comp}(L) \hookrightarrow L$ (resp. the inclusion $V\textrm{-comp}(L) \hookrightarrow L$) (cf. p. 546 \cite{El} for the general definition of induced coverage) will be denoted by $J_{can}^{L}|_{U\textrm{-comp}(L)}$ (resp. by $J_{can}^{L}|_{V\textrm{-comp}(L)}$) or simply by $J_{can}^{L}|$ when the subcategory $U\textrm{-comp}(L)$ (resp. $V\textrm{-comp}(L)$) can be unambiguously inferred from the context. Given a frame $L$, we denote by $\textbf{SC}(L)$ the set of supercompact elements of $L$.  

We distinguish between the covariant and contravariant case:

\begin{enumerate}[(i)]

\item \emph{Covariant case}

We define two categories $Ext({\cal K})$ and $Ext({\cal H})$ as follows. 

The objects of $Ext({\cal K})$ are the posets which are isomorphic to a poset in $\cal K$, while the arrows ${\cal C}\to {\cal C}'$ in $Ext({\cal K})$ are the monotone maps $f:{\cal C}\to {\cal C}'$ which induce morphisms of sites $({\cal C}, J_{\cal C})\to ({\cal C}', J_{{\cal C}'})$ such that the frame homomorphism $A(f):Id_{J_{\cal C}}({\cal C})\to Id_{J_{{\cal C}'}}({\cal C}')$ sends $V$-compact elements to $V$-compact elements and its restriction 
\[
A(f)|:V\textrm{-comp}(Id_{J_{\cal C}}({\cal C})) \to V\textrm{-comp}(Id_{J_{{\cal C}'}}({\cal C}'))
\]
to these posets yields a morphism of sites 
\[
(V\textrm{-comp}(Id_{J_{\cal C}}({\cal C})), J_{can}^{{Id_{J_{\cal C}}({\cal C})}}|) \to (V\textrm{-comp}(Id_{J_{{\cal C}'}}({\cal C}')), J_{can}^{{Id_{J_{{\cal C}'}}({\cal C}')}}|).
\]
The definition of $Ext({\cal H})$ is perfectly symmetrical to that of $Ext({\cal K})$ ($\cal H$ playing the role of $\cal K$ and $U$ playing the role of $V$). 

Let us now define two functors $D:Ext({\cal K}) \to Ext({\cal H})$ and $E:Ext({\cal H}) \to Ext({\cal K})$ which extend respectively the functions $f$ and $g$ (up to isomorphism) as follows. 

The functor $D:Ext({\cal K}) \to Ext({\cal H})$ sends a poset $\cal C$ in $Ext({\cal K})$ to the poset $V\textrm{-comp}(Id_{J_{\cal C}}({\cal C}))$ of $V$-compact elements of the frame $Id_{J_{\cal C}}({\cal C})$, and an arrow $f:{\cal C}\to {\cal C}'$ in $Ext({\cal K})$ to the restriction 
\[
A(f)|:V\textrm{-comp}(Id_{J_{\cal C}}({\cal C})) \to V\textrm{-comp}(Id_{J_{{\cal C}'}}({\cal C}'))
\]
of the frame homomorphism $A(f):Id_{J_{\cal C}}({\cal C})\to Id_{J_{{\cal C}'}}({\cal C}')$ to the subsets of $V$-compact elements of $Id_{J_{\cal C}}({\cal C})$ and of $Id_{J_{{\cal C}'}}({\cal C}')$. 

The definition of the functor $E:Ext({\cal H}) \to Ext({\cal K})$ is perfectly symmetrical to that of the functor $D$.

\begin{theorem}\label{dualcov}
Under the hypotheses specified above, the functors $D$ and $E$ are categorical inverses to each other and hence they define a categorical equivalence between $Ext({\cal K})$ and $Ext({\cal H})$.
\end{theorem}

\begin{proofs}
The fact that the functors $D$ and $E$ are well-defined, in the sense that $D$ takes values in $Ext({\cal H})$ (and, dually, $E$ takes values in $Ext({\cal K})$), and that they are categorical inverses to each other easily follows from Theorem \ref{thmcentral} by invoking the general theory of morphisms of sites. Specifically, one appeals to the fact that, given two sites $({\cal C}, J)$ and $({\cal D}, K)$ where $J$ and $K$ are subcanonical, a geometric morphism $f:\Sh({\cal D}, K)\to \Sh({\cal C}, J)$ is induced by a morphism of sites $({\cal C}, J)\to ({\cal D}, K)$ if and only if the inverse image $f^{\ast}$ sends representables to representables (Lemma C2.3.8 \cite{El}). 

\end{proofs}

From the theorem we can easily deduce the following characterization of the categories $Ext({\cal H})$ and $Ext({\cal K})$ in terms of each other, which represents the analogue of Theorem \ref{propext}.

\begin{proposition}
Let us assume that the hypotheses of Theorem \ref{dualcov} are satisfied, and that $U$ is $\cal K$-compatible relative to a property $P$. Let ${\cal F}^{U}_{P}$ denote the set of frames such that the collection of their $U$-compact elements forms a basis of them which, endowed with the induced order, has the structure of a poset in $Ext({\cal K})$ and has the property that the embedding $B_{L}\hookrightarrow L$ of it into the frame satisfies condition $P$, the property that every covering in $L$ of an element of $B_{L}$ is refined by a covering made of elements of $B_{L}$ which satisfies the invariant $U$, and the property that the $J_{B_{L}}$-covering sieves are sent by the embedding $B_{L}\hookrightarrow L$ into covering families in $L$. Then

\begin{enumerate}[(i)] 

\item The objects of $Ext({\cal H})$ are precisely the posets of $V$-compact elements of the frames in ${\cal F}^{U}_{P}$;  

\item The arrows in $Ext({\cal H})$ are the restrictions to the subsets of $V$-compact elements of the frame homomorphisms between frames in ${\cal F}^{U}_{P}$ which send $U$-compact elements to $U$-compact elements and $V$-compact elements to $V$-compact elements.
\end{enumerate}
\end{proposition}

Of course, the description of the category $Ext({\cal H})$ is perfectly symmetrical. 

\item \emph{Contravariant case}

In the contravariant case, one of the two toposes involved in the Morita-equivalence is a presheaf topos; that is, we suppose $K_{\cal D}$ to be the trivial topology for every $\cal D$ in $\cal H$. We can therefore suppose, by Theorem \ref{induced}(i) and Remark \ref{rmkprincipal}, the invariant $V$ to be the property `to be a singleton family'. 

We define two categories $Ext({\cal K})$ and $Ext({\cal H})$ as follows. 

The objects of $Ext({\cal K})$ are the posets which are isomorphic to a poset in $\cal K$, while the arrows ${\cal C}\to {\cal C}'$ in $Ext({\cal K})$ are the monotone maps $f:{\cal C}\to {\cal C}'$ such that they induce morphisms of sites $({\cal C}, J_{\cal C})\to ({\cal C}', J_{{\cal C}'})$ such that the frame homomorphism $A(f):Id_{J_{\cal C}}({\cal C})\to Id_{J_{{\cal C}'}}({\cal C}')$ is complete i.e. preserves arbitrary infima (cf. section \ref{tarski}). 

The objects of $Ext({\cal H})$ are the posets which are isomorphic to a poset in $\cal H$, while the arrows ${\cal D}\to {\cal D}'$ in $Ext({\cal H})$ are the monotone maps $f:{\cal D}\to {\cal D}'$ such that the frame homomorphism $B(f^{\textrm{op}}):Id({{\cal D}'}) \to Id({\cal D})$ sends $U$-compact elements to $U$-compact elements. 

Let us now define two functors
\[
D:Ext({\cal K})^{\textrm{op}} \to Ext({\cal H})
\]
and 
\[
E:Ext({\cal H}) \to Ext({\cal K})^{\textrm{op}}
\]
which extend respectively the functions $f$ and $g$ (up to isomorphism). 

The functor $D:Ext({\cal K})^{\textrm{op}} \to Ext({\cal H})$ sends a poset $\cal C$ in $Ext({\cal K})$ to the poset of supercompact elements of the frame $Id_{J_{\cal C}}({\cal C})$, and an arrow $f:{\cal C}\to {\cal C}'$ in $Ext({\cal K})$ to the restriction of the left adjoint $A(f)_{!}:Id_{J_{{\cal C}'}}({\cal C}')\to Id_{J_{{\cal C}}}({\cal C})$ to the functor $A(f):Id_{J_{\cal C}}({\cal C})\to Id_{J_{{\cal C}'}}({\cal C}')$ to the subsets of supercompact elements of $Id_{J_{{\cal C}'}}({\cal C}')$ and of $Id_{J_{\cal C}}({\cal C})$. 

The functor $E:Ext({\cal H}) \to Ext({\cal K})^{\textrm{op}}$ sends a poset $\cal D$ in $Ext({\cal H})$ to the poset of $U$-compact elements of the frame $Id({{\cal D}})$ and an arrow $f:{\cal D}\to {\cal D}'$ in $Ext({\cal H})$ to the restriction of the frame homomorphism $B(f^{\textrm{op}}):Id({{\cal D}'}) \to Id({\cal D})$ to the subsets of $U$-compact elements of $Id({\cal D}')$ and of $Id({\cal D})$.    
  
\begin{theorem}\label{dualabstract}
Under the hypotheses specified above, the functors $D$ and $E$ are categorical inverses to each other and hence they define a duality between $Ext({\cal K})$ and $Ext({\cal H})$.
\end{theorem}

\begin{proofs}
The fact that the functors $D$ and $E$ are well-defined, in the sense that $D$ takes values in $Ext({\cal H})$ (and, dually, $E$ takes values in $Ext({\cal K})^{\textrm{op}}$), and that $D$ and $E$ are categorical inverses to each other easily follows from Theorem \ref{thmcentral} by invoking the general theory of morphisms of sites (as in the proof of Theorem \ref{dualcov}), and from the theory of geometric morphisms induced by functors as in A4.1.4 \cite{El}. Concerning the latter, one specifically uses the fact (Lemma A4.1.5 \cite{El}) that, given two Cauchy-complete categories $\cal C$ and $\cal D$, a geometric morphism $[{\cal C}, \Set]\to [{\cal D}, \Set]$ is induced by a functor $f:{\cal C}\to {\cal D}$ as in A4.1.4 \cite{El} if and only if it is essential. 
\end{proofs}

Similarly to the covariant case, we can deduce from the theorem the following characterizations of the categories $Ext({\cal K})$ and $Ext({\cal H})$ in terms of each other. 

\begin{proposition}
Let us assume that the hypotheses of Theorem \ref{dualcov} are satisfied, and that $U$ is $\cal K$-compatible relative to a property $P$. Let ${\cal F}^{U}_{P}$ denote the set of frames such that the collection of their $U$-compact elements forms a basis of them which, endowed with the induced order, has the structure of a poset in $Ext({\cal K})$ and has the property that the embedding $B_{L}\hookrightarrow L$ of it into the frame satisfies condition $P$, the property that every covering in $L$ of an element of $B_{L}$ is refined by a covering made of elements of $B_{L}$ which satisfies the invariant $U$, and the property that the $J_{B_{L}}$-covering sieves are sent by the embedding $B_{L}\hookrightarrow L$ into covering families in $L$. Let ${\cal F}_{sc}$ denote the set of frames such that the collection of their supercompact elements forms a basis which, endowed with the induced order, has the structure of a poset in $Ext({\cal H})$. Then

\begin{enumerate}[(i)]

\item The objects of $Ext({\cal H})$ are precisely the posets of supercompact elements of the frames in ${\cal F}^{U}_{P}$;

\item The arrows in $Ext({\cal H})$ are the restrictions to the subsets of supercompact elements of the left adjoints to the frame homomorphisms between frames ${\cal F}^{U}_{P}$ which send $U$-compact elements to $U$-compact elements;

\item The objects of $Ext({\cal K})$ are precisely the posets of $U$-compact objects of the frames in ${\cal F}_{sc}$;

\item The arrows in $Ext({\cal H})$ are the restrictions to the subsets of $U$-compact elements of the frame homomorphisms between frames in ${\cal F}_{sc}$ which are complete.  

\end{enumerate}
\end{proposition}

\end{enumerate}
Note that all of this represents a clear implementation of the philosophy `toposes as bridges' of \cite{OC10}.  

The process of `lifting' of dualities with categories of locales to dualities with categories of topological spaces of section \ref{dualtop} has an analogue in this more general context. Again, we have to distinguish between the covariant and contravariant cases:

\begin{enumerate}[(i)] 

\item \emph{Covariant case}

Suppose that we have assigned to each structure $\cal C$ in $Ext({\cal K})$ a separating set of points of the topos $\Sh({\cal C}, J_{\cal C})$ functorially in $\cal C$, as in section \ref{dualtop}. We define a subcategory $\textbf{Top}_{\cal H}$ of $\textbf{Top}$ by taking as objects the topological spaces $X$ such that $V\textrm{-comp}({\cal O}(X))$ belongs to $Ext({\cal H})$ and as arrows $X\to Y$ are the continuous maps $f:X\to Y$ between spaces in $\textbf{Top}_{\cal H}$ with the property that $f^{-1}:{\cal O}(Y)\to {\cal O}(X)$ sends $V$-compact open sets to $V$-compact open sets in such a way that the restriction $f^{-1}|:V\textrm{-comp}({\cal O}(Y)) \to V\textrm{-comp}({\cal O}(X))$ is an arrow $V\textrm{-comp}({\cal O}(Y)) \to V\textrm{-comp}({\cal O}(X))$ in $Ext({\cal H})$. Thus we have, as in section \ref{dualtop}, a functor $\tilde{D}:Ext({\cal K})\to \textbf{Top}_{\cal H}^{\textrm{op}}$, and we can define an essentially surjective functor (both on the objects and on the arrows) $U_{\cal H}:\textbf{Top}_{\cal H}^{\textrm{op}} \to Ext({\cal H})$ such that $U_{{\cal H}}\circ \tilde{D}\cong D$; $U_{\cal H}$ sends a topological space $X$ to the poset of $V$-compact open subsets of $X$ and a continuous map $f:X\to Y$ in $\textbf{Top}_{\cal H}$ to the restriction of the inverse image $f^{-1}$ to the subsets of $V$-compact open sets of $X$ and $Y$. 

Of course, we can define a subcategory $\textbf{Top}_{\cal K}$ of $\textbf{Top}$ and a functor $U_{\cal K}:\textbf{Top}_{\cal K}^{\textrm{op}} \to Ext({\cal K})$ in a perfectly symmetrical way. 

\item \emph{Contravariant case}

Suppose that we have assigned to each structure $\cal C$ in $Ext({\cal K})$ a separating set of points of the topos $\Sh({\cal C}, J_{\cal C})$ functorially in $\cal C$, as in section \ref{dualtop}. Let us define a subcategory $\textbf{Top}_{\cal H}$ of $\textbf{Top}$ by taking as objects the topological spaces $X$ such that the poset $\textbf{SC}({\cal O}(X))$ of supercompact elements of ${\cal O}(X)$ belongs to $Ext({\cal H})$ and as arrows $X\to Y$ the continuous maps $f:X\to Y$ between spaces in $\textbf{Top}_{\cal H}$ such that $f^{-1}:{\cal O}(Y)\to {\cal O}(X)$ is complete and its left adjoint $f_{!}$ sends supercompact open sets to supercompact open sets in such a way that the restriction $f_{!}|:\textbf{SC}({\cal O}(X)) \to \textbf{SC}({\cal O}(Y))$ is an arrow $\textbf{SC}({\cal O}(X)) \to \textbf{SC}({\cal O}(Y))$ in $Ext({\cal H})$. Thus we have, as in section \ref{dualtop}, a functor $\tilde{D}:Ext({\cal K})^{\textrm{op}}\to \textbf{Top}_{\cal H}$, and we can define an essentially surjective functor (both on the objects and on the arrows) $U:\textbf{Top}_{\cal H} \to Ext({\cal H})$ such that $U\circ \tilde{D}\cong D$; $U$ sends a topological space $X$ to the poset of supercompact open subsets of $X$ and a continuous map $f:X\to Y$ in $\textbf{Top}_{\cal H}$ to the restriction of the left adjoint $f_{!}$ to the subsets of supercompact open sets of $X$ and $Y$. 

We can make a similar construction for the category $Ext({\cal H})$, as follows. Suppose that we have assigned to each structure $\cal D$ in $Ext({\cal H})$ a separating set of points of the topos $\Sh({\cal D}, K_{\cal D})\simeq [{\cal D}^{\textrm{op}}, \Set]$ functorially in $\cal D$, as in section \ref{dualtop}. We define a subcategory $\textbf{Top}_{\cal K}$ of $\textbf{Top}$ by taking as objects the topological spaces such that $U\textrm{-comp}({\cal O}(X))$ belongs to $Ext({\cal K})$ and as arrows $X\to Y$ the continuous maps $f:X\to Y$ between spaces in $\textbf{Top}_{\cal K}$ with the property that $f^{-1}:{\cal O}(Y)\to {\cal O}(X)$ sends $U$-compact open sets to $U$-compact open sets in such a way that the restriction $f^{-1}|:U\textrm{-comp}({\cal O}(Y)) \to U\textrm{-comp}({\cal O}(X))$ is an arrow $U\textrm{-comp}({\cal O}(Y)) \to U\textrm{-comp}({\cal O}(X))$ in $Ext({\cal K})$. Thus we have, as in section \ref{dualtop}, a functor $\tilde{E}:Ext({\cal H})\to \textbf{Top}_{\cal K}$, and we can define an essentially surjective functor (both on the objects and on the arrows) $W:\textbf{Top}_{\cal K} \to Ext({\cal K})^{\textrm{op}}$ such that $W\circ \tilde{E}\cong E$; $W$ sends a topological space $X$ to the poset $U\textrm{-comp}({\cal O}(X))$ of $U$-compact open subsets of $X$ and a continuous map $f:X\to Y$ in $\textbf{Top}_{\cal K}$ to the restriction of the inverse image $f^{-1}$ to the subsets of $U$-compact open sets of $X$ and of $Y$.     

\end{enumerate}

The topological spaces in the images of the `lifting functors' can be directly described in terms of the structures in the source category (cf. Proposition \ref{mslattice}); so the above constructions allow characterizations of the posets in the target categories as subsets consisting of the open sets of a topological space satisfying some generalized compactness condition.   

We notice that our original framework of sections \ref{locales} and \ref{dualtop} sits inside this more general setting as the particular case in which each of the toposes associated to the structures in $\cal K$ (or in $\cal H$) is of the form $\Sh({\cal C}, J_{\cal C})$, where $\cal C$ is a locale and $J_{\cal C}$ is the canonical topology on it, and the Morita-equivalences from which the dualities or equivalences are built are instances of Theorem \ref{fund}. Note that the canonical topologies $J_{\cal C}$ are all $U$-induced where $U$ is the vacuous invariant `to be a family' (note that this invariant trivially satisfies the hypotheses of Theorem \ref{construction}); in particular, the condition of $U$-compactness is vacuous in this case. 

\subsection{The role of the Comparison Lemma}

It is worth to remark the key role of Grothendieck's Comparison Lemma (cf. Theorem C2.2.3 \cite{El}) in the framework that we have developed. 

The equivalence $\Sh({\cal C}, J) \simeq \Sh(Id_{J}({\cal C}))$ of Theorem \ref{fund} is, when $J$ is subcanonical, an instance of an equivalence induced by the Comparison Lemma. Indeed, if $J$ is subcanonical then $\cal C$ can be identified with a full subcategory of $\Sh(Id_{J}({\cal C}))$ which is dense with respect to the canonical coverage, and the Grothendieck topology induced by it on $\cal C$ is precisely $J$ (cf. Proposition C2.2.16(ii) \cite{El}).    

So, all of the dualities established so far arise in fact from instances of the Comparison Lemma. Generalizing our original situation, we can expect to be able to obtain representation theorems for preordered structures, as well as new dualities or equivalences, starting from Morita-equivalences provided by the Lemma, following the method of section \ref{Mordual}. Without embarking on a comprehensive treatment of these more general situations (in fact, such an investigation would bring us far beyond the scope of the present paper), we limit ourselves to presenting a few examples in the next section.

\subsection{Examples}\label{addex}

In this section we give some examples of results obtained through the method of section \ref{Mordual} by applying it to Morita-equivalences arising from instances of the Comparison Lemma which are not particular cases of Theorem \ref{fund}. There are of course a great number of other dualities that can be established in a semi-automatic way by applying the technique of section \ref{Mordual}; these examples are just meant to give a flavor of the results that arise from the application of our machinery. 

We have discussed in section \ref{tarski} the representation of an atomic frame as the powerset on the collection of its atoms. Similarly, one can obtain a representation theorem for \emph{atomic distributive lattices} i.e. distributive lattices in which every element is a finite join of atoms: if $\cal D$ is such a lattice and $J_{\cal D}$ is the coherent topology on it then the Comparison Lemma yields an equivalence of toposes $\Sh({\cal D}, J_{\cal D})\simeq [At({\cal D}), \Set]$, from which it follows, invoking Corollary \ref{cor}(i), that $\cal D$ is isomorphic to the the lattice of compact elements of the powerset of its collection of atoms (equivalently, to the lattice of finite subsets of the set of its atoms). If we take $\cal K$ to be the collection of atomic distributive lattices, each of which equipped with the coherent topology, and $\cal H$ to be the collection of the (discrete) posets of the form $At({\cal D})$ for some atomic distributive lattice $\cal D$, each of which equipped with the trivial Grothendieck topology, then the functions $f:{\cal K} \to {\cal H}$ and $g:{\cal H}\to {\cal K}$ defined by setting $f({\cal D})=At({\cal D})$ (for any atomic distributive lattice $\cal D$) and $g(A)={\mathscr{P}}_{fin}(A)$ (for any set $A$) are inverse to each other (up to isomorphism), and the closure under isomorphisms of $\cal H$ can be identified with the collection of all the sets (since for any set $A$, ${\mathscr{P}}_{fin}(A)$ is an atomic distributive lattice whose set of atoms is isomorphic to $A$).

Now, the category $Ext({\cal K})$ is the category $\textbf{AtDLat}$ whose objects are the atomic distributive lattices and whose arrows ${\cal D}\to {\cal D}'$ are the distributive lattices homomorphisms $f:{\cal D}\to {\cal D}'$ between them such that the frame homomorphism $A(f):Id_{J_{\cal D}}({\cal D})\to Id_{J_{\cal D}'}({\cal D}')$ which sends an ideal $I$ of $\cal D$ to the ideal of ${\cal D}'$ generated by $f(I)$ preserves arbitrary infima. 
The category $Ext({\cal H})$ is the category $\Set_{f}$ whose objects are the sets and whose arrows $A\to A'$ are the functions $g:A\to A'$ such that the inverse image $g^{-1}:{\mathscr{P}}(A') \to {\mathscr{P}}(A)$ sends finite subsets of $A'$ to finite subsets of $A$.

Theorem \ref{dualabstract} thus yields two functors $D:\textbf{AtDLat}^{\textrm{op}}\to \Set_{f}$ and $E:\Set_{f} \to \textbf{AtDLat}^{\textrm{op}}$ which are categorical inverses to each other, and which therefore give a duality between $\textbf{AtDLat}$ and $\Set_{f}$. The functor 
\[
D:\textbf{AtDLat}^{\textrm{op}}\to \Set_{f}
\]
sends a poset $\cal D$ in $\textbf{AtDLat}$ to the set $At({\cal D})$ of its atoms, and an arrow $f:{\cal D}\to {\cal D}'$ in $\textbf{AtDLat}$ to the restriction $At({\cal D}')\to At({\cal D})$ of the left adjoint to $A(f):Id_{J_{\cal D}}({\cal D})\to Id_{J_{\cal D}'}({\cal D}')$ to the sets of atoms of ${\cal D}'$ and of ${\cal D}$. The functor 
\[
E:\Set_{f} \to \textbf{AtDLat}^{\textrm{op}}
\]
sends a set $A$ to the finite powerset ${\mathscr{P}}_{fin}(A)$ and an arrow $f:A\to A'$ in $\Set_{f}$ to the restriction ${\mathscr{P}}_{fin}(A')\to {\mathscr{P}}_{fin}(A)$ of the inverse image $f^{-1}:{\mathscr{P}}(A')\to {\mathscr{P}}(A)$.

Summarizing, we have the following result.

\begin{theorem}\label{atomfin}
Via the functors 
\[
D:\textbf{AtDLat}^{\textrm{op}}\to \Set_{f}
\]
and 
\[
E:\Set_{f} \to \textbf{AtDLat}^{\textrm{op}}
\]
defined above, the categories $\textbf{AtDLat}$ and $\Set_{f}$ are dual to each other. 
\end{theorem}\qed 
   
Notice that this duality restricts to a duality between the category of finite atomic distributive lattices and distributive lattice homomorphisms between them and the category of finite sets and functions between them. 

One might wonder why we have not considered, instead of the class of atomic distributive lattices, the class of weakly atomic meet-semilattices in which every element is a finite join of atoms. In fact, the two classes coincide with each other: if $\cal D$ is a weakly atomic meet-semilattices in which every element is a join of atoms then, by equipping it with the atomically generated topology $J_{\cal D}$, the Comparison Lemma yields an equivalence $\Sh({\cal D}, J_{{\cal D}})\simeq [At({\cal D}), \Set]$ from which it follows, invoking Corollary \ref{cor}(v), that $\cal D$ is isomorphic to the poset of atomically compact elements of ${\mathscr{P}}(At({\cal D}))$, in other words to ${\mathscr{P}}_{fin}(At({\cal D}))$.   

In fact, we can establish a more general duality, which extends that of Theorem \ref{atomfin} and represents the `finitary version' of the duality of Theorem \ref{proalex}. Given a distributive lattice $\cal D$, we say that an element $d$ of $\cal D$ is \emph{join-irreducible} if for any $a,b\in {\cal D}$, $a\vee b =d$ implies $a=d$ or $b=d$; we denote the set of join-irreducible elements of a distributive lattice $\cal D$ by $Irr({\cal D})$. We say that a distributive lattice $\cal D$ is \emph{irreducibly generated} if every element of $\cal D$ can be expressed as a finite join of join-irreducible elements. If $\cal D$ is such a lattice and $J_{\cal D}$ is the coherent topology on it then the Comparison Lemma yields an equivalence of toposes $\Sh({\cal D}, J_{\cal D})\simeq [Irr({\cal D})^{\textrm{op}}, \Set]$, from which it follows, invoking Corollary \ref{cor}(i), that $\cal D$ is isomorphic to the the lattice of compact elements of the frame $Id({Irr({\cal D})})$ of ideals on $Irr({\cal D})$.

Take $\cal K$ to be the collection of irreducibly generated distributive lattices, each of which equipped with the coherent topology, and $\cal H$ to be the collection of the posets of the form $Irr({\cal D})$ for some irreducibly generated distributive lattice, each of which equipped with the trivial Grothendieck topology; then the functions $f:{\cal K} \to {\cal H}$ and $g:{\cal H}\to {\cal K}$ defined by setting $f({\cal D})=Irr({\cal D})$ (for any atomic distributive lattice $\cal D$) and the image $g({\cal P})$ of a poset $\cal P$ in $\cal H$ under $g$ equal to the poset $Id_{comp}({\cal P})$ of compact elements of the frame $Id({\cal P})$ are inverse to each other (up to isomorphism), and the closure under isomorphisms of $\cal H$ can be identified with the collection of all the posets such that $Id_{comp}({\cal P})$ is closed under finite intersections in $Id({\cal P})$ (since for such poset $\cal P$, $Id_{comp}({\cal P})$ is an irreducibly generated distributive lattice whose poset of join-irreducible elements is isomorphic to $\cal P$). Notice that an ideal $I$ on $\cal P$ belongs to $Id_{comp}({\cal P})$ if and only if it is a finite union of principal ideals on $\cal P$, equivalently if there exists a finite set of elements $\{a_{k} \textrm{ | } k\in K\}$ of $I$ such that for any $p\in {\cal P}$, $p\in I$ if and only if $p\leq a_{k}$ for some $k\in K$. The condition that $Id_{comp}({\cal P})$ be closed under finite intersections in $Id({\cal P})$ is clearly equivalent to the requirement that the intersection in $Id({\cal P})$ of any two principal ideals on $\cal P$ be equal to the union of a finite family of principal ideals on $\cal P$, equivalently for any $a,b\in {\cal P}$ there exist a finite set of elements $\{c_{k} \textrm{ | } k\in K\}$ of $\cal P$ such that for any $p\in {\cal P}$, $p\leq a$ and $p\leq b$ if and only if $p\leq c_{k}$ for some $k\in K$. In particular, all the meet-semilattices are objects of the closure under isomorphisms of $\cal H$.
  
Now, the category $Ext({\cal K})$ is the category $\textbf{IrrDLat}$ whose objects are the irreducibly generated distributive lattices and whose arrows ${\cal D}\to {\cal D}'$ are the distributive lattices homomorphisms $f:{\cal D}\to {\cal D}'$ between them such that the frame homomorphism $A(f):Id_{J_{\cal D}}({\cal D})\to Id_{J_{\cal D}'}({\cal D}')$ which sends an ideal $I$ of $\cal D$ to the ideal of ${\cal D}'$ generated by $f(I)$ preserves arbitrary infima. 

The category $Ext({\cal H})$ is the category $\textbf{Pos}_{comp}$ whose objects are the posets $\cal P$ such that for any $a,b\in {\cal P}$ there exists a finite set of elements $\{c_{k} \textrm{ | } k\in K\}$ of $\cal P$ such that for any $p\in {\cal P}$, $p\leq a$ and $p\leq b$ if and only if $p\leq c_{k}$ for some $k\in K$, and whose arrows ${\cal P}\to {\cal P}'$ are the monotone maps $g:{\cal P}\to {\cal P}'$ such that the inverse image $g^{-1}:Id({\cal P}') \to Id({\cal P})$ sends ideals in $Id_{comp}({\cal P}')$ to ideals in $Id_{comp}({\cal P})$, equivalently, for any $q\in {\cal P}'$, there exists a finite family $\{a_{k} \textrm{ | } k\in K\}$ of elements of $\cal P$ such that for any $p\in {\cal P}$, $g(p)\leq q$ if and only if $p\leq a_{k}$ for some $k\in K$.   

Theorem \ref{dualabstract} thus yields two functors $D:\textbf{IrrDLat}^{\textrm{op}}\to \textbf{Pos}_{comp}$ and $E:\textbf{Pos}_{comp} \to \textbf{IrrDLat}^{\textrm{op}}$ which are categorical inverses to each other, and which therefore form a duality between $\textbf{IrrDLat}$ and $\textbf{Pos}_{comp}$. The functor $D:\textbf{IrrDLat}^{\textrm{op}}\to \textbf{Pos}_{comp} $ sends a poset $\cal D$ in $\textbf{IrrDLat}$ to the set $Irr({\cal D})$ of its join-irreducible elements, and an arrow $f:{\cal D}\to {\cal D}'$ in $\textbf{IrrDLat}$ to the restriction $Irr({\cal D}')\to Irr({\cal D})$ of the left adjoint to $A(f):Id_{J_{\cal D}}({\cal D})\to Id_{J_{\cal D}'}({\cal D}')$ to the sets of join-irreducible elements of ${\cal D}'$ and of ${\cal D}$. The functor $E:\textbf{Pos}_{comp} \to \textbf{IrrDLat}^{\textrm{op}}$ sends a poset $\cal P$ to $Id_{comp}({\cal P})$ and an arrow $g:{\cal P}\to {\cal P}'$ in $\textbf{Pos}_{comp}$ to the restriction $Id_{comp}({\cal P}')\to Id_{comp}({\cal P})$ of the inverse image $g^{-1}:Id({\cal P}') \to Id({\cal P})$.

Summarizing, we have the following result.

\begin{theorem}\label{birk}
Via the functors 
\[
D:\textbf{IrrDLat}^{\textrm{op}}\to \textbf{Pos}_{comp}
\]
and 
\[
E:\textbf{Pos}_{comp} \to \textbf{IrrDLat}^{\textrm{op}}
\]
defined above, the categories $\textbf{IrrDLat}$ and $\textbf{Pos}_{comp}$ are dual to each other.
\end{theorem}\qed

\begin{remark}
If $\cal D$ is a finite distributive lattice then \emph{a fortiori} $\cal D$ is irreducibly generated and has a finite set of join-irreducible elements; conversely, given a finite poset $\cal P$, the lattice $Id_{comp}({\cal P})$ is finite, and any monotone map ${\cal P}\to {\cal P}'$ between finite posets $\cal P$ and ${\cal P}'$ is an arrow in $\textbf{Pos}_{comp}$. We thus conclude that the duality of Theorem \ref{birk} restricts to Birkhoff's duality between finite distributive lattices and finite posets.
\end{remark}

Another example in the same style is given by disjunctively distributive frames.

Let $F$ be a disjunctively distributive frame. An element $a\in F$ is said to be \emph{indecomposable} if for any family $\{a_{i} \textrm{ | } i\in I\}$ of pairwise disjoint elements of $F$, $\mathbin{\mathop{\textrm{\huge $\vee$}}\limits_{i\in I}} a_{i}=a$ implies $a_{i}=a$ for some $i\in I$. Notice that, classically, an element $a\in F$ is indecomposable if and only if it is non-zero and \emph{connected} (i.e., for any elements $b,c \in F$ such that $b\wedge c=0$, $b\vee c=a$ implies that either $a=b$ or $a=c$). Indeed, to be indecomposable clearly implies to be non-zero and connected, while the converse can be proved as follows. Given a connected element $a\neq 0$, for any family $\{a_{i} \textrm{ | } i\in I\}$ of pairwise disjoint elements of $F$, $\mathbin{\mathop{\textrm{\huge $\vee$}}\limits_{i\in I}} a_{i}=a$ implies that for any $i\in I$, $a_{i} \vee (\mathbin{\mathop{\textrm{\huge $\vee$}}\limits_{j\in I, j\neq i}} a_{j})=a$, from which it follows, $a$ being connected, that either $a=a_{i}$ or $a=\mathbin{\mathop{\textrm{\huge $\vee$}}\limits_{j\in I, j\neq i}} a_{j}$. But, since $a\neq 0$, an equality of the latter kind cannot hold for all $i\in I$, since $a=\mathbin{\mathop{\textrm{\huge $\vee$}}\limits_{j\in I, j\neq i}} a_{j}$ implies $a\wedge a_{i}=0$; this shows that there is $i\in I$ such that $a=a_{i}$, as required. 

Let us suppose that $F$ is a \emph{disjunctive frame}, i.e. a disjunctively distributive frame in which every element can be written as a pairwise disjoint join of indecomposable elements. Then, denoted by $J_{F}$ the infinitary disjunctive topology on $F$ and by ${\cal I}_{F}$ the full subcategory of $F$ on the indecomposable elements of $F$, we have that ${\cal I}_{F}$ is $J_{F}$-dense; hence the Comparison Lemma yields an equivalence $\Sh(F, J_{F})\simeq \Sh({\cal I}_{F}, J_{F}|_{{\cal I}_{F}})$. But, since the objects of ${\cal I}_{F}$ are indecomposable, $J_{F}|_{{\cal I}_{F}}$ is the trivial Grothendieck topology on ${\cal I}_{F}$, and hence we have an equivalence $\Sh(F, J_{F})\simeq [{\cal I}_{F}^{\textrm{op}}, \Set]$, which implies by Corollary \ref{cor}(v), that $F$ is isomorphic to the poset of infinitarily disjunctively compact elements of the frame of  ideals on the poset ${\cal I}_{F}$ of indecomposable elements of $F$. 
 
Starting from the representations $\Sh(F, J_{F})\simeq [{\cal I}_{F}^{\textrm{op}}, \Set]$, holding for any disjunctive frame $F$, we can build a duality applying Theorem \ref{dualabstract}.

We take $\cal K$ to be the collection of all the disjunctive frames, each of which equipped with the infinitary disjunctive topology, and $\cal H$ to be the collection of all the posets of the form ${\cal I}_{F}$ for a disjunctive frame $F$, each of which equipped with the trivial topology. Define $f:{\cal K}\to {\cal H}$ as the function sending a disjunctive frame $F$ to the poset ${\cal I}_{F}$ and $g:{\cal H}\to {\cal K}$ as the function which sends a poset $P$ in $\cal H$ to the frame of infinitarily disjunctively compact elements of the frame $Id(P)$ of ideals on the poset $P$. Clearly, $f$ and $g$ are inverses to each other, up to isomorphism.   

Now, the category $Ext({\cal K})$ is the category $\textbf{DisFrm}$ whose objects are the disjunctive frames and whose arrows $F\to F'$ are the meet-semilattice homomorphisms $f:F\to F'$ between them which send pairwise disjoint joins to pairwise disjoint joins and have the property that the frame homomorphism $A(f):Id_{J_{F}}(F)\to Id_{J_{F'}}(F')$ sending an ideal $I$ of $F$ to the ideal of $F'$ generated by $f(I)$ preserves arbitrary infima. 

The category $Ext({\cal H})$ is the category $\textbf{Pos}_{dis}$ whose objects are the posets $\cal P$ such that their frames of ideals $Id({\cal P})$ have the property that the collection of infinitarily disjunctively compact elements forms a basis of them which is closed under finite meets, and whose arrows $P\to P'$ are the monotone maps $g:{\cal P}\to {\cal P}'$ such that the inverse image $g^{-1}:Id({\cal P}') \to Id({\cal P})$ sends infinitarily disjunctively compact elements to infinitarily disjunctively compact elements.

We can describe the category $\textbf{Pos}_{dis}$ more explicitly as follows. First, we note that for any poset $\cal P$, an ideal $I$ in $Id({\cal P})$ is infinitarily disjunctively compact in $Id({\cal P})$ if and only if it is a disjoint union of principal ideals on $\cal P$, equivalently if there exists a family $\{c_{i} \textrm{ | } i\in I\}$ of elements of $I$ such that for any $p\in {\cal P}$, $p\in I$ if and only if $p\leq c_{i}$ for a unique $i\in I$; indeed, the `if' direction follows from Proposition \ref{multicomposition}, while the `only if' direction follows from the fact that if $I$ is infinitarily disjunctively compact then the covering of $I$ formed by the principal ideals on $\cal P$ generated by elements in $I$ has a disjunctive refinement by principal ideals on $\cal P$ (since the invariant `to be a disjoint family' satisfies the condition in the statement of Theorem \ref{construction}). Now, all the principal ideals on $\cal P$ are infinitarily disjunctively compact elements of $Id({\cal P})$, and hence from the characterization above it follows that the condition that the intersection of two infinitarily disjunctively compact elements of $Id({\cal P})$ should be infinitarily disjunctively compact is equivalent to the requirement that the intersection of any two principal ideals on $\cal P$ should be equal to a disjoint union of principal ideals on $\cal P$, equivalently that for any $a,b\in {\cal P}$ there should exist a family $\{c_{i} \textrm{ | } i\in I\}$ of elements of $\cal P$ such that for any $p\in {\cal P}$, $p\leq a$ and $p\leq b$ if and only if $p\leq c_{i}$ for a unique $i\in I$. Note that if $\cal P$ is a meet-semilattice than this condition is always satisfied (take the family $\{c_{i} \textrm{ | } i\in I\}$ to be equal to the singleton family $\{a\wedge b\}$). Concerning the arrows in $\textbf{Pos}_{dis}$ we note that, given a monotone map $g:{\cal P}\to {\cal P}'$, the inverse image $g^{-1}:Id({\cal P}') \to Id({\cal P})$ sends infinitarily disjunctively compact elements to infinitarily disjunctively compact elements if and only if it sends any principal ideal on ${\cal P}'$ to an infinitarily disjunctively compact elements, equivalently if for any $b\in {\cal P}'$ there exists a family $\{c_{i} \textrm{ | } i\in I\}$ of elements of $\cal P$ such that for any $p\in {\cal P}$, $g(p)\leq b$ if and only if $p\leq c_{i}$ for a unique $i\in I$.

Summarizing, the category $\textbf{Pos}_{dis}$ has as objects the posets $\cal P$ such that for any $a,b\in {\cal P}$ there exists a family $\{c_{i} \textrm{ | } i\in I\}$ of elements of $\cal P$ such that for any $p\in {\cal P}$, $p\leq a$ and $p\leq b$ if and only if $p\leq c_{i}$ for a unique $i\in I$ and as arrows ${\cal P}\to {\cal P}'$ the monotone maps $g:{\cal P}\to {\cal P}'$ such that for any $b\in {\cal P}'$ there exists a family $\{c_{i} \textrm{ | } i\in I\}$ of elements of $\cal P$ such that for any $p\in {\cal P}$, $g(p)\leq b$ if and only if $p\leq c_{i}$ for a unique $i\in I$.

Our Theorem \ref{dualabstract} thus yields two functors $D:\textbf{DisFrm}^{\textrm{op}}\to \textbf{Pos}_{dis}$ and $E:\textbf{Pos}_{dis} \to \textbf{DisFrm}^{\textrm{op}}$ which are categorical inverses to each other and hence form a duality between $\textbf{DisFrm}$ and $\textbf{Pos}_{dis}$. The functor $D:\textbf{DisFrm}^{\textrm{op}}\to \textbf{Pos}_{dis}$ sends a poset $F$ in $\textbf{DisFrm}$ to the poset ${\cal I}_{F}$ of indecomposable elements of $F$ and an arrow $f:F\to F'$ in $\textbf{DisFrm}$ to the restriction ${\cal I}_{F'} \to {\cal I}_{F}$ of the left adjoint to $A(f):Id_{J_{F}}(F)\to Id_{J_{F'}}(F')$ to the sets of indecomposable elements of $F'$ and of $F$. The functor $E:\textbf{Pos}_{dis} \to \textbf{DisFrm}^{\textrm{op}}$ sends a poset $\cal P$ in $\textbf{Pos}_{dis}$ to the set ${\cal P}_{dis}$ of ideals $I$ on $\cal P$ with the property that there exists a family $\{c_{i} \textrm{ | } i\in I\}$ of elements of $I$ such that for any $p\in {\cal P}$, $p\in I$ if and only if $p\leq c_{i}$ for a unique $i\in I$, endowed with the subset-inclusion ordering, and an arrow $g:{\cal P}\to {\cal P}'$ in $\textbf{Pos}_{dis}$ to the restriction ${{\cal P}'}_{dis}\to {\cal P}_{dis}$ of the inverse image $g^{-1}:Id({\cal P}')\to Id({\cal P})$.

In conclusion, we have the following result.

\begin{theorem}
Via the functors 
\[
D:\textbf{DisFrm}^{\textrm{op}}\to \textbf{Pos}_{dis}
\]
and 
\[
E:\textbf{Pos}_{dis} \to \textbf{DisFrm}^{\textrm{op}}
\]
defined above, the categories $\textbf{DisFrm}$ and $\textbf{Pos}_{dis}$ are dual to each other. 
\end{theorem}\qed 

Note that this duality restricts to a duality between the full subcategory of $\textbf{DisFrm}^{\textrm{op}}$ on the posets $F$ such that the meet of two indecomposable elements of $F$ is indecomposable and the full subcategory of $\textbf{Pos}_{dis}$ on the meet-semilattices. 

By using the finitary version of the disjunctive topology, one can obtain a similar duality for disjunctively distributive lattices in which every element can be written as a finite pairwise disjoint join of connected elements.

For disjunctive frames with the property that the meet of any two indecomposable elements is indecomposable we can `functorialize' the representation $\Sh(F, J_{F})\simeq [{\cal I}_{F}^{\textrm{op}}, \Set]$ in a covariant way, by using Theorem \ref{dualcov}. Let $\cal K$ be the collection of such frames, each of which equipped with the infinitary disjunctive topology, and let $\cal H$ be the collection of all the posets (in fact, meet-semilattices) of the form ${\cal I}_{F}$ for a disjunctive frame $F$ in $\cal K$, each of which equipped with the trivial topology. Let $f:{\cal K}\to {\cal H}$ and $g:{\cal H}\to {\cal K}$ be defined exactly as in the example above. 

Using the covariant method for `functorializing' Morita-equivalences, we obtain the categories $Ext({\cal K})$ and $Ext({\cal H})$ defined as follows. 

The category $Ext({\cal K})$ coincides with the category $\textbf{DIFrm}$ having as objects the disjunctive frames $F$ such that the meet in $F$ of any two indecomposable elements of $F$ is indecomposable and as arrows $F\to F'$ the meet-semilattice homomorphisms $f:F\to F'$ between them which send pairwise disjoint joins to pairwise disjoint joins and have the property that the frame homomorphism $A(f):Id_{J_{F}}(F)\to Id_{J_{F'}}(F')$ sending an ideal $I$ in $Id_{J_{F}}(F)$ to the $J_{F'}$ ideal on $F'$ generated by $f(I)$ sends supercompact elements to supercompact elements. We note that, for any disjunctive frame $F$, since $0_{F}$ is not indecomposable in $F$, if the meet in $F$ of two indecomposable elements is indecomposable then it is non-zero; this implies that the disjunctive frames in $F$ can be identified with the meet-semilattices $F$ with a bottom element $0_{F}$ which have the property that for any $a, b\in F$, $a\wedge b=0$ implies that either $a=0$ or $b=0$; under this identification, the arrows $F\to F'$ in $\textbf{DIFrm}$ correspond precisely to the meet-semilattice homomorphisms $F\to F'$ which send $0_{F}$ to $0_{F'}$ and any non-zero element of $F$ to a non-zero element of $F'$.   

Let us denote by $\textbf{MSLat}^{\ast}$ the category having as objects the meet-semilattices $\cal P$ which have a bottom element $0_{\cal P}$ and satisfy the property that for any $a, b\in F$, $a\wedge b=0$ implies that either $a=0$ or $b=0$, and as arrows ${\cal P}\to {\cal P}'$ the meet-semilattice homomorphisms ${\cal P}\to {\cal P}'$ which send $0_{{\cal P}}$ to $0_{{\cal P}'}$ and any non-zero element of ${\cal P}$ to a non-zero element of ${\cal P}'$.     

The category $Ext({\cal H})$ coincides with the category $\textbf{MPos}$ which has as objects the meet-semilattices $\cal P$ such that their frames of ideals $Id({\cal P})$ have the property that the collection of infinitarily disjunctively compact elements forms a basis of them which is closed under finite meets, and whose arrows $P\to P'$ are the meet-semilattice homomorphisms $g:{\cal P}\to {\cal P}'$ such that the homomorphism $A(g):Id({\cal P}) \to Id({\cal P}')$ sends infinitarily disjunctively compact elements to infinitarily disjunctively compact elements. From the characterization of the category $\textbf{Pos}_{dis}$ obtained above, we thus conclude that the category $\textbf{MPos}$ coincides with the category $\textbf{MSLat}$ of meet-semilattices and meet-semilattice homomorphisms between them.  

Theorem \ref{dualcov} thus yields two functors 
\[
D:\textbf{MSLat}^{\ast}\to \textbf{MSLat} 
\]
and 
\[
E:\textbf{MSLat} \to \textbf{MSLat}^{\ast}
\]
which are categorical inverses to each other. The functor $D:\textbf{MSLat}^{\ast}\to \textbf{MSLat}$ sends a meet-semilattice ${\cal P}$ in $\textbf{MSLat}^{\ast}$ to the meet-semilattice ${\cal P}^{\ast}$ of its non-zero elements, and an arrow $f:{\cal P}\to {\cal P}'$ in $\textbf{MSLat}^{\ast}$ to its restriction ${\cal P}^{\ast} \to {{\cal P}'}^{\ast}$. The functor $E:\textbf{MSLat} \to \textbf{MSLat}^{\ast}$ sends a meet-semilattice $\cal P$ to the set $Id_{\ast}({\cal P})$ of the ideals on $\cal P$ which are either principal or empty, endowed with the subset-inclusion ordering, and a meet-semilattice homomorphism $g:{\cal P}\to {\cal P}'$ to the restriction $Id_{\ast}({\cal P}) \to Id_{\ast}({\cal P}')$ of $A(g):Id({\cal P})\to Id({\cal P}')$.  

Summarizing, we have the following result.

\begin{theorem}
Via the functors 
\[
D:\textbf{MSLat}^{\ast}\to \textbf{MSLat} 
\]
and 
\[
E:\textbf{MSLat} \to \textbf{MSLat}^{\ast}
\]
defined above, the categories $\textbf{MSLat}^{\ast}$ and $\textbf{MSLat}$ are equivalent to each other. 
\end{theorem}\qed 

As another example, we construct a duality for a natural class of preframes including all the algebraic lattices.

Given a preframe $\cal D$, we say that an element $d\in {\cal D}$ is \emph{directedly irreducible} if any directed sieve on $d$ is maximal, i.e. if for any directed family $\{a_{i} \textrm{ | } i\in I\}$ of elements of $\cal D$ such that $d=\mathbin{\mathop{\textrm{\huge $\vee$}}\limits_{i\in I}}a_{i}$ there exists $i\in I$ such that $d=a_{i}$; given a preframe $\cal D$, we denote by $DirIrr({\cal D})$ the poset of directedly irreducible elements of $\cal D$. 

We shall call the preframes in which every element is a directed join of directedly irreducible elements the \emph{directedly generated} preframes. For any directedly generated preframe $\cal D$, the Comparison Lemma yields an equivalence of toposes $\Sh({\cal D}, J_{\cal D})\simeq [DirIrr({\cal D})^{\textrm{op}}, \Set]$, where $J_{\cal D}$ is the directed topology on $\cal D$ (cf. section \ref{charinv} above). Let us now show that the invariant `to be a directed family' satisfies the hypotheses of Theorem \ref{construction} with respect to the class of directedly generated preframes $\cal D$, each of which equipped with the directed topology $J_{\cal D}$. For any directedly generated preframe $\cal D$, the directed topology $J_{\cal D}$ is easily seen to be $C$-induced relative to $P$ where $P$ is the following property of embeddings $i:B_{L}\hookrightarrow L$: `$i$ preserves finite meets' or equivalently (since the canonical embedding ${\cal C}\hookrightarrow Id_{J}({\cal C})$ satisfies the property that ${\cal C}$ is closed under finite meets in $Id_{J}({\cal C})$) the property `$i(B_{L})$ is closed in $L$ under finite meets'. Also, since $Id_{J_{\cal D}}({\cal D})\cong Id(DirIrr({\cal D}))$ via the Comparison Lemma, any family $\cal F$ of principal $J_{\cal D}$-ideals on $\cal D$, $\cal F$ has a directed refinement (if and) only if it has a directed refinement made of principal $J_{\cal D}$-ideals on $\cal D$. Indeed, any principal $J_{\cal D}$-ideal $(d)\downarrow_{J_{\cal D}}$ on $\cal D$ corresponds via the isomorphism $Id_{J_{\cal D}}({\cal D})\cong Id(DirIrr({\cal D}))$ to the set $(d)\downarrow=\{a\in DirIrr({\cal D}) \textrm{ | } a\leq d\}$, and it is immediate to see that if a family $\cal F$ of ideals on $DirIrr({\cal D})$ of the form $(d)\downarrow$ for $d\in {\cal D}$ admits a refinement in $Id(DirIrr({\cal D}))$ by a directed family then $\cal F$ itself must be directed. Thus, if we denote by $\cal K$ the collection of all the directedly generated preframes, Theorem \ref{construction} ensures that the invariant `to be directed' is both $\cal K$-adequate and $\cal K$-compatible; this implies in particular that any directedly generated preframe $\cal D$ is isomorphic to the poset of directedly compact elements (cf. section \ref{charinv} above) of the frame $Id_{J_{\cal D}}({\cal D}) \cong Id(DirIrr({\cal D}))$. 

Now, take $\cal K$ to be the collection of directedly generated preframes, each of which equipped with the directed topology, and $\cal H$ to be the collection of the posets of the form $DirIrr({\cal D})$ for some directedly generated preframe, each of which equipped with the trivial Grothendieck topology. The functions $f:{\cal K} \to {\cal H}$ and $g:{\cal H}\to {\cal K}$ defined by setting $f({\cal D})=DirIrr({\cal D})$ (for any atomic distributive lattice $\cal D$) and $g({\cal P})$ equal to the poset $Id_{dir}({\cal P})$ of directedly compact elements of the frame $Id({\cal P})$ are inverse to each other (up to isomorphism), and that the closure of $\cal H$ under isomorphisms can be identified with the collection of all the posets $\cal P$ such that the poset $Id_{dir}({\cal P})$ is closed in $Id({\cal P})$ under finite meets (since for any poset $\cal P$, $Id_{dir}({\cal P})$ is an irreducibly generated distributive lattice whose poset of join-irreducible elements is isomorphic to $\cal P$).
  
The category $Ext({\cal K})$ is the category $\textbf{DirIrrPFrm}$ whose objects are the directedly generated preframes and whose arrows ${\cal D}\to {\cal D}'$ are the preframe homomorphisms $f:{\cal D}\to {\cal D}'$ between them such that the frame homomorphism $A(f):Id_{J_{\cal D}}({\cal D})\to Id_{J_{{\cal D}'}}({\cal D}')$ which sends an ideal $I$ of $\cal D$ to the ideal on ${\cal D}'$ generated by $f(I)$ preserves arbitrary infima. 

The category $Ext({\cal H})$ is the category $\textbf{Pos}_{dir}$ whose objects are the posets $\cal P$ such that $Id_{dir}({\cal P})$ is closed in $Id({\cal P})$ under finite meets and whose arrows ${\cal P}\to {\cal P}'$ are the monotone maps $g:{\cal P}\to {\cal P}'$ between them such that the inverse image $g^{-1}:Id({\cal P}') \to Id({\cal P})$ sends ideals in $Id_{dir}({\cal P}')$ to ideals in $Id_{dir}({\cal P})$.

We can describe the category $\textbf{Pos}_{dir}$ in more concrete terms, as follows. First, we note that for any poset $\cal P$, an ideal $I$ in $Id({\cal P})$ is directedly compact in $Id({\cal P})$ if and only if it is a directed union of principal ideals on $\cal P$, equivalently if $I$ is directed, i.e. it is non-empty and for any $a,b\in I$ there exists $c\in I$ such that $a\leq c$ and $b\leq c$. Indeed, the `if' direction follows from Proposition \ref{multicomposition}, while the `only if' direction follows from the fact that if $I$ is directedly compact then the covering of $I$ formed by the principal ideals on $\cal P$ generated by an element of $I$ has a directed refinement by principal ideals on $\cal P$ (since the invariant `to be directed' satisfies the condition in the statement of Theorem \ref{construction}, cf. above) and hence $I$ itself must be directed. Now, all the principal ideals on $\cal P$ are directedly compact elements of $Id({\cal P})$, and hence from the characterization above it follows that the condition that the intersection of two directedly compact elements of $Id({\cal P})$ should be directedly compact is equivalent to the requirement that the intersection of any two principal ideals on $\cal P$ should be a directed ideal, equivalently that for any $a,b\in {\cal P}$ there should be $c\in {\cal P}$ such that $c\leq a$ and $c\leq b$, and for any elements $d,e\in {\cal P}$ such that $d,e\leq a$ and $d,e\leq b$ there should exist $z\in {\cal P}$ such that $z\leq a$, $z\leq b$, $d,e \leq z$. Note that if $\cal P$ is a meet-semilattice than this condition is always satisfied (the intersection of $(a)\downarrow$ and $(b)\downarrow$ is equal to $(a\wedge b)\downarrow$). 

Concerning the arrows in $\textbf{Pos}_{dir}$ we note that, given a monotone map $g:{\cal P}\to {\cal P}'$, the inverse image $g^{-1}:Id({\cal P}') \to Id({\cal P})$ sends directedly compact elements to directedly compact elements if and only if for any $b\in {\cal P}'$, $g^{-1}((b)\downarrow)$ is a directed ideal on $\cal P$, i.e. there exists $a\in {\cal P}$ such that $g(a)\leq b$ and for any $u,v\in {\cal P}$ such that $g(u)\leq b$ and $g(v)\leq b$ there exists $z\in {\cal P}$ such that $u,v\leq z$ and $g(z)\leq b$. 

Summarizing, the category $\textbf{Pos}_{dir}$ has as objects the posets $\cal P$ such that for any $a,b\in {\cal P}$ there is $c\in {\cal P}$ such that $c\leq a$ and $c\leq b$ and for any elements $d,e\in {\cal P}$ such that $d,e\leq a$ and $d,e\leq b$ there exists $z\in {\cal P}$ such that $z\leq a$, $z\leq b$, $d,e \leq z$, and as arrows ${\cal P}\to {\cal P}'$ the monotone maps $g:{\cal P}\to {\cal P}'$ with the property that for any $b\in {\cal P}'$ there exists $a\in {\cal P}$ such that $g(a)\leq p$ and for any two $u,v\in {\cal P}$ such that $g(u)\leq b$ and $g(v)\leq b$ there exists $z\in {\cal P}$ such that $u,v\leq z$ and $g(z)\leq b$.

Our Theorem \ref{dualabstract} thus yields two functors $D:\textbf{DirIrrPFrm}^{\textrm{op}}\to \textbf{Pos}_{dir}$ and $E:\textbf{Pos}_{dir} \to \textbf{DirIrrPFrm}^{\textrm{op}}$ which are categorical inverses to each other and hence form a duality between $\textbf{DirIrrPFrm}$ and $\textbf{Pos}_{dis}$. The functor $D:\textbf{DirIrrPFrm}^{\textrm{op}}\to \textbf{Pos}_{dir}$ sends a preframe ${\cal D}$ in $\textbf{DirIrrPFrm}$ to the poset $DirIrr({\cal D})$ of directedly irreducible elements of ${\cal D}$ and an arrow $f:{\cal D}\to {\cal D}'$ in $\textbf{DirIrrPFrm}$ to the restriction $DirIrr({\cal D}) \to DirIrr({\cal D})$ of the left adjoint to $A(f):Id_{J_{{\cal D}}}({\cal D})\to Id_{J_{{\cal D}'}}({\cal D}')$ to the sets of directedly irreducible elements of ${\cal D}'$ and of ${\cal D}$. The functor $E:\textbf{Pos}_{dir} \to \textbf{DirIrrPFrm}^{\textrm{op}}$ sends a poset $\cal P$ in $\textbf{Pos}_{dir}$ to the poset $Id_{dir}({\cal P})$ of directed ideals on $\cal P$, and an arrow $g:{\cal P}\to {\cal P}'$ in $\textbf{Pos}_{dir}$ to the restriction $Id_{dir}({\cal P}')\to Id_{dir}({\cal P})$ of the inverse image $g^{-1}:Id({\cal P}')\to Id({\cal P})$.

In conclusion, we have the following result.

\begin{theorem}\label{dirirr}
Via the functors 
\[
D:\textbf{DirIrrPFrm}^{\textrm{op}}\to \textbf{Pos}_{dir}
\]
and 
\[
E:\textbf{Pos}_{dir} \to \textbf{DirIrrPFrm}^{\textrm{op}}
\]
defined above, the categories $\textbf{DirIrrPFrm}$ and $\textbf{Pos}_{dir}$ are dual to each other. 
\end{theorem}\qed 

Given a poset $\cal P$ in $\textbf{Pos}_{dir}$, if $\cal P$ has binary joins and a bottom element $0_{\cal P}$ then the condition that an ideal $I$ on $\cal P$ be directed can be reformulated as the requirement that $I$ be non-empty and that for any $a,b\in I$, $a\vee b\in I$. Moreover, it is easy to verify that for any monotone map $g:{\cal P}\to {\cal P}'$ between posets in $\textbf{Pos}_{dir}$ having a bottom element and binary joins, $g(0_{\cal P})=0_{{\cal P}'}$ and $g$ sends binary joins in $\cal P$ to binary joins in ${\cal P}'$ if and only if $g$ is an arrow ${\cal P}\to {\cal P}'$ in the category $\textbf{Pos}_{dir}$. Therefore, the category $\textbf{SSLat}$ of sup-semilattices and sup-semilattice homomorphisms between them can be identified with a full subcategory of $\textbf{Pos}_{dir}$, and hence the duality of Theorem \ref{dirirr} restricts to a duality between the full subcategory $\textbf{DirIrrPFrm}_{s}$ of $\textbf{DirIrrPFrm}$ on the preframes $\cal D$ such that $DirIrr{\cal D}$ is a sup-semilattice and the category $\textbf{SSLat}$. That is, we have the following result.

\begin{theorem}
Via the restrictions 
\[
D|:\textbf{DirIrrPFrm}_{s}^{\textrm{op}}\to \textbf{SSLat}
\]
and 
\[
E|:\textbf{SSLat} \to \textbf{DirIrrPFrm}_{s}^{\textrm{op}}
\]
of the functors 
\[
D:\textbf{DirIrrPFrm}^{\textrm{op}}\to \textbf{Pos}_{dir}
\]
and 
\[
E:\textbf{Pos}_{dir} \to \textbf{DirIrrPFrm}^{\textrm{op}}
\]
defined above, the categories $\textbf{DirIrrPFrm}_{s}$ and $\textbf{SSLat}$ are dual to each other. 
\end{theorem}

\begin{remark}
The restriction $E|:\textbf{SSLat}\to \textbf{Pos}^{\textrm{op}}$ of the functor $E$ coincides with the functor from $\textbf{SSLat}$ to $\textbf{Pos}^{\textrm{op}}$ giving one half of the duality of Theorem 1.5 \cite{duality} between algebraic lattices and sup-semilattices; from this it follows that the category $\textbf{DirIrrPFrm}_{s}$ coincides with the category of algebraic lattices defined in \cite{duality} (since both categories coincide, by Theorem \ref{dirirr} and the duality theorem of \cite{duality}, with the extended image of the same functor, namely $E|:\textbf{SSLat}\to \textbf{Pos}^{\textrm{op}}$).     
\end{remark}

As a final example, consider frames in which every element is a join of compact elements; by the Comparison Lemma, they can be represented as frames of ideals (i.e. lower subsets which are closed under finitary joins) on the join-semilattice of their compact elements. 

Of course, these examples do not exhaust at all the possibilities of application of our techniques; they are just meant to show that a great variety of insights on preordered structures and topological spaces can be easily obtained by applying our methods. The reader is invited to build his or her favorite dualities or representation theorems applying our techniques to his or her cases of interest.

\section{Adjunctions}\label{adj}

Our approach to Stone-type dualities described so far provides us with a natural way of building adjunctions between categories of preorders and categories of locales or topological spaces which restrict, on appropriate subcategories, to the categorical equivalences established in the previous sections. Again, this represents an application of the philosophy `toposes as bridges' of \cite{OC10} in the context of the Morita-equivalence given by Theorem \ref{fund}.
 
\subsection{Frames presented by sites} 
  
Starting from the equivalence 
\[
\Sh({\cal C}, J)\simeq \Sh(Id_{J}({\cal C}))
\]
of Theorem \ref{fund}, we can investigate the behaviour of a particular class of invariants in relation to such equivalence, namely geometric morphisms from localic toposes to the given topos.

First, let us analyze the behavior of this invariant with respect to the first site of definition $({\cal C}, J)$ of the topos. For any locale $L$, the (isomorphism classes of) geometric morphisms $\Sh(L)\to \Sh({\cal C}, J)$ correspond exactly, by Diaconescu's equivalence, to the flat $J$-continuous functors from $\cal C$ to $\Sh(L)$; note that these latter functors always take values in the frame $L$ of subterminals in $\Sh(L)$ (cf. section \ref{logical}).

Next, let us describe the behavior of the invariant with respect to the second site of definition for the topos. For any locale $L$, the geometric morphisms $\Sh(L)\to \Sh(Id_{J}({\cal C}))$ correspond to the frame homomorphisms $Id_{J}({\cal C}) \to L$ (cf. Proposition C1.4.5 \cite{El}).

Therefore, we can conclude that the frame homomorphisms $Id_{J}({\cal C}) \to L$ correspond bijectively to the flat $J$-continuous functors ${\cal C}\to \Sh(L)$, under the bijection sending a frame homomorphism $Id_{J}({\cal C})\to L$ to the composite of it with the canonical map $l:{\cal C}\to Id_{J}({\cal C})$. The maps $f:{\cal C}\to L$ such that the composite $y\circ f:{\cal C}\to \Sh(L)$ of $f$ with the Yoneda embedding $y:L\to \Sh(L)$ is a flat functor ${\cal C}\to \Sh(L)$, which we call the \emph{filtering} maps, can be characterized explicitly.

\begin{proposition}\label{filtering}
Let $\cal C$ be a preorder, $L$ be a frame and $f:{\cal C}\to L$ be a monotone map. Then $f$ is filtering if and only if
\begin{enumerate}[(i)]
\item $1_{L}=\mathbin{\mathop{\textrm{\huge $\vee$}}\limits_{c\in {\cal C}}}f(c)$;

\item For any $c,c'\in {\cal C}$, $f(c)\wedge f(c')=\mathbin{\mathop{\textrm{\huge $\vee$}}\limits_{b\in B_{c, c'}}}f(b)$ where $B_{c, c'}$ is the set
\[
\{b\in {\cal C} \textrm{ | } b\leq c \textrm{ and } b\leq c'\}.
\]  
\end{enumerate}
In particular, if $\cal C$ is a meet-semilattice, $f:{\cal C}\to L$ is filtering if and only if it is a meet-semilattice homomorphism ${\cal C}\to L$.
\end{proposition} 

\begin{proofs}
The proposition immediately follows from the characterization of flat functors as filtering functors given by Theorem VII 10.1 \cite{MM}.
\end{proofs}

Note that a filtering map $f:{\cal C}\to L$ corresponds to a \emph{$J$-continuous} flat functor ${\cal C}\to \Sh(L)$ as above if and only if it sends $J$-covering sieves to covering families in $L$.

The argument above, combined with Proposition \ref{filtering}, thus yields the following result.

\begin{theorem}\label{freeframes}
Let $\cal C$ be a preorder and $J$ be a Grothendieck topology on $\cal C$. Then the frame $Id_{J}({\cal C})$, together with the map $\eta:{\cal C}\to Id_{J}({\cal C})$ sending an element $c\in {\cal C}$ to the principal ideal $(c)\downarrow_{J}$, satisfies the following universal property: for any  map $f:{\cal C}\to L$ to a frame $L$, $f$ is filtering and sends every $J$-covering sieve to a covering family in $L$ if and only if there is a (unique) frame homomorphism $\tilde{f}:Id_{J}({\cal C})\to L$ such that $\tilde{f}\circ \eta=f$ (given by the formula $\tilde{f}(I)=\mathbin{\mathop{\textrm{\huge $\vee$}}\limits_{c\in I}}f(c)$ for any $I\in Id_{J}({\cal C})$).
\end{theorem}\qed

\begin{remark}
The particular case of Theorem \ref{freeframes} when $\cal C$ is a meet-semilat-\\tice appears as Proposition II.2.11 of \cite{stone}; the reader might find it interesting to compare the simple topos-theoretic argument which proves Theorem \ref{freeframes} with the proof of the result for meet-semilattices in \cite{stone}, which consists of direct but elaborate technical arguments involving frames (cf. also section \ref{appendix} for an elementary proof of Theorem \ref{freeframes}).  
\end{remark}

Provided that some natural hypotheses are satisfied, we can read the universal property of $Id_{J}({\cal C})$ stated in Theorem \ref{freeframes} as an adjunction between the category \textbf{Frm} of frames and a category having as objects the preorders $\cal C$. To this end, suppose to have, as in section \ref{general} above, a category $\cal K$ of posets $\cal C$, each of which equipped with a subcanonical Grothendieck topology $J_{{\cal C}}$. If $\cal K$ contains, among its objects, all the frames, and for any structure $\cal C$ in $\cal K$ and any frame $L$ in $\cal K$ the arrows ${\cal C}\to L$ in $\cal K$ coincide with the filtering maps ${\cal C}\to L$ which send $J_{\cal C}$-covering sieves to a covering families in $L$, then we have a forgetful functor $U_{\cal K}:\textbf{Frm}\to {\cal K}$, and the following result holds.

\begin{theorem}
Let $\cal K$ satisfy the hypotheses above. Then the functor $A:\cal K \to \textbf{Frm}$ sending an object $\cal C$ of $\cal K$ to $Id_{J}({\cal C})$ and an arrow $f:{\cal C}\to {\cal D}$ in $\cal K$ to the frame homomorphism $A(f):Id_{J_{\cal C}}({\cal C})\to Id_{J_{\cal D}}({\cal D})$ sending any ideal $I$ in $Id_{J_{\cal C}}({\cal C})$ to the $J_{\cal D}$-ideal on $\cal D$ generated by the image $f(I)$ of $I$ under $f$ is left adjoint to the forgetful functor $U_{\cal K}:\textbf{Frm}\to {\cal K}$.   
\end{theorem}

\begin{proofs}
Under the hypotheses of the proposition, the arrows $f:{\cal C} \to U_{\cal K}(L)$ in $\cal K$ can be identified with the filtering maps ${\cal C}\to L$ which send $J_{\cal C}$-covering sieves to covering families in $L$, and these maps in turn correspond to the frame homomorphisms $Id_{J}({\cal C}) \to L$ (by Theorem \ref{freeframes}). It is easy to check that this bijection is natural in ${\cal C}\in {\cal K}$ and in $L\in \textbf{Frm}$ and hence yields an adjunction between $A$ and $U_{\cal K}$ in which $U_{\cal K}$ is the left adjoint. The naturality in $L$ is obvious, while the naturality in $\cal C$ follows from the fact that any arrow $f$ in $\cal K$ is naturally isomorphic to the restriction of $A(f)$ to the subsets of principal ideals on $\cal C$ and on ${\cal C}'$. 
\end{proofs}

As particular cases of the theorem we recover:

\begin{enumerate}[(i)]
\item The reflection from the category of meet-semilattices to the category of frames (take $\cal K$ to be the category of meet-semilattices, each of which equipped with the trivial Grothendieck topology);

\item The reflection from the category of distributive lattices to the category of frames (take $\cal K$ to be the category of distributive lattices, each of which equipped with the coherent topology);

\item The reflection from the category of preframes to the category of frames (take $\cal K$ to be the category of preframes, each of which equipped with the directed topology, cf. section \ref{charinv}).   
\end{enumerate}

As a novel application of the theorem we obtain a reflection from the category of disjunctively distributive frames to the category of frames (take $\cal K$ to be the category of disjunctively distributive frames, each of which equipped with the infinitary disjunctive topology). Also, for any regular cardinal $k$, we obtain a reflection from the category of $k$-frames (cf. section \ref{charinv} above) to the category of frames.

As another application of our usual method (of transferring invariants, in this case geometric morphisms involving localic toposes, across different sites of definition of a given topos), we establish an adjunction between the opposite of the category $\textbf{Bool}$ of Boolean algebras and the category $\textbf{Loc}$ of locales, which restricts to the equivalence between $\textbf{Bool}^{\textrm{op}}$ and the full subcategory of the category of coherent locales on the locales which have a basis of complemented elements (cf. section \ref{stonedist} above).

Our starting point is, as above, the Morita-equivalence   
\[
\Sh({\cal C}, J)\simeq \Sh(Id_{J}({\cal C}))
\]
of Theorem \ref{fund}, where we suppose $\cal C$ to be a Boolean algebra (regarded as a preorder coherent category) and $J$ to be the coherent topology on it. 

For any locale $L$, the (isomorphism classes of) geometric morphisms
\[
\Sh(L)\to \Sh(Id_{J}({\cal C}))
\]
correspond exactly to the frame homomorphisms $Id_{J}({\cal C})\to L$, while the (isomorphism classes of) geometric morphisms 
\[
\Sh(L)\to \Sh({\cal C}, J)
\]
correspond precisely to the meet-semilattice homomorphisms ${\cal C}\to L$ which send $J$-covering sieves to covering families, in other words to the lattice homomorphisms ${\cal C}\to L$. 

Now, since $\cal C$ is a Boolean algebra, every lattice homomorphism ${\cal C}\to L$ takes values in the sublattice $L_{c}$ of $L$ consisting of the complemented elements of $L$ and hence the arrows $L \to Id_{J}({\cal C})$ in \textbf{Loc} are in bijection with the Boolean algebra homomorphisms ${\cal C}\to L_{c}$; it is immediate to verify that this bijection is natural in ${\cal C}\in \textbf{Bool}$ and in $L\in \textbf{Loc}$. Thus, define $c:\textbf{Loc}\to \textbf{Bool}^{\textrm{op}}$ as the functor sending a locale $L$ to the lattice $L_{c}$ of complemented elements of $L$ and a frame homomorphism $L\to L'$ to its restriction $L_{c}\to L'_{c}$, and $Id:\textbf{Bool}^{\textrm{op}}\to \textbf{Loc}$ as the functor sending a Boolean algebra $\cal C$ to the frame $Id_{J_{\cal C}}({\cal C})$ of ideals of $\cal C$ (where $J_{\cal C}$ is the coherent topology on $\cal C$) and a morphism $f:{\cal C}\to {\cal C}'$ in $\textbf{Bool}$ to the frame homomorphism $Id_{J_{\cal C}}({\cal C})\to Id_{J_{\cal C}'}({\cal C}')$ sending an ideal $I$ in $Id_{J_{\cal C}}({\cal C})$ to the ideal in $Id_{J_{\cal C}'}({\cal C}')$ generated by the image $f(I)$ of $I$ under $f$.
  
We have thus established the following result.

\begin{theorem}
The functors $Id:\textbf{Bool}^{\textrm{op}}\to \textbf{Loc}$ and $c:\textbf{Loc}\to \textbf{Bool}^{\textrm{op}}$ defined above are adjoint to each other, where $c$ is the left adjoint and $Id$ is the right adjoint.
\end{theorem}\qed

\begin{remark}
Composing this adjunction with the well-known adjunction between locales and spaces yields the usual Stone adjunction between the opposite of the category of Boolean algebras and the category of topological spaces. 
\end{remark}

\subsection{Disjunctive and atomic frames}

Below, we shall adopt the terminology of section \ref{addex}.

Given a disjunctively distributive frame $F$, and denoted by ${\cal I}_{F}$ the collection of its indecomposable elements, we have a map $\phi_{F}: F\to Id({\cal I}_{F})$ sending an element $a\in F$ to the set $\{b\in {\cal I}_{F} \textrm{ | } b\leq a\}$. Notice that $\phi_{F}$ takes values in the subset of infinitarily disjunctively compact elements of $Id({\cal I}_{F})$; indeed, for any $a \in F$, $\phi_{F}(a)$ is supercompact and hence in particular infinitarily disjunctively compact.    

\begin{theorem}\label{disthm}

\begin{enumerate}[(i)]
\item For any disjunctively distributive frame $F$ such that the subset $Dis(Id({\cal I}_{F}))$ of infinitarily disjunctively compact elements of the frame $Id({\cal I}_{F})$ is closed in $Id({\cal I}_{F})$ under finite meets, $Dis(Id({\cal I}_{F}))$ is a disjunctive frame with the induced order, and the map $\phi_{F}:F\to Dis(Id({\cal I}_{F}))$ is a disjunctively distributive frame homomorphism which restricts to an isomorphism from ${\cal I}_{F} \subseteq F$ to the subset ${\cal I}_{Dis(Id({\cal I}_{F}))}$ of indecomposable elements of $Dis(Id({\cal I}_{F}))$;

\item A disjunctively distributive frame $F$ is disjunctive if and only if the subset $Dis(Id({\cal I}_{F}))$ of infinitarily disjunctively compact elements of the frame $Id({\cal I}_{F})$ is closed in $Id({\cal I}_{F})$ under finite meets and the map\\ $\phi_{F}:F\to Dis(Id({\cal I}_{F}))$ is an isomorphism (of disjunctively distributive frames);

\item A topological space $X$ is locally connected if and only if the frame ${\cal O}(X)$ of open sets of $X$ is isomorphic via the map $\phi_{{\cal O}(X)}$ to the poset $Dis(Id({\cal I}_{{\cal O}(X)}))$ of infinitarily disjunctively compact elements of the\\ frame of ideals on the poset of non-empty connected open subsets of $X$ and the (set-theoretic) intersection of any two infinitarily disjunctively compact elements of $Id({\cal I}_{{\cal O}(X)})$ is infinitarily disjunctively compact.    
\end{enumerate}

\end{theorem}

\begin{proofs}
$(i)$ By Proposition \ref{multicomposition}, $Dis(Id({\cal I}_{F}))$ is closed in $Id({\cal I}_{F})$ under pairwise disjoint joins and, by our hypothesis, it is closed in $Id({\cal I}_{F})$ under finite meets. Therefore $Dis(Id({\cal I}_{F}))$ is, with the order induced by that of $Id({\cal I}_{F})$, a disjunctively distributive frame. The fact that $\phi_{F}:F\to Dis(Id({\cal I}_{F}))$ is a disjunctively distributive frame homomorphism is immediate to verify. Now, let us prove that $Dis(Id({\cal I}_{F}))$ is disjunctive. If $I$ is an infinitarily disjunctively compact element of the frame $Id({\cal I}_{F})$ then the covering of $I$ in $Id({\cal I}_{F})$ formed by the principal ideals generated by the elements of $I$ has a refinement consisting of pairwise disjoint elements; clearly, the elements of this refining family must all be principal ideals, from which it follows, these ideals being supercompact, and hence indecomposable, elements of $Id({\cal I}_{F})$, that $I$ can be expressed as a pairwise disjoint join in $Id({\cal I}_{F})$ of indecomposable elements of $Id({\cal I}_{F})$. Now, by Proposition \ref{multicomposition}, $Dis(Id({\cal I}_{F}))$ is closed in $Id({\cal I}_{F})$ under pairwise disjoint joins, and from the fact that $Dis(Id({\cal I}_{F}))$ is closed in $Id({\cal I}_{F})$ under finite meets it follows that two elements of $Dis(Id({\cal I}_{F}))$ are pairwise disjoint in $Dis(Id({\cal I}_{F}))$ if and only if they are pairwise disjoint in $Id({\cal I}_{F})$; therefore $I$ can be expressed as a pairwise disjoint join in $Dis(Id({\cal I}_{F}))$ of indecomposable elements of $Dis(Id({\cal I}_{F}))$. This proves that $Dis(Id({\cal I}_{F}))$ is a disjunctive frame; moreover, the argument shows that the indecomposable elements of $Dis(Id({\cal I}_{F}))$ are precisely those of the form $\phi_{F}(a)$ for $a\in Id({\cal I}_{F})$. 

$(ii)$ We have already proved the `only if' direction (cf. the discussion above), so it remains to prove the `if' one. Let us suppose that $\phi_{F}:F\to Dis(Id({\cal I}_{F}))$ is an isomorphism of posets. By part $(i)$ $Dis(Id({\cal I}_{F}))$ is disjunctive and hence $F$, being isomorphic to it, is disjunctive as well, as required.     

$(iii)$ This follows immediately from part $(ii)$ by recalling that a topological space is (classically) locally connected if and only if every open set is a disjoint union of non-empty connected open sets.
\end{proofs}

We now proceed to show that the construction above of a disjunctive frame consisting of the infinitarily disjunctively compact elements of a frame $Id({\cal I}_{F})$ can be naturally made into an adjunction.

Given a disjunctively distributive frames $F$ such that the subset of infinitarily disjunctively compact elements of $Id({\cal I}_{F})$ is closed in $Id({\cal I}_{F})$ under finite meets, the composite $y \circ \phi_{F}:F\to [{{\cal I}_{F}}^{\textrm{op}}, \Set]$ of the map $\phi_{F}:F\to Id({\cal I}_{F})$ with the Yoneda embedding $y:Id({\cal I}_{F}) \to \Sh(Id({\cal I}_{F}))\simeq [{{\cal I}_{F}}^{\textrm{op}}, \Set]$ preserves finite meets and sends $J_{F}$-covering sieves to covering families. Indeed, given a $J_{F}$-covering sieve $\{a_{i} \leq a \textrm{ | } \in I\}$ on an object $a$ of $F$, $\phi_{F}(\mathbin{\mathop{\textrm{\huge $\vee$}}\limits_{i\in I}}a_{i})=\mathbin{\mathop{\textrm{\huge $\cup$}}\limits_{i\in I}}\phi_{F}(a_{i})$ since for any $b\in \phi_{F}(a)$, $b=b\wedge a=\mathbin{\mathop{\textrm{\huge $\vee$}}\limits_{i\in I}}(b\wedge a_{i})$ implies, $b$ being indecomposable, $b=b\wedge a_{i}$ for some $i\in I$, equivalently $b\in \phi_{F}(a_{i})$.

So, by Diaconescu's equivalence, the flat $J_{F}$-continuous functor 
\[
y \circ \phi_{F}:F\to [{{\cal I}_{F}}^{\textrm{op}}, \Set]
\]
corresponds to a geometric morphism 
\[
\chi_{F}:[{{\cal I}_{F}}^{\textrm{op}}, \Set] \to \Sh(F, J_{F}).
\]   
   
Let $\textbf{DJFrm}_{dis}$ be the category whose objects are the disjunctively distributive frames $F$ such that the subset of infinitarily disjunctively compact elements of $Id({\cal I}_{F})$ is closed in $Id({\cal I}_{F})$ under finite meets, and whose arrows $F\to F'$ are the disjunctively distributive frames homomorphisms $f:F\to F'$ such that there exists a monotone map $u:{\cal I}_{F'}\to {\cal I}_{F}$ such that the diagram  
\[  
\xymatrix {
[{{\cal I}_{F'}}^{\textrm{op}}, \Set] \ar[r]^{E_{u}} \ar[d]^{\chi_{F'}} & [{{\cal I}_{F'}}^{\textrm{op}}, \Set] \ar[d]^{\chi_{F}} \\
\Sh(F', J_{F'}) \ar[r]^{\dot{f}} & \Sh(F, J_{F})}
\]
commutes (up to isomorphism), where $E_{u}:[{{\cal I}_{F'}}^{\textrm{op}}, \Set] \to [{{\cal I}_{F}}^{\textrm{op}}, \Set]$ is the geometric morphism induced by $u^{\textrm{op}}:{{\cal I}_{F'}}^{\textrm{op}}\to {{\cal I}_{F}}^{\textrm{op}}$ as in Example A4.1.4 \cite{El}, and the inverse image $u^{-1}:Id({\cal I}_{F})\to Id({\cal I}_{F'})$ sends the infinitarily disjunctively compact elements of $Id({\cal I}_{F})$ to infinitarily disjunctively compact elements of $Id({\cal I}_{F'})$.

Given an arrow $f:F\to F'$ in $\textbf{DJFrm}_{dis}$ there is exactly one monotone function $u:{\cal I}_{F'}\to {\cal I}_{F}$ such that the diagram 
\[  
\xymatrix {
[{{\cal I}_{F'}}^{\textrm{op}}, \Set] \ar[r]^{E_{u}} \ar[d]^{\chi_{F'}} & [{{\cal I}_{F}}^{\textrm{op}}, \Set] \ar[d]^{\chi_{F}} \\
\Sh(F', J_{F'}) \ar[r]^{\dot{f}} & \Sh(F, J_{F})}
\]
commutes (up to isomorphism). Indeed, there is at most one such map by definition of the category $\textbf{DJFrm}_{dis}$, while the uniqueness can be proved as follows. The commutativity of the square above (up to isomorphism) is equivalent to the (strict) commutativity of the diagram
\[  
\xymatrix {
F \ar[r]^{f} \ar[d]^{\phi_{F}} & F' \ar[d]^{\phi_{F'}} \\
Id({\cal I}_{F}) \ar[r]^{u^{-1}} & Id({\cal I}_{F'})},
\]
which, in light of the fact that $u^{-1}$ preserves unions of ideals, forces $u^{-1}:Id({\cal I}_{F}) \to Id({\cal I}_{F'})$ to be equal to the map $\xi_{f}:Id({\cal I}_{F}) \to Id({\cal I}_{F'})$ sending an ideal $I$ in $Id({\cal I}_{F})$ to the union $\mathbin{\mathop{\textrm{\huge $\vee$}}\limits_{p\in I}}\phi_{F'}(f(p))$ in $Id({\cal I}_{F'})$. Our claim then follows from the fact that $u$ is uniquely determined by $u^{-1}$ (cf. Theorem \ref{proalex}).

Note that the map $\xi_{f}:Id({\cal I}_{F}) \to Id({\cal I}_{F'})$ is a frame homomorphism, which implies that given a disjunctively distributive frame homomorphism $f:F\to F'$, there exists a monotone map $u:{\cal I}_{F'}\to {\cal I}_{F}$ such that the square above commutes, equivalently $u^{-1}=\xi_{f}:Id({\cal I}_{F}) \to Id({\cal I}_{F'})$, if and only if $\xi_{f}$ is complete, i.e. it preserves arbitrary infima (cf. section \ref{prealex} above).

Therefore the category $\textbf{DJFrm}_{dis}$ can be alternatively described as the category whose objects are the disjunctively distributive frames $F$ such that the subset of infinitarily disjunctively compact elements of $Id({\cal I}_{F})$ is closed in $Id({\cal I}_{F})$ under finite meets, and whose arrows $F\to F'$ are the disjunctively distributive frames homomorphisms $f:F\to F'$ such that the map $\xi_{f}:Id({\cal I}_{F}) \to Id({\cal I}_{F'})$ is complete and sends infinitarily disjunctively compact elements to infinitarily disjunctively compact elements.

Let us consider the category $\textbf{DisFrm}$ of disjunctive frames defined in section \ref{addex}; recall that the objects of $\textbf{DisFrm}$ are the disjunctive frames and the arrows $F\to F'$ of $\textbf{DisFrm}$ are the meet-semilattice homomorphisms $f:F\to F'$ between them which send pairwise disjoint joins to pairwise disjoint joins and have the property that the frame homomorphism $A(f):Id_{J_{F}}(F)\to Id_{J_{F'}}(F')$ sending an ideal $I$ of $F$ to the ideal of $F'$ generated by $f(I)$ preserves arbitrary infima. 

We can define two functors 
\[
i:\textbf{DisFrm}\to \textbf{DJFrm}_{dis}
\]
and 
\[ 
L:\textbf{DJFrm}_{dis} \to \textbf{DisFrm}
\]
as follows; $i$ is the inclusion functor of $\textbf{DisFrm}$ into $\textbf{DJFrm}_{dis}$ while $L:\textbf{DJFrm}_{dis} \to \textbf{DisFrm}$ is the functor sending a poset $F$ in $\textbf{DJFrm}_{dis}$ to the poset $Dis(Id({\cal I}_{F}))$ of infinitarily disjunctively compact elements of the frame $Id({\cal I}_{F})$ and an arrow $f:F\to F'$ in $\textbf{DJFrm}_{dis}$ to the restriction of the map $\xi_{f}:Id({\cal I}_{F})\to Id({\cal I}_{F'})$ to the subsets $Dis(Id({\cal I}_{F}))$ and $Dis(Id({\cal I}_{F'}))$ of infinitarily disjunctively compact elements of $Id({\cal I}_{F})$ and of $Id({\cal I}_{F'})$ (note that the functors $i$ and $L$ are well-defined by Theorem \ref{disthm}). 

\begin{theorem}\label{disadj}
The embedding
\[
i:\textbf{DisFrm}\hookrightarrow \textbf{DJFrm}_{dis}
\] 
identifies $\textbf{DisFrm}$ a full reflective subcategory of the category $\textbf{DJFrm}_{dis}$, with reflector $L:\textbf{DJFrm}_{dis} \to \textbf{DisFrm}$.  
\end{theorem}       

\begin{proofs}
The fact that $L\circ i:\textbf{DisFrm} \to \textbf{DisFrm}$ is naturally isomorphic to the identity functor follows immediately from Theorem \ref{disthm} and the results of section \ref{Mordual}.

To prove that $L$ is left adjoint to $i$, we observe that we have the following natural correspondences, natural in $F\in \textbf{DJFrm}_{dis}$ and in $G\in \textbf{DisFrm}$, between the arrows $L(F)\to G$ in $\textbf{DisFrm}$ and the arrows $F\to i(G)$ in $\textbf{DJFrm}_{dis}$; an arrow $\alpha:L(F)\to G$ in $\textbf{DisFrm}$ corresponds to the arrow $\alpha\circ \phi_{F}:F\to i(G)$, while an arrow $\beta:F\to i(G)$ corresponds to the arrow $\phi_{G}^{-1}\circ L(\beta)$ given by the composite of $\phi_{G}^{-1}:Dis(Id({\cal I}_{G}))\to G$ with the arrow $L(\beta):L(F)\to L(i(G))$.   

Note that the unit of the adjunction at an object $F$ of $\textbf{DJFrm}_{dis}$ is given by the map $\phi_{F}:F\to Dis(Id({\cal I}_{F}))$, while the counit at an object $G$ of $\textbf{DisFrm}$ is given by the isomorphism $\phi_{G}^{-1}:Dis(Id({\cal I}_{G}))\to G$.    
\end{proofs}

Similarly, we can establish an adjunction based on the Lindenbaum-Tarski representation of atomic frames. The analogue of Theorem \ref{disthm} for atomic frames reads as follows.

\begin{theorem}\label{atthm}

\begin{enumerate}[(i)]
\item For any frame $F$, ${\mathscr{P}}(At(F))$ is an atomic frame (with the subset inclusion order), and the map $\psi_{F}:F\to {\mathscr{P}}(At(F))$ sending an element $p\in F$ to the subset $\{a\in At(F) \textrm{ | } a\leq p\}$ is a frame homomorphism which restricts to an isomorphism from $At(F) \subseteq F$ to the subset $At({\mathscr{P}}(At(F)))$ of atoms of ${\mathscr{P}}(At(F))$;

\item A frame $F$ is atomic if and only if the map $\psi_{F}:F\to {\mathscr{P}}(At(F))$ is a frame isomorphism;

\item A topological space $X$ is discrete if and only if it is sober and the frame ${\cal O}(X)$ of open sets of $X$ is isomorphic via the map $\psi_{{\cal O}(X)}$ to the powerset ${\mathscr{P}}(At({\cal O}(X)))$ of the set $At({\cal O}(X))$ of atoms of the frame ${\cal O}(X)$ of open sets of $X$. 
\end{enumerate}
\end{theorem}

\begin{proofs}
$(i)$ It is clear that ${\mathscr{P}}(At(F))$ is an atomic frame, and it is immediate to verify that $\psi_{F}:F\to {\mathscr{P}}(At(F))$ is a frame homomorphism which restricts to an isomorphism from $At(F) \subseteq F$ to the subset $At({\mathscr{P}}(At(F)))$ of atoms of ${\mathscr{P}}(At(F))$. 

$(ii)$ The `only if' direction follows from the Comparison Lemma (cf. section \ref{tarski}), so it remains to prove the `if' one. Let us suppose that $\psi_{F}:F\to {\mathscr{P}}(At(F))$ is an isomorphism. By part $(i)$ ${\mathscr{P}}(At(F))$ is atomic and hence $F$, being isomorphic to it, is atomic as well, as required.     

$(iii)$ This follows immediately from part $(ii)$ by recalling that a topological space is discrete if and only if it is sober and has a basis of atomic open subsets (cf. Proposition \ref{discrete}).
\end{proofs}

Similarly to the case of disjunctively distributive frames, for any frame $F$ we have a geometric morphism $\chi_{F}:[At(F), \Set] \to \Sh(F, J_{F})$ (where $J_{F}$ is the canonical topology on $F$) corresponding to the map $\psi_{F}:F\to {\mathscr{P}}(At(F))$.

Let us define $\textbf{Frm}_{at}$ as the category having as objects the frames and as arrows $F\to F'$ the frame homomorphisms $f:F\to F'$ such that there exists a monotone map $u:At(F')\to At(F)$ making the diagram  
\[  
\xymatrix {
[At(F'), \Set] \ar[r]^{E_{u}} \ar[d]^{\chi_{F'}} & [At(F), \Set] \ar[d]^{\chi_{F}} \\
\Sh(F', J_{F'}) \ar[r]^{\dot{f}} & \Sh(F, J_{F})}
\]
commute (up to isomorphism), where $E_{u}:[At(F'), \Set] \to [At(F), \Set]$ the geometric morphism induced by the map $u: At(F')\to At(F)$ as in Example A4.1.4 \cite{El}. 

One can easily see that the category $\textbf{Frm}_{at}$ can be alternatively described as the category whose objects are the frames $F$ and whose arrows $F\to F'$ are the frame homomorphisms $f:F\to F'$ such that the map $\xi_{f}:{\mathscr{P}}(At(F)) \to {\mathscr{P}}(At(F'))$ which sends a subset $S$ in ${\mathscr{P}}(At(F))$ to the union $\mathbin{\mathop{\textrm{\huge $\vee$}}\limits_{p\in S}}\psi_{F'}(f(p))$ in ${\mathscr{P}}(At(F'))$ is complete. Note that if $f$ is complete then $\xi_{f}$ is complete; so the category of frames and complete frame homomorphisms between them sits as a subcategory of $\textbf{Frm}_{at}$. 

Let $\textbf{AtFrm}$ be the category of atomic frames and complete frame homomorphisms between them, as defined in section \ref{tarski}.

We can define two functors $i':\textbf{AtFrm}\to \textbf{Frm}_{at}$ and $L':\textbf{Frm}_{dis} \to \textbf{AtFrm}$ as follows; $i'$ is the inclusion functor of $\textbf{AtFrm}$ into $\textbf{Frm}_{at}$ (note that this functor is well-defined by Theorem \ref{atthm}) while $L':\textbf{Frm}_{at} \to \textbf{AtFrm}$ sends a poset $F$ in $\textbf{Frm}_{at}$ to the poset ${\mathscr{P}}(At(F))$ and an arrow $f:F\to F'$ in $\textbf{Frm}_{at}$ to the map $\xi_{f}:{\mathscr{P}}(At(F))\to {\mathscr{P}}(At(F'))$. 

\begin{theorem}
The embedding
\[
i':\textbf{AtFrm} \hookrightarrow \textbf{Frm}_{at}
\] 
identifies $\textbf{AtFrm}$ a full reflective subcategory of the category $\textbf{Frm}_{at}$, with reflector $L':\textbf{Frm}_{at} \to \textbf{AtFrm}$.  
\end{theorem}       

\begin{proofs}
The proof is entirely analogous to that of Theorem \ref{disadj}.

The fact that $L\circ i:\textbf{AtFrm} \to \textbf{AtFrm}$ is naturally isomorphic to the identity functor follows immediately from Theorem \ref{disthm} and the results of section \ref{tarski}.

To prove that $L'$ is left adjoint to $i'$, we observe that we have the following natural correspondences, natural in $F\in \textbf{Frm}_{at}$ and in $G\in \textbf{AtFrm}$, between the arrows $L(F)\to G$ in $\textbf{DisFrm}$ and the arrows $F\to i'(G)$ in $\textbf{DJFrm}_{dis}$: an arrow $\alpha:L'(F)\to G$ in $\textbf{AtFrm}$ corresponds to the arrow $\alpha\circ \psi_{F}:F\to i'(G)$, while an arrow $\beta:F\to i'(G)$ corresponds to the arrow $\psi_{G}^{-1}\circ L'(\beta)$ given by the composite of $\psi_{G}^{-1}:{\mathscr{P}}(At(G))\to G$ with the arrow $L(\beta):L'(F)\to L'(i'(G))$.   

Note that the unit of the adjunction at an object $F$ of $\textbf{Frm}_{at}$ is given by the map $\psi_{F}:F\to {\mathscr{P}}(At(F))$, while the counit at an object $G$ of $\textbf{AtFrm}$ is given by the isomorphism $\psi_{G}^{-1}:{\mathscr{P}}(At(G))\to G$.    
\end{proofs}

Of course, by using the same technique of these examples, one can establish similar reflections or, more generally, adjunctions between subcategories of locales and categories of locales satisfying topological properties which can be expressed in terms of the existence of bases satisfying particular conditions; note that the topological properties of locales which we have addressed in the two examples above are local connectedness and atomicity.

\subsection{Toposes paired with points and topological spaces}

In this section we establish, as a further application of the philosophy `toposes as bridges' of \cite{OC10}, adjunctions between full subcategories of the category $\textbf{Top}$ of topological spaces and full subcategories of the category $\mathfrak{Top}_{p}$ of toposes paired with points defined in section \ref{subterminal}.

First, we note that to any topological space $X$ we can associate an object $(\Sh(X), \xi_{X})$ of $\mathfrak{Top}_{p}$ obtained by equipping $\Sh(X)$ with the set of points of $\Sh(X)$ indexed by the set $X_{0}$ of the points of $X$ (see Remark \ref{topolinvariant}(c) above). In fact, this assignment defines a functor $P:\textbf{Top}\to {\mathfrak{Top}_{p}}$, which sends a space $X$ to the object $(\Sh(X), \xi_{X})$ of $\mathfrak{Top}_{p}$ and a continuous map $f:X\to Y$ of topological spaces to the arrow $(P(f), f):(\Sh(X), \xi_{X}) \to (\Sh(Y), \xi_{Y})$ in $\mathfrak{Top}_{p}$, where $P(f)$ is the geometric morphism $\Sh(X)\to \Sh(Y)$ whose inverse image acts, at the level of \'{e}tale bundles, as the pullback functor along $f$ (note that the pair $(P(f), f)$ defines an arrow $(\Sh(X), \xi_{X}) \to (\Sh(Y), \xi_{Y})$ in ${\mathfrak{Top}_{p}}$ since for any $x\in X_{0}$, $f\circ s(x)\cong t(f(x))$, where $s(x):\Set \to \Sh(X)$ is the point of $\Sh(X)$ corresponding to the point $x\in X_{0}$ and $t(f(x)):\Set\to \Sh(Y)$ is the point of $\Sh(Y)$  corresponding to the point $f(x)\in Y_{0}$).

Let $\textbf{\cal U}$ be a full subcategory of $\textbf{Top}$ and $\mathfrak{V}$ be a full subcategory of ${\mathfrak{Top}_{p}}$ which contains all the objects of the form $P(X)$ for $X\in \textbf{\cal U}$. Suppose to have an assignment $\Gamma$ sending each pair $({\cal E}, \xi)$ in $\mathfrak{V}$ to a subframe $\Gamma_{({\cal E}, \xi)}$ of $\Sub_{\cal E}(1)$ in such a way that for any $\cal E$ in $\mathfrak{V}$, $X_{{\tau}^{\cal E}_{\Gamma_{\cal E}}, \xi}$ belongs to $\cal E$, and for any arrow  $(f, l):({\cal E}, \xi)\to ({\cal F}, \xi')$ in $\mathfrak{V}$, $f^{\ast}$ sends $\Gamma_{({\cal F}, \xi')}$ to $\Gamma_{({\cal E}, \xi)}$. Then, by the results of section \ref{subterminal}, we have a functor $Q^{\Gamma}:\mathfrak{V}\to \textbf{\cal U}$ sending each pair $({\cal E}, \xi)$ in $\mathfrak{V}$ to the topological space $X_{{\tau}^{\cal E}_{\Gamma_{\cal E}}, \xi}$.

We want to investigate under which conditions we have an adjunction between (restrictions of) the functors $P$ and $Q^{\Gamma}$. The notation below is borrowed from section \ref{subterminal}.

First, we recall that an arrow $({\cal E}, \chi)\to (\Sh(X), \xi_{X})$ in $\mathfrak{Top}_{p}$ consists of a pair $(f, l)$, where $f$ is a geometric morphism ${\cal E} \to \Sh(X)$ and $l:Z\to X_{0}$ (where $Z$ is the domain of $\chi$) is a function such that the diagram

\[  
\xymatrix {
[Z, \Set] \ar[r]^{E(l)} \ar[d]^{\chi} & [X_{0}, \Set] \ar[d]^{\xi_{X}} \\
{\cal E} \ar[r]^{f} & \Sh(X)}
\]

commutes (up to isomorphism). 

By Diaconescu's equivalence, the geometric morphism $f:{\cal E} \to \Sh(X)$ corresponds to a flat $J$-continuous functor $F(f):{\cal O}(X)\to {\cal E}$, where ${\cal O}(X)$ is the frame of open sets of $X$ and $J$ is the canonical topology on it; note that this functor takes values in the frame of subterminals $\Sub_{\cal E}(1)$ of $\cal E$ and hence it can be equivalently regarded as a frame homomorphism ${\cal O}(X)\to \Sub_{\cal E}(1)$. 

The commutativity condition of the square above amounts precisely to the requirement that for every open set $U$ of $X$, $\phi_{\Sub_{\cal F}(1), {\cal E}}(f^{\ast}(U))=l^{-1}(U)$. If ${\cal O}(X)=\Gamma_{(\Sh(X), \xi_{X})}$ (equivalently, $X=Q^{\Gamma}(P(X))$) then this condition implies that $F(f)$ takes values in $\Gamma_{(\cal E, \chi)}$ (we shall denote this image restriction of $F(f)$ to $\Gamma_{(\cal E, \chi)}$ by $F(f)|$) and $l:Z\to X_{0}$ is the underlying function of a continuous map of topological spaces $Q^{\Gamma}(({\cal E}, \chi))\to X$; in fact, $l^{-1}:{\cal O}(X)\to {\cal O}(Q^{\Gamma}(({\cal E}, \chi)))$ can be identified with the composite of $F(f)|$ with the canonical surjection $\phi_{\Gamma_{(\cal E, \chi)}, {\cal E}}:\Gamma_{(\cal E, \chi)} \to {\cal O}(Q^{\Gamma}(({\cal E}, \chi)))$. 

Conversely, given a continuous map $l:Q^{\Gamma}(({\cal E}, \chi))\to X$, $l^{-1}$ defines a frame homomorphism ${\cal O}(X) \to {\cal O}(Q^{\Gamma}(({\cal E}, \chi)))$. So, if the canonical surjection $\phi_{\Gamma_{(\cal E, \chi)}, {\cal E}}:\Gamma_{(\cal E, \chi)} \to {\cal O}(Q^{\Gamma}({\cal E}, \xi))$ is injective (equivalently, bijective), we can lift this homomorphism to a frame homomorphism ${\cal O}(X) \to \Gamma_{(\cal E, \chi)} \hookrightarrow \Sub_{\cal E}(1)$. This in turn corresponds to a flat $J$-continuous functor ${\cal O}(X)\to {\cal E}$ and hence to a geometric morphism $f:{\cal E} \to \Sh(X)$ such that $(f,l)$ is an arrow $({\cal E}, \chi)\to (\Sh(X), \xi_{X})$ in $\mathfrak{Top}_{p}$. 

So far we have identified two conditions:
\begin{enumerate}[(1)]
\item $X=Q^{\Gamma}(P(X))$ and

\item the injectivity (equivalently, bijectivity) of the canonical surjection\\ $\phi_{\Gamma_{(\cal E, \chi)}, {\cal E}}:\Gamma_{(\cal E, \chi)} \to {\cal O}(Q^{\Gamma}(({\cal E}, \chi)))$,
\end{enumerate}

under which we can establish a bijective correspondence between the arrows $({\cal E}, \chi)\to P(X)$ in $\mathfrak{Top}_{p}$ and the arrows $l:Q^{\Gamma}(({\cal E},\chi))\to X$ in $\textbf{Top}$, natural in $X\in \textbf{{\cal U}}$ and $({\cal E}, \chi)\in \mathfrak{V}$ .  

In order to make this correspondence into an adjunction between a subcategory of $\mathfrak{V}$ and a subcategory of $\textbf{{\cal U}}$, we consider the restriction $P|:\textbf{{\cal U}}_{\Gamma} \to \mathfrak{V}$ of $P$ to the full subcategory $\textbf{{\cal U}}_{\Gamma}$ of $\textbf{{\cal U}}$ on the objects $X$ such that $X=Q^{\Gamma}(P(X))$ and the restriction $Q^{\Gamma}|:\mathfrak{V}' \to \textbf{{\cal U}}$ of $Q^{\Gamma}$ to the full subcategory $\mathfrak{V}'$ of $\mathfrak{V}$ on the objects $({\cal E}, \chi)$ such that the canonical surjection $\phi_{\Gamma_{({\cal E}, \chi)}, {\cal E}}:\Gamma_{({\cal E}, \chi)} \to {\cal O}(Q^{\Gamma}(({\cal E}, \chi)))$ is an isomorphism. 

Now, since the correspondences established above are inverse to each other and natural in $({\cal E}, \chi) \in \mathfrak{V}'$ and $X\in \textbf{{\cal U}}_{\Gamma}$, to obtain an adjunction between $P|$ and $Q^{\Gamma}|$ it would be enough to show that the functor $P|$ takes values in $\mathfrak{V}'$ and the functor $Q^{\Gamma}|$ takes values in $\textbf{{\cal U}}_{\Gamma}$. The second of these conditions is always satisfied, since for any topological space $X$, $\xi_{X}$ is a separating set of points of $\Sh(X)$ (cf. Theorem \ref{topol}(ii) and Remark \ref{topolinvariant}(c)). On the other hand, for any $({\cal E}, \chi)$ in $\mathfrak{V}$, $Q^{\Gamma}(({\cal E}, \chi))$ belongs to $\textbf{{\cal U}}_{\Gamma}$ if and only if ${\cal O}(Q^{\Gamma}(({\cal E}, \chi)))=\Gamma_{\Sh(Q^{\Gamma}({\cal E}, \chi)), \xi_{Q^{\Gamma}({\cal E}, \chi)}}$ (equivalently, $Q^{\Gamma}(({\cal E}, \chi))=Q^{\Gamma}(P(Q^{\Gamma}(({\cal E}, \chi))))$). Let us denote by $(2')$ the conjunction of this latter condition with condition $(2)$ above. Clearly, for any $X\in \textbf{{\cal U}}_{\Gamma}$, $P(X)$ satisfies condition $(2')$; so, if we denote by ${\mathfrak{V}_{p}}_{\Gamma}$ the full subcategory of $\mathfrak{V}$ on the objects $({\cal E}, \chi)$ which satisfy condition $(2')$, we have the following result.

\begin{theorem} 
With the above notation, the restrictions $P_{\Gamma}:\textbf{{\cal U}}_{\Gamma} \to {\mathfrak{V}_{p}}_{\Gamma}$ and $Q_{\Gamma}:{\mathfrak{V}_{p}}_{\Gamma} \to \textbf{{\cal U}}_{\Gamma}$ respectively of the functors $P$ and $Q^{\Gamma}$ define a pair of adjoint functors, with $Q_{\Gamma}$ being the left adjoint and $P_{\Gamma}$ being the right adjoint.   
\end{theorem}\qed 

Note that the adjunction of the theorem is in fact a reflection, since for any arrow $f:X\to Y$ in $\textbf{Top}$, $f=Q(P(f))$ (cf. condition (1) above), and hence $\textbf{{\cal U}}_{\Gamma}$ can be regarded, via $P_{\Gamma}$, as a full subcategory of ${\mathfrak{V}_{p}}_{\Gamma}$.  

Finally, let us discuss a couple of applications of this result.

\begin{enumerate}[(i)] 

\item If $\mathfrak{V}=\mathfrak{Top}_{p}$, $\textbf{{\cal U}}=\textbf{Top}$ and $\Gamma_{({\cal E}, \chi)}$ is the subframe $\{0,1\}$ of $\Sub_{\cal E}(1)$ for any $({\cal E}, \chi)$ in $\mathfrak{V}$ then the spaces in $\textbf{{\cal U}}_{\Gamma}$ are exactly the trivial topological spaces (so $\textbf{Top}_{\Gamma}$ is isomorphic to the category $\Set$ of sets), while the objects in ${\mathfrak{V}_{p}}_{\Gamma}$ are exactly the pairs $({\cal E}, \chi)$ such that $\cal E$ is trivial (i.e., $0\cong 1$ in $\cal E$) if and only if $dom(\chi)=\emptyset$. The functor $P_{\Gamma}:\textbf{{\cal U}}_{\Gamma} \to {\mathfrak{V}_{p}}_{\Gamma}$ sends a topological space $X$ in $\textbf{{\cal U}}_{\Gamma}$ to $(\textbf{1}, \xi_{X})$ (where $\textbf{1}\simeq \Sh(X)$ is the trivial topos) if $X$ is empty and to $({\cal S}, \xi_{X})$ if $X$ is non-empty, where ${\cal S}\simeq \Sh(X)$ is the Sierpinski topos i.e. the category of sheaves on the Sierpinski space (equivalently the topos $[\textbf{2}, \Set]$ where $\textbf{2}$ is the preorder category on the natural number $2$), while the functor $Q_{\Gamma}:{\mathfrak{Top}_{p}}_{\Gamma} \to \textbf{Top}_{\Gamma}$ sends an object $({\cal E}, \chi)$ of ${\mathfrak{Top}_{p}}_{\Gamma}$ to the topological space obtained by equipping the domain of $\chi$ with the trivial topology.

\item If $\mathfrak{V}=\mathfrak{Top}_{p}$, $\textbf{{\cal U}}=\textbf{Top}$ and $\Gamma_{({\cal E}, \chi)}$ is the whole frame $\Sub_{\cal E}(1)$ for any $({\cal E}, \chi)$ in $\mathfrak{V}$ then $\textbf{{\cal U}}_{\Gamma}=\textbf{Top}$, while a pair $({\cal E}, \chi)$ in ${\mathfrak{Top}_{p}}$ belongs to ${\mathfrak{V}_{p}}_{\Gamma}$ if and only if the canonical surjection $\phi_{\Sub_{\cal E}(1) , {\cal E}}:\Sub_{\cal E}(1)\to {\cal O}(Q(({\cal E}, \chi))$ is a bijection (that is, if and only if $\chi$ separates the subterminals in $\cal E$ i.e. for any two subterminals $U$ and $V$ in $\cal E$, if $U\ncong V$ then there exists $z\in dom(\chi)$ such that $\chi(z)^{\ast}(U)\ncong \chi(z)^{\ast}(V)$). The functor $P_{\Gamma}:\textbf{{\cal U}}_{\Gamma} \to {\mathfrak{V}_{p}}_{\Gamma}$ sends a topological space $X$ in $\textbf{{\cal U}}_{\Gamma}$ to $(\Sh(X), \xi_{X})$ while the functor $Q_{\Gamma}:{\mathfrak{V}_{p}}_{\Gamma} \to \textbf{{\cal U}}_{\Gamma}$ sends an object $({\cal E}, \chi)$ of ${\mathfrak{V}_{p}}_{\Gamma}$ to the topological space obtained by equipping the domain of $\chi$ with the subterminal topology.  

\end{enumerate}

\section{Insights obtained by using invariants}\label{insights}

In this section our aim is to show that the technique of \cite{OC10} of using toposes as `bridges' for transferring information between Morita-equivalent theories (in the form of different sites of definition for the same topos), can be profitably applied in the context of our topos-theoretic interpretation of Stone-type dualities, for translating properties of preordered structures into properties of the corresponding locales or topological spaces (or, more generally, for translating properties between preordered structures $\cal C$ and $\cal D$ related by Morita-equivalences of the form $\Sh({\cal C}, J)\simeq \Sh({\cal D}, K)$).

\subsection{The logical interpretation}\label{logical}

The equivalence $\Sh({\cal C}, J)\simeq \Sh(Id_{J}({\cal C}))$ of Theorem \ref{fund} (for a preorder category ${\cal C}$ and a Grothendieck topology $J$ on ${\cal C}$) can be read as a Morita-equivalence between two distinct geometric theories: the theory of $J$-contin-\\uous flat functors on $\cal C$ and the theory of $J_{can}^{Id_{J}({\cal C})}$-continuous flat functors on $Id_{J}({\cal C})$, where $J_{can}^{Id_{J}({\cal C})}$ is the canonical topology on the frame $Id_{J}({\cal C})$ (cf. \cite{OC10}).

On the other hand, any $J$-continuous flat functor $F:{\cal C}\to {\cal E}$ from a preorder category $\cal C$ to a Grothendieck topos $\cal E$ sends every object in $\cal C$ to a subterminal object of $\cal E$. Indeed, by Diaconescu's equivalence, $F$ is isomorphic to a functor of the form $f^{\ast}\circ a_{J}\circ y$, where $f:{\cal E}\to \Sh({\cal C}, J)$ is a geometric morphism, $y:{\cal C}\to [{\cal C}^{\textrm{op}}, \Set]$ is the Yoneda embedding and $a_{J}:[{\cal C}^{\textrm{op}}, \Set]\to \Sh({\cal C}, J)$ is the associated sheaf functor; and, $\cal C$ being a preorder, the objects of $\cal C$ are sent by $y$ to subterminals in $[{\cal C}^{\textrm{op}}, \Set]$, which are in turn sent by $f^{\ast}\circ a_{J}$ to subterminals in $\cal E$. 

From the characterization of flat functors ${\cal C}\to {\cal E}$ as filtering functors given in section VII.9 of \cite{MM} it thus follows that the $J$-continuous flat functors ${\cal C}\to {\cal E}$ and natural transformations between them can be identified, naturally in $\cal E$, respectively with the models in $\cal E$ and model homomorphisms between them of the propositional theory ${\mathbb T}^{{\cal C}}_{J}$ defined as follows: the signature of ${\mathbb T}^{{\cal C}}_{J}$ has no sorts and one atomic proposition $F_{a}$ for each element $a\in {\cal C}$, and the axioms of ${\mathbb T}^{{\cal C}}_{J}$ are the following:
\[
(\top \vdash \mathbin{\mathop{\textrm{\huge $\vee$}}\limits_{c\in {\cal C}}} F_{c});
\]
\[
(F_{a} \vdash F_{b})
\] 
for any $a\leq b$ in $\cal C$; 
\[
(F_{a}\wedge F_{b} \vdash \mathbin{\mathop{\textrm{\huge $\vee$}}\limits_{c\in K_{a,b}}} F_{c})
\]
for any $a, b \in {\cal C}$, where $K_{a,b}$ is the collection of all the elements $c\in {\cal C}$ such that $c\leq a$ and $c\leq b$;
\[
(F_{a} \vdash \mathbin{\mathop{\textrm{\huge $\vee$}}\limits_{i\in I}}F_{a_{i}})
\]
for any $J$-covering sieve $\{a_{i} \to a \textrm{ | } i\in I\}$ in $\cal C$.

We note that the models of the theory ${\mathbb T}^{{\cal C}}_{J}$ in $\Set$ are precisely the $J$-prime filters on $\cal C$, as defined in section \ref{exsub}; accordingly, we call the theory ${\mathbb T}^{{\cal C}}_{J}$ the \emph{theory of $J$-prime filters}. In particular, if $L$ is a frame and $J$ is the canonical topology on it then the theory ${\mathbb T}^{{\cal C}}_{J}$ specializes to the propositional theory of completely prime filters on $L$ defined in section D1.1 of \cite{El}. Given a topological space $X$, we define the \emph{theory of completely prime filters on $X$} as the theory of completely prime filters on the corresponding frame ${\cal O}(X)$ of open sets of $X$.

Summarizing, we have the following result.

\begin{theorem}\label{classif}
Let $\cal C$ be a preorder category and $J$ be a Grothendieck topology on $\cal C$. Then the theory ${\mathbb T}^{{\cal C}}_{J}$ of $J$-prime filters is classified by the topos $\Sh({\cal C}, J)$.
\end{theorem}\qed

\begin{remarks}
\begin{enumerate}[(a)]

\item If ${\cal C}={\cal P}^{\textrm{op}}$, where $({\cal P}, \leq)$ is a preorder category, then the models of ${\mathbb T}^{{\cal C}}_{J}$ in $\Set$, i.e. the $J$-prime filters on $\cal C$, are precisely the non-empty directed ideals on $\cal P$;   

\item If $\cal C$ is cartesian (resp. coherent, geometric) and $J$ is the trivial (resp. coherent, geometric) topology on $\cal C$, the theory ${\mathbb T}^{{\cal C}}_{j}$ specializes to the theory of filters (resp. prime filters, completely prime filters) on $\cal C$ described in section D1.1 \cite{El}.

\end{enumerate}
\end{remarks}

The equivalence $\Sh({\cal C}, J)\simeq \Sh(Id_{J}({\cal C}))$ of Theorem \ref{fund} can thus be read logically as a Morita-equivalence between the theory of $J$-prime filters on $\cal C$ and the theory of completely prime filters on $Id_{J}({\cal C})$. In particular, Stone duality for distributive lattices admits the following logical interpretation: given a locale $L$, $L$ is the Stone locale associated to a distributive lattice $\cal D$ if and only if the theory of prime filters on $\cal D$ is Morita-equivalent to the theory of completely prime filters on $L$. Of course, the other dualities of sections \ref{ex} and \ref{addex} admit similar interpretations. 

As we have already remarked in section \ref{general}, the `Stone-type dualities' can be seen as arising from the process of `functorializing' a bunch of Morita-equivalences; we have one Morita-equivalence for each of the preordered structures, and these Morita-equivalences can be `merged together' to produce a `global' duality or equivalence of categories. In fact, this method of generating dualities or equivalences of categories starting from `parametrised' Morita-equivalences, which we have exploited in section \ref{Mordual}, is likely to find, because of its generality, new applications in a variety of different mathematical contexts in the near future.

From Theorem \ref{classif} it follows, by recalling the syntactic construction of classifying toposes, that for any site $({\cal C}, J)$ whose underlying category $\cal C$ is a preorder, $Id_{J}({\cal C})$ can be identified with the geometric syntactic category of the theory ${\mathbb T}^{{\cal C}}_{J}$, and hence the theory of completely prime filters on $Id_{J}({\cal C})$ is Morita-equivalent to the theory of completely prime filters on the logical space of ${\mathbb T}^{{\cal C}}_{J}$, as defined in section \ref{exsub}. 

The localic dualities of section \ref{ex} can thus be read logically as `functorialized' Morita-equivalences (in the sense explained above) between propositional geometric theories and the theories of completely prime filters on their logical spaces.       

Recall that the technique `toposes as bridges' introduced in \cite{OC10} consists, broadly speaking, in expressing a given topos-theoretic invariant in terms of two distinct sites of definition of a given topos; the topos acts as a `bridge' enabling one to transfer properties from one site into the other. Of course, the feasibility of this method heavily depends on how `natural' is, for the given invariant, the relationship between the topos and its different sites of definition. There is an important class of invariants for which this transfer of information is always feasible and, in a sense, even \emph{automatic}. Indeed, for logically motivated invariants such as the property of a topos to be two-valued (resp. Boolean, De Morgan), one has bijective characterizations (holding `uniformly' for any site $({\cal C}, J)$) of the kind `$\Sh({\cal C}, J)$ satisfies the invariant if and only if the site $({\cal C}, J)$ satisfies some `tractable' categorical property'. In fact, any first-order sequent $\sigma$ written in the algebraic language of the theory of Heyting algebras can be interpreted in the internal Heyting algebra $\Omega_{\cal E}$ of a topos $\cal E$ given by its subobject classifier, and as such it gives rise to a topos-theoretic invariant (namely, the interpretation of $\sigma$ in the algebra $\Omega_{{\cal E}}$) admitting a site characterization of the above kind (such characterizations can be obtained by using the canonical description of the algebra $\Omega_{\cal E}$ in terms of the given site of definition of $\cal E$, and the explicit descriptions of the interpretation of the first-order connectives and quantifiers in a topos of sheaves on a site, as given for example in Chapter III of \cite{MM}). 

One might naturally wonder whether, given a sequent $\sigma$ as above, its (internal) validity in the algebra $\Omega_{\Sh(L)}$ of the topos $\Sh(L)$ is equivalent to its (external) validity in the locale $L$ (regarded as a model of the theory of Heyting algebras). This is not always the case in general, but there is an important class of sequents $\sigma$ for which these two notions of validity coincide.

\begin{proposition}\label{cartesian}
Let $\sigma$ be a cartesian (in particular, Horn) sequent in the theory of Heyting algebras. Then for any locale $L$, $\sigma$ is valid in the internal algebra $\Omega_{\Sh(L)}$ of the topos $\Sh(L)$ if and only if it is valid in $L$ (regarded as a model of the theory of Heyting algebras).
\end{proposition}

\begin{proofs}
For any locally small topos $\cal E$, from the fact the Yoneda embedding $y:{\cal E}\to [{\cal E}^{\textrm{op}}, \Set]$ is a cartesian functor, it follows that $y$ preserves and the interpretation of all the cartesian formulae and hence, given any cartesian sequent $\sigma$ in the theory of Heyting algebras, the internal Heyting algebra $\Omega_{{\cal E}}$ satisfies $\sigma$ if and only if every frame $\Sub_{\cal E}(e)\cong Hom_{\cal E}(e, \Omega)$ in $\cal E$ satisfies $\sigma$. Now, given a Grothendieck topos $\cal E$ and an object $e\in \cal E$, if $\cal C$ is a separating set for $\cal E$ then $e$ can be expressed as a quotient of a coproduct of objects in $\cal C$, that is there exists a set-indexed family $\{c_{i} \textrm{ | } i\in I\}$ of objects in $\cal C$ and an epimorphism $p:\coprod_{i\in I}c_{i}\epi e$. Since $p$ is an epimorphism, the pullback functor $p^{\ast}:\Sub_{\cal E}(e) \rightarrow \Sub_{\cal E}(\coprod_{i\in I}c_{i})\cong \prod_{i\in I}\Sub_{\cal E}(c_{i})$ is logical and conservative (cf. Example A4.2.7(a) \cite{El}); so $\Sub_{\cal E}(e)$ satisfies a first-order sequent $\sigma$ if all the $\Sub_{\cal E}(c_{i})$ do. Also, if $m:b\mono a$ is a monomorphism in $\cal E$ then the pullback functor $m^{\ast}:\Sub_{\cal E}(a)\to \Sub_{\cal E}(b)$ is logical and essentially surjective; so, if $\Sub_{\cal E}(a)$ satisfies $\sigma$ then $\Sub_{\cal E}(b)$ satisfies $\sigma$. Our thesis now follows from the combination of all these facts by recalling that, for any locale $L$, the collection of all the subterminals of $\Sh(L)$ forms a separating set of $\Sh(L)$, and $L\cong \Sub_{\Sh(L)}(1)$.    
\end{proofs}

An example of this kind of invariants, which we shall consider below, is the interpretation of the sequent 
\[
(\top \vdash_{x,y} (x\imp y) \vee (y\imp x)),
\] 
whose validity in the algebra $\Omega_{\cal E}$ of a topos $\cal E$ amounts to saying that $\cal E$ satisfies G\"{o}del-Dummett's law.

Recall that the subobject classifier $\Omega_{\Sh({\cal C}, J)}:{\cal C}^{\textrm{op}}\to \Set$ of a topos $\Sh({\cal C}, J)$ of sheaves on a site $({\cal C}, J)$ is (isomorphic to) the functor sending any object $c$ of $\cal C$ to the collection $\Omega_{\Sh({\cal C}, J)}(c)$ of $J$-closed sieves on $c$, and an arrow $f:d\to c$ to the pullback operation $f^{\ast}:\Omega_{\Sh({\cal C}, J)}(c)\to \Omega_{\Sh({\cal C}, J)}(d)$ of sieves along $f$. Proposition \ref{cartesian} thus enables us to express `external' properties of locales $L$ formulated as cartesian sequents in the theory of Heyting algebras as `internal' properties of the corresponding localic topos $\Sh(L)$, by means of a topos-theoretic invariant which can in turn be reformulated in terms of any other site of definition $({\cal C}, J)$ of the same topos, leading to categorical characterizations involving $J$-closed sieves on the category $\cal C$. Note that this technique is bound to bring more substantial insights than a straightforward reformulation of the given property of $L$ as a topos-theoretic invariant on the frame $\Sub_{\cal E}(1)$ of subterminals of the topos ${\cal E}\simeq \Sh(L)$, since a formulation of this latter kind would merely consist in a translation of the given property of $L$ across the isomorphism $L\simeq \Sub_{\cal E}(1)$. 

In the next sections, we shall provide applications of the methodology just described by using, as invariants, the property of a topos to be Boolean, to be De Morgan, to be two-valued, to satisfy G\"{o}del-Dummett logic; anyway, the reader should bear in mind that these particular examples are just meant to show the effectiveness of our technique, and notice that a whole range of new results can be `automatically' generated by applying the same method to different invariants.  

In addition to the invariants discussed above, there are of course many other ones which behave `naturally' with respect to sites (cf. \cite{OC10} for an extensive discussion of these aspects) and hence are appropriate for an application of the methodology `toposes as bridges'. Two of them, which we will consider in the next section in relation to the Morita-equivalence of Theorem \ref{fund}, are the notions of point and subtopos of a given topos. 

\subsection{$J$-ideals and $J$-prime filters}

Recall that, for any site $({\cal C}, J)$, the points of a topos $\Sh({\cal C}, J)$ (i.e., the geometric morphisms $\Set \to \Sh({\cal C}, J)$), can be naturally identified with the flat $J$-continuous functors on $\cal C$, while the subtoposes of $\Sh({\cal C}, J)$ (i.e. the equivalence classes of geometric inclusions into $\Sh({\cal C}, J)$) correspond precisely to the Grothendieck topologies on $\cal C$ which contain $J$. 

The method `toposes as bridges' applied to the invariant notion of point and to the Morita-equivalence of Theorem \ref{fund} produces an identification between the $J$-prime filters on $\cal C$ and the completely prime filters on $Id_{J}({\cal C})$.

\begin{theorem}\label{bijectionpoints}
Let $\cal C$ be a preorder and $J$ be a Grothendieck topology on $\cal C$. Then the assignment sending a filter $F$ on $Id_{J}({\cal C})$ to the $J$-prime filter $\{c\in {\cal C} \textrm{ | } (c)\downarrow_{J}\in F\}$ on $\cal C$ defines a bijection between the completely prime filters on the frame $Id_{J}({\cal C})$ of $J$-ideals of $\cal C$ and the $J$-prime filters on $\cal C$. 

In particular, 
\begin{enumerate}[(i)]
\item if $\cal C$ is a meet-semilattice, we have a natural bijection between the filters on $\cal C$ and the completely prime filters on the frame of ideals (i.e. lower sets) of $\cal C$, and 

\item if $\cal C$ is a distributive lattice we have a natural bijection between the prime filters on $\cal C$ and the completely prime filters on the frame of ideals (i.e. lower sets which are closed under finite joins) of $\cal C$.
\end{enumerate} 
\end{theorem}

\begin{proofs}
Starting from the equivalence $\Sh({\cal C}, J)\simeq \Sh(Id_{J}({\cal C}))$ of Theorem \ref{fund}, it suffices to observe that the bijection between the points of the two toposes induced by such equivalence yields, at the level of filters, the assignment sending a filter $F$ on $Id_{J}({\cal C})$ to the filter $\{c\in {\cal C} \textrm{ | } (c)\downarrow_{J}\in F\}$ on $\cal C$.
\end{proofs}

We note that, while this result follows an a natural and immediate application of our method `toposes as bridges', lengthier and less conceptual arguments would be necessary to prove this result directly, that is without appealing to the topos-theoretic machinery (the reader might find instructive to compare our proof of the thoerem with the topological arguments used to establish the particular case of the result for distributive lattices given in Chapter II of \cite{stone}, or with the direct proof of the theorem given in Appendix \ref{appendix} below).

A further application of the philosophy `toposes as bridges' to the Morita-equivalence of Theorem \ref{fund} is the following result, obtained by considering as topos-theoretic invariant the notion of subtopos.

Recall that, given a site $({\cal C}, J)$ whose underlying category is a preorder, and any lower-set $I$ on $\cal C$, its $J$-closure $cl_{J}(I)$ is the smallest $J$-ideal on $\cal C$ which contains $I$, i.e. the ideal $cl_{J}(I)=\{c\in {\cal C} \textrm{ | } \{f:d\to c \textrm{ in $\cal C$ | } d\in I\} \in J(c)\}$. 

\begin{theorem}\label{existence}
Let $\cal C$ be a preorder and $J$ be a Grothendieck topology on $\cal C$. For any surjective frame homomorphism $f:Id_{J}({\cal C})\to F$ onto a frame $F$ there exists a Grothendieck topology $J'\supseteq J$ on $\cal C$ such that $F\cong Id_{J}({\cal C})$ and $f$ corresponds, under this isomorphism, to the frame homomorphism $cl_{J'}:Id_{J}({\cal C}) \to Id_{J'}({\cal C})$ sending a $J$-ideal $I$ on $\cal C$ to its $J'$-closure. In particular, any surjective frame homomorphism whose domain is the frame $Id({\cal C})$ of lower sets on $\cal C$ is, up to isomorphism, of the form $cl_{J}:Id({\cal C}) \to Id_{J}({\cal C}) $ for some Grothendieck topology $J$.  
\end{theorem}

\begin{proofs}
A surjective frame homomorphism $f:Id_{J}({\cal C})\to F$ corresponds to an embedding of the corresponding locales and hence to a geometric inclusion $\Sh(f):\Sh(F)\hookrightarrow \Sh(Id_{J}({\cal C}))$. Through the equivalence
\[
\Sh({\cal C}, J)\simeq \Sh(Id_{J}({\cal C}))
\]
of Theorem \ref{fund}, this inclusion transfers to a subtopos of $\Sh({\cal C}, J)$; and this subtopos must necessarily be, up to equivalence, of the form $i_{J'}:\Sh({\cal C}, J')\hookrightarrow \Sh({\cal C}, J)$ for a Grothendieck topology $J'$ on $\cal C$ which contains $J$, where $i_{J'}$ is the canonical geometric inclusion. 

Since the equivalence of Theorem \ref{fund} is natural with respect to inclusions of Grothendieck topologies, we have a commutative diagram
\[  
\xymatrix {
\Sh(Id_{J}({\cal C})) \ar[r]^{\simeq} & \Sh({\cal C}, J) \ar[r]^{\simeq}  & \Sh(Id_{J}({\cal C})) \\
\Sh(Id_{J'}({\cal C})) \ar[u]^{\Sh(cl_{J'})} \ar[r]^{\simeq} & \Sh({\cal C}, J') \ar[u]^{i_{J'}} \ar[r]^{\simeq} & \Sh(F) \ar[u]^{\Sh(f)}}
\]
from which it follows that $\Sh(f)$ and $\Sh(cl_{J'})$ are equivalent as subtoposes of $\Sh(Id_{J}({\cal C}))$, equivalently $f$ is isomorphic to $cl_{J'}$, as required. 
\end{proofs}

Considering as topos-theoretic invariant the notion of equivalence of sub-\\toposes, applied to the Morita-equivalence   
\[
\Sh({\cal C}, J)\simeq \Sh(Id_{J}({\cal C}))
\]
of Theorem \ref{fund}, we obtain the following result.

\begin{theorem}\label{unique}
Let $\cal C$ be a preorder and let $J_{1}$ and $J_{2}$ be two Grothendieck topologies on $\cal C$. If for any lower set $I$ in $\cal C$, $I$ is a $J_{1}$-ideal if and only if it is a $J_{2}$-ideal then $J_{1}=J_{2}$. 
\end{theorem}

\begin{proofs}
By Theorem \ref{fund}, for any Grothendieck topology $J$ on $\cal C$ we have a commutative diagram
\[  
\xymatrix {
\Sh({\cal C}, J) \ar[r]^{\simeq} \ar[d]^{i_{J}} & \Sh(Id_{J}({\cal C})) \ar[d]^{\Sh(cl_{J})} \\
[{\cal C}^{\textrm{op}}, \Set] \ar[r]^{\simeq} & \Sh(Id({\cal C}))}
\]
where $i_{J}:\Sh({\cal C}, J)\to [{\cal C}^{\textrm{op}}, \Set]$ is the canonical inclusion, $Id({\cal C})$ is the frame of lower sets in $\cal C$ and $\Sh(cl_{J}):\Sh(Id_{J}({\cal C})) \to \Sh(Id({\cal C}))$ is the geometric morphism induced by the frame homomorphism $cl_{J}:Id({\cal C}) \to Id_{J}({\cal C})$.

Now, the condition `for any lower set $I$ in $\cal C$, $I$ is a $J_{1}$-ideal if and only if it is a $J_{2}$-ideal' can be expressed topos-theoretically by saying that there exists an equivalence $e:\Sh(Id_{J_{1}}({\cal C}))\to \Sh(Id_{J_{2}}({\cal C}))$ such that $\Sh(cl_{J_{2}}) \circ e \simeq \Sh(cl_{J_{1}})$. But transferring such an equivalence through the Morita-equivalence of Theorem \ref{fund} yields an equivalence $e':\Sh({\cal C}, J_{1}) \to \Sh({\cal C}, J_{2})$ such that $i_{J_{2}} \circ e' \simeq i_{J_{1}}$. To conclude the proof of the theorem, it suffices to recall the standard topos-theoretic fact that such an equivalence exists if and only if $J_{1}=J_{2}$.
\end{proofs}

Note that Theorem \ref{unique} assures the uniqueness of the Grothendieck topology $J'$ in the statement of Theorem \ref{existence}. 

\subsection{Topological properties as topos-theoretic invariants}

In this section we apply our usual technique `toposes as bridges' (cf. \cite{OC10}) to obtain insights concerning the relationship between preordered structures and the locales or topological spaces which correspond to them under the dualities or equivalences of section \ref{ex}.

\subsubsection{Almost discreteness via the law of excluded middle}

Let us start by investigating the almost discreteness of the locale corresponding to a distributive lattice (resp. to a meet-semilattice, to a preorder) via Stone duality (resp. the duality for meet-semilattices of Theorem \ref{meetsm}, the functor $B:\textbf{Pro}\to \textbf{AlexLoc}$ of section \ref{prealex}).

If $\cal D$ is a distributive lattice, equipped with the coherent topology $J_{\cal D}$, we have two different sites of definition for the topos associated to it via the technique of section \ref{locales}: $\Sh({\cal D}, J_{\cal D}) \simeq \Sh(Id_{J_{\cal D}}({\cal D}))$. We call $Id_{J_{\cal D}}({\cal D})$ the \emph{Stone locale} associated to $\cal D$ (as in section \ref{stonedist}), and we denote it by $L_{\cal D}$.  

Moreover, the Comparison Lemma yields a third site of definition of the topos, obtained by cutting down the site $({\cal D}, J_{\cal D})$ to the full subcategory ${\cal D}^{\ast}$ on the non-zero objects and equipping it with the induced Grothendieck topology $J_{\cal D}^{\ast}=J_{\cal D}|_{{\cal D}^{\ast}}$. 

As we shall see below, rephrasing a given topos-theoretic invariant in terms of these three different sites of definition leads to three different expressions of it, each of them written in the `language' of the corresponding site, which are nonetheless equivalent to each other. 

The following proposition provides several alternative characterizations of the property of the Stone locale associated to a distributive lattice to be almost discrete, obtained by applying this technique. The relevant topos-theoretic invariant is the property of a topos to be Boolean; in fact, it is well-known that for any locale $L$, $\Sh(L)$ is Boolean if and only if $L$ is almost discrete (i.e. every element of $L$ is complemented).

In the proofs below, the symbol $\neg$, applied to an element $l$ of a locale $L$, denotes the Heyting pseudocomplement of $l$ in $L$ (where $L$ is regarded as a Heyting algebra). 

\begin{proposition}\label{boolean}
Let $\cal D$ be a distributive lattice, and let $L_{\cal D}$ (resp. $X_{\cal D}$) be its associated Stone locale (resp. Stone space). Then the following conditions are equivalent:
\begin{enumerate}[(i)]

\item $L_{\cal D}$ (equivalently, $X_{\cal D}$) is almost discrete;

\item For any ideal $I$ of $\cal D$, if $a\in {\cal D}$ is an element satisfying the property that for any non-zero element $b\leq a$ there exists a non-zero element $c\in I$ such that $c\leq b$ then $a\in I$; 

\item For any collection $\{a_{i} \textrm{ | } i\in I\}$ of non-zero elements of $\cal D$ with the property that for any non-zero element $a$ there exists $i\in I$ such that $a_{i}\wedge a\neq 0$, there exists a finite subset $J\subseteq I$ such that $1=\mathbin{\mathop{\textrm{\huge $\vee$}}\limits_{i\in J}}a_{i}$;

\item Every element of $\cal D$ is complemented in $\cal D$, $\cal D$ is complete and the supremum in $\cal D$ of any subset $S$ of $\cal D$ is a finite join of elements of $S$;   

\item $\cal D$ is a finite Boolean algebra.
\end{enumerate}
\end{proposition}

\begin{proofs}
It is well-known that a locale $L$ is almost discrete if and only if the topos $\Sh(L)$ is Boolean. 

Now, condition $(ii)$ represents the expression of the invariant property of a topos to be Boolean in terms of the representation $\Sh(Id_{J_{\cal D}}({\cal D}))$ of the topos $\Sh(L_{\cal D})$, it being the assertion that every ideal $I$ in $Id_{J_{\cal D}}({\cal D})$ satisfies $\neg\neg I=I$ (equivalently, $\neg\neg I\subseteq I$). 

Condition $(iii)$ is the expression of the invariant in terms of the representation $\Sh({\cal D}^{\ast}, J_{\cal D}^{\ast})$ of the topos $\Sh(L_{\cal D})$, as established in \cite{OC3}.

Starting from the assumption that $L_{\cal D}$ is almost discrete, condition $(iv)$ can be deduced as follows. The top element of $L_{\cal D}=Id_{J_{\cal D}}({\cal D})$ being compact, every (complemented) element of $Id_{J_{\cal D}}({\cal D})$ is compact, i.e. it is a principal ideal (cf. Corollary \ref{cor}(i)). Since every element of $L_{\cal D}$ is complemented, every element of $\cal D$ is complemented in $\cal D$. Also, since every ideal in $Id_{J_{\cal D}}({\cal D})$ is principal, $\cal D$ is isomorphic to $Id_{J_{\cal D}}({\cal D})$ and hence in particular $\cal D$ is complete. For a given subset $S$ of $\cal D$, consider the $J_{\cal D}$-ideal $I$ on $\cal D$ generated by $S$; since $I$ is principal then $I=(a)\downarrow$ for some element $a\in {\cal D}$. Clearly, by construction of $I$, $a$ is the supremum of $S$ and can be expressed as a finite join of elements in $S$. Conversely, let us suppose that condition $(iv)$ holds. Since $\cal D$ is complete and the supremum of any subset $S$ of $\cal D$ is a finite join of elements in $S$ then every ideal in $Id_{J_{\cal D}}({\cal D})$ is principal. From the fact that every element of $\cal D$ is complemented we can thus conclude that every ideal in $Id_{J_{\cal D}}({\cal D})$ is complemented, in other words that $L_{\cal D}$ is almost discrete. 

Condition $(iv)$ is clearly equivalent to condition $(v)$; indeed, $(iv)$ implies $(v)$ since, the top element of $L_{\cal D}$ being compact, if every element in $\cal D$ is complemented then $\cal D$ is a finite Boolean algebra, while the fact that $(v)$ implies $(iv)$ is obvious.         
\end{proofs}

The next proposition provides characterizations of the meet-semilattices (resp. preorders) whose corresponding locales via the duality for meet-semilattices of Theorem \ref{meetsm} (resp. the functor $B:\textbf{Pro}\to \textbf{AlexLoc}$ of section \ref{prealex}) are almost discrete. Again, the technique consists in using the topos as a bridge for transferring a given invariant (in this case, the property of a topos to be Boolean) across two different sites of definition for it.

\begin{proposition}\label{Boole2}
\begin{enumerate}[(i)]

\item For any meet-semilattice $\cal M$, the ideal locale $S_{\cal M}$ (equivalently, the ideal topological space) associated to $\cal M$ is almost discrete if and only if $\cal M$ is a singleton.

\item For any preorder $\cal P$, the Alexandrov locale $A_{\cal P}$ (equivalently, the Alexandov space) associated to $\cal P$ is almost discrete if and only if for any $p,q\in {\cal P}$, $p\leq q$ implies $q\leq p$.
   
\end{enumerate}
\end{proposition}

\begin{proofs}
We use the following well-known site characterizations for the invariant property of a topos to be Boolean:

\begin{enumerate}[(a)]
\item A presheaf topos $[{\cal C}^{\textrm{op}}, \Set]$ is Boolean if and only if the category $\cal C$ is a groupoid;

\item A localic topos $\Sh(L)$ is Boolean if and only if the locale $L$ is almost discrete.
\end{enumerate}

The thesis follows immediately from expressing the invariant property of the topos $\Sh(S_{\cal M})\simeq [{\cal M}^{\textrm{op}}, \Set]$ (resp. $\Sh(A_{\cal P})\simeq [{\cal P}, \Set]$) to be Boolean in terms of the two different sites of definition of it, according to the site characterizations reported above.
\end{proofs}

\subsubsection{Extremal disconnectedness via De Morgan's law}

The following proposition represents the analogue of Proposition \ref{boolean} for the property of a Stone locale to be extremally disconnected. Again, the proof is based on the consideration of a topos-theoretic invariant, namely the property of a topos to be De Morgan, in relation to different sites of definition of a given topos.

\begin{proposition}\label{morgan}
Let $\cal D$ be a distributive lattice, and let $L_{\cal D}$ (resp. $X_{\cal D}$) be its associated Stone locale (resp. Stone space). Then the following conditions are equivalent:
\begin{enumerate}[(i)]

\item $L_{\cal D}$ (equivalently, $X_{\cal D}$) is extremally disconnected;

\item For every ideal $I$ of $\cal D$ there exists a finite covering $1=\mathbin{\mathop{\textrm{\huge $\vee$}}\limits_{i\in I}}a_{i}$ of $1$ such that each $a_{i}$ either belongs to $\neg I$ (i.e., for any non-zero element $b\leq a_{i}$, $b\notin I$) or to $\neg\neg I$ (i.e., for any non-zero element $b\leq a_{i}$ there exists a non-zero element $c\in I$ such that $c\leq b$); 

\item For any collection $\{a_{i} \textrm{ | } i\in I\}$ of non-zero elements of $\cal D$ there exists a finite family $\{b_{j} \textrm{ | } j\in J\}$ of non-zero elements of $\cal D$ such that $1=\mathbin{\mathop{\textrm{\huge $\vee$}}\limits_{j\in J}}b_{j}$ and for each $j\in J$ either $b_{j}\wedge a_{i}=0$ for all $i\in I$ or for every non-zero element $x\leq b_{j}$ there exists $i\in I$ such that $x\wedge a_{i}\neq 0$; 

\item For any ideal $I$ of $\cal D$ with the property that for any element $a\in {\cal D}$ such that for any non-zero element $b\leq a$ there exists a non-zero element $c\in I$ such that $c\leq b$, $a\in I$, there exists a complemented element $x\in {\cal D}$ such that $I=(x)\downarrow$. In fact, the lattice of complemented elements of $\cal D$ is isomorphic, via the map sending any complemented elements to the principal ideal which it generates, to the frame of ideals $I$ such that $\neg\neg I=I$.

\end{enumerate}
\end{proposition}

\begin{proofs}
It is well-known that a locale $L$ is extremally disconnected if and only if the topos $\Sh(L)$ is De Morgan.

Condition $(ii)$ represents the expression of this invariant in terms of the representation $\Sh(Id_{J_{\cal D}}({\cal D}))$ of the topos, it being the assertion that every ideal $I$ in $Id_{J_{\cal D}}({\cal D})$ satisfies $\neg I \vee \neg\neg I=L_{\cal D}$.

Condition $(iii)$ is the expression of the invariant in terms of the representation $\Sh({\cal D}^{\ast}, J_{\cal D}^{\ast})$ of the topos $\Sh(L_{\cal D})$, as established in \cite{OC3}.
 
The equivalence $(iv)\biimp (i)$ can be proved as follows. It is well-known that $L_{\cal D}$ is extremally disconnected if and only if every $\neg\neg$-stable element is complemented. But the complemented ideals, being compact, must all be principal, from which our thesis follows immediately.
\end{proofs}

If $\cal D$ is a Boolean algebra then the property of completeness of $\cal D$ is sufficient (as well as necessary) to ensure that condition $(iv)$ holds. We can show this as follows. If $\cal D$ is complete and $I$ is an ideal on $\cal D$ satisfying the property in condition $(iv)$, $I=(x)\downarrow$, where $x$ is the supremum of $I$ in $\cal D$. Indeed, it is clear that $I\subseteq (x)\downarrow$, while $x\in I$ because for any non-zero element $b\leq x$, since $b=b\wedge \mathbin{\mathop{\textrm{\huge $\vee$}}\limits_{a\in I}}a=\mathbin{\mathop{\textrm{\huge $\vee$}}\limits_{a\in I}}(b\wedge a)$ (note that, $\cal D$ being a complete Boolean algebra, $\cal D$ is a frame and hence the infinite distributive law holds in $\cal D$), there is an element $a\in I$ such that $b\wedge a\leq b$ is non-zero and belongs to $I$ ($I$ being an ideal). We have thus recovered Lemma III.3.5 \cite{stone}.

The next proposition represents the analogue of Proposition \ref{Boole2} for meet-semilattices and preorders. 

\begin{proposition}
\begin{enumerate}[(i)]

\item For any meet-semilattice $\cal M$, the ideal locale $S_{\cal M}$ (equivalently, the ideal topological space) associated to $\cal M$ is extremally disconnected;

\item For any preorder $({\cal P}, \leq)$, the Alexandrov locale $A_{\cal P}$ (equivalently, the Alexandrov space) associated to $\cal P$ is extremally disconnected if and only if $\cal P$ satisfies the amalgamation property (i.e. for any elements $a,b,c \in {\cal P}$ such that $c \leq a, b$ there exists $d\in {\cal P}$ such that $a,b \leq d$).
   
\end{enumerate}
\end{proposition}

\begin{proofs}
The thesis follows from a similar argument to that in the proof of Proposition \ref{Boole2}, by using the following site characterizations for the invariant property of a topos to be De Morgan (cf. \cite{El}):

\begin{enumerate}[(a)]
\item A presheaf topos $[{\cal C}^{\textrm{op}}, \Set]$ is De Morgan if and only if the category $\cal C$ satisfies the right Ore condition (i.e. the dual of the amalgamation property);

\item A localic topos $\Sh(L)$ is De Morgan if and only if the locale $L$ is extremally disconnected.
\end{enumerate}
\end{proofs}

\subsubsection{Triviality via two-valuedness}

By a two-valued locale we mean a locale $L$ such that the only two elements of $L$ are $0$ and $1$, and they are distinct from each other. Note that, for a topological space $X$, ${\cal O}(X)$ is two-valued if and only if the underlying set $X_{0}$ of $X$ is non-empty and the topology on $X$ is trivial.

Considering the invariant property of a topos to be two-valued, we obtain the following result. 

\begin{proposition}
\begin{enumerate}[(i)]

\item If $\cal D$ is a distributive lattice, its associated Stone locale $L_{\cal D}$ is two-valued (equivalently, its associated Stone space is trivial and non-empty) if and only if $\cal D$ is two-valued (i.e., the only two elements of ${\cal D}$ are $0$ and $1$, and they are distinct from each other);

\item If $\cal M$ is a meet-semilattice and $S_{\cal M}$ is the ideal locale associated to $\cal M$ then $S_{\cal M}$ is two-valued (equivalently, its associated ideal topological space is trivial and non-empty) if and only if $M$ is a singleton;

\item If $({\cal P}, \leq)$ is a preorder and $A_{\cal P}$ is the Alexandrov locale associated to $\cal P$ then $A_{\cal P}$ is two-valued (equivalently, its associated Alexandrov space is trivial and non-empty) if and only if for every $p,q\in {\cal P}$, $p\leq q$ and $q\leq p$.   
\end{enumerate}
\end{proposition}

\begin{proofs}
The method of proof is always the same as that employed in the proofs of the previous propositions. In this case, the invariant is the property of a topos to be two-valued, while the two different site representations for the topos are $\Sh(L_{\cal D}) \simeq \Sh({\cal D}, J_{\cal D})$ (resp. $\Sh(S_{\cal M})\simeq [{\cal M}^{\textrm{op}}, \Set]$, $\Sh(A_{\cal P})\simeq [{\cal P}, \Set]$).

The characterization of the property of two-valuedness in terms of the site $({\cal D}, J_{\cal D})$ is easily seen to yield the property of $\cal D$ to be two-valued (i.e., the condition that the only two elements of ${\cal D}$ are $0$ and $1$, and they are distinct from each other).
 
The other site characterizations, leading to $(i)$, $(ii)$ and $(iii)$, are the following:

\begin{enumerate}[(a)]
\item A presheaf topos $[{\cal C}^{\textrm{op}}, \Set]$ is two-valued if and only if the category $\cal C$ is strongly connected (i.e. for any two objects $a, b \in {\cal C}$, there exist arrows $a\to b$ and $b\to a$);

\item A localic topos $\Sh(L)$ is two-valued if and only if the locale $L$ is two-valued.
\end{enumerate}
\end{proofs}

\subsubsection{G\"{o}del-Dummett's law as an invariant}

Another logically-motivated topos-theoretic invariant which admits natural site characterizations is the property of a topos to satisfy G\"{o}del-Dummett's law, in the sense that the sequent 
\[
(\top \vdash_{x,y} (x\imp y) \vee (y\imp x))
\]
holds in the internal Heyting algebra of the topos given by the subobject classifier. One easily calculates, by using the well-known explicit descriptions of the internal Heyting algebra operations on the subobject classifier $\Omega_{\Sh({\cal C}, J)}$ of a topos $\Sh({\cal C}, J)$ of sheaves on a site (as for example given in chapter III of \cite{MM}), that a Grothendieck topos $\Sh({\cal C}, J)$ satisfies G\"{o}del-Dummett's law  if and only if for any $J$-closed sieves $R$ and $S$ on an object $c\in {\cal C}$, the sieve $\{f:d\to c \textrm{ | } f^{\ast}(R)\subseteq f^{\ast}(S) \textrm{ or } f^{\ast}(S)\subseteq f^{\ast}(R)\}$ is $J$-covering. In particular, a presheaf topos $[{\cal C}^{\textrm{op}}, \Set]$ satisfies G\"{o}del-Dummett's law if and only if $\cal C$ satisfies the following property: for any arrows $f:b\to a$ and $g:c\to a$ with common codomain, either $f$ factors through $g$ or $g$ factors through $f$ (cf. also Proposition 3.1 \cite{morgan}). On the other hand, in \cite{morgan} Johnstone established the following result: for any topological space $X$, the topos $\Sh(X)$ satisfies G\"{o}del-Dummett logic if and only if every closed subspace of $X$ is extremally disconnected.  

Let us unravel the property of the topos $\Sh({\cal D}, J_{\cal D})$ of coherent sheaves on a distributive lattice $\cal D$ to satisfy G\"{o}del-Dummett's law in terms of the site $({\cal D}, J_{\cal D})$. Given two sets $A:=\{a_{i}\leq x \textrm{ | } i\in I\}$ and $B:=\{b_{j}\leq x \textrm{ | } j\in J\}$ of elements of $\cal D$, we say that $A$ \emph{refines} $B$ if the sieve in $\cal D$ on $x$ generated by $A$ is contained in the sieve in $\cal D$ on $x$ generated by $B$, equivalently if for every $i\in I$ there exists $j\in J$ such that $a_{i}\leq b_{j}$. We define a set of elements $A:=\{a_{i}\leq x \textrm{ | } i\in I\}$ to be \emph{finitely closed} if the sieve on $x$ generated by it is $J_{\cal D}$-closed (i.e., for every $b\leq a$, if there is a finite subset $I'$ of $I$ such that $b=\mathbin{\mathop{\textrm{\huge $\vee$}}\limits_{i\in I'}}(a_{i}\wedge b)$ then $b\leq a_{i}$ for some $i\in I$). In these terms, the condition that $\Sh({\cal D}, J_{\cal D})$ satisfies G\"{o}del-Dummett's law rephrases as follows: for any finitely closed sets of elements $A:=\{a_{i}\leq x \textrm{ | } i\in I\}$ and $B:=\{b_{j}\leq x \textrm{ | } j\in J\}$ in $\cal D$, there exists a finite collection $\{c_{k} \textrm{ | } k\in K\}$ of elements satisfying $x=\mathbin{\mathop{\textrm{\huge $\vee$}}\limits_{k\in K}}c_{k}$ such that for any $k\in K$ either $A_{c_{k}}:=\{a_{i}\wedge c_{k} \leq c_{k} \textrm{ | } i\in I\}$ refines $B_{c_{k}}:=\{b_{j}\wedge c_{k}\leq c_{k} \textrm{ | } j\in J\}$ or $B_{c_{k}}$ refines $A_{c_{k}}$. We call a distributive lattice $\cal D$ satisfying this property a G\"{o}del-Dummett distributive lattice.  

Combining together all the site characterizations discussed above, the usual method `toposes as bridges' applied to the invariant `to satisfy G\"{o}del-Dumm-\\ett's law' yields the following result.   

\begin{theorem}
With the notation above, we have:
\begin{enumerate}[(i)]

\item If $\cal D$ is a distributive lattice and $X_{\cal D}$ is its associated Stone space then every closed subspace of $X_{\cal D}$ is extremally disconnected if and only if $\cal D$ is a G\"{o}del-Dummett distributive lattice;  

\item If $({\cal M}, \leq)$ is a meet-semilattice and $X_{\cal M}$ is the ideal topological space associated to $\cal M$ then every closed subspace of $X_{\cal M}$ is extremally disconnected if and only if for every $p,q\in {\cal M}$ such that $p,q\leq r$ for some $r$ in $\cal M$, either $p\leq q$ or $q\leq p$ (i.e., if and only if ${\cal M}$ is a forest);   

\item If $({\cal P}, \leq)$ is a preorder and $A_{\cal P}$ is the Alexandrov space associated to $\cal P$ then $A_{\cal P}$ has the property that every closed subspace of $X$ is extremally disconnected if and only if for every $p,q\in {\cal P}$ such that $r\leq p,q$ for some $r$ in $\cal P$, either $p\leq q$ or $q\leq p$ (i.e., if and only if ${\cal P}^{\textrm{op}}$ is a forest).   
\end{enumerate}
\end{theorem}\qed

\section{Spaces of models of propositional theories}\label{spacesprop}

We have seen in Example \ref{exa}(f) that, for any geometric theory $\mathbb T$, the space of models of $\mathbb T$ in $\Set$ can be endowed with a natural topology, of logical nature, which is a particular case of the subterminal topology introduced in section \ref{subterminaltop}. 

In the following sections we focus our attention on the logical spaces of propositional geometric theories. In section \ref{subs} we give an explicit description of the logical topology on a set of models of a propositional geometric theory. Then we characterize the frames of open sets of these logical topological spaces in terms of preordered structures presented by generators and relations. To this end, we introduce in section \ref{gensyn} an abstract notion of first-order mathematical theory, and a corresponding notion of syntactic category, which is seen in sections \ref{genrel} and \ref{clex} to subsume the known ones and to provide a uniform way for building structures presented by generators and relations for certain `ordered algebraic theories'.

\subsection{Subsets and propositional theories}\label{subs}

Propositional theories are particularly convenient means for describing subsets of a given set having particular properties, for example ideals of a commutative ring (cf. section \ref{zariski}), filters on a meet-semilattice (cf. section \ref{logical}), etc.

A propositional geometric theory $\mathbb T$ can be formally defined as a geometric theory over a signature $\Sigma_{\mathbb T}$ with no sorts (cf. Part D of \cite{El}). Thus $\Sigma_{\mathbb T}$ has no function symbols or constants, and consists only of a set of $0$-ary relation symbols ${\cal R}_{\mathbb T}$.  Every symbol $R$ in ${\cal R}_{\mathbb T}$ gives therefore rise to an atomic formula, and every atomic formula is of this form.  Note that $\Sigma_{\mathbb T}$-structures in $\Set$ can be identified with functions $f_{M}:{\cal R}_{\mathbb T}\to \{0,1\}$, or equivalently with subsets $f^{-1}(\{1\})\subseteq {\cal R}_{\mathbb T}$ of ${\cal R}_{\mathbb T}$.

Let $\mathbb T$ be a propositional geometric theory over a signature $\Sigma_{\mathbb T}$. 

Recall from Example \ref{exa}(f) that the collection $X_{\mathbb T}$ of the (isomorphism classes of) models of $\mathbb T$ in $\Set$, endowed with the subterminal topology, has as open sets exactly those of the form $F_{\phi}=\{M \in X_{\mathbb T} \textrm{ | } M \vDash \phi \}$, where $\phi$ ranges among all the geometric sentences over $\Sigma_{\mathbb T}$. In particular, since (by Lemma D1.3.8(ii)) every geometric formula over $\Sigma_{\mathbb T}$ is provably equivalent to a disjunction of conjunctions of atomic formulae in the same context, the collection of subsets of the form $F_{R}=\{M \in X_{\mathbb T} \textrm{ | } M \vDash R\}$ for an atomic formula $R$ over $\Sigma_{\mathbb T}$ forms a subbasis of the topological space $X_{\mathbb T}$; that is, the open sets of $X_{\mathbb T}$ are exactly the unions of finite intersections of sets of the form $F_{R}$ for a $0$-ary relation symbol $R$ over $\Sigma_{\mathbb T}$. Clearly, if the conjunction of any two atomic formulae over $\Sigma_{\mathbb T}$ is $\mathbb T$-provably equivalent to an atomic formula then the sets of the form $F_{R}$ (for $R\in {\cal R}_{\mathbb T}$) form a basis of the topological space $X_{\mathbb T}$ (cf. section \ref{zariski} below for a concrete example of this situation).
  
Using the identification of $\Sigma_{\mathbb T}$-structures in $\Set$ with subsets of ${\cal R}_{\mathbb T}$, the space $X_{\mathbb T}$ acquires the following description: its points are the subsets of ${\cal R}_{\mathbb T}$ which correspond to models of $\mathbb T$, and its open sets are the unions of subsets of the form $\{L\subseteq {\cal R}_{\mathbb T} \textrm{ | } R_{1}, \ldots, R_{n} \in L\}$ for a finite number $R_{1}, \ldots, R_{n}$ of $0$-ary relation symbols over $\Sigma_{\mathbb T}$.   

As an example, take $\mathbb T$ to be the empty theory over a signature $\Sigma$ consisting of a set $A$ of $0$-ary relation symbols.
The $\Sigma$-structures in $\Set$ can be identified with the subsets of $A$, and the resulting topology on the powerset $\mathscr{P}(A)$ is given by the collection of subsets which are unions of subsets of the form $\{L\subseteq A \textrm{ | } a_{1}, \ldots, a_{n} \in L\}$ for a finite number $a_{1}, \ldots, a_{n}$ of elements of $A$. This topology is an interesting one; in fact, as we shall see in section \ref{zariski}, it specializes, under a natural bijection, precisely to the Zariski topology on the prime spectrum of a commutative ring with unit, and the frame of open sets of the resulting topological space with underlying set $A$ can be characterized as the free frame on $A$ (by the results of sections \ref{genrel} and \ref{clex} below). We will refer below to this topology as to the \emph{elemental topology}. 

Given a propositional geometric theory $\mathbb T$ over a signature $\Sigma_{\mathbb T}$, we have an embedding $\xi:{\cal R}_{\mathbb T} \to {\cal C}_{\mathbb T}$ of ${\cal R}_{\mathbb T}$ into the (underlying set of the) geometric syntactic category ${\cal C}_{\mathbb T}$ of $\mathbb T$, sending a symbol $R$ in ${\cal R}_{\mathbb T}$ to the corresponding atomic formula over $\Sigma_{\mathbb T}$. The operation of taking the inverse image of subsets under $\xi$ defines a bijection between the collection of completely prime filters on ${\cal C}_{\mathbb T}$ and the collection of the subsets of ${\cal R}_{\mathbb T}$ which correspond to models of $\mathbb T$; in fact, this bijection becomes a homeomorphism of topological spaces if we endow the two sets respectively with the topology on the set of points of a locale (cf. Example \ref{exa}(e)) and with the elemental topology.  

Of course, analogous results hold if we replace the geometric syntactic category ${\cal C}_{\mathbb T}$ of $\mathbb T$ with a (cartesian) syntactic category appropriate to a given fragment of logic in which the theory $\mathbb T$ lies (for example, if we replace ${\cal C}_{\mathbb T}$ with the coherent syntactic category of $\mathbb T$ in case $\mathbb T$ is coherent), and we endow the resulting space of ($J$-)filters on the corresponding syntactic category with the subterminal topology as in Proposition \ref{mslattice} above.  
 
Note that these results arise from the expression of a particular invariant, namely the subterminal topology, in terms of different representations of the locales (equivalently, of the toposes) involved.

\subsection{General first-order theories}\label{gensyn}

In this section we introduce a general notion of first-order mathematical theory and a corresponding notion of syntactic category, which encompass all the instances of the concepts given in chapter D1 of \cite{El}. We shall present these notions in a sketchy form, since a completely formal and detailed treatment of them would bring us far beyond the scope of this paper; anyway, we are confident that the interested reader will have no trouble in filling in the details by generalizing the framework of chapter D1 of \cite{El} according to our indications. 

Let $\Sigma$ be a first-order signature. Let us suppose to have a collection $\cal S$ of symbols $s$ (to be thought of as generalized connectives) each of which equipped with an arity $ar(s)$ given by a cardinal number; we treat constants as $0$-ary function symbols, and we also allow the symbols in ${\cal S}$ to have arity $\infty$, i.e. to take as inputs sets of arbitrary cardinality. Starting with atomic formulae over $\Sigma$ we can inductively (and freely) build a collection $F_{({\cal S}, \Sigma)}$ of words, which we call ${\cal S}$-formulae over $\Sigma$, as follows: $F_{({\cal S}, \Sigma)}$ is the smallest set of words such that all the atomic formulae over $\Sigma$ belong to $F_{({\cal S}, \Sigma)}$, and for any $s\in {\cal S}$ and any family of words $\{F_{i} \textrm{ | } i\in ar(s)\}$ in $F_{({\cal S}, \Sigma)}$ the word $s(\{F_{i} \textrm{ | } i\in ar(s)\})$ belongs to $F_{({\cal S}, \Sigma)}$; the notion of free variable in a word of $F_{({\cal S},\Sigma)}$ is defined by simultaneous recursion as in classical cases, by defining the set of free variables of the word $s(F_{1}, \ldots F_{a(s)})$ equal to the union of the sets of free variables of the $F_{i}$. Of course, one could define a larger fragment of words by requiring $F_{({\cal S}, \Sigma)}$ to be closed also under existential or universal quantifications, but we shall not pursue these generalizations in this paper, our interest being primarily focused on the propositional fragments of these generalized logics. We say that a ${\cal S}$-theory over a signature $\Sigma$ is \emph{propositional} if $\Sigma$ has no sorts.  

Given a generalized set of connectives $\cal S$ and a first-order signature $\Sigma$, we define a \emph{${\cal S}$-theory} over $\Sigma$ to be a collection of axioms and inference rules involving sequents of the form $\phi \vdash_{\vec{x}} \psi$ where $\phi$ and $\psi$ are formulae in $F_{(\cal S, \Sigma)}$ in the context $\vec{x}$ . Note that we have a notion of provability in a $\cal S$-theory of sequents involving ${\cal S}$-formulae over its signature. 

Any generalized signature $({\cal S}, \Sigma)$ determines a class of categories in which the formulae in $F_{({\cal S}, \Sigma)}$ are interpretable. Specifically, given a category $\cal C$, we interpret the $0$-ary relation symbols over $\Sigma$ as objects $c$ of $\cal C$ such that for any object $c'\in {\cal C}$ there is at most one arrow $c'\to c$ in $\cal C$; we call these objects the \emph{subterminal objects} of $\cal C$ even if $\cal C$ does not have a terminal object (in fact, in the latter case, these objects are precisely the subterminal objects in the classical sense, that is the objects $c$ such that the unique arrow $c\to 1_{\cal C}$ is monic). If for any connective $s$ in $\cal S$, there is an operation $f_{s}$ of arity $ar(s)$ on the set of subterminal objects of $\cal C$ (if $ar(s)=0$, a subterminal object of $\cal C$), the interpretation of the formulae in $F_{({\cal S}, \Sigma)}$ with no free variables can be defined inductively induction on their structure by setting the interpretation of a formula $s(\{F_{i} \textrm{ | } i\in ar(s)\})$ equal to the result of applying the operation $f_{s}$ to the subterminals given by the interpretation of the formulae $F_{i}$. Provided that $\cal C$ has at least finite limits and if for any object $c$ of $\cal C$ and any connective $s$ in $\cal S$ there is an operation $f_{s}$ of arity $ar(s)$ on the set of subobjects of $c$ in $\cal C$, the formulae with at least a free variable can be similarly interpreted in a $\Sigma$-structure in $\cal C$ (as subobjects of (products of) the underlying object(s) of the structure) by interpreting the atomic formulae over $\Sigma$ as usual, and setting the interpretation of a formula $s(\{F_{i} \textrm{ | } i\in ar(s)\})$ equal to the result of applying the operation $f_{s}$ to the subobjects given by the interpretations of the formulae $F_{i}$. We shall call a category equipped with operations on its subobjects (if $\Sigma$ has at least one sort) or subterminals (if $\Sigma$ has no sorts) which interpret the connectives in $\cal S$ a \emph{$({\cal S}, \Sigma)$-category}.  

Given a functor $F:{\cal C}\to {\cal D}$ between two $({\cal S}, \Sigma)$-categories, we say that $F$ is \emph{$({\cal S}, \Sigma)$-preserving} if either $\Sigma$ has no sorts and $F$ sends subterminal objects of $\cal C$ to subterminal objects of $\cal D$ and commutes with the operations on subterminals in the two categories which interpret the connectives in $\cal S$, or, if $\Sigma$ has at least one sort, $F$ preserves finite limits and commutes with the operations on subobjects in the two categories which interpret the connectives in $\cal S$. The category of $({\cal S}, \Sigma)$-preserving functors ${\cal C} \to {\cal D}$ and natural transformations between them will be denoted by $({\cal S}, \Sigma)\textrm{-}\textbf{Fun}({\cal C}, {\cal D})$.   

Given a $\cal S$-theory $\mathbb T$ over a signature $\Sigma$, and a $({\cal S}, \Sigma)$-category $\cal C$, we define a \emph{model} of $\mathbb T$ in $\cal C$ to be a $\Sigma$-structure in $\cal C$ in which all the axioms of $\mathbb T$ are valid and the inference rules of $\mathbb T$ are sound (as in \cite{El}, we say that a sequent $\phi \vdash_{\vec{x}} \psi$ is valid in the structure if the interpretation of $\phi(\vec{x})$ factors, as a subobject or subterminal, through the interpretation of $\psi(\vec{x})$, while the soundness of inference rules has the usual meaning). We remark that, unlike in Definition D1.2.1 of \cite{El}, we also consider $\Sigma$-structures in categories lacking finite products or a terminal object, provided that $\Sigma$ has no sorts (in fact, as explained above, we interpret $0$-ary relation symbols as subterminal objects in the category); this extra-generality will be crucial for our purposes in constructing syntactic categories of propositional theories which lie in fragments which are weaker than Horn logic.

We define a homomorphism of models of a $\cal S$-theory $\mathbb T$ over a signature $\Sigma$ in a $({\cal S}, \Sigma)$-category $\cal C$ to be a homomorphism of the underlying $\Sigma$-structures. Clearly, models of $\mathbb T$ in $\cal C$ and homomorphisms between them form a category, which we denote by ${\mathbb T}\textrm{-mod}_{\cal S}({\cal C})$. 

If $\mathbb T$ is a propositional ${\cal S}$-theory over a signature $\Sigma$ then a model of $\mathbb T$ in a $({\cal S}, \Sigma)$-preorder category $({\cal P}, \leq)$ can be identified with a function sending every $0$-ary relation symbol over $\Sigma$ to an element of $\cal P$ in such a way its extension $f:F_{({\cal S}, \Sigma)}\to {\cal P}$ to the set of all $\cal S$-formulae over $\Sigma$ satisfies the property that whenever $F \vdash_{[]} G$ is provable in $\mathbb T$, $f(F)\leq f(G)$ in $\cal P$. Note that, if $\cal S$ consists of the usual connectives $\wedge$ and $\vee$, that is, in the case of usual first-order logic, and $\cal P$ has a top element, this notion specializes to the classical one; indeed, the elements of such a preorder $\cal P$ can be clearly identified with the subobjects of the terminal object of $\cal P$, when the latter is regarded as a preorder category. 

For a $\cal S$-theory $\mathbb T$ over a signature $\Sigma$, we define the \emph{syntactic category} ${\cal C}^{({\cal S}, \Sigma)}_{\mathbb T}$ of $\mathbb T$ as follows. The objects of ${\cal C}^{({\cal S}, \Sigma)}_{\mathbb T}$ are the $\mathbb T$-provable equivalence classes of $\cal S$-formulae-in-context over $\Sigma$ (considered up to `renaming' equivalence), and the arrows of ${\cal C}^{({\cal S}, \Sigma)}_{\mathbb T}$ are the $\mathbb T$-provable equivalence classes of $\cal S$-formulae over $\Sigma$ which are $\mathbb T$-provably functional from the domain to the codomain (this is a condition which ensures that their interpretation in any model of $\mathbb T$ in a $({\cal S}, \Sigma)$-theory is the graph of an arrow from the interpretation of the $\cal S$-formula in the domain to the interpretation of the $\cal S$-formula in the codomain, cf. p. 841 \cite{El} for the details). If the signature $\Sigma$ has at least one sort and the fragment $F_{({\cal S}, \Sigma)}$ contains cartesian logic, we have a model $M_{\mathbb T}$ of $\mathbb T$ in ${\cal C}^{({\cal S}, \Sigma)}_{\mathbb T}$, defined exactly as in p. 844 \cite{El}, in which the $\cal S$-sequents over $\Sigma$ which are valid coincide precisely with those which are provable in $\mathbb T$, and, as in \cite{El} (Theorem D1.4.7 etc.), we have a categorical equivalence
\[
({\cal S}, \Sigma)\textrm{-}\textbf{Fun}({\cal C}^{({\cal S}, \Sigma)}_{\mathbb T}, {\cal C})\simeq {\mathbb T}\textrm{-mod}_{\cal S}({\cal C})
\]
for any $({\cal S}, \Sigma)$-category $\cal C$.

Note that the syntactic category ${\cal C}^{({\cal S}, \Sigma)}_{\mathbb T}$ of a propositional $\cal S$-theory $\mathbb T$ over a signature $\Sigma$ can be identified with the poset obtained by equipping the set of provable-equivalence classes of $\cal S$-formulae over $\Sigma$ with the order given by the relation of provable entailment; the model $M_{\mathbb T}$ is defined as the function sending a $0$-ary relation symbol over $\Sigma$ to the $\mathbb T$-provable equivalence class of the corresponding atomic formula, and we have a categorical equivalence $({\cal S}, \Sigma)\textrm{-}\textbf{Fun}({\cal C}^{({\cal S}, \Sigma)}_{\mathbb T}, {\cal C})\simeq {\mathbb T}\textrm{-mod}_{\cal S}({\cal C})$ exactly as above. In particular, if the axioms of $\mathbb T$ are all bisequents then the preorder ${\cal C}^{({\cal S}, \Sigma)}_{\mathbb T}$ is a discrete category (since if $F \vdash_{[]} G$ is provable in $\mathbb T$ then necessarily $F$ is provably equivalent to $G$ in $\mathbb T$).     

All the notions introduced in this section are meant to provide a general notion of fragment of first-order logic which encompasses the (quantifier-free versions of the) classical ones, including Horn, regular, coherent and geometric logic (cf. Part D of \cite{El}). The usual logical connectives, as well as the inference rules which govern them, appear as instances of an abstract notion of connective (respectively, of inference rule) within the unifying framework that we have developed. Specifically, quantifier-free Horn (resp. coherent, geometric, finitary first-order) can be seen as a $\cal S$-theory, where $\cal S$ consists of the binary connective $\wedge$ and the constant $\top$ (resp. of the binary connectives $\wedge$ and $\vee$ and the constants $\top$ and $\bot$, of the binary connective $\wedge$, the infinitary connective $\mathbin{\mathop{\textrm{\huge $\vee$}}\limits_{}}$ and the constants $\top$ and $\bot$, of the binary connectives $\wedge$, $\vee$, $\neg$, $\imp$ and the constants $\top$ and $\bot$), and the corresponding notion of model of a regular (resp. coherent, geometric, finitary first-order logic) theory in a regular (resp. coherent, geometric, Heyting) category is subsumed by our notion of model in a $({\cal S}, \Sigma)$-category. In fact, if we regard a (quantifier-free) Horn (resp. coherent, geometric, first-order) theory as a $\cal S$-theory then our notion of syntactic category yields a category which is equivalent to the cartesian (resp. coherent, geometric, first-order) syntactic category of the theory, as defined in chapter D1 of \cite{El}.     

We note that the main technical differences between our approach and that of \cite{El} consist in allowing arbitrary inference rules in the definition of a ($\cal S$-)theory, and in extending the notion of structure in a category in such a way as to allow the interpretation in categories of fragments of logic which are weaker than Horn logic (in particular, to allow the interpretation of propositional formulae in preorders which do not possess meets or a top element).  

\subsection{Generators and relations}\label{genrel}

In this section we show that the syntactic categories of propositional $\cal S$-theories can be naturally realized as structures presented by generators and relations (in the sense of universal algebra).

The following definition will be central for our purposes. 

\begin{definition}
\begin{enumerate}[(a)]

\item An \emph{infinitary Horn formula} over a first-order signature $\Sigma$ is a formula built from atomic formulae by only using possibly infinitary conjunctions (including the conjunction over the empty set, which we identify with the truth formula $\top$);

\item An \emph{infinitary Horn theory} over a signature $\Sigma$ is an infinitary first-order theory over $\Sigma$ whose axioms are \emph{infinitary Horn sequents}, that is sequents of the form $\phi \vdash_{\vec{x}} \psi$ where $\phi$ and $\psi$ are infinitary Horn formulae in the same context $\vec{x}$;

\item By \emph{infinitary Horn logic} we mean the logic of infinitary Horn theories, i.e. the fragment of infinitary first-order logic in which only the structural rules and the rules for infinitary conjunctions are present;

\item An \emph{ordered algebraic theory} is an infinitary Horn theory over a one-sorted signature $\Sigma$ consisting of a set of (possibly infinitary) function symbols and a binary relation symbol $\leq$, in which the following three Horn sequents, expressing the idea that $\leq$ is a partial order, are provable (in Horn logic):
\[
(\top \vdash_{x} x\leq x),
\]
\[
(x\leq y \vdash_{x,y} y\leq x),
\]
\[
((x\leq y) \wedge (y\leq z) \vdash_{x,y,z} x\leq z).
\]

\item Following \cite{freestr}, we say that an ordered algebraic theory $\mathbb T$ over a signature $\Sigma$ is \emph{small} if for every cardinal number $k$, there is only a \emph{set} of $k$-ary terms over $\Sigma$ up to $\mathbb T$-provable equivalence.
\end{enumerate}

\end{definition}

\begin{remarks}\label{rem}
\begin{enumerate}[(a)]

\item We can suppose, without loss of generality, all the axioms of an infinitary Horn theory to be of the form $F_{1} \vdash_{\vec{x}} F_{2}$, where $F_{1}(\vec{x})$ is a (possibly infinitary) conjunction of atomic formulae and $F_{2}(\vec{x})$ is an atomic formula; 

\item Any (possibly infinitary) algebraic theory $\mathbb T$ over a signature $\Sigma$ can be considered as an ordered algebraic theory whose axioms are of the form $\top \vdash_{\vec{x}} t_{1}=t_{2}$ where $t_{1}$ and $t_{2}$ are terms over $\Sigma$, and in which the sequent 
\[
(x\leq y \vdash_{x,y} x=y)
\]
is provable. 
\end{enumerate}
\end{remarks}

Let us now introduce the notion of model of an ordered algebraic theory presented by generators and relations.

Let $M$ be a model (in $\Set$) of an ordered algebraic theory $\mathbb A$ over a signature $\Sigma_{\mathbb A}$. Given a set $A$ of \emph{generators}, let $\Sigma_{\mathbb A}^{A}$ be the signature having, in addition to the function symbols of $\Sigma_{\mathbb A}$, a constant symbol $C_{a}$ for each element $a\in A$. Note that, given a function $\xi:A\to M$ from $A$ to (the underlying set of) a $\Sigma_{\mathbb A}$-structure $M$, $M$ can be made into a $\Sigma_{\mathbb A}^{A}$-structure, where for any $a\in A$ the interpretation of the constant symbol $C_{a}$ is given by the element $\xi(a)\in M$. Given a set $R$ of \emph{relations}, i.e. infinitary Horn sequents of the form $\phi \vdash_{[]} \psi$ where $\phi$ and $\psi$ are closed infinitary Horn formulae over the signature $\Sigma_{\mathbb A}^{A}$, we denote by ${\mathbb T}_{A, R}$ the infinitary Horn theory over $\Sigma_{\mathbb A}$ having as axioms all the axioms of $\mathbb A$ and the sequents in $R$. Clearly, ${\mathbb T}_{A, R}$ is an ordered algebraic theory over the signature $\Sigma_{\mathbb A}^{A}$. 

\begin{definition}
With the notation above, given a function $\xi:A\to M$ and a set $R$ of (infinitary) Horn sequents of the form $\phi \vdash_{[]} \psi$ where $\phi$ and $\psi$ are closed infinitary Horn formulae over the signature $\Sigma_{\mathbb A}^{A}$, we say that $M$ is presented via $\xi$ by the set of \emph{generators} $A$ subject to the \emph{relations} $R$, briefly that $M$ is presented by $(A, R)$, if $M$, regarded as a $\Sigma_{\mathbb A}^{A}$-structure as specified above, is an initial object in the category ${\mathbb T}_{A, R}\textrm{-mod}(\Set)$ of ${\mathbb T}_{A, R}$-models in $\Set$ and homomorphisms between them.   
\end{definition} 

Intuitively, a model $M$ of $\mathbb A$ in $\Set$ is presented by the set of generators $A$ via $\xi:A\to M$ subject to relations $R$ if all the `relations' in $R$ are satisfied in $M$ when evaluated in $M$ via $\xi$ and for any function $f:A\to N$ to a model $N$ of $\mathbb A$ in $\Set$ such that all the relations in $R$ are satisfied in $N$ when evaluated in $N$ via $f$ there exists a unique $\mathbb A$-model homomorphism $g:M\to N$ such that $g\circ \xi=f$.   

The following result is a natural generalization of the classical construction of free models of algebraic theories. 

\begin{theorem}\label{oat}
Let $\mathbb A$ be an small ordered algebraic theory, $A$ be a set and $R$ be a set of relations over $\Sigma_{\mathbb A}^{A}$ (in the sense specified above).
 
Let $M_{(A, R)}$ be the set of equivalence classes $[t]$ of closed terms $t$ over $\Sigma^{A}_{\mathbb A}$ with respect to the equivalence relation $E$ defined by
\[
(t_{1}, t_{2})\in E \textrm{ if and only if } \top \vdash_{[]} t_{1}=t_{2} \textrm{ is provable in } {\mathbb T}_{A, R}.
\]
If we define the interpretation in $M_{(A, R)}$ of the symbol $\leq$ as the relation 
\begin{center}
$[t_{1}]\leq [t_{2}]$ if and only if $t_{1} \vdash_{[]} t_{2}$ is provable in ${\mathbb T}_{A, R}$.
\end{center}
and we set the interpretation in $M_{(A, R)}$ of a function symbol $f$ of arity $k$ over $\Sigma_{\mathbb A}$ equal to the function $M_{(A, R)}^{n}\to M_{(A, R)}$ sending any $k$-tuple $\{[t_{i}] \textrm{ | } i\in k\}$ in $M_{(A, R)}^{k}$ to the element $[f(\{[t_{i}] \textrm{ | } i\in k\})]\in M_{(A, R)}$, $M_{(A, R)}$ becomes a model of $\mathbb A$ which is presented by the set of generators $A$ and relations $R$ via the function $\xi:A\to M_{(A, R)}$ sending any $a\in A$ to $\xi(a):=[C_{a}]$.
\end{theorem}

\begin{proofs}
First, we note that, since $\mathbb A$ is small by hypothesis, $M_{(A, R)}$ is actually \emph{set}, rather than a proper class. 

It is immediate to see, by using the deduction theorem in first-order logic, that for any formula $\phi(\vec{x})$ over $\Sigma_{\mathbb A}$ and any closed terms $\vec{t}$ over $\Sigma_{\mathbb A}^{A}$ (substitutable in place of $\vec{x}$), the formula $\phi$ is valid in $M_{(A, R)}$ when evaluated in $\vec{t}$ (regarded as a string of elements of $M_{(A, R)}$) if and only if the sequent $\top \vdash_{[]} \phi([\vec{t}\slash \vec{x}])$ is provable in the theory ${\mathbb T}_{A, R}$. From this it easily follows that all the axioms of $\mathbb A$ are valid in $M_{(A, R)}$. Clearly, $M_{(A, R)}$ is a model of ${\mathbb T}_{A, R}$ in $\Set$, and is presented by $A$ via $\xi$ subject to relations $R$, since for any ${\mathbb T}_{A, R}$-model $N$, regarded a function $f:A\to N$ to a $\mathbb A$-model $N$, there exists a unique $\mathbb A$-model homomorphism $g:M_{(A, R)} \to N$ such that $g\circ \xi=f$.    
\end{proofs}        

Let us now show how the construction of general syntactic categories provides an alternative way for building models of small ordered algebraic theories presented by generators and relations.

Let $\mathbb A$ be an ordered algebraic theory over a signature $\Sigma_{\mathbb A}$. For each function symbol $f$ over $\Sigma_{\mathbb A}$ of arity $k$ we define a corresponding connective on formulae $s_{f}$ of arity $k$. Given a set $A$, we define a propositional signature $\Sigma_{A}$ having exactly one $0$-ary relation symbol $R_{a}$ for each element $a\in A$. We denote by $F_{{\cal S}, A}$ the collection of $\cal S$-formulae over $\Sigma_{A}$, as defined in section \ref{gensyn}.  

We have a natural bijection between the closed terms over $\Sigma_{\mathbb A}^{A}$ and the $\cal S$-formulae over $\Sigma_{A}$. Indeed, we can define a correspondence between them inductively, as follows: to a constant $C_{a}$ (for $a\in A$) over $\Sigma_{\mathbb A}^{A}$ we associate the atomic formula $R_{a}$, and to a term of the form $f(t_{1}, \ldots t_{n})$, for a function symbol $f$ over $\Sigma_{\mathbb A}$ of arity $k$ and a collection of terms $\{t_{i} \textrm{ | } i\in k\}$ corresponding to $\cal S$-formulae $\{F_{i} \textrm{ | } i\in k\}$, we associate the formula $s_{f}(\{F_{i} \textrm{ | } i\in k\})$. Clearly, this correspondence is a bijection; we denote by $F_{t}$ the $\cal S$-formula over $\Sigma_{A}$ corresponding to a closed term $t$ over $\Sigma_{\mathbb A}^{A}$.

Let us build a $\cal S$-theory ${\mathbb S}_{A, R}$ over $\Sigma_{A}$, as follows.

\begin{itemize} 

\item First, we put as axioms of ${\mathbb S}_{A, R}$ the usual structural rules (as described in section D1.3 \cite{El}). 

\item Second, to each axiom $\phi \vdash_{\vec{x}} \psi$ of the theory $\mathbb A$, where $\phi$ is either $\top$ or a conjunction of atomic formulae over $\Sigma_{\mathbb A}$ and $\psi$ is an atomic formula (cf. Remark \ref{rem2}(b)), we associate a scheme of inference rules of ${\mathbb S}_{A, R}$, each obtained by substituting arbitrary formulae $\vec{F}$ in $F_{\cal S}$ in place of the variables $\vec{x}$ in the following way. To each atomic formula $\chi$ over $\Sigma_{\mathbb A}$ we associate a sequent $\Gamma_{\chi}$ involving $\cal S$-formulae over $\Sigma_{A}$, obtained by putting the formulae $\vec{F}$ in place of the variables of the terms occurring in the atomic formulae, replacing the function symbols $f$ over $\Sigma_{\mathbb A}$ with the corresponding connectives $s_{f}$ in $\cal S$ and replacing the relation $\leq$ (resp. the equality $=$) between terms with the implication $\vdash$ (resp. the biimplication $\dashv \vdash$) (in case $\chi$ is $\top$ then we define $\Gamma_{\chi}$ to be a tautological sequent, as for example provided by the structural rules); we then define the inference rule corresponding to an axiom $\phi \vdash_{\vec{x}} \psi$ of $\mathbb A$ to be the rule having as premises all the sequents of the form $\Gamma_{\chi}$ where $\chi$ ranges among the atomic subformulae of $\phi$, and as conclusion the sequent $\Gamma_{\psi}$. So, for example, the sequent 
\[
((x \cdot (y + z) \leq (x \cdot y) + (x \cdot z) \wedge (x=y)) \vdash_{x,y,z} x=y\cdot z)
\]
in the signature of rings corresponds to the inference rule whose premis-\\es are the sequents
\[
(s_{\cdot}(F_{1}, s_{+}(F_{2}, F_{3})) \vdash s_{+}(s_{\cdot}(F_{1}, F_{2}), s_{\cdot}(F_{1}, F_{2})))
\]
and 
\[
(F_{1} \dashv\vdash F_{2})
\]
and whose conclusion is the sequent
\[
(F_{1} \dashv\vdash s_{\cdot}(F_{2}, F_{3})).
\]

\item Third, to each of the (possibly infinitary) Horn sequents over the signature $\Sigma_{\mathbb A}^{A}$ corresponding to the relations in $R$, we associate an inference rule obtained by the method of the last paragraph, the only difference being that the constants $C_{a}$ (rather than the variables) arising in the terms are replaced by the corresponding $\cal S$-formulae $R_{a}$ (for each $a \in A$). 

\item Finally, we add inference rules asserting the invariance under provable equivalence of the connectives; that is, for each connective $s\in {\cal S}$ of arity $k$, we add an inference rule scheme which enables to derive the bisequent $(s_{f}(\{F_{i} \textrm{ | } i\in k\}) \dashv \vdash s_{f}(\{F_{i}' \textrm{ | } i\in k\}))$ from the set of bisequents $(F_{i} \dashv \vdash F_{i}')$ (for any $\cal S$-formulae $F_{i}$ and $F_{i}'$ (for $i\in k$) over $\Sigma_{A}$).      
\end{itemize}
 
This completes the description of the $\cal S$-theory ${\mathbb S}_{A, R}$.

Let us now show that the (underlying poset of the) syntactic category ${\cal C}^{({\cal S}, \Sigma_{A})}_{{\mathbb S}_{A, R}}$ of the $\cal S$-theory ${\mathbb S}_{A, R}$ is a model of the theory ${\mathbb T}_{A, R}$ presented by $(A, R)$. 

First, we note that, since our theory ${\mathbb S}$ is small by hypothesis, the underlying preorder of the category ${\cal C}^{({\cal S}, \Sigma_{A})}_{{\mathbb S}_{A, R}}$ is actually a \emph{set}.

Next, we observe that ${\cal C}_{{\mathbb S}_{A, R}}$ can be made into a ${\mathbb T}_{A, R}$-model structure, by defining the interpretation of any generator $a\in A$ to be the ${\mathbb S}_{A, R}$-provable equivalence class $[R_{a}]$ of the atomic formula $R_{a}$. One can easily verify by induction on the structure of $\cal S$-formulae over $\Sigma_{A}$ that for any term $t$ over $\Sigma_{\mathbb A}$, the interpretation of $t$ in ${\cal C}^{({\cal S}, \Sigma_{A})}_{{\mathbb S}_{A, R}}$ is equal to $[F_{t}]$; from this it easily follows, by definition of the theory ${\mathbb S}_{A, R}$, that all the axioms of ${\mathbb T}_{A, R}$ are satisfied in ${\cal C}_{{\mathbb S}_{A, R}}$. 

To prove that, with the structure just defined, ${\cal C}^{({\cal S}, \Sigma_{A})}_{{\mathbb S}_{A, R}}$ is actually an initial object of ${\mathbb T}_{A, R}\textrm{-mod}(\Set)$, we appeal to the universal property of the syntactic category ${\cal C}^{({\cal S}, \Sigma_{A})}_{{\mathbb S}_{A, R}}$ (cf. section \ref{gensyn}): for any $({\cal S}, \Sigma_{A})$-category $\cal D$, the category ${{\mathbb S}_{A, R}}\textrm{-mod}_{\cal S}({\cal D})$ of $\cal S$-models of ${\mathbb S}_{A, R}$ in $\cal D$ is equivalent to the category $({\cal S}, \Sigma_{A})\textrm{-}\textbf{Fun}({\cal C}_{{\mathbb S}_{A, R}}, {\cal D})$ of $({\cal S}, \Sigma_{A})$-preserving functors from ${\cal C}^{({\cal S}, \Sigma_{A})}_{{\mathbb S}_{A, R}}$ to ${\cal D}$, via the equivalence sending a functor $F:{\cal C}_{{\mathbb S}_{A, R}} \to {\cal D}$ to the model of ${\mathbb S}_{A, R}$ in $\cal D$ in which a generator $a\in A$ is interpreted as the image $F([R_{a}])$ of the ${\mathbb S}_{A, R}$-provable $R_{a}$ of the atomic formula $R_{a}$. 

Let us apply this property to poset categories $\cal D$. Note that in this case, all of our categories being posets, the equivalence 
\[
{{\mathbb S}_{A, R}}\textrm{-mod}_{\cal S}({\cal D})\simeq ({\cal S}, \Sigma_{A})\textrm{-}\textbf{Fun}({\cal C}^{({\cal S}, \Sigma_{A})}_{{\mathbb S}_{A, R}}, {\cal D})
\]
is in fact an isomorphism of categories (equivalently, of posets). 

Given a ${\mathbb T}_{A, R}$-model $N$ in $\Set$, we can regard $N$ as a poset $({\cal S}, \Sigma_{A})$-category ${\cal D}_{N}$ (in which the connectives in $\cal S$ are interpreted precisely as the operations on $N$ which interpret the function symbols over $\Sigma_{\mathbb A}$) containing a ${\mathbb S}_{A, R}$-model, obtained by interpreting any $0$-ary relation symbol $R_{a}$ (for $a\in A$) as the element of $N$ which interprets the constant $C_{a}$ over $\Sigma_{\mathbb A}^{A}$ in $N$ regarded as a ${\mathbb T}_{A, R}$-model. Now, the ${\mathbb T}_{A, R}$-model homomorphisms ${\cal C}_{{\mathbb S}_{A, R}} \to N$ can be identified exactly with the $({\cal S}, \Sigma_{A})$-preserving functors ${\cal C}_{{\mathbb S}_{A, R}} \to {\cal D}_{N}$ which send any $[R_{a}]$ (for $a\in A$) to the interpretation of $C_{a}$ in $N$ (regarded as a ${\mathbb T}_{A, R}$-model); but these correspond, via the isomorphism $({\cal S}, \Sigma_{A})\textrm{-}\textbf{Fun}({\cal C}^{({\cal S}, \Sigma_{A})}_{{\mathbb S}_{A, R}}, {\cal D})\cong {\mathbb T}\textrm{-mod}_{\cal S}({\cal D})$, precisely to the ${\mathbb S}_{A, R}$-models in ${\cal D}_{N}$ in which the interpretation of the formula $R_{a}$ (for any $a\in A$) is equal to the interpretation of $C_{a}$ in $N$, and there is obviously just one such model, namely $N$. This proves our claim.

Summarizing, we have obtained the following result.

\begin{theorem}
Let $\mathbb A$ be an small ordered algebraic theory and $(A, R)$ a set of generators and relations for it. With the above notation, the poset ${\cal C}_{{\mathbb S}_{A, R}}$ is `the' $\mathbb A$-model presented by the set of generators $A$ and relations $R$ via the function $\xi:A\to {\cal C}_{{\mathbb S}_{A, R}}$ sending any $a\in A$ to the ${{\mathbb S}_{A, R}}$-provable equivalence class $[R_{a}]$ of the formula $R_{a}$.

In particular, models of small ordered algebraic theories presented by arbitrary generators and relations always exist.
\end{theorem}\qed

In fact, in light of the identification between the closed terms over $\Sigma_{\mathbb A}^{A}$ and the $\cal S$-formulae over $\Sigma_{A}$ observed above, the construction given by this theorem is isomorphic to that of Theorem \ref{oat}. 
 
From the proof given above it is clear that the poset ${\cal C}_{{\mathbb S}_{A, R}}$ does not only satisfy the universal property of the $\mathbb A$-model presented by $(A, R)$ with respect to the \emph{$\mathbb A$-models}, but also with respect to every structure $N$ over the signature of $\mathbb A$ equipped with interpretations of the constant symbols $C_{a}$ (for $a\in A$) in such a way that the substructure of $N$ generated by these interpretations is a model of $\mathbb A$.   
  
We have used the syntactic categories of propositional $\cal S$-theories to build models presented by generators and relations of (small) ordered algebraic theories; the question thus naturally arises whether every such syntactic category is of this form, that is if it can be regarded, in a natural way, as a model of an ordered algebraic theory presented by generators and relations. The answer to this question is positive. Given a propositional $\cal S$-theory $\mathbb T$, we take one function symbol of arity $k$ for each connective in $\cal S$ of arity $k$; this, together with the binary relation symbol $\leq$, defines a signature $\Sigma$. Take $\mathbb H$ to be the theory over $\Sigma$ having as axioms the infinitary Horn sequents over $\Sigma$ whose associated inference rule schemes (via the method above) are provably valid in $\mathbb T$. Define $A$ to be equal to the set of $0$-ary relation symbols of the signature of $\cal S$, and the set $R$ of relations involving elements in $A$ as the set of Horn sequents over $\Sigma_{\mathbb A}^{A}$ whose associated inference rule scheme (via the method above) is provably valid in $\mathbb T$. It is then clear that the $\cal S$-syntactic category of $\mathbb T$ is precisely `the' $\mathbb H$-model presented by the set of generators $A$ subject to relations $R$. 
 
\subsection{Classical examples}\label{clex}

\begin{itemize}

\item \emph{Meet-semilattices}

The concept of meet-semilattice can be formalized as a small ordered algebraic theory, as follows. Take $\Sigma_{\mathbb A}$ to be the signature consisting of a binary relation symbol $\leq$, one constant $1$ and one binary function symbol $\wedge$. Consider the following Horn axioms over $\Sigma$:
\[
(\top \vdash_{x} (x \leq 1)),
\]
\[
(x\leq (y\wedge z) \dashv\vdash_{x,y,z} (x\leq y) \wedge (x\leq z)).
\]
Clearly, the Horn theory $\mathbb A$ over $\Sigma_{\mathbb A}$ obtained by adding to these axioms the sequents which express the fact that $\leq$ is a partial order has the property that its models in $\Set$ and homomorphisms between them can be identified with the meet-semilattices and meet-semilattice homomorphisms between them. 

It is immediate to see that, for any choice of generators and relations $(A, R)$ for $\mathbb A$ such that the relations in $R$ are all of the form $\top \vdash t_{1} T t_{2}$ where $t_{1}$ and $t_{2}$ are closed terms over $\Sigma_{\mathbb A}^{A}$ and $T$ is either the relation $\leq$ or the equality relation $=$, the $\cal S$-theory $\mathbb L$ associated to the theory $\mathbb A$ via the method of the last section is precisely the Horn propositional theory $\mathbb M$ over the signature $\Gamma$ consisting of one $0$-ary relation symbol for each of the elements of $A$ (note that we do not have to add the constant $\top$ since it is already present in Horn logic), with axioms (in addition to those of Horn logic):
\[
(F_{t_{1}} \vdash_{[]} F_{t_{2}})
\]
(respectively, 
\[
(F_{t_{1}} \dashv \vdash_{[]} F_{t_{2}}))
\]
for each relation in $R$ of the form $\top \vdash t_{1} \leq t_{2}$ (resp. of the form $\top \vdash t_{1}=t_{2}$), where $t_{1}$ and $t_{2}$ are closed terms over $\Sigma_{\mathbb A}^{A}$ and $F_{t_{1}}$ and $F_{t_{2}}$ are the atomic formulae over $\Gamma$ corresponding to them via the method of section \ref{genrel}. In particular, the syntactic category of the $\cal S$-theory $\mathbb L$ is equivalent to the cartesian syntactic category of the cartesian theory $\mathbb M$ (as defined in section D1 of \cite{El}). 

\item \emph{Distributive lattices}

The concept of distributive lattice can also be formalized as a small ordered algebraic theory, as follows. Let $\Sigma_{\mathbb A}$ be the signature consisting of a binary relation symbol $\leq$, two constants $0, 1$ and two binary function symbols $\wedge, \vee$. Consider the following Horn axioms over $\Sigma_{\mathbb A}$:
\[
(\top \vdash_{x} (x \leq 1)),
\]
\[
(\top \vdash_{x} (0 \leq x)),
\]
\[
(x\leq (y\wedge z) \dashv\vdash_{x,y,z} (x\leq y) \wedge (x\leq z)),
\]
\[
((x\vee y) \leq z \dashv\vdash_{x,y,z} (x\leq z) \wedge (y\leq z)),
\]
\[
(\top \vdash_{x,y,z} (x\wedge (y\vee z) = (x\wedge y) \wedge (x\wedge z))).
\]
These axioms, in addition to the sequents which express the fact that $\leq$ is a partial order, define a small ordered algebraic theory $\mathbb A$ whose models in $\Set$ and homomorphisms between them can be identified with the distributive lattices and homomorphisms between them. 

As with meet-semilattices, it is immediate to see that, for any choice of generators and relations $(A, R)$ for $\mathbb A$ such that the relations in $R$ are all of the form $\top \vdash t_{1} T t_{2}$ where $t_{1}$ and $t_{2}$ are closed terms over $\Sigma_{\mathbb A}^{A}$ and $T$ is either the relation $\leq$ or the equality relation $=$, the $\cal S$-theory $\mathbb L$ associated to the theory $\mathbb A$ via the method of the last section can be identified with the coherent propositional theory $\mathbb M$ over the signature $\Gamma$ consisting of one $0$-ary relation symbol for each of the elements of $A$ (note that we do not have to add the constants $\top$ and $\bot$ since they are already present in coherent logic), having as axioms (in addition to those of coherent logic) are the following the sequents of the form
\[
(F_{t_{1}} \vdash_{[]} F_{t_{2}})
\]
(respectively, 
\[
(F_{t_{1}} \dashv \vdash_{[]} F_{t_{2}}))
\]
for each relation in $R$ of the form $\top \vdash t_{1} \leq t_{2}$ (resp. of the form $\top \vdash t_{1}=t_{2}$), where $t_{1}$ and $t_{2}$ are closed terms over $\Sigma_{\mathbb A}^{A}$ and $F_{t_{1}}$ and $F_{t_{2}}$ are the atomic formulae over $\Gamma$ corresponding to them via the method of section \ref{genrel}. In particular, the syntactic category of the $\cal S$-theory $\mathbb L$ is equivalent to the coherent syntactic category of the coherent theory $\mathbb M$. 

\item \emph{Frames}

The notion of frame is the infinitary analogue of that of distributive lattice. Since our definition of ordered algebraic theory allows infinitary function symbols and infinitary conjunctions of atomic formulae, the arguments given above for distributive lattices straightforwardly extend to frames. Specifically, the small ordered algebraic theory formalizing the concept of frame is obtained by taking a signature $\Sigma_{\mathbb A}$ consisting of a binary relation symbol $\leq$, two constants $0, 1$, one binary function symbol $\wedge$ and an infinitary function symbol $D$. Consider the following Horn axioms over $\Sigma$:
\[
(\top \vdash_{x} (x \leq 1)),
\]
\[
(\top \vdash_{x} (0 \leq x)),
\]
\[
(x\leq (y\wedge z) \dashv\vdash_{x,y,z} (x\leq y) \wedge (x\leq z)),
\]
\[
(D(x_{i} \textrm{ | } i\in I) \leq z \dashv\vdash_{x_{i}, z} \mathbin{\mathop{\textrm{\huge $\wedge$}}\limits_{i\in I}}(x_{i}\leq z)),
\]
\[
(\top \vdash_{x_{i},y} (D(x_{i} \textrm{ | } i\in I) \wedge y = D(x_{i}\wedge y \textrm{ | } i\in I)).
\]

These axioms, in addition to the sequents which express the fact that $\leq$ is a partial order, define a small ordered algebraic theory $\mathbb A$ whose models in $\Set$ and homomorphisms between them can be identified with the frames and frame homomorphisms between them. 

As in the case of distributive lattices, for any choice of generators and relations $(A, R)$ for $\mathbb A$ such that the relations in $R$ are all of the form $\top \vdash t_{1} T t_{2}$ where $t_{1}$ and $t_{2}$ are closed terms over $\Sigma_{\mathbb A}^{A}$ and $T$ is either the relation $\leq$ or the equality relation $=$, the $\cal S$-theory $\mathbb L$ associated to the theory $\mathbb T$ via the method of section \ref{genrel} can be identified with the geometric propositional theory $\mathbb M$ over the signature $\Gamma$ consisting of one $0$-ary relation symbol for each of the elements of $A$, and whose axioms (in addition to those of geometric logic) are the following:
\[
(F_{t_{1}} \vdash_{[]} F_{t_{2}})
\]
(respectively, 
\[
(F_{t_{1}} \dashv \vdash_{[]} F_{t_{2}})
\]
for each relation in $R$ of the form $\top \vdash t_{1} \leq t_{2}$ (resp. of the form $\top \vdash t_{1}=t_{2}$), where $t_{1}$ and $t_{2}$ are closed terms over $\Sigma_{\mathbb A}^{A}$ and $F_{t_{1}}$ and $F_{t_{2}}$ are the atomic formulae over $\Gamma$ corresponding to them via the method of section \ref{genrel}.

\end{itemize} 

\subsection{Structures presented by generators and relations}\label{concreteness}

We have seen that models of small ordered algebraic theories presented by generators and relations can always be constructed as preordered syntactic categories of generalized propositional theories. 

A problem which frequently arises in practice is that of finding concrete descriptions for poset structures presented by generators and relations. The techniques elaborated in this paper, combined with the philosophy `toposes as bridges' of \cite{OC10}, provide flexible and effective tools for addressing this kind of problems. The key idea is to equip such a poset structure $\cal P$ with a Grothendieck topology $J$ such that $\cal P$ can be recovered (up to isomorphism) from the topos $\Sh({\cal P}, J)$ by means of a topos-theoretic invariant $U$ which has a `natural behaviour' with respect to sites (as we did for example in section \ref{charinv}); indeed, under these hypotheses, any other site of definition $({\cal C}, K)$ of the topos $\Sh({\cal P}, J)$ leads to a different representation for $\cal P$ as the poset of subterminals in the topos $\Sh({\cal C}, K)$ which satisfy the invariant $U$. For example, given a commutative ring with unit $(A, +, \cdot, 0_{A}, 1_{A})$, the distributive lattice generated by symbols $D(a)$, $a \in A$, subject to the relations $D(1_{A})=1_{L(A)}$, $D(a\cdot b) = D(a) \wedge D(b)$, $D(0_{A})=0_{L(A)}$, and $D(a+ b) \leq D(a) \vee D(b)$ can be characterized (up to isomorphism) as the lattice of compact elements of the frame of open sets of the prime spectrum of $A$ endowed with the Zariski topology (cf. section \ref{zariski} below).

The notion of subterminal topology can be profitably used to give `topological descriptions' of poset structures $\cal P$ (that is, descriptions in terms of frames of open sets of topological spaces) whenever the corresponding toposes $\Sh({\cal P}, J)$ have enough points. Indeed, if $U$ is a topos-theoretic invariant of families of subterminals in a topos which is $({\cal P}, J)$-adequate (in the sense of Definition \ref{ade}), so that $\cal P$ can be recovered from $\Sh({\cal P}, J)$ as the poset of $U$-compact subterminals of $\Sh({\cal P}, J)$ then, provided that $\Sh({\cal P}, J)$ has enough points, $\cal C$ is isomorphic (as a poset) to the poset of $U$-compact open sets of any topological space obtained by endowing a separating set of points of $\Sh({\cal P}, J)$ with the subterminal topology (cf. section \ref{subterminaltop}).   

For example, the free frame $F(A)$ on a set $A$ can be described (up to isomorphism) as the frame of open sets of the topological space obtained by endowing the powerset of $A$ with the elemental topology (cf. section \ref{subs}). Similarly, the free meet-semilattice ${\cal M}(A)$ on a set $A$ can be regarded as the cartesian syntactic category of the empty propositional theory ${\mathbb T}_{A}$ over the signature having exactly one $0$-ary relation symbol for each element $a$ of $A$, and hence it can be described, by Corollary \ref{cor}(iii), as the poset of supercompact open sets of the topological space obtained by equipping the powerset of $A$ with the elemental topology (cf. section \ref{subs}). 

In connection with the method described above, it is worth to remark that the theory of classifying toposes can be profitably exploited to obtain different representations for a given topos, regarded as a classifying topos of a geometric theory. For example, by regarding the theory ${\mathbb T}_{A}$ defined above as a cartesian theory, we immediately obtain a description of its classifying topos as the presheaf topos $[\textrm{f.p.} {\mathbb T}_{A}\textrm{-mod}(\Set), \Set]$, where $\textrm{f.p.} {\mathbb T}_{A}\textrm{-mod}(\Set)$ is the category of (representatives of isomorphism classes of) finitely presentable models of ${\mathbb T}_{A}$ in $\Set$, from which it follows that the free meet-semilattice ${\cal M}(A)$ on a set $A$ is isomorphic to the opposite of the poset $\textrm{f.p.} {\mathbb T}_{A}\textrm{-mod}(\Set)$; this poset can be clearly identified with the poset ${\mathscr{P}}_{fin}(A)$ of finite subsets of $A$ (with the subset-inclusion ordering), and from this description we can recover the well-known characterization of ${\cal M}(A)$ as ${\mathscr{P}}_{fin}(A)^{\textrm{op}}$. Similarly, one can recover the classical description of the free frame $F(A)$ on a set $A$ by regarding $F(A)$ as the geometric syntactic category of the theory ${\mathbb T}_{A}$; indeed, the latter can be identified with the frame of subterminals of the classifying topos $\Set[{\mathbb T}_{A}]\simeq [{\mathscr{P}}_{fin}(A), \Set]$ of ${\mathbb T}_{A}$, from which it follows that the free frame on a set $A$ is (isomorphic to) the frame of upper sets on ${\mathscr{P}}_{fin}(A)$. 

\subsection{The free frame on a complete join-semilattice}\label{freecomjoin}

A particularly significant application of the general method introduced in the last section for explicitly describing structures presented by generators and relations is the solution to the problem of describing the free frame $L(A)$ on a complete join-semilattice $A$. 

If we denote by $\mathbin{\mathop{\textrm{\huge $\vee$}}\limits_{i\in I}}a_{i}$ the join of a family of elements $\{a_{i} \textrm{ | } i\in I\}$ in a complete join-semilattice $A$ then, by the results of sections \ref{genrel} and \ref{clex}, $L(A)$ can be identified with the geometric syntactic category ${\cal C}_{{\mathbb L}_{A}}$ of the geometric theory ${\mathbb L}_{A}$ defined as follows: the signature $\Sigma_{A}$ of ${\mathbb L}_{A}$ consists of one $0$-ary relation symbol $F_{a}$ for each element $a\in A$, and the axioms of ${\mathbb L}_{A}$ are, besides those of geometric logic, all the sequents over $\Sigma_{A}$ of the form
\[
\mathbin{\mathop{\textrm{\huge $\vee$}}\limits_{i\in I}}F_{a_{i}} \dashv \vdash F_{a} 
\]    
for any family of elements $\{a_{i} \textrm{ | } i\in I\}$ in $A$ such that $a=\mathbin{\mathop{\textrm{\huge $\vee$}}\limits_{i\in I}}a_{i}$ in $A$.

Now, the classifying topos of ${\mathbb L}_{A}$ can be realized as a subtopos 
\[
\Sh(\textrm{f.p.} {\mathbb T}_{A}\textrm{-mod}(\Set)^{\textrm{op}}, J)
\]
of the classifying topos of the (cartesian) empty theory ${\mathbb T}_{A}$ over $\Sigma_{A}$, which, as we have seen above, is the presheaf topos  
\[
[{\cal C}_{{\mathbb T}_{A}}^{\textrm{op}}, \Set]\simeq [\textrm{f.p.} {\mathbb T}_{A}\textrm{-mod}(\Set), \Set]\simeq [{\mathscr{P}}_{fin}(A), \Set],
\]
where ${\cal C}_{{\mathbb T}_{A}}$ is the cartesian syntactic category of ${\mathbb T}_{A}$.        

The Grothendieck topology $J$ on $\textrm{f.p.} {\mathbb T}_{A}\textrm{-mod}(\Set)^{\textrm{op}}\simeq {\cal C}_{{\mathbb T}_{A}}$ is defined as follows: denoted by $[F_{a}]$ the object in ${\cal C}_{{\mathbb T}_{A}}$ corresponding to the atomic formula $F_{a}$ (for an element $a\in A$), $J$ is generated by the families of sieves which contain sinks of the form     
\[
\{[F_{a_{i}} \wedge F_{a}] \leq [F_{a}] \textrm{ | } i\in I\} 
\]    
for some family of elements $\{a_{i} \textrm{ | } i\in I\}$ in $A$ such that $a=\mathbin{\mathop{\textrm{\huge $\vee$}}\limits_{i\in I}}a_{i}$ in $A$.

Under the equivalence ${\cal C}_{{\mathbb T}_{A}} \simeq \textrm{f.p.} {\mathbb T}_{A}\textrm{-mod}(\Set)^{\textrm{op}}$, the closure under pullbacks of this family of sieves corresponds to the coverage $K$ on the category $\textrm{f.p.} {\mathbb T}_{A}\textrm{-mod}(\Set)^{\textrm{op}}\simeq {\mathscr{P}}_{fin}(A)^{\textrm{op}}$ described as follows: the $K$-covering cosieves on ${\mathscr{P}}_{fin}(A)$ are those generated by families of the form 
\[
\{U \cup \{a_{i}\} \supseteq U \textrm{ | } i\in I\} 
\]  
for any finite subset $U$ of $A$ and any family $\{a_{i} \textrm{ | } i\in I\}$ of elements of $A$ such that $\mathbin{\mathop{\textrm{\huge $\vee$}}\limits_{i\in I}}a_{i}\in U$.

Now, the geometric syntactic category ${\cal C}_{{\mathbb L}_{A}}$ of the theory ${\mathbb L}_{A}$ is isomorphic, as a poset, to the frame of subterminals of the classifying topos $\Set[{\mathbb L}_{A}]$ of ${\mathbb L}_{A}$ and hence, from the representation 
\[
\Set[{\mathbb L}_{A}] \simeq \Sh(\textrm{f.p.} {\mathbb T}_{A}\textrm{-mod}(\Set)^{\textrm{op}}, J)
\]
it follows that ${\cal C}_{{\mathbb L}_{A}}$ is isomorphic to the frame $Id_{J}(\textrm{f.p.} {\mathbb T}_{A}\textrm{-mod}(\Set)^{\textrm{op}})$ of $J$-ideals on $\textrm{f.p.} {\mathbb T}_{A}\textrm{-mod}(\Set)^{\textrm{op}}$. Note that this argument represents an application of the philosophy `toposes as bridges' of \cite{OC10}. 
 
Since $J$ is generated by the coverage $K$ then the $J$-ideals on the category $\textrm{f.p.} {\mathbb T}_{A}\textrm{-mod}(\Set)^{\textrm{op}}$ coincide exactly with the $K$-ideals on $\textrm{f.p.} {\mathbb T}_{A}\textrm{-mod}(\Set)^{\textrm{op}}$. This can be proved directly (by showing that, for any $K$-ideal $I$ the collection of all the sieves $S$ on $\cal C$ such that the implication `$(dom(f)\in I$ for all $f\in S)$ entails $cod(f)\in I$' holds defines a Grothendieck topology which contains $K$) or, more conceptually, can be deduced as a consequence of the general well-known fact (cf. Proposition C2.1.9 \cite{El}) that if $J$ is the Grothendieck topology generated by a coverage $K$ on a small category $\cal C$ (i.e., the smallest Grothendieck topology $J$ on $\cal C$ such that every $K$-covering family generates a $J$-covering sieve) then the toposes $\Sh({\cal C}, J)$ and $\Sh({\cal C}, K)$ are equal, by observing that for any site $({\cal E}, W)$, the $W$-ideals on ${\cal E}$ can be (canonically) identified with the subterminals of the topos $\Sh({\cal E}, W)$.

The $K$-ideals on $\textrm{f.p.} {\mathbb T}_{A}\textrm{-mod}(\Set)^{\textrm{op}}\simeq {\mathscr{P}}_{fin}(A)^{\textrm{op}}$ can be identified precisely with the upward closed subsets ${\cal I}$ of $\mathscr{P}_{fin}(A)$ such that for every $U\in \mathscr{P}_{fin}(A)$ and every family of elements $\{a_{i} \textrm{ | } i\in I\}$ of $A$ such that $\mathbin{\mathop{\textrm{\huge $\vee$}}\limits_{i\in I}}a_{i}\in U$, if $U\cup \{a_{i}\}\in {\cal I}$ for all $i\in I$ then $U\in {\cal I}$. We thus conclude that $L(A)$ is given by the collection of all such subsets, endowed with the natural subset-inclusion ordering. The universal complete join-semilattice homomorphism $\eta_{A}:A\to L(A)$ is the map sending an element $a\in A$ to the $J$-closure of the set $I_{a}$ of all the finite subsets of $A$ which contain the element $a$. We can easily see that $I_{a}$ is $J$-closed (equivalently, $K$-closed); indeed, given $U\in \mathscr{P}_{fin}(A)$ and any family of elements $\{a_{i} \textrm{ | } i\in I\}$ of $A$ such that $\mathbin{\mathop{\textrm{\huge $\vee$}}\limits_{i\in I}}a_{i}\in U$, $a\in U\cup \{a_{i}\}$ for all $i\in I$ implies $a\in U$. The universal map $\eta_{A}:A\to L(A)$ thus sends any element $a\in A$ to the set $I_{a}$, and is therefore injective; this solves an open problem of Resende and Vickers (cf. the remarks preceding Proposition 3.7 \cite{RV}). Summarizing, we have the following result.   

\begin{theorem}\label{solution}
The free frame on a complete join semilattice $A$ is (up to isomorphism) the set $L(A)$ of upward closed (with respect to the subset-inclusion ordering on $\mathscr{P}_{fin}(A)$) subsets ${\cal I}$ of $\mathscr{P}_{fin}(A)$ with the property that for every $U\in \mathscr{P}_{fin}(A)$ and every family of elements $\{a_{i} \textrm{ | } i\in I\}$ of $A$ such that $\mathbin{\mathop{\textrm{\huge $\vee$}}\limits_{i\in I}}a_{i}\in U$, if $U\cup \{a_{i}\}\in {\cal I}$ for all $i\in I$ then $U\in {\cal I}$, endowed with the subset-inclusion ordering, while the universal complete join-semilattice homomorphism $\eta_{A}:A\to L(A)$ is the map sending any element $a\in A$ to the set $I_{a}\in L(A)$ of all the finite subsets of $A$ which contain the element $a$.   
\end{theorem}\qed

We can describe the frame operations on $L(A)$ explicitly as follows: given any $I_{1}, I_{2}\in L(A)$, the meet $I_{1}\wedge I_{2}$ in $L(A)$ is simply the set-theoretic intersection $I_{1}\cap I_{2}$, while the join $\mathbin{\mathop{\textrm{\huge $\vee$}}\limits_{j\in J}}I_{j}$ in $L(A)$ of a family $\{I_{j} \textrm{ | } j\in J\}$ of elements $I_{j}\in L(A)$ is equal to the $J$-closure of the set-theoretic union $\mathbin{\mathop{\textrm{\huge $\cup$}}\limits_{j\in J}}I_{j}$, i.e. to the collection of all the finite subsets $U$ of $A$ with the property that there exists a family $\{a_{h} \textrm{ | } h\in H\}$ of elements of $A$ such that $\mathbin{\mathop{\textrm{\huge $\vee$}}\limits_{h\in H}}a_{h}\in U$ and for every $h\in H$, $U\cup \{a_{h}\}\in I_{j}$ for some $j \in J$.

Given a frame $L$, and a homomorphism $f:A\to L$ of complete join-semilattices, the unique frame homomorphism $g:L(A)\to L$ such that $g\circ \eta_{A}=f$ is defined by the formula:
\[
g(I)=\mathbin{\mathop{\textrm{\huge $\vee$}}\limits_{U\in I}} (\mathbin{\mathop{\textrm{\huge $\wedge$}}\limits_{a\in U}} f(a)),
\]
for any $I\in L(A)$.
 
Now that we have obtained a solution to our problem, it is natural to look for an elementary proof of the result which does not involve topos-theoretic concepts. Such a proof clearly exists, but we do not bother to write it down because, although being elementary, it involves a good amount of tedious and intricate verifications; on the other hand, we are confident that the interested reader will have no trouble in carrying out the necessary calculations by himself or herself. Rather, we emphasize that this is just an example of what can be achieved by using our methods; many similar problems can be solved by applying the same techniques, and the reader is invited to try them out in the context of his or her questions of interest.

\begin{remark}
Let $\textbf{CJSLat}$ denote the category having as objects the complete join-semilattices and as arrows the maps between them which preserve arbitrary joins, and let $\textbf{Frm}$ be the category of frames. The map from complete join-semilattices to frames which sends a complete join-semilattice $A$ to the complete frame $L(A)$ on it can be made into a functor $L:\textbf{CJSLat}\to \textbf{Frm}$ which is left adjoint to the inclusion functor $\textbf{Frm} \hookrightarrow \textbf{CJSLat}$. It is easy to see that the functor $L$ sends a complete join-semilattice $A$ to the frame $L(A)$ and an arrow $f:A\to B$ in $\textbf{CJSLat}$ to the unique frame homomorphism $L(f):L(A)\to L(B)$ such that $\xi\circ \eta_{A}=\eta_{B}\circ f$; concretely, $L(f)$ sends a subset $\cal I$ in $\mathscr{P}_{fin}(A)$ to the smallest upward closed subsets $\cal J$ of $\mathscr{P}_{fin}(B)$ (with respect to the subset-inclusion ordering on $\mathscr{P}_{fin}(B)$) which contains $f({\cal I})$ and has the property that for every $U\in \mathscr{P}_{fin}(B)$ and every family of elements $\{b_{i} \textrm{ | } i\in I\}$ of $B$ such that $\mathbin{\mathop{\textrm{\huge $\vee$}}\limits_{i\in I}}b_{i}\in U$, if $U\cup \{b_{i}\}\in {\cal J}$ for all $i\in I$ then $U\in {\cal J}$.      
\end{remark}

We conclude this section by presenting a couple of other examples. First, we note that the finitary analogue of the construction (that is, the version of the construction obtained by requiring all the sets $I$ to be finite) provides a way for building the free frame $F(A)$ on a join-semilattice $A$. Note that the free distributive lattice on a join-semilattice $A$ can be characterized as the lattice of compact elements of $F(A)$ (cf. Corollary \ref{cor}(i)). 

The free frame $F(A)$ on a join-semilattice $(A, \vee)$ can also be described in topological terms by using the method of section \ref{concreteness}. Specifically, we get that $F(A)$ is isomorphic to the frame of open sets of the space obtained by endowing the collection of all the subsets $U$ of $A$ such that $0\notin U$ and for any pair of elements $a,b\in A$, $a\vee b \in U$ if and only if $a\in U$ or $b\in U$ with the elemental topology. Note that this topological space is homeomorphic to the space obtained by endowing the set of subsets $V$ of $A$ such that $0\in V$ and for any pair of elements $a,b\in A$, $a\vee b \in V$ if and only if $a\notin V$ and $b\notin V$ with the topology of which the sets of the form $D_{a}:=\{V \textrm{ | } a\notin A \}$ for $a\in A$ form a basis (notice the similarity between this topology and the Zariski topology, cf. section \ref{zariski} below).

\subsection{The Zariski topology}\label{zariski}

In section V.III of \cite{stone}, two different ways of building the Zariski spectrum of a commutative ring with unit are presented. In this section, we provide a topos-theoretic interpretation of these two constructions as a Morita-equivalence, in the form of two different sites of definition for the same topos. This is particularly relevant in relation to the philosophy `toposes as bridges' of \cite{OC10}; in fact, as we shall see below, one can effectively transfer information between the different sites of definition of the topos through topos-theoretic invariants, according to the indications given in \cite{OC10}. 

Let us begin by discussing the constructions given in \cite{stone}. 

Let $A$ be an arbitrary commutative ring with unit. Let $S(A)$ be the quotient of the multiplicative monoid of $A$ by the smallest monoid congruence $\simeq$ which identifies $a$ with $a^{2}$ for each $a \in A$; we denote by $\pi_{S}:A\to S(A)$ the canonical quotient map. Note that $S(A)$ is a commutative monoid in which every element is idempotent. 
For any monoid $M$ in which every element is idempotent, we can define a partial order relation $\leq_{M}$ on $M$ as follows: for any $a,b\in M$, $a \leq_{M} b$ if and only if $a\cdot b=a$; in fact, for any $a,b\in M$, $a\cdot b$ is the greatest lower bound of $a$ and $b$ with respect to $\leq_{M}$, and hence $(M, \leq_{M})$ is a meet-semilattice.
So $S(A)$ is a meet-semilattice and hence, regarded as a category, it is a cartesian category. We can define a coverage $C$ on $S(A)$ as follows: $\emptyset \in C(\pi_{S}(0_{A}))$, and the $C$-covering sieves on an object $x\in S(A)$ are those which contain finite families of the form $\pi_{S}(a_{i})\to \pi_{S}(a)$ for $i=1, \ldots, n$, where $\pi_{S}(a_{i}) \leq_{S_{A}} \pi_{S}(a)$ and $\pi_{S}(a_{1}+ \cdots + a_{n})=\pi_{S}(a)=x$. 

The Zariski spectrum ${\cal Z}_{A}$ of the ring $A$ is defined in \cite{stone} as the space of points of the locale $Id_{C}({S(A)})$ of $C$-ideals in $S(A)$. 

The second approach to the construction of ${\cal Z}_{A}$ defines ${\cal Z}_{A}$ as the prime spectrum (in the sense of section II3.4 \cite{stone}) of the distributive lattice $L(A)$ generated by symbols $D(a)$, $a \in A$, subject to the relations $D(1_{A})=1_{L(A)}$, $D(a\cdot b) = D(a) \wedge D(b)$, $D(0_{A})=0_{L(A)}$, and $D(a+ b) \leq D(a) \vee D(b)$.

By the results of section \ref{clex}, $L(A)$ can be seen as the coherent syntactic category of the coherent propositional theory over the signature having a propositional symbol $P(a)$ for each element $a\in A$, whose axioms are the following:
\[
(\top \vdash P_{1_{A}});
\]
\[
(P_{0_{A}} \vdash \bot);
\]
\[
(P_{a\cdot b} \dashv\vdash P_{a} \wedge P_{b})
\] 
for any $a, b$ in $A$; 
\[
(P_{a + b} \vdash P_{a} \vee P_{b})
\]
for any $a, b \in A$. 
 
The prime spectrum of $L(A)$ can thus be identified with the space of points of the locale $Id_{J}(L(A))$, where $J$ is the coherent topology on $L(A)$. 

To prove that this second definition of the Zariski spectrum of $A$ coincides with the first one, one observes that the universal map $\pi_{L}: A\to L(A)$ factors through the quotient map $\pi_{S}: A\to S(A)$, and the resulting factorization is universal among the meet-semilattice homomorphisms from $S(A)$ to distributive lattices which carry the covering families in $C$ to joinly covering families. Indeed, this implies, by Theorem \ref{freeframes}, that the frames $Id_{C}({S(A)})$ and $Id_{J}({L(A)})$ are isomorphic. Note that, by the results of sections \ref{genrel} and \ref{clex}, these two isomorphic frames can also be identified with the frame $F$ generated by symbols $F(a)$, $a \in A$, subject to the relations $F(1_{A})=1_{F}$, $F(a \cdot b) = F(a) \wedge F(b)$, $F(0_{A})=0_{F}$, and $F(a+ b) \leq F(a) \vee F(b)$.  

Now, since the locales $Id_{C}({S(A)})$ and $Id_{J}({L(A)})$ are isomorphic, Theorem \ref{fund} yields an equivalence of toposes 
\[
\Sh(S(A), C)\simeq \Sh(L(A), J)
\]
The fact that the two definitions of the Zariski spectrum ${\cal Z}_{A}$ given above coincide then follows at once from the fact that the subterminal topology is a topos-theoretic invariant and that on a localic topos $\Sh(L)$ equipped with the collection of all its points it yields the space of points of $L$ (cf. Example \ref{exa}(e) above). 

To give an explicit description of ${\cal Z}_{A}$ as a topological space, we can use either of the two topos-theoretic representations. For example, using the representation $\Sh(S(A), C)$, we see that the points of ${\cal Z}_{A}\cong X_{{\tau}^{\Sh(S(A), C)}}$ are precisely the points of the topos $\Sh(S(A), C)$. These correspond bijectively, by definition of $S(A)$, to the functions $f:A\to \{0,1\}$ such that they factor through $\pi_{S}:A\to S(A)$ and this factorization is a meet-semilattice homomorphism $S(A)\to \{0,1\}$ which sends every $C$-covering sieve to a covering family in $\{0,1\}$. Identifying such functions $f$ with the subsets $f^{-1}(1)$ of $A$, we see that the these functions correspond bijectively to the subsets $S\subseteq A$ with the property that $1\in S$, $0_{A}\notin S$, $a\cdot b\in S$ if and only if $a\in S$ and $b\in S$, and $a+b\in S$ implies that $a\in S$ or $b\in S$. These subsets are called in \cite{stone} the \emph{prime filters} on $A$, and we denote their collection by ${\cal F}_{A}$. 
On the other hand, identifying the functions $f$ above with the (complementary) subsets $f^{-1}(0)$ of $A$, we get a bijection between the set of such functions $f$ and the set $Spec(A)$ of prime ideals of $A$.          

By Proposition \ref{mslattice}, the subterminal topology on the set ${\cal F}_{A}\cong X_{{\tau}^{\Sh(S(A), C)}}$ has as a sub-basis of open sets the collection of sets of the form 
\[
{\cal F}_{a}=\{S\in {\cal F}_{A} \textrm{ | } a\in S\},
\] 
where $a$ varies among the elements of $A$. 

Now, under the bijection ${\cal F}_{A}\cong Spec(A)$ sending each subset to its complement in $A$, this topology on ${\cal F}_{A}$ corresponds precisely to the \emph{Zariski topology} on $Spec(A)$, that is to the topology whose closed sets are those of the form
\[
{\cal P}_{I}=\{P\in Spec(A) \textrm{ | } P\supseteq I\},
\]
where $I$ varies among the ideals of $A$. Indeed, $Spec(A)\setminus {\cal P}_{I}=\mathbin{\mathop{\textrm{\huge $\cup$}}\limits_{a\in I}}\{P\in Spec(A) \textrm{ | } a\in A\setminus P\}$.

Note that, the ideals in $Spec(A)$ being prime, the open subsets of $Spec(A)$ corresponding to the sub-basic open sets ${\cal F}_{a}$ of ${\cal F}_{A}$ form actually a basis of $Spec(A)$ (since for any $a,b\in A$, $\{P\in Spec(A) \textrm{ | } a\notin P\}\cap \{P\in Spec(A) \textrm{ | } b\notin P\}=\{P\in Spec(A) \textrm{ | } a\cdot b\notin P\}$).   

So the topological space obtained by equipping the set $Spec(A)$ with the Zariski topology is homeomorphic to the space ${\cal Z}_{A}\simeq X_{{\tau}^{\Sh(S(A), C)}}$, equivalently to the topological space $X_{{\tau}^{\Sh(L(A), J)}}$. Note that from this latter characterization it immediately follows that ${\cal Z}_{A}$ is a spectral space.
    
We conclude this section by showing that the technical lemma V3.2 \cite{stone} used to establish the bijection between the ideals of the distributive lattice $L(A)$ and the radical ideals of $A$ admits an alternative (and possibly more transparent) proof obtained by translating the thesis, naturally formulated in terms of the site $(L(A), J)$, into a property of the site $(S(A), C)$. This is done by transferring an appropriate topos-theoretic invariant across the Morita-equivalence
\[
\Sh(S(A), C)\simeq \Sh(L(A), J),
\]
according to the philosophy `toposes as bridges' of \cite{OC10}. 

To this end, we need the following result.

\begin{lemma}\label{msites}
Let $F:({\cal C}, J) \to ({\cal D}, K)$ be a morphism of sites. Then we have a commutative diagram (in the $2$-category of toposes)
\[  
\xymatrix {
[{\cal D}^{\textrm{op}}, \Set] \ar[r]^{f'} & [{\cal C}^{\textrm{op}}, \Set] \\
\Sh({\cal D}, K) \ar[u]^{i_{K}} \ar[r]^{f} & \ar[u]_{i_{J}} \Sh({\cal C}, J),}
\]
where $i_{J}$ and $i_{K}$ are the canonical geometric inclusions, $f$ is the geometric morphism induced by $F$ as in in Corollary C2.3.4 \cite{El}, and $f'$ is the geometric morphism induced by the functor $F:{\cal C}\to {\cal D}$ as in A4.1.10 \cite{El}.
\end{lemma}

\begin{proofs}
This follows immediately from the proof of Corollary C2.3.4 \cite{El}, observing that the direct image functor of $f$ is the restriction of the direct image functor of $f'$. Indeed, by the uniqueness (up to isomorphism) of adjoints, in order to prove that two geometric morphisms are isomorphic, it suffices to establish an isomorphism between their direct image functors.  
\end{proofs}

Given a commutative ring with unit $A$, let us denote by $\xi:S(A)\to L(A)$ the factorization of $\pi_{L}:A\to L(A)$ through $\pi_{S}:A\to S(A)$. Clearly, $\xi$ is a morphism of sites $(S(A), C) \to (L(A), J)$ and hence it induces a geometric morphism $f:\Sh(S(A), C) \to \Sh(L(A), J)$. It is immediate to verify that this geometric morphism is precisely the equivalence $\Sh(S(A), C) \simeq \Sh(L(A), J)$ induced by Theorem \ref{fund} as indicated above. So, by Lemma \ref{msites}, we have a commutative diagram
\[  
\xymatrix {
[L(A)^{\textrm{op}}, \Set] \ar[r]^{f'} & [S(A)^{\textrm{op}}, \Set] \\
\Sh(L(A), J) \ar[u]^{i_{J}} \ar[r]^{f} & \ar[u]^{i_{S}} \Sh(S(A), C),}
\]
where $i_{S}$ and $i_{J}$ are the canonical inclusions and $f'$ is the geometric morphism induced by $\xi:S(A)\to L(A)$ as in A4.1.10 \cite{El}.

Now, the statement of Lemma V3.2 \cite{stone} reads as follows:

\emph{Let $a, b_{1}, \ldots, b_{r}$ be elements of a commutative ring $A$ with unit. Then $D(a)\leq D(b_{1})\vee \cdots \vee D(b_{r})$ in $L(A)$ if and only if there exists an integer $n \geq 1$ and elements $b_{1}, \ldots, b_{r}$ such that $a^{n}=b_{1}\cdot c_{1}+ \cdots + b_{r}\cdot c_{r}$.}

To prove the lemma, we can clearly suppose, without loss of generality, that $b_{i}$ is a multiple of $a$ in $A$ (i.e., $b_{i}=a\cdot b_{i}'$ for some $b_{i}'\in A$) for all $i=1, \ldots, r$; this condition ensures that for all $i$, $\pi_{S}(b_{i}) \leq \pi_{S}(a)$ and $\pi_{L}(b_{i}) \leq \pi_{L}(a)$.

Consider the sieve $S_{B}$ in $S(A)$ on $\pi_{S}(a)$ generated by the family 
\[
\{\pi_{S}(b_{i}) \to \pi_{S}(a) \textrm{ | } i=1, \ldots, r \}
\]
in $S(A)$. 

Clearly, $\xi$ sends this family of arrows to the family $\{\pi_{L}(b_{i}) \to \pi_{L}(a) \textrm{ | } i=1, \ldots, r \}$ in $L(A)$; we denote by $L_{B}$ the sieve in $L(A)$ on $\pi_{L}(a)$ generated by this family. So, if $Y_{S}:S(A)\to [S(A)^{\textrm{op}}, \Set]$ and $Y_{L}:L(A)\to [L(A)^{\textrm{op}}, \Set]$ are the Yoneda embeddings, the inverse image functor of $f'$ sends the mono-\\morphism $B_{S}\mono Y_{S}(\pi_{S}(a))$ in $[S(A)^{\textrm{op}}, \Set]$ corresponding to the sieve $S_{B}$ to the monomorphism $B_{L}\mono Y_{L}(\pi_{L}(a))$ in $[L(A)^{\textrm{op}}, \Set]$ corresponding to the sieve $L_{B}$. Now, by the commutativity of the square above and the fact that $f$ is an equivalence, the associated sheaf functor $a_{S}:[S(A)^{\textrm{op}}, \Set] \to \Sh(S(A), C)$ sends $B_{S}$ to an isomorphism (equivalently, $S_{B}$ is $\overline{C}$-covering, where $\overline{C}$ is the Grothendieck topology on $S(A)$ generated by $C$) if and only if the associated sheaf functor $a_{L}:[L(A)^{\textrm{op}}, \Set] \to \Sh(L(A), J)$ sends $B_{L}$ to an isomorphism (equivalently, $L_{B}$ is $J$-covering). We thus deduce that $S_{B}$ is $\overline{C}$-covering if and only if $L_{B}$ is $J$-covering. 

Now, the hypothesis `$D(a)\leq D(b_{1})\vee \cdots \vee D(b_{r})$ in $L(A)$' in the statement of the lemma is equivalent to the statement that $L_{B}$ is $J$-covering. From this, our argument enables us to conclude that $S_{B}$ is $\overline{C}$-covering. We now show that this condition implies the second condition of the lemma. To this end, we give an explicit description of the monoid congruence $\simeq$ involved in the definition of $S(A)$. It is easy to see that $\simeq$ is the transitive closure of the binary relation $R$ on $A$ defined by saying that for any $a,b\in A$, $a R b$ if and only if there exist $c,d\in A$ and non-zero integers $n,m\geq 1$ such that $a=c^{n}\cdot d$ and $b=c^{m}\cdot d$. From this it easily follows that if $a\simeq b$ then there exists a positive integer $k\geq 1$ such that $b=a^{k}\cdot e$ for some element $e$ of $A$; below, we will refer to this condition on $a$ and $b$ as to $(\ast)$. This remark allows us to achieve an explicit description of the Grothendieck topology $\overline{C}$ on $S(A)$.

\begin{lemma}\label{topC}
Let $A$ be a commutative ring with unit. With the above notation, for any sieve $S$ in $S(A)$ on an element $x\in S(A)$, $S\in \overline{C}(x)$ if and only if $S$ contains a finite family of arrows in $S(A)$ of the form $\{\pi_{S}(a_{i})\to \pi_{S}(a) \textrm{ | } i=1, \ldots, n\}$, where $\pi_{S}(a_{i}) \leq_{S_{A}} \pi_{S}(a)$ and $a^{k}=a_{1}\cdot b_{1}+ \cdots + a_{n}\cdot b_{n}$ for some positive integer $k\geq 1$ and elements $b_{1}, \ldots, b_{n}$ in $A$.
\end{lemma}
 
\begin{proofs}
The `if' direction is clear from the definition of the coverage $C$; indeed, if $a^{k}=a_{1}\cdot b_{1}+ \cdots + a_{n}\cdot b_{n}$ with $\pi_{S}(a_{i}) \leq_{S_{A}} \pi_{S}(a)$ then the family of arrows $\{\pi_{S}(a_{i}\cdot b_{i})\to \pi_{S}(a) \textrm{ | } i=1, \ldots, n\}$ is $C$-covering and hence the sieve generated by the family $\{\pi_{S}(a_{i})\to \pi_{S}(a) \textrm{ | } i=1, \ldots, n\}$ is $\overline{C}$-covering, since it contains a $C$-covering family.

The `only if' direction can be deduced as a consequence of the following two facts:

\begin{enumerate}[(i)]

\item Every $C$-covering sieve satisfies the condition in the statement of the lemma; 

\item The collection of the sieves on $S(A)$ which satisfy the condition in the statement of the lemma forms a Grothendieck topology;
\end{enumerate}    

Fact $(i)$ follows easily from the definition of the coverage $C$ and the explicit description of the monoid congruence $\simeq$ obtained above, while Fact $(ii)$ is immediately verified.
\end{proofs} 

Coming back to our original problem, the `if' direction of Lemma V3.2 \cite{stone} is trivial, so we only care to prove the `only if' one. If $S_{B}$ is $\overline{C}$-covering then, by Lemma \ref{topC}, $S_{B}$ contains a finite family of arrows of the form $\{\pi_{S}(a_{i})\to \pi_{S}(a) \textrm{ | } i=1, \ldots, n\}$, where $\pi_{S}(a_{i}) \leq_{S_{A}} \pi_{S}(a)$ and $a^{k}=a_{1}\cdot c_{1}+ \cdots + a_{n}\cdot c_{n}$ for some positive integer $k\geq 1$ and elements $c_{1}, \ldots, c_{n}$ in $A$. Now, by definition of $S_{B}$, for every $i\in \{1,\ldots, n\}$ there exists $j_{i}\in \{1, \ldots, r\}$ such that $\pi_{S}(a_{i}) \leq \pi_{S}(b_{j_{i}})$, equivalently $a_{i}\cdot b_{j_{i}}\simeq a_{i}$. From this we deduce, by invoking the characterization of the congruence $\simeq$ on $A$ obtained above, that for any sufficiently large natural number $k_{i}$, $a_{i}^{k_{i}}$ is a multiple of $b_{j_{i}}$ in $A$, which in turn implies that for a sufficiently large natural number $k$, $a^{k}$ is a sum of multiples of the $b_{i}$. This completes our proof of Lemma V3.2 \cite{stone}.  

It is natural to wonder whether there is a natural way to generalize the Zariski topology on the collection of prime ideals of a ring to the collection of all the proper ideals of the ring, or to some other more general class of subsets of the ring. Thanks to the techniques developed in the paper, we are able to give a positive answer to this question. We can use propositional geometric theories to describe subsets of a ring with particular properties, such as the class of ideals of the ring; the subterminal topology then provides a way of endowing the collection of models of such a theory with a topology such that the topos of sheaves on the resulting topological space is equivalent to the classifying topos of the propositional theory; also, as we have seen above, the results of section \ref{genrel} enable us to achieve, in many cases of interest, an explicit semantic description of such classifying topos as a topos of sheaves on a poset structure presented by generators and relations with respect to some Grothendieck topology on it.  

For example, consider the case of (proper) ideals of a commutative ring with unit $A$. In order to obtain the Zariski topology, we characterized the prime ideals of $A$ as the complements in $A$ of the models of a particular propositional theory, namely the theory of prime filters on $A$. Therefore, a first natural step to generalize the Zariski topology to a topology on the set of proper ideals of $A$ is to try to find a propositional theory whose models are exactly the complements in $A$ of proper ideals; we will refer to these subsets as to the \emph{op-ideals} on $A$. Explicitly, an op-ideal on $A$ is a subset $S\subseteq A$ such that for all $a, b\in A$, $a\cdot b\in S$ implies $a\in S$, $a+b\in I$ implies $a\in I$ or $b\in I$, $0_{A}\notin I$ and $1_{A}\in I$.       

Op-ideals can be described as the models of the propositional coherent theory $\mathbb Q$ over the signature consisting of one $0$-ary relation symbol $P_{a}$ for each element $a\in A$, with axioms:
\[
P_{0_{A}} \vdash \bot;
\]

\[
\top \vdash P_{1_{A}};
\]

\[
(P_{a\cdot b} \vdash P_{a} \wedge P_{b})
\] 
for any $a, b$ in $A$, and 
\[
(P_{a + b} \vdash P_{a} \vee P_{b})
\]
for any $a, b \in A$. 

The classifying topos of $\mathbb Q$, which is equivalent to the category of sheaves on the coherent syntactic category ${\cal C}_{\mathbb Q}$ of $\mathbb Q$ with respect to the coherent topology on it, can be represented as the topos of sheaves on the topological space obtained by equipping the set $Spec_{id}(A)$ of ideals on $A$ with the elemental topology (cf. section \ref{subs}). Note that ${\cal C}_{\mathbb Q}$ can be characterized as the distributive lattice $D(A)$ generated by symbols $D(a)$, $a \in A$, subject to the relations $D(0_{A})=0_{D(A)}$, $D(1_{A})=1_{D(A)}$, $D(a\cdot b) \leq D(a) \wedge D(b)$, $D(a+ b) \leq D(a) \vee D(b)$. The frame of open sets of $Spec_{id}(A)$ is thus isomorphic to the frame of ideals of the distributive lattice ${\cal C}_{\mathbb Q}$, and, by the results of sections \ref{clex} and \ref{genrel}, can also be characterized as the geometric syntactic category of the theory $\mathbb Q$, that is as the frame $F(A)$ generated by symbols $F(a)$, $a \in A$, subject to the relations $F(0_{A})=0_{F(A)}$, $F(1_{A})=1_{F(A)}$, $F(a\cdot b) \leq F(a) \wedge F(b)$ and $F(a+ b) \leq F(a) \vee F(b)$.

\section{Conclusions}\label{conclusions}

The theory developed in this paper paves the way for a vast world of new possibilities. Three natural directions that one could immediately pursue are the following. 

First, one can investigate particular dualities already generated through our machinery, and described in sections \ref{ex} and \ref{addex}; this should be interesting because, as it is clear from our topos-theoretic interpretation, all of these dualities have essentially the same level of `mathematical depth' as the classical Stone duality. 

Second, one can generate new dualities for particular classes of preordered structures by using our machinery. As we have already remarked in the course of the paper, the process of generation of new `Stone-type' dualities is essentially automatic, and relies on the choice of two (in the case of dualities with locales) or three (in the case of dualities with topological spaces) ingredients to give our `machine' as `inputs': the initial category $\cal K$ of preordered structures, the Grothendieck topologies $J_{\cal C}$ associated to the structures in $\cal K$, and the sets of points of the toposes $\Sh({\cal C}, J_{\cal C})$ corresponding to the structures. These choices are, although in a relation of sequential dependence each on the previous ones, essentially independent from each other, in the sense that the choice of a certain ingredient at one stage does not uniquely determine the choices of the ingredients at later stages. These two (resp. three) `degrees of freedom' in the choice of ingredients confer to our machinery a great level of techical flexibility, which allows us in particular to identify a given category of preorders with a subcategory of a category of locales or topological spaces in several non-equivalent ways (cf. section \ref{ex} for some examples of this phenomenon). In fact, the more general methodology of section \ref{Mordual} shows that there are essentially \emph{three} degrees of freedom in generating dualities (resp. equivalences) between general categories of preordered structures, with \emph{two} additional degrees of freedom if one wants to `lift' these dualities (resp. equivalences) to `topological' ones.   

Third, one can perform topos-theoretic translations between between properties of preordered structures and properties of the locales or topological spaces corresponding to them via both the `Stone-type dualities', according to the technique `toposes as bridges' of \cite{OC10}; as we remarked in section \ref{insights}, this translation can be performed automatically in many cases, and, as the examples provided in section \ref{insights} show, the results generated in this way are non-trivial in general.

Overall, this paper can be read as an assay of what can be achieved by applying the methodologies of \cite{OC10}. Indeed, in \cite{OC10}, we emphasized the importance of taking Morita-equivalences as starting points of topos-theoretic investigations, and of using topos-theoretic invariants to extract information about them which is relevant for classical mathematics; the common classifying topos acts in this context as a `bridge' for transferring properties and constructions between its two different sites of definition related to each other by the Morita-equivalence.  Now, in this paper, starting from the Morita-equivalence of Theorem \ref{fund}, we have proceeded to extract information about it through the use of topos-theoretic invariants of various kinds. Through the lenses of the unifying framework that we have developed, the different Stone-type dualities, as well as the other results that we have obtained in the paper, appear to be just `variations on the same theme', this theme being precisely the original Morita-equivalence (or, more generally, any equivalence of the form $\Sh({\cal C}, J)\simeq \Sh({\cal D}, K)$ as considered in section \ref{generalization}).

In order to keep our treatment self-contained, we have chosen not to discuss in detail in the paper any examples involving categories which are not preorder, but in fact there are many ideas in the paper, motivated by the principles introduced in \cite{OC10}, which can be profitably extended from the preorder context to the level of arbitrary (small) categories. For example, the idea of recovering a category $\cal C$ from a topos $\Sh({\cal C}, J)$, where $J$ is a subcanonical topology on $\cal C$, through a topos-theoretic invariant is clearly very general and as such it is liable to be applied in a variety of different contexts. In fact, whenever we have a Morita-equivalence $\Sh({\cal C}, J)\simeq \Sh({\cal D}, K)$ such that $\cal C$ can be recovered from the topos $\Sh({\cal C}, J)$ through a topos-theoretic invariant $U$, we can represent $\cal C$ in terms of the site $({\cal D}, K)$ as the category of objects of the topos $\Sh({\cal D}, K)$ which satisfy the invariant $U$, and if this equivalence holds `naturally' for $\cal C$ varying in a (large) category $\cal K$ of small categories, we can expect to be able to `functorialize' this representation of $\cal C$ in terms of $({\cal D}, K)$ so to yield a duality between $\cal K$ and a category consisting of objects from which the toposes $\Sh({\cal D}, K)$ can be `directly built' (so, for example, if the topologies $K$ on $\cal D$ are `uniformly defined' in terms of $\cal D$ by means of a topos-theoretic invariant then specifying $\cal D$ is enough to `build' the site $({\cal D}, K)$ and hence the topos $\Sh({\cal D}, K)$). 

While we are on the topic of Morita-equivalences holding for large classes of toposes, we remark that a `first-order analogue' of the Morita-equivalence of Theorem \ref{fund} - which represents the topos of sheaves on a site whose underlying category is a preorder as the topos of sheaves on a locale - is the representation theorem \cite{AT} of Joyal and Tierney, which gives a representation of a general Grothendieck topos as the topos of sheaves on a localic groupoid, while an analogue of the representation of a localic topos with enough points as a topos of sheaves on a topological space is provided by Butz and Moerdijk's representation theorem \cite{BM} for Grothendieck toposes with enough points as toposes of sheaves on topological groupoids. As we have already remarked, we can expect the general methodology `toposes as bridges' recalled above to yield interesting results also in the context of these Morita-equivalences, in the form of representation theorems, adjunctions, dualities for various classes of categories or topos-theoretic translations of properties (or constructions) from one side of such a Morita-equivalence to the other. Incidentally, a result which already witnesses the effectiveness of our abstract technique also in the first-order context is the duality between the category of ($k$-small) Boolean pretoposes and Stone topological groupoids obtained in \cite{AF}, which can be seen as arising from the processes of functorializing (an adaptation of) Butz and Moerdijk's representation theorem and recovering a Boolean pretopos from the corresponding topos of coherent sheaves though a topos-theoretic invariant. 

Finally, as already remarked in section \ref{concreteness}, the methods of this paper, combined with the philosophy `toposes as bridges' of \cite{OC10}, provide a topos-theoretic interpretation of the problem of finding `concrete' descriptions of models of ordered algebraic theories (in the sense of section \ref{spacesprop}) presented by generators and relations. Indeed, given such a model $\cal C$, regarded as a preorder category, for any Grothendieck topology $J$ such that $\cal C$ can be recovered, up to isomorphism, from $\Sh({\cal C}, J)$ through a topos-theoretic invariant, any alternative representation of the topos $\Sh({\cal C}, J)$ in terms of a different site $({\cal D}, K)$ yields an alternative description of $\cal C$ in terms of this latter site (applications of this methodology are given in section \ref{concreteness}); conversely, any explicit description of $\cal C$ as a preorder isomorphic to a certain structure $\cal D$ yields a Morita-equivalence $\Sh({\cal C}, J)\simeq \Sh({\cal D}, K)$, where $K$ is the Grothendieck topology on $\cal D$ given by the image of $J$ under the isomorphism ${\cal C}\simeq {\cal D}$. The problem of finding alternative descriptions of structures presented by generators and relations is thus strictly connected to the problem of finding alternative representations for toposes of sheaves on the structures with respect to appropriate subcanonical topologies. This interpretation paves the way for the use of topos-theoretic methods for addressing this kind of problems, which are in general rather hard to solve (cf. for example the problem discussed in section \ref{freecomjoin}); compelling examples of the validity of these methods are already given in the paper (specifically in section \ref{concreteness}), and the reader can easily generate new ones by following the same principles.
 
\newpage
       
\section{Appendix}\label{appendix}

In this appendix, in order to make (parts of) the theory developed in the paper intelligible to the reader who is not familiar with Topos Theory, we provide elementary proofs (i.e. proofs consisting of direct, non-topos-theoretic arguments) of some of the key results in the paper, and indicate how the arguments which constitute the general machinery for building dualities of sections \ref{locales}, \ref{dualtop} and \ref{Mordual} can be rewritten in `elementary' terms.  

Let us start from the construction of frames from preorders (in general, arbitrary categories) equipped with Grothendieck topologies.  

The following notions are the specializations to preorders of the well-known topos-theoretic notions of coverage and of site. For an element $c\in {\cal C}$ of a preorder $({\cal C}, \leq)$, we define the principal ideal $(c)\downarrow$ on $c$ as the set of all the elements $d$ of $\cal C$ such that $d\leq c$ in $\cal C$. 

\begin{definition}
Let $({\cal C}, \leq)$ be a preorder. 

\begin{enumerate}[(i)]

\item A \emph{coverage} on $\cal C$ is a function $J$ which assigns to every element $c\in {\cal C}$ a family $J(c)$ of subsets of $(c)\downarrow$ such that for any $S\in J(c)$ and any $c'\leq c$ the subset $S_{c'}=\{d\leq c' \textrm{ | } d\in S\}$ belongs to $J(c')$; 

\item A coverage on $\cal C$ is said to be a \emph{Grothendieck coverage} if for any $c\in {\cal C}$, $(c)\downarrow \in J(c)$ and for any subset $S\subseteq (c)\downarrow$, if $S_{c'}\in J(c')$ for every $c'\in T$ where $T\in J(c)$ then $S\in J(c)$; 

\item A \emph{site} is a pair $({\cal C}, J)$, where $\cal C$ is a preorder and $J$ is a coverage on $\cal C$;

\item A coverage $J$ on $\cal C$ is \emph{subcanonical} if for every $c\in {\cal C}$ and any subset $S\in J(c)$, $c$ is the supremum in $\cal C$ of the elements $d\in S$ (i.e., for any element $c'$ in $\cal C$ such that for every $d\in S$ $d\leq c'$, we have $c\leq c'$). 
\end{enumerate}
\end{definition}

For example, on a distributive lattice $\cal C$ we can put the \emph{coherent coverage} $J$, defined by saying that for any $c\in {\cal C}$, the subsets $S$ in $J(c)$ are precisely the subsets of $(c)\downarrow$ which contain a finite subset whose join is $c$. Given a frame $L$, we define the \emph{canonical coverage} on $L$ as the coverage $J_{can}^{L}$ such that for any $c\in L$ and $S\subseteq (c)\downarrow$, $S\in J_{can}^{L}(c)$ if and only if $\mathbin{\mathop{\textrm{\huge $\vee$}}\limits_{d\in S}}d=c$.

\begin{remarks}

\begin{enumerate}[(a)]
\item If $({\cal C}, \wedge)$ is a meet-semilattice and $J$ is a function which assigns to every element $c\in {\cal C}$ a family $J(c)$ of subsets of $(c)\downarrow$ such that for any $S\in J(c)$ and any $c'\leq c$, the subset $\{a\wedge c' \textrm{ | } a\in S\}$ belongs to $J(c')$ then the function $J'$ which sends an element $c$ of $\cal C$ to the collection of the subsets of $(c)\downarrow$ which contain a subset in $J(c)$ is a coverage on ${\cal C}$; 

\item Given a coverage $J$ on a preorder $\cal C$, there is a smallest Grothendieck coverage $J'$ on $\cal C$ which contains $J$ (i.e. such that for any $c\in {\cal C}$, $S\in J(c)$ implies $S\in J'(c)$), called the \emph{Grothendieck coverage generated by $J$}: $J'$ can be defined by saying that for any $c\in {\cal C}$, the subsets in $J(c)$ are precisely the subsets which belong to $J'(c)$ for every Grothendieck coverage $J'$ containing $J$.  
\end{enumerate}
\end{remarks}

Let us recall from section \ref{general} the following notions.

\begin{definition}

Let $({\cal C}, \leq)$ be a preorder and $J$ be a coverage on $\cal C$.
\begin{enumerate}

\item  An \emph{ideal} on $\cal C$ is a subset $I\subseteq {\cal C}$ such that for any $a,b \in {\cal C}$ such that $b\leq a$ in $\cal C$, $a\in I$ implies $b\in I$. A \emph{$J$-ideal} on $\cal C$ is an ideal on $\cal C$ such that for any $R\in J(c)$, if $a\in I$ for every $a\in R$ then $c\in I$. 

\item Given a subset $I$ on $\cal C$, the \emph{$J$-closure $\overline{I}^{J}$} is the smallest $J$-ideal on $\cal C$ containing $I$, equivalently the intersection of all the $J$-ideals which contain $I$. 

\item Given an object $c$ of $\cal C$, we define the \emph{principal $J$-ideal} $(c)\downarrow_{J}$ generated by $c$ as the $J$-closure of the subset $(c)\downarrow$ of $\cal C$.
\end{enumerate}
\end{definition}

\begin{remark}
Note that, if $J$ is the trivial coverage on $\cal C$ (i.e., the coverage in which for any $c\in {\cal C}$, $S\in J(c)$ if and only if $S=(c)\downarrow$) then the $J$-ideals on $\cal C$ are precisely the ideals on $\cal C$.
\end{remark}

Given a preorder ${\cal C}$, equipped with a coverage $J$, we define \emph{$Id_{J}({\cal C})$} to be the set of all the $J$-ideals on $\cal C$, equipped with the subset-inclusion ordering $\subseteq$. The set of all the ideals on $\cal C$ will be denoted by $Id({\cal C})$.

\begin{remark}\label{grotcov}

\begin{enumerate}[(a)]
\item
If $J$ is a Grothendieck coverage on a preorder $({\cal C}, \leq)$ then for any ideal $I$ on $\cal C$, the $J$-closure $\overline{I}^{J}$ of $I$ is equal to the set of elements $d\in {\cal C}$ such that there exists $R\in J(d)$ with the property that for every $a\in R$, $a\in I$;  

\item If $J$ is a coverage on a preorder $({\cal C}, \leq)$ and $J'$ is the Grothendieck coverage on $\cal C$ generated by $J$ then for any ideal $I$ on $\cal C$, $I$ is $J$-closed if and only if $I$ is $J'$-closed; in particular, the $J$-closure of $I$ coincides with the $J'$-closure of $I$. 
\end{enumerate}
\end{remark}

\begin{theorem}
Let $\cal C$ be a preorder and $J$ be a coverage on $\cal C$. Then $(Id_{J}({\cal C}), \subseteq)$ is a frame.
\end{theorem}

\begin{proofs}
Clearly, the intersection of any two $J$-ideals on $\cal C$ is a $J$-ideal; this gives the meet in $Id_{J}({\cal C})$. The join $\mathbin{\mathop{\textrm{\huge $\vee$}}\limits_{k\in K}}I_{k}$ of a family $\{I_{k} \textrm{ | } k\in K \}$ of $J$-ideals on $\cal C$ is given by the $J$-closure of the set-theoretic union of the $I_{k}$; explicitly, $\mathbin{\mathop{\textrm{\huge $\vee$}}\limits_{k\in K}}I_{k}=\{d\in {\cal C} \textrm{ | } \{e\leq d \textrm{ | } e\in \mathbin{\mathop{\textrm{\huge $\cup$}}\limits_{k\in K}}I_{k}\}\supseteq S \textrm{ for some } S\in J'(d)\}$, where $J'$ is the Grothendieck topology on $\cal C$ generated by $J$ (cf. Remark \ref{grotcov}(a)).

The infinite distributive law of joins with respect to finite meets is easily verified, by using the infinite distributive law of unions with respect to finite intersections and the fact that the operation of $J'$-closure commutes with finite intersections.   
\end{proofs}

Let us recall the following notion from section \ref{locales}.

\begin{definition}
Let $\cal C$ and $\cal D$ be preorders. A monotone map $f:{\cal C}\to {\cal D}$ is said to be \emph{flat} if the following two conditions hold:
\begin{enumerate}[(i)]
\item For any $d\in {\cal D}$ there exists $c\in {\cal C}$ such that $d\leq f(c)$;

\item For any element $d\in {\cal D}$ and any elements $c,c'\in {\cal C}$ such that $d\leq f(c)$ and $d\leq f(c')$ there exists $c''\in {\cal C}$ such that $c''\leq c$, $c''\leq c'$ and $d\leq f(c'')$.
\end{enumerate}
\end{definition}

\begin{remark}
If $\cal C$ and $\cal D$ are meet-semilattices then a monotone map $f:{\cal C}\to {\cal D}$ is flat if and only if it is a meet-semilattice homomorphism.
\end{remark}

If $\cal C$ and $\cal D$ are two preorders equipped respectively with coverages $J$ and $K$, we say that a flat map $f:{\cal C}\to {\cal D}$ is \emph{cover-preserving} if for any element $c\in {\cal C}$ and any $R\in J(c)$, the set of elements $d$ of $\cal D$ such that $d\leq f(c')$ for some $c'\in R$ contains a subset belonging to $K(f(c))$.

\begin{theorem}\label{concretecontrov}
Let $\cal C$ and $\cal D$ be two preorders equipped respectively with coverages $J$ and $K$, and let $f:{\cal C}\to {\cal D}$  be a flat cover-preserving map between them. Then the map 
\[
A_{f}:Id_{J}({\cal C}) \to Id_{K}({\cal D})
\]
sending any $J$-ideal $I$ on $\cal C$ to the smallest $K$-ideal on $\cal D$ containing $f(I)$ is a frame homomorphism (where $Id_{J}({\cal C})$ and $Id_{K}({\cal D})$ are endowed with the subset-inclusion ordering).  
\end{theorem}

\begin{proofs}
Let us begin by proving that $A_{f}$ preserves arbitrary joins. Given a family $\{I_{h} \textrm{ | } h\in H\}$ of $J$-ideals on $\cal C$, we have that  $A_{f}(\mathbin{\mathop{\textrm{\huge $\vee$}}\limits_{h\in H}}I_{h})=\overline{f(\mathbin{\mathop{\textrm{\huge $\vee$}}\limits_{h\in H}}I_{h})}^{K}=\overline{f(\overline{\mathbin{\mathop{\textrm{\huge $\cup$}}\limits_{h\in H}}I_{h}}^{J})}^{K}=\overline{\mathbin{\mathop{\textrm{\huge $\cup$}}\limits_{h\in H}}f(I_{h})}^{K}=\mathbin{\mathop{\textrm{\huge $\vee$}}\limits_{h\in H}}\overline{f(I_{h})}^{K}=\mathbin{\mathop{\textrm{\huge $\vee$}}\limits_{h\in H}}A_{f}(I_{h})$, where the third equality follows from the fact that $f$ is cover-preserving.

To prove that $A_{f}$ is a frame homomorphism, it remains to show that it preserves the top element and binary meets.

The fact that $A_{f}$ preserves the top element follows immediately from condition $(i)$ in the definition of flat map.

To prove that $A_{f}$ preserves finite meets, suppose that $I_{1}$ and $I_{2}$ are $J$-ideals on $\cal C$. Then, condition $(ii)$ in the definition of flat map ensures that the ideal on $\cal D$ generated by $f(I_{1})\cap f(I_{2})$ coincides with the ideal on $\cal D$ generated by $f(I_{1}\cap I_{2})$; therefore $A_{f}(I_{1}\cap I_{2})=\overline{f(I_{1}\cap I_{2})}^{K}=\overline{f(I_{1})\cap f(I_{2})}^{K}=\overline{f(I_{1})}^{K}\cap\overline{f(I_{2})}^{K}=A_{f}(I_{1})\cap A_{f}(I_{2})$.
\end{proofs}

\begin{definition}
Let $\cal C$ be a preorder and $J$ be a coverage on $\cal C$. We say that $J$ is \emph{subcanonical} if for every $c\in {\cal C}$ and any $S\in J(c)$, for any $c'\in {\cal C}$, if $d\leq c'$ for every $d\in S$ then $c\leq c'$).
\end{definition}

\begin{remark}
A coverage $J$ on a preorder $\cal C$ is subcanonical if and only if $(c)\downarrow_{J}=(c)\downarrow$ for any $c\in {\cal C}$.
\end{remark}

Given a flat cover-preserving map $f:{\cal C}\to {\cal D}$ as in the statement of Theorem \ref{concretecontrov}, and denoted by $i_{\cal C}:{\cal C}\to Id_{J}({\cal C})$ (resp. $i_{\cal D}:{\cal D}\to Id_{K}({\cal D})$) the map sending an element $c\in {\cal C}$ to the principal $J$-ideal $(c)\downarrow_{J}\in Id_{J}({\cal C})$ (resp. an element $d\in {\cal D}$ to the principal $K$-ideal $(d)\downarrow_{K}\in Id_{K}({\cal D})$), $A_{f}\circ i_{\cal C}=i_{\cal D}\circ f$. In particular, if the coverages $J$ and $K$ satisfy $(c)\downarrow_{J}=(c)\downarrow$ for any $c\in {\cal C}$ and $(d)\downarrow_{K}=(d)\downarrow$ for any $d\in {\cal D}$ (that is, if $J$ and $K$ are subcanonical) and the preorders $\cal C$ and $\cal D$ are posets then $f$ can be recovered from $A_{f}$ as its restriction to the subsets of principal ideals. 

Theorem \ref{concretecontrov} provides a way for obtaining frame homomorphisms between frames of ideals on preordered structures starting from monotone maps between the structures. On the other hand, if the structures are not equipped with any coverage, we can build frame homomorphisms between the frames of ideals on them, as follows. 

For any monotone map $f:{\cal C}\to {\cal D}$ between preorders, the map $B_{f}:Id({\cal D})\to Id({\cal C})$ sending an ideal $I$ on $\cal D$ to the inverse image $f^{-1}(I)$ of $I$ under $f$ is a frame homomorphism. It is natural to wonder if it is possible to recover $f$ from $B_{f}$ and characterize `intrinsically' the frame homomorphisms $Id({\cal D})\to Id({\cal C})$ of the form $B_{f}$ for some monotone map $f:{\cal C}\to {\cal D}$. The following result provides an answer to these questions. Below, for a preorder $\cal C$, we denote by $i_{\cal C}:{\cal C}\to Id({\cal C})$ the map sending an element $c\in {\cal C}$ to the principal ideal $(c)\downarrow \in Id({\cal C})$.

\begin{theorem}\label{concretecovariant}
\begin{enumerate}[(i)]
\item Let $\cal C$ and $\cal D$ be two posets. A frame homomorphism $F:Id({\cal D})\to Id({\cal C})$ is of the form $B_{f}$ for some monotone map $f:{\cal C}\to {\cal D}$, where $B_{f}:Id({\cal D})\to Id({\cal C})$ is the frame homomorphism sending an ideal $I$ on $\cal D$ to the inverse image $f^{-1}(I)$ of $I$ under $f$, if and only if $F$ preserves arbitrary infima, equivalently it has a left adjoint $F_{!}:Id({\cal C})\to Id({\cal D})$, given by the formula $F_{!}(I)=\mathbin{\mathop{\textrm{\huge $\cap$}}\limits_{I \subseteq F(I')}I'}$ (for any $I\in Id({\cal C})$);

\item For any monotone map $f:{\cal C}\to {\cal D}$, ${(B_{f})}_{!}\circ i_{\cal C}=i_{\cal D}\circ f$; in particular, $f$ can be recovered from $B_{f}$ as the restriction of its left adjoint ${(B_{f})}_{!}$ to the subsets of principal ideals.
\end{enumerate}
\end{theorem}

\begin{proofs}
$(i)$ It is clear that for any monotone map $f:{\cal C}\to {\cal D}$, $B_{f}:Id({\cal D})\to Id({\cal C})$ preserves arbitrary intersections (i.e., arbitrary infima). Conversely, suppose that $F:Id({\cal D})\to Id({\cal C})$ preserves arbitrary intersections, equivalently (by the Special Adjoint Functor Theorem) that it has a left adjoint $F_{!}:Id({\cal C})\to Id({\cal D})$ given by the formula $F_{!}(I)=\mathbin{\mathop{\textrm{\huge $\cap$}}\limits_{I \subseteq F(I')}I'}$ (for any $I\in Id({\cal C})$). Let us show that $F_{!}$ sends principal ideals to principal ideals. Given $c\in {\cal C}$, the formula for $F_{!}$ yields $F_{!}((c)\downarrow)=\mathbin{\mathop{\textrm{\huge $\cap$}}\limits_{c \in F(I')}I'}$. Now, since $F$ preserves arbitrary intersections, $T:=\mathbin{\mathop{\textrm{\huge $\cap$}}\limits_{c \in F(I')}I'}$ satisfies the property that $c\in F(T)$ and is the smallest ideal on $\cal D$ with this property. This minimality condition implies that $F_{!}((c)\downarrow)=T$ is a principal ideal on $\cal D$, since $c\in F(T)$ implies that $c\in F((t)\downarrow)$ for some $t\in T$. Now, since $\cal C$ and $\cal D$ are posets, the principal ideals on them can be identified with the elements that generate them, and hence we have a monotone map $f:{\cal C}\to {\cal D}$ such that for every $c\in {\cal C}$, $F_{!}((c)\downarrow)= (f(c))\downarrow$. This implies, since $F_{!}$ preserves unions, that for any ideal $I$ on $\cal C$, $F_{!}(I)=f(I)$, which in turn implies that the right adjoint $F$ to $F_{!}$ is equal to the inverse image map $f^{-1}:Id({\cal D})\to Id({\cal C})$, in other words to $B_{f}$. 

$(ii)$ This follows immediately from the proof of part $(i)$.
\end{proofs}

The method of section \ref{locales} for constructing dualities or equivalences between a given category of posets and a category of locales can be reformulated in `concrete' terms as follows. 

Let us start with the controvariant case. Given a collection of preorders $\cal C$, each of which equipped with a coverage $J_{\cal C}$, consider the category $\cal K$ having as objects the posets $\cal C$ and as arrows ${\cal C} \to {\cal D}$ the flat cover-preserving maps $({\cal C}, J_{\cal C})\to ({\cal D}, J_{\cal D})$. Then we have a functor $A:{\cal K}\to \textbf{Frm}$ sending any preorder $\cal C$ in $\cal K$ to $Id_{J_{\cal C}}({\cal C})$ and any arrow $f:{\cal C}\to {\cal D}$ in $\cal K$ to the frame homomorphism $A_{f}:Id_{J_{\cal C}}({\cal C}) \to Id_{J_{\cal D}}({\cal D})$ given by Theorem \ref{concretecontrov}. 

In the covariant case, one defines a functor $B:\textbf{Pro}^{\textrm{op}}\to \textbf{Frm}$, where $\textbf{Pro}$ is the category of preorders and monotone maps between them, by sending a preorder $\cal C$ to the frame $Id({\cal C})$ of ideals on $\cal C$, and a monotone map $f$ to the frame homomorphism $B_{f}:Id({\cal D})\to Id({\cal C})$ of Theorem \ref{concretecovariant}.

In the controvariant case, if the topologies $J_{\cal C}$ are subcanonical and the preorders $\cal C$ are posets then the map $i_{\cal C}:{\cal C}\to Id_{J}({\cal C})$ is an embedding and hence we can hope to find a categorical inverse to the functor $A$, by using the technique of section \ref{charinv}. 

In the covariant case, we can find a categorical inverse to the restriction to the opposite of the category $\textbf{Pos}$ of posets of the functor $B:\textbf{Pro}^{\textrm{op}}\to \textbf{Frm}$, by using the technique of section \ref{charinv}. 

The proofs of the main results of section \ref{charinv} are already direct and all the concepts which appear in them can be straightforwardly reformulated in terms of frame-theoretic invariants, rather than in terms of topos-theoretic invariants of families of subterminals in a locally small cocomplete topos (cf. Remark \ref{loctop}). 

Summarizing the results obtained in section \ref{charinv}, we have that if all the (Grothendieck topologies generated by the) coverages $J_{\cal C}$ are $C$-induced (in the sense of Definition \ref{inducedtop}) for a frame-theoretic invariant $C$ satisfying the property that for any structure $\cal C$ in $\cal K$ and for any family $\cal F$ of principal $J_{\cal C}$-ideals on $\cal C$, $\cal F$ has a refinement satisfying $C$ (if and) only if it has a refinement made of principal $J_{\cal C}$-ideals on $\cal C$ satisfying $C$, then we can define a functor $I_{A}$ (resp. a functor $I_{B}$) on the extended image of the functor $A$ (resp. of the functor $B$) which is a categorical inverse to the functor $A$ (resp. to the functor $B$). The notion of $C$-compactness plays a central role in obtaining `intrinsic' characterization of the extended image of the functor $A$ (resp. of the functor $B$), and in defining the functor $I_{A}$ (resp. the functor $I_{B}$). Recall that an element $l$ of a frame $L$ is said to be $C$-compact if every covering of $l$ in $L$ admits a refinement which satisfies $C$. The functor $I_{A}$ (resp. the functor $I_{B}$) sends a frame $L$ in the extended image of $A$ (resp. of $B$) to the collection of the elements of $L$ which are $C$-compact and acts on the arrows accordingly. We refer to section \ref{charinv} for the rigorous definitions and the details of this technique.

The more general methodology of section \ref{Mordual} takes as starting point a bunch of equivalences of the form $Id_{J_{\cal C}}({\cal C})\cong Id_{K_{\cal D}}({\cal D})$ and relies on the existence of invariants $C$ and $D$ of families of elements of frames which satisfy the property that for any family  $\cal F$ of principal $J_{\cal C}$-ideals on $\cal C$ (resp. of $K_{\cal D}$-ideals on $\cal D$), $\cal F$ has a refinement satisfying $C$ (resp. $D$) (if and) only if it has a refinement made of principal $J_{\cal C}$-ideals on $\cal C$ (resp. of principal $K_{\cal D}$-ideals on $\cal D$)) satisfying $C$ (resp. $D$), and such that all the topologies $J_{\cal C}$ (resp. of $K_{\cal D}$) are $C$-induced (resp. $D$-induced). The technique of section \ref{Mordual} then produces equivalences between categories of posets consisting of the $C$-compact elements of some frame $L$ and categories of posets consisting of the $D$-compact elements of $L$ (cf. section \ref{Mordual} for the details). The technique of section \ref{locales} is recovered as a particular case of this more general technique when the original equivalences are of the form $Id_{J_{\cal C}}({\cal C})\cong Id_{Id_{J_{\cal C}}({\cal C})}(Id_{J_{\cal C}}({\cal C}))$, where $J^{Id_{J_{\cal C}}({\cal C})}_{can}$ is the canonical topology on the frame $Id_{J_{\cal C}}({\cal C})$.

The method of section \ref{Mordual} works for arbitrary bunches of equivalences 
\[
Id_{J_{\cal C}}({\cal C})\cong Id_{K_{\cal D}}({\cal D}).
\] 
An effective way for generating such equivalences is to use the following result, which is the particular case of the Comparison Lemma in Topos Theory for preorder categories (in fact, all the known Stone-type dualities considered in this paper, as well as the new ones that we generate through our methodology, arise from equivalences which are instances of the Comparison Lemma).

Let $\cal C$ be a preorder equipped with a coverage $J$, and let $\cal D$ be a subset of $\cal C$ (endowed with the induced order) which is \emph{$J$-dense}, i.e. such that for every $c\in {\cal C}$ there exists $R\in J(c)$ such that for every $b\in R$, $b\in {\cal D}$. We can define the \emph{induced coverage} $J|_{{\cal D}}$ on $\cal D$ as follows: for any $d\in {\cal D}$, $R\in J|_{{\cal D}}(d)$ if and only if $R=T\cap {\cal D}$ for some $T\in J(d)$.

\begin{theorem}

\begin{enumerate}[(i)]
\item Let $\cal C$ be a preorder equipped with a coverage $J$, and let $\cal D$ be a $J$-dense subset of $\cal C$. Then the map $\phi:Id_{J}({\cal C})\to Id_{J|_{{\cal D}}}({\cal D})$ sending a $J$-ideal $I$ on $\cal C$ to its intersection with $\cal D$ is a frame isomorphism, with inverse the map $\psi:Id_{J|_{{\cal D}}}({\cal D}) \to Id_{J}({\cal C})$ sending a subset in $Id_{J|_{{\cal D}}}({\cal D})$ to its $J$-closure in $\cal C$.  

\item Any locale $L$ with a basis $B_{L}$ is isomorphic to the frame $Id_{J^{L}_{can}|_{B_{L}}}(B_{L})$ of $J_{can}|_{B_{L}}$-ideals on $B_{L}$, via the map 
\[
\phi:L\to Id_{J^{L}_{can}|_{B_{L}}}(B_{L})
\]
sending an element $l\in L$ to the subset given by the intersection $B_{L}\cap (l)\downarrow$ and the map 
\[
\psi:Id_{J^{L}_{can}|_{B_{L}}}(B_{L})\to L
\]
sending an ideal $I$ in $Id_{J^{L}_{can}|_{B_{L}}}$ to the supremum $\mathbin{\mathop{\textrm{\huge $\vee$}}\limits_{l\in I}l}$ of $I$ in $L$.
\end{enumerate}
\end{theorem} 

\begin{proofs}

$(i)$ It is immediate to verify that $\phi$ is a frame homomorphism with inverse $\psi$. Indeed, given a $J$-ideal $I$ on $\cal C$, the fact that $\cal D$ is $J$-dense ensures that $I$ is the $J$-closure of its intersection with $\cal D$; conversely, given a $J|_{{\cal D}}$-ideal $I'$ on $\cal D$, it is clear that $I'$ is equal to the intersection of $\cal D$ with its $J$-closure in $\cal C$. 

$(ii)$ This can be deduced from part $(i)$ of the theorem by using the identification of $L$ with the frame of $J^{L}_{can}$-ideals on $L$. Anyway, we give a direct proof of this result for the reader's convenience.
 
Clearly, $\psi \circ \phi=1_{L}$. To prove that $\phi\circ \psi=1_{B_{L}}$, we have to verify that for any ideal $I$ in $Id_{J^{L}_{can}|_{B_{L}}}(B_{L})$, $I=B_{L}\cap (\mathbin{\mathop{\textrm{\huge $\vee$}}\limits_{}I})\downarrow$. The inclusion $\subseteq$ is obvious, so it remains to prove the converse one. Given $s\in B_{L}\cap (\mathbin{\mathop{\textrm{\huge $\vee$}}\limits_{}I})\downarrow$, $s=s\wedge \mathbin{\mathop{\textrm{\huge $\vee$}}\limits_{}I}=\mathbin{\mathop{\textrm{\huge $\vee$}}\limits_{a\in I}s\wedge a}$; but, $B_{L}$ being a basis of $L$, each $s\wedge a$ can be written as a join of elements of $B_{L}$ which belong to $I$ (they being less then or equal to $a$, which belongs to $I$). The fact that $I$ is a $J^{L}_{can}|_{B_{L}}$-ideal on $B_{L}$ then implies that $s$ belongs to $I$, as required.  
\end{proofs}

The method of section \ref{dualtop} for `lifting' dualities or equivalences of categories of preorders with categories of frames (equivalently, with categories of locales) admits a natural frame-theoretic interpretation. Let us begin by specializing the notions of subterminal topology and the construction of the category of toposes paired with points of section \ref{subterminal} in the context of locales.

Recall that the points of the locale $L$ are the frame homomorphisms $L\to \{0,1\}$ from $L$ (regarded as a frame) to the two-element frame $\{0,1\}$.

The indexing functions $\xi:X\to \textsc{P}_{l}$ of a set of points $\textsc{P}_{l}$ of a locale $L$ by a set $X$ correspond precisely to the frame homomorphisms $L\to {\mathscr{P}}(X)$. Indeed, ${\mathscr{P}}(X)$ is the product in the category $\textbf{Frm}$ of frames of $X$-times the frame $\{0,1\}$; for any $x\in X$, we have a product projection $\overline{x}:{\mathscr{P}}(X)\to \{0,1\}$ sending a subset $S\in {\mathscr{P}}(X)$ to $1$ if $x\in S$ and to $0$ otherwise, and for any indexing function $\xi$ sending any element $x\in X$ to a frame homomorphism $\xi(x):L\to \{0,1\}$, there exists a unique frame homomorphism $g_{\xi}:L\to {\mathscr{P}}(X)$ such that for every $x\in X$, $\xi(x)=\overline{x}\circ g_{\xi}$, defined by the formula
\[
g_{\xi}(l)=\{x\in X \textrm{ | } \xi(x)(l)=1\}
\] 
for any $l\in L$.

Given an indexing $\xi:X\to \textsc{P}$ of a set of points of a locale $L$, the construction of the subterminal topology on the set $X$ (cf. section \ref{subterminal}) reformulates as follows; the underlying set of the topological space is $X$, while the open sets are the subsets in the image of the frame homomorphism $g_{\xi}:L\to {\mathscr{P}}(X)$ corresponding to the indexing $\xi$ as specified above, in other words the subsets of the form $g_{\xi}(l)=\{x\in X \textrm{ | } \xi(x)(l)=1\}$ where $l$ ranges among the elements of $L$. We denote the topological space induced via this construction by an indexing $\xi:X\to \textsc{P}$ of a set of points of a locale $L$ by $X_{{\tau}^{L}_{\xi}}$.  

The construction of the category of toposes paired with points of section \ref{subterminal} specializes, when restricted to the context of locales, to the following construction, of the category $\textbf{Loc}_{p}$ of \emph{locales paired with points}: the objects of $\textbf{Loc}_{p}$ are the pairs $(L, \xi)$ where $L$ is a locale and $\xi:X\to \textsc{P}$ is an indexing of a set of points $\textsc{P}$ of $L$, while the arrows $(L, \xi)\to (L', \xi')$ in $\textbf{Loc}_{p}$, where $\xi:X\to \textsc{P}$ and $\xi':Y\to \textsc{P}'$, are the pairs $(f, l)$ where $f$ is a frame homomorphism $f:L'\to L$ and $l:X\to Y$ is a function such that the diagram
\[  
\xymatrix {
{\mathscr{P}}(Y) \ar[r]^{l^{-1}}  & {\mathscr{P}}(X)  \\
L' \ar[u]^{g_{\tilde{\xi'}}} \ar[r]^{f} & L \ar[u]^{g_{\tilde{\xi}}}}
\]
commutes.

Identities and composition in $\textbf{Loc}_{p}$ are defined in the obvious way.   

Given an indexing $\xi:X\to \textsc{P}$ of a set of points of a locale $L$, we say that $\xi$ is \emph{separating} if for any $l,l'\in L$, $\xi(x)(l)=\xi(x)(l')$ for every $x\in X$ implies that $l=l'$.

In these terms, the method of \ref{dualtop} for `lifting' dualities or equivalences of categories of preorders with categories of frames reformulates as follows. 

Let us first consider the contravariant case. 

Let as assume to start with a duality $A:{\cal K}^{\textrm{op}}\to \textbf{Loc}$ obtained by the method of section \ref{locales}. Suppose that we have assigned to every structure $\cal C$ in $\cal K$ a separating indexing $\xi_{\cal C}:X_{\cal C}\to \textsc{P}_{\cal C}$ of a set of points $\textsc{P}_{\cal C}$ of the locale $Id_{J_{\cal C}}({\cal C})$, and to each arrow $f:{\cal C}\to {\cal D}$ in $\cal K$ a function $l_{f}:X_{\cal D} \to X_{\cal C}$ such that the pair $(A_{f}, l_{f})$ (see the statement of Theorem \ref{concretecontrov} above for the definition of the frame homomorphism $A_{f}:Id_{J_{\cal C}}({\cal C}) \to Id_{J_{\cal D}}({\cal D})$) defines an arrow $(Id_{J_{\cal D}}({\cal D}), \xi_{\cal D}) \to (Id_{J_{\cal C}}({\cal C}), \xi_{\cal C})$ in the category $\textbf{Loc}_{p}$ of locales paired with points. Then we have a functor $\tilde{A}:{\cal K}^{\textrm{op}}\to \textbf{Top}$ such that $\tilde{A}({\cal C})={X_{\cal C}}_{{{\tau}^{Id_{J_{\cal C}}({\cal C})}_{\xi_{\cal C}}}}$ for any ${\cal C}\in {\cal K}$ and $\tilde{A}(f)=l_{f}:X_{\cal D} \to X_{\cal C}$ for any arrow $f:{\cal C}\to {\cal D}$ in $\cal K$.

In the covariant case, let us assume to start with a functor $B:{\cal K} \to \textbf{Loc}$ obtained by the method of section \ref{locales}. Similarly as above, suppose that we have assigned to each structure $\cal C$ in $\cal K$ a separating indexing $\xi_{\cal C}:X_{\cal C}\to \textsc{P}_{\cal C}$ of a set $\textsc{P}_{\cal C}$ of points of the locale $Id({\cal C}^{\textrm{op}})$, and to each arrow $f:{\cal C}\to {\cal D}$ in $\cal K$ a function $l_{f}:X_{\cal C} \to X_{\cal D}$ such that the pair $(B_{f^{\textrm{op}}}, l_{f})$ (see the statement of Theorem \ref{concretecovariant} above for the definition of $B_{f^{\textrm{op}}}$) defines an arrow $(Id({\cal C}^{\textrm{op}}), \xi_{\cal C}) \to (Id({\cal D}^{\textrm{op}}), \xi_{\cal D})$ in the category $\textbf{Loc}_{p}$. Then we have a functor $\tilde{B}:{\cal K}\to \textbf{Top}$ such that $\tilde{B}({\cal C})={X_{\cal C}}_{{{\tau}^{Id({\cal C}^{\textrm{op}})}_{\xi_{\cal C}}}}$ for ${\cal C}\in {\cal K}$ and $\tilde{B}(f)=l_{f}:X_{\cal C} \to X_{\cal D}$ for any arrow $f:{\cal C}\to {\cal D}$ in $\cal K$.

We refer the reader to section \ref{dualtop} for a discussion of the properties of the functors $\tilde{A}$ and $\tilde{B}$. 

Finally, let us give elementary proofs of two results in the paper concerning frames of ideals on preorders, namely Theorem \ref{freeframes} and \ref{bijectionpoints}.

Let us start with Theorem \ref{freeframes}. This result is based on the notion of filtering map, which we recall below.

\begin{definition}[cf. Proposition \ref{filtering}]
Let $\cal C$ be a preorder, $L$ be a frame and $f:{\cal C}\to L$ be a monotone map. We say that $f$ is \emph{filtering} if the following conditions hold:
\begin{enumerate}[(i)]
\item $1_{L}=\mathbin{\mathop{\textrm{\huge $\vee$}}\limits_{c\in {\cal C}}}f(c)$;

\item For any $c,c'\in {\cal C}$, $f(c)\wedge f(c')=\mathbin{\mathop{\textrm{\huge $\vee$}}\limits_{b\in B_{c, c'}}}f(b)$ where $B_{c, c'}$ is the set
\[
\{b\in {\cal C} \textrm{ | } b\leq c \textrm{ and } b\leq c'\}.
\]  
\end{enumerate}
\end{definition}

Given a map $f:{\cal C}\to L$ and a coverage $J$ on $\cal C$, we say that $f$ is \emph{$J$-filtering} if $f$ is filtering and satisfies the property that for any $c\in {\cal C}$ and any $S\in J(c)$, $f(c)=\mathbin{\mathop{\textrm{\huge $\vee$}}\limits_{d\in S}}f(d)$.

Let us now give a direct proof of Theorem \ref{freeframes}.

\begin{theorem}[cf. Theorem \ref{freeframes}]\label{freeframeselem}
Let $\cal C$ be a preorder and $J$ be a coverage on $\cal C$. Then the frame $Id_{J}({\cal C})$, together with the map $\eta:{\cal C}\to Id_{J}({\cal C})$ sending an element $c\in {\cal C}$ to the principal $J$-ideal $(c)\downarrow_{J}$, satisfies the following universal property: for any  map $f:{\cal C}\to L$ to a frame $L$, $f$ is $J$-filtering if and only if there exists a (necessarily unique) frame homomorphism $\tilde{f}:Id_{J}({\cal C})\to L$ such that $\tilde{f}\circ \eta=f$ (given by the formula $\tilde{f}(I)=\mathbin{\mathop{\textrm{\huge $\vee$}}\limits_{c\in I}}f(c)$ for any $I\in Id_{J}({\cal C})$).
\end{theorem}

\begin{proofs}
Given a frame homomorphism $g:Id_{J}({\cal C})\to L$, it is clear that the composite $g\circ \eta:{\cal C}\to L$ is a filtering map. Indeed, condition $(i)$ in the definition of filtering map follows from the fact that $\mathbin{\mathop{\textrm{\huge $\vee$}}\limits_{c\in {\cal C}}}(g\circ \eta)(c)=g(\mathbin{\mathop{\textrm{\huge $\vee$}}\limits_{c\in {\cal C}}}(c)\downarrow_{J})=g({\cal C})=1_{L}$, while condition $(ii)$ follows from the fact that for any $c, c'\in {\cal C}$, $(g\circ \eta)(c) \wedge (g\circ \eta)(c')=g((c)\downarrow_{J} \cap (c')\downarrow_{J})=g(\mathbin{\mathop{\textrm{\huge $\vee$}}\limits_{b\in B_{c, c'}}}(b)\downarrow_{J}))=\mathbin{\mathop{\textrm{\huge $\vee$}}\limits_{b\in B_{c, c'}}}g((b)\downarrow_{J})=\mathbin{\mathop{\textrm{\huge $\vee$}}\limits_{b\in B_{c, c'}}}(g\circ \eta)(b)$. The fact that $g\circ \eta$ satisfies the property that for any $c\in {\cal C}$ and any $S\in J(c)$, $(g\circ \eta)(c)=\mathbin{\mathop{\textrm{\huge $\vee$}}\limits_{d\in S}}(g\circ \eta)(d)$ is obvious.

Conversely, given a filtering map $f:{\cal C}\to L$, there is exactly one frame homomorphism $\tilde{f}:Id_{J}({\cal C})\to L$ such that $\tilde{f}\circ \eta=f$. Indeed, $\tilde{f}$ is forced (by the fact that it must preserve arbitrary joins) to be equal to the map sending an ideal $I$ in $Id_{J}({\cal C})$ to the join $\mathbin{\mathop{\textrm{\huge $\vee$}}\limits_{c\in I}}f(c)$. The fact that $\tilde{f}$ preserves the top element (resp. binary meets) easily follows from condition $(i)$ (resp. condition $(ii)$) in the definition of filtering map, while the fact that $\tilde{f}$ preserves arbitrary joins follows from the fact that for any $c\in {\cal C}$ and any $S\in J(c)$, $f(c)=\mathbin{\mathop{\textrm{\huge $\vee$}}\limits_{d\in S}}f(d)$. 
\end{proofs}

Let us now turn our attention to Theorem \ref{bijectionpoints}. To this end, we recall from section \ref{exsub} the notion of $J$-prime filter on a preorder $\cal C$ equipped with a coverage $J$: a \emph{$J$-prime filter} on $\cal C$ is a subset $F\subseteq {\cal C}$ such that $F$ is non-empty, $a\in F$ implies $b\in F$ whenever $a\leq b$ in $\cal C$, for any $a, b\in F$ there exists $c\in F$ such that $c\leq a$ and $c\leq b$, and for any $J$-covering sieve $\{a_{i} \to a \textrm{ | } i\in I\}$ in $\cal C$ if $a\in F$ then there exists $i\in I$ such that $a_{i}\in F$. If $\cal C$ is a frame and $J$ is the canonical coverage on $\cal C$, the notion of $J$-prime filter on $\cal C$ specializes to the notion of \emph{completely prime filter} on $\cal C$. 

Theorem \ref{bijectionpoints} asserts that the assignment sending a filter $F$ on $Id_{J}({\cal C})$ to the $J$-prime filter $\{c\in {\cal C} \textrm{ | } (c)\downarrow_{J}\in F\}$ on $\cal C$ defines a bijection between the completely prime filters on the frame $Id_{J}({\cal C})$ of $J$-ideals of $\cal C$ and the $J$-prime filters on $\cal C$. This result can be deduced as the particular case of Theorem \ref{freeframeselem} for $L=\{0,1\}$. Indeed, the $J$-filtering maps ${\cal C} \to \{0,1\}$ correspond exactly to the $J$-prime filters on $\cal C$, while the frame homomorphisms $Id_{J}({\cal C})\to \{0,1\}$ correspond exactly to the completely prime filters on $Id_{J}({\cal C})$.

\vspace{0.5cm}
\textbf{Acknowledgements:} I am grateful to Peter Johnstone for useful comments on a preliminary version of this paper.

\end{document}